\documentclass[11pt]{article}
\usepackage{graphicx}
\usepackage{amssymb,amsfonts,amsthm}
\usepackage{enumerate}
\usepackage[T2B]{fontenc}
\usepackage[cp1251]{inputenc}
\usepackage{enumitem}
\usepackage{amsmath,amsthm,amssymb}
\usepackage{amscd}
\usepackage{amsfonts}
\usepackage{color}

\newcommand{\pr}{\pageref}

\newtheorem{theorem}{Theorem}[section]
\newtheorem{lemma}[theorem]{Lemma}

\newtheorem{cy}[theorem]{Corollary}

\theoremstyle{definition}
\newtheorem{df}[theorem]{Definition}

\newtheorem{rk}[theorem]{Remark}

\newcommand{\area}{\mathrm{Area}}

\newcommand{\me}{\medskip}

\newcommand{\Lab}{{\mathrm{Lab}}}

\newcommand{\tool}{\stackrel{\ell}{\too} }

\newcommand{\ttt}{{\cal T}}

\newcommand{\bb}{{\cal B}}
\newcommand{\topp}{{\bf top}}
\newcommand{\ttopp}{{\bf ttop}}
\newcommand{\tbott}{{\bf tbot}}
\newcommand{\bott}{{\bf bot}}
\newcommand{\vk}{van Kampen }

\newcommand{\iv}{^{-1}}

\newcommand{\too}{\to }

\newcommand{\qq}{{\cal Q} }

\newcommand{\sss}{{\cal S} }

\begin{document}
\renewcommand{\theequation}{\thesection.\arabic{equation}}
\bigskip

\title{Polynomially-bounded Dehn functions of groups}

 \author{A.Yu. Ol'shanskii\thanks{The
author was supported in part by the NSF grant DMS-1500180 and by
the Russian Fund for Basic Research,
grant 15-01-05823}}
\date{}
\maketitle

\begin{abstract} It is well known that every subqadratic Dehn function is linear. A question by Bridson asked to describe the isoperimetric spectrum of groups, that is the set of all numbers $\alpha$ such that $n^\alpha$ is equivalent to the Dehn function of a finitely presented group.
The goal of this paper is to give a description of the isoperimetric spectrum. Earlier a similar description was given by Sapir, Birget and Rips for the intersection of the isoperimetric spectrum with $[4,\infty]$. Lowering the bound from 4 to 2 required significant new ideas and tools.
\end{abstract}


\medskip


\tableofcontents

\section{Introduction}

\subsection{Formulation of the theorem and corollaries}

The minimal non-decreasing function $d\colon \mathbb{N}\to \mathbb{N}$ such that
every word $w$ vanishing in a group $G=\langle A\mid R\rangle$ and having length \label{lengthw||} $||w||\le n$ is
freely equal to a product of at most $d(n)$ conjugates of relators
from $R^{\pm 1}$,
is called the \label{Dehnf}{\em Dehn function} of the presentation
$G=\langle A\mid R\rangle$ \cite{GrHyp}. In other words, the the value $d(n)$ is the smallest integer that bounds from above the areas of loops of length $\le n$ in the Cayley complex $Cay(G)$, and so
by van Kampen's Lemma, $d(n)$ is equal to the maximal
area of minimal filling diagrams $\Delta$ with perimeter $\le n.$ (See Subsection \ref{md} for the definitions.)

The values $d(n)$ are defined if the set of generators $A$ is finite. For a {\it finitely presented} group (i.e., both sets $A$ and $R$ are finite), the
Dehn function exists and it is usually taken up to equivalence to get rid of the dependence  on
a finite presentation of $G$ (see \cite{MO}). To introduce this \label{equivf}{\em equivalence} $\sim$, let $f$ and $g$ be
non-decreasing functions $\mathbb{N}\to \mathbb{R}_+$.
We write $f\preceq g$
if there is a positive integer $c$ such that
$f(n)\le cg(cn)+cn\;\;\; for \;\; every \;\;n\in \mathbb{N}.$
Two non-decreasing functions $f$ and $g$ on $\mathbb{N}$ are called equivalent if
$f\preceq g$ and $g\preceq f.$

Note that for many functions (for example, for $n^{\alpha}$, $n^{\alpha}(\log n)^{\beta}(\log\log n)^{\gamma}$,
 and so on),
their $\sim$-equivalence classes coincide with their \label{Tequiv}$\Theta$-{\it equivalence} classes,
where the symbol $\Theta$ is borrowed from the theory of computational complexity: one says that $f(n)=\Theta(g(n))$ if both properties $f(n)=O(g(n))$ and $g(n)=O(f(n))$ hold.

The Dehn function $d=d_G$ of a finitely presented group $G$ is also called an isoperimetric function of $G$ since it is equivalent to the usual isoperimetric function of a simply connected Riemannian manifold $M$, provided $G$ acts properly and co-compactly on $M$
by isometries. So the concept of Dehn function is derived from geometry, and one can find much more regarding this connection in \cite{Gr}.

Another connection is to Computational Complexity. The algorithmic word problem in a finitely presented group is decidable if and only if the Dehn function is recursive,
and the Dehn function of a group bounds the computational complexity of the word problem. It was shown in \cite{BORS} that conversely, every recursively presented group $G$ with decidable word problem embeds into a finitely presented group whose Dehn function is only polynomially larger than the computational complexity (the time function) of the word problem in $G$. In particular, groups with word problem in {\bf NP} are precisely the subgroups of finitely presented groups $G$ with at most polynomial Dehn functions $d_G$.

For every positive integer $\alpha$, there are (nilpotent) groups with Dehn function $n^{\alpha}$ \cite{BMS}. The first examples
of Dehn functions $n^{\alpha}$ with non-integer $\alpha$ can be found in \cite{B}, where the description of possible exponents $\alpha$, forming the {\it isoperimetric spectrum},  was called the most fundamental question concerning isoperimetric functions. (Obviously some conditions on the real exponent $\alpha$ are inevitable since the set of real numbers is uncountable, while the set of non-isomorphic finitely presented groups is countable.)

Almost all possible Dehn functions $F(n)\ge n^4$ of finitely presented groups were described in \cite{SBR} in terms of time functions of non-deterministic Turing machines. By Theorem 1.2 \cite{SBR}, to obtain a group with Dehn function $F$, it suffices to assume that  the function $F$ is super-additive (i.e. $F(m+n)\ge F(m)+F(n)$ for $m,n\in \mathbb{N}$) and the integral part of $\sqrt[4]{F(n)}$ is a time function of a non-deterministic Turing machine (see the definition in \cite{DK}).
As a corollary, it was proved in \cite{SBR} that if $\alpha\ge 4$ and the real number $\alpha$ is  computable in time $\le 2^{2^m},$
then there is a finitely presented group with Dehn function equivalent to $n^{\alpha}$ (One should use the integral part sign
for functions on $\mathbb{N}$, but we omit this sign speaking on
asymptotic behavior.).
The computability can be defined as follows.

\begin{df}\label{co}
 Let $T\colon \mathbb{N}\to\mathbb{N}$. A real number $\alpha$ \label{compT} {\it is computable in time} $\le T(m)$ if there exists a Turing machine which, given a natural number $m,$ computes a binary rational approximation of $\alpha$
with an error  $O(2^{-m}),$ and the time of this computation $\le T(m).$
\end{df}

The algebraic and many transcendental
numbers are computable much faster, and so there are examples of groups with
Dehn functions equivalent to $n^{\pi+e}$, and so on.
If, conversely, for a real $\alpha$, the function $n^\alpha$ is equivalent to the Dehn function of a finitely presented group, then for some $c$, the exponent $\alpha$ is computable in time $2^{2^{2^{cm}}}$ \cite{SBR}.



Still the class of Dehn functions $< n^4$ was unclear even though it has drawn
attention (see, for example, \cite{B}, \cite{BB}, \cite{BBFS}). We note the paper \cite{BB} (also the references there), where the Dehn functions equivalent to $n^{\alpha}$ were constructed for special non-integer exponents of the form $\alpha=2\log_2\frac {2p}{q}$, where $p$ and $q$ are integers, $p>q>0$.

Since all finitely presented groups with subquadratic Dehn functions are hyperbolic \cite{GrHyp, Ol91, Bow}, i.e. their Dehn functions are in fact linear, the only interval
$2<\alpha< 4$ remained misty. To formulate a theorem that fills this gap, we need

\begin{df}\label{sui}
We say that a non-decreasing function $f: {\mathbb N}\to {\mathbb N}$ is {\it suitable} if the following properties hold.
\begin{itemize}

\item{$f(n)^3= O(n)$}

\item{For every integer $c>0$, there is $C>0$ such that $f(cn)\le Cf(n)$ , i.e., $$f(O(n))=O(f(n)).$$}

\item{There is a (non-deterministic) Turing machine ${\bf M}_0$
 recognizing the values of the function $f(n)$ with computation time $O(n^{1/3})$.

 It works as follows. An integer $k\ge 0$ is an input of ${\bf M}_0$
in the form $a^k$ for a fixed letter $a$. The Turing machine ${\bf M}_0$ produces a value $f(n)$ for some $n\ge 1$, i.e. it obtains the word $c^{f(n)}$
on a special
tape. Then ${\bf M}_0$ compares $f(n)$ and $k$, accepting $k$ if $k=f(n)$.
It can accept the input word $a^k$
if and only if $k=f(n)$ for some natural number $n$.}
\end{itemize}
\end{df}

\begin{theorem} \label{beta} For every suitable function $f(n)$ and every integer $s\ge 2$,
the function $F(n)=n^sf(n)^3$ is equivalent to the Dehn function of
a finitely presented group.
\end{theorem}

Given a suitable function $f(n)$, we denote $g(n)=f(n)^3$.

\begin{cy} \label{22m} If $\alpha\ge 2$ and the real number $\alpha$ is  computable in time $O(2^{2^m}),$
then there is a finitely presented group with Dehn function equivalent to $n^{\alpha}$.
\end{cy}

\begin{rk} It is easy to see that we  have  an equivalent statement when replacing $2$ with any integer $d>1$ both in the statement of Corollary \ref{22m} and in Definition
\ref{co} (resp., "binary" with "d-ary").
The formulation and the proof of Corollary \ref{22m} are close to those for
Corollary 1.4 \cite{SBR} (up to a minor inaccuracy in the formulation of Corollary 1.4  \cite{SBR}). However below we give a proof of our Corollary 1.4 since one should draw it from  different assumptions of Theorem \ref{beta}.
\end{rk}

{\it Proof} of Corollary \ref{22m}. Assume that $\alpha$ is computable in time $O(2^{2^m})$.
Clearly, the same property holds for the number $\beta=\frac13 (\alpha - s)$, where $s=[\alpha]$.

Let $M_0$ be a (non-deterministic) Turing machine producing $r=[n^{1/3}]\ge 1$ (in
unary) with time $O(r)$.
Then it computes $[\log_2 r]$ and $m= [\log_2\log_2 r]$ (in binary) with time $O(r)$ using divisions by $2$.

It follows from the assumption of the corollary that  one can recursively compute
binary rational numbers $\beta_m$ such that
\begin{equation}\label{nalpha}
|\beta-\beta_m|=O(2^{-m})= O((\log_2 r)^{-1})
\end{equation}
and the time of the computation of $\beta_m$ is $O(r).$ Let
${\bf M}_0$ accomplish this computation.
In addition, one may assume that the number of digits in the
binary expansion of $\beta_m$ is $O(m).$ Therefore the computation
of the product $[\beta_m[\log_2 r]]$ needs  time at most $O((\log_2 n)^2).$
Since $r=[n^{1/3}]$, the Turing machine can now obtain the binary presentation $\gamma_m$ of $[\beta_m[\log_2 n]]$
with time $O((\log_2 n))$ and error $O(1)$.

Next, let ${\bf M}_0$ rewrite the binary presentation of $\gamma_m$
in unary (as a sequence of $1$-s). This well-known rewriting (e.g., see p.352 in \cite{SBR}) has
time complexity of the form $O(\gamma_m)=O(\beta_m\log_2 n).$ One more ${\bf M}_0$-rewriting
of this type applied to the unary presentation of $\gamma_m$
(considered now as binary one),
will have time complexity of the form $$\Theta(2^{\gamma_m})=\Theta(2^{[\beta_m[\log_2 n]]+O(1)})=\Theta(2^{[\beta_m[\log_2 n]]}).$$
One can rewrite as
$\Theta(2^{\beta\log_2 n})=\Theta(n^{\beta})=O(n^{1/3})$ by inequalities  (\ref{nalpha}) and
$\beta<1/3$, because $\log_2 n=O(2^m)$.

During the last (deterministic) rewriting, one can count the number $f(n)$ of
commands and obtain a word $b^{f(n)}$ on a special tape, where
$f(n)=\Theta(n^{\beta})$. It is easy to see that the rewriting can be defined
so that the function $f(n)$ is non-decreasing.

One more tape of the Turing machine ${\bf M}_0$ under construction
has the input word $a^k$. It remains to check whether the lengths of the words
$b^{f(n)}$ and $a^k$ are equal or not. This takes the time $O(f(n))=O(n^{1/3}).$
Since the time of the entire procedure is $O(n^{1/3})$, the function $f(n)$ is suitable.
Now by Theorem \ref{beta}, the function $n^sf(n)^3=\Theta(n^s n^{3\beta})=\Theta(n^{\alpha})$
is equivalent to the Dehn function of a finitely presented group.
$\Box$

In particular, Corollary \ref{22m} implies

\begin{cy} The functions $n^{\alpha}$ for every real algebraic $\alpha\ge 2$, the functions $n^{\pi-1}$, $n^{\sqrt{e}+1}, \dots$
are equivalent to Dehn functions of finitely presented groups.
\end{cy} $\Box$

 As we mentioned above, the analog of Corollary \ref{22m} for $\alpha\ge 4$ was proved in \cite{SBR}. But weakening the restriction to $\alpha\ge 2$ now, although uses $S$-machines, as in \cite{SBR}, it requires substantially new ideas.

Theorem \ref{beta} gives a tremendous class of new Dehn functions  of the form $O(n^4)$. The following examples can be validated in absolutely similar way as Corollary \ref{22m}.

\begin{cy} The functions $n^{\alpha}(\log n)^{\beta}$, $n^{\alpha}(\log n)^{\beta}(\log\log n)^{\gamma},\dots $ are equivalent to the Dehn functions of finitely
presented groups, provided the real $\alpha,\beta,\gamma,\dots $ are computable in time $O(2^{2^m})$ and $\alpha>2$ or $\alpha=2$ and $\beta>0$, or $\alpha=2, \beta=0,\gamma>0$ ...
\end{cy}$\Box$

Since every finitely presented group is
a fundamental group of a connected closed Riemannian manifold $X$
and therefore acts properly and
co-compactly by isometries on its universal cover $\tilde X$, one can use Theorem
\ref{beta} and Corollary \ref{22m} to formulate one more

\begin{cy} \label{Riem} For every function $F(n)$ satisfying the assumption of Theorem \ref{beta},
there exists a closed connected Riemannian manifold $X$ such that the isoperimetric function of the universal cover $\tilde X$ is equivalent to $F(n)$.

In particular,
if  a real number $\alpha\ge 2$ is computable in time $O( 2^{2^m}),$
then there exists such a universal cover $\tilde X$ with  isoperimetric function equivalent to $n^{\alpha}.$
\end{cy}
$\Box$

Since the condition $\alpha \ge 2$ is the best possible and the obtained upper bound for the Dehn function must be equivalent to the lower one, all inequalities throughout this paper should  be uniformly sharp, up to multiplicative constants.

We collect all the definitions and  terms at the end of the paper (see Subject Index). The next subsection presents a short outline of the plan.

\subsection{Brief description of the proof of Theorem \ref{beta}}

The idea of simulation of the commands of a Turing machine by group relations goes back
to the works of P.Novikov, W.Boone and many other authors (see \cite{R}, \cite{S}). However one has to
properly code the work of a Turing machine in terms of group relations, and the interpretation problem for groups remains much harder than for semigroups, because the group theoretic simulation can execute unforeseen computations with non-positive
words. Boone and Novikov secured the positiveness of admissible configurations  with the help of an
additional `quadratic letter' (see \cite{R}, Ch.12). However this old trick implies that the constructed
group $G$ contains Baumslag - Solitar subgroups $B_{1,2}$ and has at least exponential Dehn function. Since we want to obtain at most polynomial Dehn functions, we use the S-machines introduced in \cite{SBR}. Those S-machines invented by M. Sapir can work with non-positive words on the tapes and they are polynomially equivalent to classical Turing machines.

According to the original version, $S$-machines are special
rewriting systems. All necessary definitions are given in Subsection \ref{SM}. On the other hand, the state, tape and command
letters of an $S$-machine can be regarded as group generators,
and the commands can be interpreted as defining relations
(see Subsection \ref{MG}). The obtained group $M$ is a multiple
HNN-extension of a free group. Every computation of the $S$-machine
is simulated by  van Kampen diagram over this group called trapezia (Subsection \ref{md}).

To construct a finitely presented groups $G$ with desired Dehn functions, one needs to add to $M$ a special relation called the \emph{hub}. It consists of state letters. There are very particular
van Kampen diagrams, called disks, built of the hub and many trapezia attached around.

It is proved in \cite{SBR} that for every  Turing machine ${\bf M}_0$ with  time function $T(n)$, there is an equivalent $S$-machine ${\bf M}_1$ with
time function (and the generalized time function) $\sim T(n)^3$.
It follows that if an accepting computation starts with an input of length $\sim n$, it has  length $\sim T(n)^3$, and the computational disk
 has perimeter $\sim n$, and so its area $\sim nT(n)^3$. Since time functions are at least linear, this approach gave the lower bound $\ge n^4$ for the Dehn function of $G$.

 To get the linear time of accepting computation, in comparison with the length of the initial configuration for an $S$-machine ${\bf M}_2$, one can add an additional tape, where the whole history of the forthcoming computation is written. Then every command will erase one letter on this tape. This trick gives disks with quadratic area with respect to their perimeters.

 However we want to construct disks with prescribed area $F(n)$, as in Theorem \ref{beta}.
 In this paper, we first prove this theorem for $s=2$, i.e. for $F(n)=O(n^3)$, and in the final Subsection \ref{sc},
we show that the value of $s$ can be increased, since a non-difficult modification of the main S-machine $\bf M$ constructed in
Subsection \ref{M6} linearly slows down the work of $\bf M$.

The main S-machine is composed of the S-machine ${\bf M}_2$ repeating the same computation many times and another S-machine that can stop the computations of ${\bf M}_2$ after $\sim g(n)$ such cycles with subsequent erasure of all the tapes and acceptance.
This gives the lower bound $\sim ng(n)$ for the (general) time
function and the lower bound $\sim F(n)$ for the areas of computational disks.

The obtainment of the upper bounds is the major job in this paper. First of all, to obtain quadratic upper bound for the areas of trapezia, one needs a linear bound of the space of every computation
(i.e. the maximal lengths of all admissible words of it) in terms
of the lengths of the first and the last word. This task
is aggravated by inaccurate simulation of the work of Turing machines
by S-machines and so by group relations. Standard trapezia correspond to the prescribed work of S-machines, but there are non-standard ones simulating undesired computations when the same tape is simultaneously changed
at both ends.   The  features of standard (accurate) computations of the main S-machine  and of non-standard ones
are considered in Subsections \ref{SC} and \ref{fau}, respectively.
To reduce the effect of non-standard computations, the basic steps
of the work alternate with control steps in the definition of the main S-machine given in Subsection \ref{M6}.

Whereas a non-standard computation has linearly bounded space
in terms of the lengths of the first and the last words
(and so the width of corresponding trapezia is linearly bounded too), there exist much wider trapezia in the standard case.
Hence standard trapezia can have too large areas in the
group $M$. The new idea is decrease the area of their boundary
labels in the quotient group $G$. We do this in   Subsection \ref{end1} applying the properties of long computations
obtained in Subsection \ref{long}. It turned out that wide trapezia
with super-quadratic areas in $M$ can be replaced, preserving the boundary label, by diagrams of quadratic areas over $G$, i.e. by diagrams containing hubs.

However before that, we prescribe artificial (but quadratic !) $G$-areas to special
`big' subtrapezia (Subsection \ref{qub}), which leads to the definition of  $G$-area for an arbitrary
diagram over $M$ or over $G$. A quadratic upper bound for the
$G$-areas of diagrams over $M$ is given in Subsection \ref{qub}. It turns out later, that
the are diagrams of quadratic area over $G$
with the same boundary labels. The induction on the perimeter is based on a non-trivial surgery. (See
the proof of Lemma \ref{main} in Subsection \ref{qub}, where we cut
and paste diagrams.) In fact, we give an upper bound $G$-for area
in terms of the perimeter and the mixture of $q$- and $\theta$-letters
in the boundary label, because different types of surgeries
decrease either perimeter or the mixture. The mixture is defined
for arbitrary necklace with beads of two colors (see Subsection
\ref{mix}), and it is bounded by the square of the number of beads, and so we finally obtain a quadratic upper bound in terms
of perimeter only.

Note that instead of the combinatorial length of words (length of paths, perimeter) we consider a modified length, where different
letters and syllabi have different lengths (Subsection \ref{lf}).
With respect to this modified length $|\cdot|$, the length of the top/bottom of every $q$-band $Q$ is just the number of 2-cells in $Q$, and the rim $\theta$-bands with bounded number of $(\theta,q)$-cells
can be removed from a diagram with decrease of perimeter. Such properties are exploited in the paper many times. It is easy
to reformulate the final results in terms of the combinatorial length $||\cdot||$ since it is $\Theta$-equivalent to $|\cdot|$.

Our presentation of the group $G$ (Subsection \ref{MG}) is highly non-aspherical, and
so van Kampen diagrams with the same boundary label can differ
widely. We choose minimal diagrams in Section \ref{midi}, i.e. reduced diagrams with
minimal number of disks and, for given number of disks, with minimal number of $(\theta,q)$-cells. We do not claim
that a minimal diagram
has minimal area or minimal $G$-area  for fixed boundary label, but to obtain the upper bound for the Dehn function, it suffices to bound from above the $G$-areas of minimal diagrams.

However, even one has quadratic estimates  for $G$-areas of disk-free diagrams (i.e. diagrams over $M$) and the required upper bound $\sim F(n)$ for the areas of disks (see the definition in Section 7), the `snowman' decomposition of diagrams in the union of subdiagrams with single disk defined in \cite{SBR}, would give at least cubic upper bound for the $G$-area of the entire minimal diagram. Thus, without new tools one could only hope to weaken the restriction from \cite{SBR} to $\alpha\ge 3$.

The helpful property is that a minimal diagram cannot contain
subdiagram formed by a disk  and a very special trapezium
connected by a `shaft' (Section \ref{midi}). At first sight, it seems that
such subdiagrams are extremely rare. But they become ordinary
if the work of an S-machine has sufficiently many control steps (Steps $1^- - 5^-$ in Subsection 4.1).
Hence the absence of these exotic subdiagrams can help. And
it helps indeed provided the sum $\sigma$ of the lengths of all `shafts' (see Definition \ref{shf})
linearly bounded in terms of the perimeter.

To obtain such a linear estimate, in Section \ref{des}, we introduce designs formed by two finite sets of segments and prove a pure combinatorial proposition.
So there are neither machines nor groups, nor van Kampen diagrams
in Section \ref{des}, and the reader can start with that short section. (Is there even shorter proof or a reference to a known property? Althogh the author did believe that the linear estimate took place, he wasted time devising a proof.)

Since $\sigma\le  c n$ for some constant $c>0$, where $n$ is the perimeter, one can estimate the $G$-area of a diagram in terms of the sum $n+\sigma$ instead of $n$ (Section $\ref{ub}$), and this is another new tool for obtaining  the required upper bound for the area.

\section{General properties of S-machines}

\subsection{S-machines as rewriting systems}\label{SM}

There are several interpretations of $S$-machines in groups.
In particular, one can define an $S$-machine
as a group that is a multiple HNN extension of a free group.
Here we  slightly modify the original
definition \cite{SBR} using \cite{OS06} and  define $S$-machines as rewriting systems
working with words in group alphabets.
The precise definition
of an \label{Smachine}$S$-machine $\cal S$ is as follows.

The hardware of an $S$-machine $\cal S$ is a pair $(Y,Q),$ where $Q=\sqcup_{i=0}^NQ_i$ and $Y= \sqcup_{i=1}^N Y_i$ (for convenience we always set $Y_0=Y_{N+1}=\emptyset$).
The elements from $Q$ are called \label{statel}{\it state letters}, the elements from $Y$ are \label{tapel}{\it tape letters}. The sets $Q_i$ (resp.
$Y_i$) are called {\em parts} of $Q$ (resp. $Y$).

The {\it language of \label{admissiblew} admissible words}
consists of all reduced words $W$ of the form
\begin{equation}\label{admiss}
q_1^{\pm 1}u_1q_2^{\pm 1}\dots u_k q_{k+1}^{\pm 1},
\end{equation}
where  every subword $q_i^{\pm 1}u_iq_{i+1}^{\pm 1}$ either
\begin{itemize}
\item belongs to $(Q_jF(Y_{j+1})Q_{j+1})^{\pm 1}$ for some $j$ and $u_i\in F(Y_{j+1}),$ where $F(Y_i)$ is the set of reduced group words in the alphabet $Y_i^{\pm 1},$
or
\item has the form $quq^{-1}$ for some $q\in Q_j$ and $u\in F(Y_{j+1}),$
or
\item is of the form $q^{-1}uq$ for $q\in Q_j$ and $u\in F(Y_{j}).$
\end{itemize}

(The 2d and 3d items extend the definition of admissible words in comparison with
\cite{SBR}, and the language of admissible words is equal to the language from \cite{OS06}.)

We shall follow the tradition of calling state letters \label{qletter}{\em
$q$-letters} and tape letters \label{aletter}{\em $a$-letters}, even though we
shall sometimes use  letters different from $q$ and $a$ as state or tape letters. The number
of $a$-letters in a word $W$ is the \label{alen} $a$-{\it length} $|W|_a$ of $W$. Usually parts of the set $Q$ of state letters are denoted by capital letters. They may differ from $Q_i$ for some S-machines. For example, a set $P$ would consist of letters $p$ with various indices. Then we shall say that letters in $P$ are $p$-letters or $P$-letters.) The \label{lengthc} length
of a word $W$, i.e. the number of all letters in $W$, is denoted by $||W||$.

If a group word $W$ over $Q\cup Y$
has the form $u_0q_1u_1q_2u_2...q_su_s,$
and $q_i\in Q_{j(i)}^{\pm 1},$
$i=1,...,s$, $u_i$  are
group words in $Y$, then we shall say that the \label{basew}{\em base} of $W$ is
$base(W)\equiv Q_{j(1)}^{\pm 1}Q_{j(2)}^{\pm 1}...Q_{j(s)}^{\pm 1}$ Here $Q_i$ are just letters, denoting the parts of the set of state letters. Note that the base is not necessarily a reduced word, and the sign \label{sequiv}
$\equiv$ is used for letter-by-letter equality of words.
The subword of $W$ between the $Q_{j(i)}^{\pm 1}$-letter and the $Q_{j(i+1)}^{\pm 1}$-letter
will be called a \label{sectorw}
$Q_{j(i)}^{\pm 1}Q_{j(i+1)}^{\pm 1}$-{\em sector} of $W$. A word can have many $Q_{j(i)}^{\pm 1}Q_{j(i+1)}^{\pm 1}$-sectors.

Instead of specifying the names of the parts of $Q$
and their order as in \\$Q=Q_0\sqcup Q_2\sqcup ... \sqcup Q_N$, we say that \label{standardb}{\em the standard} base
of the $S$-machine
is $Q_0...Q_N$. An admissible word with standard base is called a \label{config}{\em configuration} of the S-machine.

An $S$-machine also has a set of rewriting \label{rule}{\em rules} $\Theta$.  To every $\theta\in \Theta$,
two sequences of
reduced words are asigned:
$[U_1,...,U_m]$,
$[V_1,...,V_m]$, and a subset $Y(\theta)=\cup
Y_j(\theta)$ of $Y$, where \label{Yit} $Y_j(\theta)\subseteq Y_j$.
A rule
has the form:
$$[U_1\to V_1,...,U_m\to V_m],$$
where the following conditions hold:
\begin{itemize}
\item Each of $U_i$ and $V_i$ is a subword of an admissible word,
both  $U_i$ and $V_i$ have base $Q_{\ell}Q_{\ell+1}\dots Q_r$  ($\ell=\ell(i)\le
r=r(i)$) and have $a$-letters from $Y(\theta)$,
\item  $\ell(i+1)=r(i) +1$ for $i=1,\dots m-1$.

\item $U_1$ and $V_1$ must start with a $Q_0$-letter and
 $U_m$ and $V_m$ must end with a $Q_{N}$-letter.
\end{itemize}

The pair of words $U_i,V_i$ is called a \label{part}{\em part} of the rule, and is denoted $[U_i\to V_i]$.

The notation $\theta:\;[U_1\to
V_1,...,U_m\to V_m]$ contains all the necessary
information about the rule except for the sets $Y_j(\theta)$. In
most cases it will be clear what these sets are, and very often the sets $Y_j(\theta)$ will be equal to
either $Y_j$ or $\emptyset$. By default $Y_j(\theta)=Y_j$.

Every $S$-rule $\theta=[U_1\to V_1,...,U_m\to V_m]$ has an inverse
$\theta\iv=[V_1\to U_1,...,V_m\to U_m]$ which is also a rule of $\sss$; we set
$Y_i(\theta\iv)=Y_i(\theta)$. We always divide the set of rules
$\Theta$ of an $S$-machine into two disjoint parts, $\Theta^+$ and
$\Theta^-$ such that for every $\theta\in \Theta^+$, $\theta\iv\in
\Theta^-$ and for every $\theta\in\Theta^-$, $\theta\iv\in\Theta^+$ (in particular $\Theta\iv=\Theta$, that is any $S$-machine is symmetric).
The rules from $\Theta^+$ (resp. $\Theta^-$) are called \label{positiver}{\em
positive} (resp. {\em negative}). In particular,
$[U_1\to U_1; \dots; U_m\to U_m]$ is never an S-rule. It is always the case that
$Y_i(\theta\iv)=Y_i(\theta)$ for every $i$.

For every word $U_i\equiv u_0q_lu_1q_{l+1}\dots q_ru_{r-l+1}$ from the definition of the rule $\theta$, we denote by $\bar U_i$ its trimmed
form $q_lu_1q_{l+1}\dots q_r$ starting and ending
with state letters.
 To \label{application} apply an $S$-rule $\theta$ to an admissible word $W$ (\ref{admiss})
means

\begin{itemize}
\item{to check if all tape letters of $W$ belong to the alphabet $Y(\theta)$ and every state letter of $W$ is contained in some subword $\bar U_i^{\pm 1}$ of $W$,}

\item{if $W$ satisfies
 this condition, then to replace simultaneously every subword $\bar U_i^{\pm 1}\equiv (q_lu_1q_{l+1}\dots q_r)^{\pm 1}$ by subword $(u_0^{-1}V_iu_{r-l+1}^{-1})^{\pm 1}\equiv(u_0^{-1}v_0q'_1v_1\dots q'_rv_{r-l+1}u_{r-l+1}^{-1})^{\pm 1}$ ($i=1,\dots,m$)},

 \item{to  trim a few first and last $a$-letters (to obtain an admissible
 word starting and ending with $q$-letters)
 followed by the reduction of
 the resulted word.}

\end{itemize}

 The following convention is important in the definition of
$S$-machine: {\it After every application of a rewriting rule, the word is automatically
reduced. The reduction is not regarded as a separate step of an $S$-machine.}

For example, applying the rule $\theta\colon [q_1\to a q_1'b^{-1}, q_2\to cq_2'd]$ to the admissible word $W\equiv q_1b q_2dq_2^{-1} q_1^{-1}$ one first obtains the word
$aq_1'b^{-1}b cq_2'ddd^{-1} (q_2')^{-1} c^{-1} b(q_1')^{-1} a^{-1},$ then after trimming and reduction one has  $q_1'cq_2'd (q_2')^{-1} c^{-1} b(q_1')^{-1}.$  But the rule $\theta$ would not
be applicable to $W$ if $Y_2(\theta)=\emptyset$
or $Y_2(\theta)=\{a'\}$, where $a'\ne b$.

 If a rule $\theta$ is applicable to an admissible word $W$ (i.e., $W$ belongs to the \label{domain} domain
 of $\theta$) then the word $W$ is called
 $\theta$-{\it admissible}, and $W\cdot \theta$ denotes the word obtained after the
 application of $\theta$.

 A \label{computation}{\em computation} of {\it length} or {\it time} $t\ge 0$ is a sequence of admissible words $W_0\to \dots\to W_t$ such that for
every $i=0,..., t-1$ the S-machine passes from $W_i$ to $W_{i+1}$ by applying one
of the rules $\theta_i$ from $\Theta$.  The word $H=\theta_1\dots\theta_t$ is called the \label{historyc} {\it history}
of the computation. Since $W_t$ is determined by $W_0$ and the history $H$, we use the notation $W_t=W_0\cdot H$.

A computation is called \label{reducedc} {\em reduced} if its history is a reduced word. Clearly, every computation can be made reduced (without changing the initial and final words of the computation) by removing consecutive mutually inverse rules.

An S-machine is called {\it recognizing} if it
has the following attributes. There are admissible words with the standard base called \label{inpt} {\it input configurations} and \label{acceptc} {\it accept (stop) configuration.} There are {\it input sectors} (at least one) and other sectors are empty for input configurations, and all sectors are empty for the accept one. (However in this paper, some S-machines have no input or accept configurations.)
The state letters of the input (of accept) configuration form a special vector $\vec s_1$ (vector $\vec s_0$)
whose letters are involved in one rule only and are completely changed by this rule.

A configuration $W$ is said to be \label{acceptedc} {\em
accepted} by an S-machine $M$ if there exists at least one computation, called {\it accepting computation}, which starts
with $W$
and ends with the accept configuration.

Assume that $M$ is an $S$-machine,
and there is (only one) accept configuration.
Then the $a$-length $|W|_a$ of an input configuration $W$ is the number of tape letters in $W$.
If the configuration $W$ is accepted, denote by $T(W)$ the minimal time of computations accepting it.
Then the \label{timef} time function $T(n)=T_M(n)$ is defined as $\max\{T(W)\}$ over all accepted input configurations $W$ with $|W|_a\le n$.

The \label{gentf} {\it generalized time function} $T'(n)$ is defined for every S-machine having a unique
accept configuration. The definition is similar to the above definition of time function
but one should consider all accepted configurations $W$, not just input ones. Therefore
$T(n)\le T'(n)$. (Presumably, $n$, in the definition of $T'(n)$, corresponds to the a-length of the accepted configurations.)

Time functions and generalized time functions are taken up to $\Theta$-equivalence.

\subsection{Simplifying the rules of $S$-machines}\label{modif}

We say that two recognizing S-machines are \label{equivM} equivalent if they have the same
language of acceptable words and $\Theta$-equivalent time functions. Next lemma simplifies the rules of S-machine. In particular, one needs Property (1) to define
trapezia (Definition \ref{dftrap}).

\begin{lemma} \label{simp} Every S-machine $\cal S$ is equivalent to an S-machine ${\cal S}'$, where

(1) Every
part $U_i\to V_i$ of every rule $\theta$ has $1$-letter base:
$U_i\equiv v_iq_iu_{i+1}, \quad V_i\equiv v_{i}'q_{i}'u_{i+1}',$ where $q_{i}, q_{i}'$ are state letters in
$Q_{i}$

(2) In every part
$v_{i}q_iu_{i+1}\to v_{i}'q_i'u_{i+1}'$, we have that $||v_{i}||+ ||v'_{i}||\le 1$ and , $||u_{i+1}||+ ||u'_{i+1}||\le 1$.

(3) Moreover, one can construct ${\cal S}'$ so that for every rule, we have $\sum_i (||v_i||+||v'_i||+||u_i||+||u'_i||) \le 1.$
\end{lemma}

\proof (1) Property (1) can be obtained after adding auxiliary state letters and splittings the rules of $\cal S$. Assume, for example, that the part $U_1\to V_1$ has $2$-letter base: $q_1 a q_2 \to q'_1a'q'_2 $. Then we introduce auxiliary state letters $q_j(1)$, $q_j(2)$  ($j=0,\dots,N$) and replace the rule $\theta$ by the product of three rules $\theta_1$, $\theta_2$ and $\theta_3$, where

 ($\theta_1$) For $\theta_1$ and $j>1$, we replace base letters $q$ in $V_j$ by their $q(1)$-copies and obtain the
parts $U_j\to V_j(1)$, while
the part $U_1\to V_1$ is replaced by the two parts
$q_1a\to q_1(1)$, $q_2\to q_2(1)$;

  ($\theta_2$) For $\theta_2$ and $j>1$, we have now the parts $V_j(1)\to V_j(2)$, where $V_j(2)$ is a copy
	of $V_j(1)$ after replacement $q(1)\to q(2)$, while the part
	$U_1\to V_1$ of $\theta$ is replaced by two parts $q_1(1)\to q_1(2)$ with $Y_2(\theta_2)= \emptyset$ and $q_2(1)\to q_2(2)$;
	
	($\theta_3$) For $\theta_3$ and $j>1$, we have $V_j(2)\to V_j$, while the first part splits now
	as $q_1(2) \to q_1'a'$ and $q_2(2)\to q_2'$.
	
	The key feature of the new S-machine $\tilde{\cal S}$ is in the following obvious property.

{\it There is a one-to-one correspondence between computations $w_0\to...\to w_t$ of $\tilde\sss$ (with any base) such that $w_0$, $w_t$ do not have auxiliary $q$-letters
and computations of $\cal S$ connecting the same words. For every history $H$ of such computation of $\cal S$, the corresponding history of computation of $\tilde S$ is obtained from $H$ by replacing every occurrence of the rule $\theta^{\pm 1}$ by the 3-letter word $(\theta_1\theta_2\theta_3)^{\pm 1}$.}

Clearly, by applying this transformation to an S-machine $\cal S$ several times, we obtain an equivalent S-machine satisfying Property (1).

(2) Suppose Property (2) is not satisfied for a part $U_i\to V_i$.
For example, suppose a  rule $\theta$ of an $S$-machine $\cal S$ has the $i$-th part of the form $av_{i}q_iu_{i+1}\to v_{i}'q_i'u_{i+1}',$ where $u_{i+1}, v_i, u_{i+1}', v_i'$ are words in the appropriate parts of the alphabet of $a$-letters, $v_{i}$ is not empty, $a$ is an $a$-letter, $q_i,q_i'$ are $q$-letters (a very similar procedure can be done in all other cases).

We want to replace $\theta$ with two rules with smaller sums of lengths of their parts. For this aim, we create a new $S$-machine $\tilde\sss$ with the same standard base and the same $a$-letters as $\sss$. In order to build $\tilde\sss$, we add one new ({\em auxiliary}) $q$-letter $\tilde q_i$ to each part of the set of $q$-letters, and replace the rule $\theta$ by two rules $\theta'$ and $\theta''$. The first rule $\theta'$ is obtained from $\theta$ by replacing the part $v_{i}aq_iu_{i+1}\to v_{i}'q_i'u_{i+1}'$ by $aq_iu_{i+1}\to \tilde q_iu_{i+1}'$, and all other parts $U_j\to V_j$ by $U_j\to \tilde q_j$ (here $\tilde q_j$ is the auxiliary $q$-letter in the corresponding part of the set of $q$-letters). The second rule $\theta''$ is obtained from $\theta$ by replacing the part $v_{i}aq_iu_{i+1}\to v_{i}'q_i'u_{i+1}'$ by $v_{i}\tilde q_i\to v'_{i}q_i'$, and all other parts $U_j\to V_j$ by $\tilde q_j\to V_j$.

Note that the sum of lengths of words in all parts of $\theta'$ (resp. $\theta''$) in $\tilde\sss$ is smaller than the similar sum for $\theta$. Therefore, applying this transformation to an S-machine $\sss$ several times, we obtain an equivalent S-machine satisfying conditions (1) and (2).

(3) Similarly, one can obtain Property (3).
\endproof

If $Y_{i+1}(\theta)=\emptyset$ for an S-machine with Property (1), then the corresponding
component $U_i\to V_i$ will be denoted \label{tool}$U_i\tool V_i$ and we shall say that the rule $\theta$ \label{locks}{\em locks} the
$Q_{i}Q_{i+1}$-sectors.
In that case we always assume that $U_i, V_i$ do not have tape letters to the right of the state letters, i.e.,
it has the form $v_iq_i\tool v_i'q_i'$. Similarly, these words have no tape letters to the
left of the state letters if the $Q_{i-1}Q_i$-sector is locked by the rule.

\medskip

\begin{rk} The definitions of an admissible word and a rule application given in Subsection \ref{SM} coincide with the the definitions from \cite{SBR} in case of standard base. However computations do not change the base. So to obtain the statement of Lemma \ref{S}, we may use the main property of computations with standard base obtained in \cite{SBR}.
\end{rk}

S-machines resembles  multi-tape Turing machines (or algorithms).
(The main difference is that a Turing machine does not deal with negative letters.) We do not
give an accurate definition of Turing machines here since from now on we will not use them in this paper (see, for example \cite{DK} or \cite{SBR} for the definition).
However, it is important that the S-machine $S(M)$ constructed in \cite{SBR} simulates the work of  a Turing machine $M$ with time function $T(n)$ as follows.
(See \cite{SBR}, Lemma 3.1 and Proposition 4.1, though we use simpler notation below.)

Let $M$ have one input sector $YZ$, and the input configurations have the form $W\equiv yvz\dots$, where
$v$ is a positive word in an alphabet $A$. Then there is an S-machine $S(M)$ with input
configurations of the form $$\sigma(W)=y_1\alpha^nz_1\dots xy_2vz_2\dots y_3\delta^n z_3\dots y_4\omega^n z_4,$$
where $n=||v||$
(so $S(M)$ has four input sectors with tape words $\alpha^n$, $v$, $\delta^n$ and $\omega^n$, resp.); the S-machine $S(M)$ has time function $\Theta(T(n)^3)$, and it accepts
the configuration $\sigma(W)$ if and only if the configuration $W$ is accepted by $M$.
Moreover, the S-machine $S(M)$ can be constructed so
that for every configuration $W\in\cal L$ accepted
by $M$ with time $T(W)$, the S-machine $S(M)$ accepts
this word with time $\Theta(T(W)^3)$.

\begin{rk} The part "Moreover" is not formulated in \cite{SBR}
explicitly, but it follows from Proposition 4.1.3 (b) since every
Turing machine can be easily modified so that the length of
every accepting computation is $\Theta$-equivalent to the space
of this computation.
\end{rk}

In the present article, we will assume that the basic S-machine ${\bf M}_1$ has only one input sector. Therefore
${\bf M}_1$ has to have a few more rules in comparison with $S(M)$. The input configurations of ${\bf M}_1$
have the form $$\bar\sigma(W)=\bar y_1\bar z_1\dots\bar x \bar v\bar y_2 \bar z_2\dots\bar y_3\bar z_3\dots\bar y_4\bar z_4,$$
where $\bar v$ is a word in an alphabet $\bar A$, which copies the alphabet $A$ (so the only input word is $\bar v$). The following rules of ${\bf M}_1$ are added to the rules of $S(M)$.

For every (positive) letter $a\in A$, there is a rule
$$\rho_a: [ \bar y_1\to \bar y_1\alpha, \;\bar y_2 \to \bar a^{-1}\bar y_2 a,\; \bar y_3\to \bar y_3\delta,\;
\bar y_4\to \bar y_4\omega ],$$
where $\bar a$ is a copy of $a$ in the alphabet $\bar a$ and all other sectors are locked by $\rho_a$.
They also are locked by the connecting rule

$$\rho: [\bar y_1 \to y_1,\dots  \bar x\tool x, \bar y_2\to y_2,\dots \bar z_4\to z_4],$$
which switches on the S-machine $S(M)$. If an input word $\sigma(W)$ is accepted by $S(M)$, then $\bar\sigma(W)$
is accepted by ${\bf M}_1$, where $\bar v$ is a copy of $v$. Indeed, if $\bar a$ is the last letter of the word $\bar v$, then the application
of the rule $\rho_a$ moves the state letter $\bar y_2$  leftward and replaces $\bar a$ by $a$ from the right of $\bar y_2$.
After $n$ rules of this type, one can obtain the word $v$ between $y_2$ and $z_2$. These rules
will also insert $\alpha^n$, $\delta^n$ and $\omega^n$ in the sectors
$Y_1Z_1$, $Y_3Z_3$ and $Y_4Z_4$, resp. So it remains to apply the  rule $\rho$ to obtain the
configuration $\sigma(W)$ accepted by $S(M)$. Also for the times of computations, we see that $T_{{\bf M}_1}(\bar\sigma(W))=O(T_{S(M)}(\sigma(W))$.

Conversely, assume that the configuration $\bar\sigma(W)$ is accepted by ${\bf M}_1$. Then the history
of this accepting reduced computation has to be of the form $H\equiv H_0\rho H_1\rho^{-1} H_2\rho\dots H_{2s-1}$,
where $H_0, H_2,\dots$ contain $\rho_a^{\pm 1}$-rules only and $H_1,H_3,\dots$ are histories of $S(M)$.
However $H$ cannot have subwords of the form $\rho^{-1}H_i\rho$. Indeed, if here $H_i$ is a reduced word $\rho_{a_1}^{\pm 1}\dots\rho_{a_t}^{\pm 1}$, then the tape word $u'$ in the $XY_2$ sector at the end of the computation
${\cal C}_i$ with history $H_i$ is obtained from the word $u$ written there in the beginning of ${\cal C}_i$,
after free multiplication from the right by the reduced word $\bar a_1^{\mp 1}\dots \bar a_t^{\mp 1}$. But
both $u$ and $u'$ are empty since the rule $\rho$ locks the $XY_2$-sector. Hence $t=0$ and so $H_i$ is
empty, a contradiction.

Thus, $H\equiv H_0\rho H_1$, and the $XY_2$-sector of the word $\bar\sigma(W)\cdot H_0$ is empty
being locked by $\rho$. Hence the history $H_0$ has to be $\rho_{a_t}\dots \rho_{a_1}$ if $v\equiv a_1\dots a_t$,
and the tape word in the $Y_2Z_2$ sector of $\bar\sigma(W)\cdot H_0$ has to be $v$ while
the tape words in the sectors $Y_1Z_1$, $Y_3Z_3$ and $Y_4Z_4$ become $\alpha^n$, $\delta^n$ and
$\omega^n$ as it follows from the definition of $\rho_a$-rules. Hence $\bar\sigma(W)\cdot H_0\rho\equiv
\sigma(W)$, and so $\sigma(W)$ is accepted my $S(M)$ and $T_{{\bf M}_1}(\bar\sigma(W))>T_{S(M)}(\sigma(W)$.

It follows that the S-machine ${\bf M}_1$ with one input sector enjoys the properties of $S(M)$ from \cite{SBR}:

\begin{lemma} \label{S}  Let ${\bf M}_0$ be a non-deterministic Turing machine
accepting the language $\cal L$ with a time function $T(n)$.
Then there is an S-machine ${\bf M}_1$ with a single input sector accepting the
language $\cal L$ with time function $\Theta$-equivalent to $T(n)^3$.

Moreover, the S-machine ${\bf M}_1$ can be constructed so
that for every word $W\in\cal L$ accepted
by ${\bf M}_0$ with time $T(W)$, the S-machine ${\bf M}_1$ accepts
this word with time $\Theta(T(W)^3)$.
\end{lemma}

\begin{rk} \label{fn} Later we will assume that the Turing machine ${\bf M}_0$ recognizes the values of some suitable function according to Definition \ref{sui}.
\end{rk}

\subsection{Some elementary properties of $S$-machines}

\label{gpsm}
The base of an admissible word is not always a reduced word. However the following is an immediate corollary of the definition of admissible word.

\begin{lemma}\label{qqiv}(\cite{O12}, Lemma 3.4)
If the $i$-th component of the rule $\theta$ has the form $v_iq_i\tool v'_iq_i',$
i.e. $Y_{i+1}(\theta)=\emptyset$, then the
base of any $\theta$-admissible word cannot have subwords $Q_iQ_i\iv$ or $Q_{i+1}^{-1}Q_{i+1}.$
\end{lemma} $\Box$

In this paper we are often using copies of words.
If $A$ is an alphabet and $W$ is a word involving no letters from $A^{\pm 1}$, then to obtain a \label{copyw}{\em copy}
of $W$ in the alphabet $A$ we substitute letters from $A$ for letters in $W$ so that different letters from $A$
substitute for different letters. Note that if $U'$ and $V'$ are copies of $U$ and $V$ respectively corresponding to
the same substitution, and $U'\equiv V'$, then $U\equiv V.$

\begin{lemma}\label{gen} Suppose that the base of an admissible word $W$ is $Q_{i}Q_{i+1}$. Suppose that each rule of
a reduced computation starting with $W\equiv q_iuq_{i+1}$ and ending with $W'\equiv q_i'u'q_{i+1}'$ multiplies the
$Q_iQ_{i+1}$-sector by a letter on the left (resp. right).
And suppose that different rules multiply that sector by
different letters.
Then

(a) the history of computation is a copy
of the reduced form of the word $u'u\iv$ read from right to left
(resp. of the word $u\iv u'$ read from left to right). In particular,
 if $u\equiv u'$, then the computation is empty;

(b) the length of the history $H$ of the computation does not exceed $||u||+||u'||$;

(c) for every admissible word  $q_i''u''q_{i+1}''$ of the computation, we have
$||u''||\le \max (||u||, ||u'||)$.
\end{lemma}

\proof Part (a) is obvious. To prove part (b), we choose a word $W_i$ of the computation with
shortest tape word $u_i$. This factorizes the history as $H\equiv H_1H_2$, where $H_2$ is the history
of the subcomputation $W_i\to W_{i+1}\to\dots \to W'$. It follows that $||W_{i+1}||= ||W_{i}||+1$.
The next rule increases the length of admissible word again since the computation is reduced
and different rules multiply the sector by different letters, i.e.  $||W_{i+2}||= ||W_{i+1}||+1$.
By induction, we have $||u'||\ge ||u_i||+||H_2||\ge||H_2||.$ To obtain the inequality
$||u||\ge ||H_1||$, we consider the inverse computation $W'\to\dots\to W$. Hence
$||H||=||H_1||+||H_2||\le||u||+||u'||$.

The same argument proves Statement (c) since the length of $u''$
is either between $||u_i||$ and $||u'||$ or between $||u_i||$ and $||u||$.
\endproof

\begin{lemma}\label{gen1} Suppose the base of an admissible word $W$ is $Q_{i}Q_{i+1}$. Assume that each rule of a reduced computation starting with $W\equiv q_iuq_{i+1}$ and ending with $W'\equiv q_i'u'q_{i+1}'$ multiplies the
$Q_iQ_{i+1}$-sector by a letter on the left and by a letter from the right.
Suppose different rules multiply that sector by
different letters and the left and right
letters are taken from disjoint alphabets.
Then

(a) for every intermediate admissible word $W_j$ of the computation, we have $||W_j||\le \max(||W||, ||W'||)$

(b) the length of the history $H$ of the computation does not exceed $\frac 12 (||u||+||u'||)$.
\end{lemma}

\proof (a) If we choose the word $W_i$ of minimal length, then after multiplications
of the form $u_i\to u_{i+1}=au_{i}b$ we have no cancellation from the left or from the right.
If we have the former option, then we will have no cancellation from the left after the
transition $u_{i+1}\to u_{i+2}$, and therefore $||u_{i+1}||\le ||u_{i+2}||\le\dots\le||u'||$.
Hence $||u_j||\le ||u'||$ if $j\ge i$. Analogously, $||u_j||\le ||u||$ if $j\le i$.

(b) The word $u'$ results from $u$ after multiplication from the left and from the right
by reduced words of length $||H||$: i.e., $u'$ is freely equal to $AuB$, where $||A||=||B||=||H||$.
There can be cancellations in the products $Au$ and $uB$ but afterwards there are no cancellations
since the words $A$ and $B$ are written in disjoint alphabets. Hence the reduced length of $u'$
is at least $||A||+||B||-||u||=2||H||-||u||$, whence $2||H|| \le ||u||+||u'||$, as required.
\endproof

The following lemma is proved in \cite{O12} (Lemma 3.7).

\begin{lemma} \label{gen2}
Suppose the base of an admissible word $W$ is $Q_{i}Q_{i}\iv$
(resp., $Q_i\iv Q_i$). Suppose each rule $\theta$ of a reduced
computation $\cal C$ starting with $W\equiv q_iuq_i\iv$ (resp., $q_i\iv uq_i$), where
$u\ne 1$, and ending with $W'\equiv q_i'u'(q_i')\iv$ (resp., $W'\equiv (q_i')\iv
u'q_i')$
has a part $q_i\to a_\theta
q_i'b_\theta,$ where $b_{\theta}$ (resp., $a_{\theta}$) is a letter, and  for different $\theta$-s the $b_\theta$-s
(resp., $a_\theta$-s) are different. Then the history of the computation
has the form $H_1H_2^kH_3,$ where $k\ge0$, $||H_2||\le \min(||u||, ||u'||),$ $||H_1||\le ||u||/2,$
and $||H_3||\le ||u'||/2.$
\end{lemma} $\Box$

\begin{lemma} \label{gen3} Under the assumptions
of Lemma \ref{gen2}, we have $|W_i|_a\le\max(||u||,||u'||)$ for every admissible word $W_i$ of the computation $\cal C$.
\end{lemma}
\proof It suffices to repeat the argument from the proof of Lemma \ref{gen} (c). \endproof

\section{Auxiliary S-machines and constructions}\label{am}

\subsection{Primitive S-machines} \label{pm}

Here we define a very simple S-machine $\bf Pr$, which has neither input nor accept
configurations. As a part of other S-machines, it will be used to read the  tape words and to recognize a computation by its history and also to check the order of state letters in the bases of computations.

The standard base of $\bf Pr$ has three letters $Q^1PQ^2$, where $Q^1=\{q^1\}$,
$P=\{p^1,p^2\}$ and $Q^2=\{q^2\}$. The alphabet $Y$ is $Y^1\sqcup Y^2$, where $Y^2$
is a copy of $Y^1$.
The positive rules of $\bf Pr$ are defined
as follows.

\begin{itemize}

\item $\zeta^1(a)=[q^1\to q^1, p^1\to a\iv p^1a', q^2\to q^2]$,

where $a$ is a positive letter from $Y^1$ and $a'$ is
its copy from $Y^2$.
\me

{\em Comment.} The state letter $p^1$ moves left replacing letters $a$ from $Y^1$ by their copies $a'$
from $Y^2$.

\me

\item $\zeta^{12}=[q^1\tool q^1, p^1\to p^2, q^2\to q^2]$.

\me

{\em Comment.} When $p^1$ meets $q^1$, it  turns into $p^2$.

\me

\item $\zeta^2(a) =[q^1\to q^1, p^2\to ap^2(a')^{-1}, q^2\to q^2]$

{\em Comment.} The state letter $p^2$ moves right towards $q^2$ replacing letters $a'$ from $Y^2$ by their copies $a$
from $Y^1$.

\begin{lemma}\label{prim} Let ${\cal C}:\;W_0\to\dots \to W_t$ be a reduced computation of the S-machine $\bf Pr$ with the standard
base and  with $t\ge 1$. Then

(1) if $|W_i|_a>|W_{i-1}|_a$ for some $i=1,\dots,t-1$, then $|W_{i+1}|_a>|W_i|_a$;

(2) $|W_i|_a\le\max(|W_0|_a, |W_t|_a)$ for every $i=0,1,\dots t$;

(3) if $W_0\equiv q^1up^{1}q^2$ and $W_t\equiv q^1vp^{2}q^2$ for some words $u,v$, then $u\equiv v$,
$|W_i|_a=|W_0|_a$ for every $i=0,\dots,t$, $t=2k+1$, where $k=|W_0|_a$, and $p^1$ (resp., $p^2$)
meets $q^1$ in $W_k$ (in $W_{k+1}$) and the sector $Q^1P$ is empty in $W_k$ and in $W_{k+1}$;
moreover, the history $H$ of $\cal C$ is uniquely determined by $W_0$ (by $W_t$), provided $W_0$ and $W_t$ have the form  $q^1up^{1}q^2$ and $q^1vp^{2}q^2$; vice versa, the
history $H$ uniquely determines words $u$ and $v$ under this
assumption.

(4) it is not possible that $W_0\equiv q^1up^{1}q^2$ and $W_t\equiv q^1vp^{1}q^2$ for some $u,v$,
and it is not possible that  $W_0\equiv q^1up^{2}q^2$ and $W_t\equiv q^1vp^{2}q^2$;

(5) if $W_0\equiv q^1up^{1}q^2$ or $W_0\equiv q^1p^{1}uq^2$, or $W_0\equiv q^1up^{2}q^2$, or $W_0\equiv q^1p^{2}uq^2$ for some word $u$, then $|W_i|_a\ge |W_0|_a$ for every $i=0,\dots,t$.
\end{lemma}

\proof Note that each of the rules $(\zeta^j)^{\pm 1}(a)$, ($j=1,2$) either moves the state letter left or moves
it right, or deletes one letter from left and one letter from right, or insert letters from both sides. In the later
case, the next rule of a computation must be again $\zeta(j)^{\pm 1}(b)$ for some $b$, and if the computation is
reduced, it again must increase the length of the configuration by two. Therefore Statement (1) is true and (2)
is also true since one can choose a shortest $W_j$ and consider the subcomputation $W_j\to\dots\to  W_t$ and the inverse
subcomputation $W_j\to\dots\to W_0$.

Since $\zeta^{12}$ locks $Q^1P$-sector, the $p$-letter must reach $q_1$ moving always left to change $p^1$ by $p^2$,
 and so $W_k\equiv q^1p^1\dots$. The next rule of the form $\zeta^1(a)^{\pm 1}$ could increase the length of the configuration, which would imply that
all consecutive  rules have to have the same type and $p^1$ would never been replaced by $p^2$, a contradiction.
Hence the next rule is $\zeta^{12}$, and arguing in this way, one uniquely reconstructs the whole computation
in case (3) for given $W_0$ or $W_t$, and vice versa, the history $H$ determines both $u$ and $v$. Propery (4) holds for same reasons.

Note that no rule of $\bf Pr$ changes the projection of a word onto the free group with basis $Y_1$ if the state letter are mapped to $1$ and the letters from $Y_2$ are
maped to their copies from $Y_1$. Since the word $u$ is mapped
to itself, we have $|q^1up^{1}q^2|_a =||u||\le |W_i|_a$.
The other cases of (5) are similar.
\endproof

\begin{rk} \label{proj} Similar
tricks will later be referred to as {\it projection arguments}.
\end{rk}

\begin{lemma}\label{ewe} If $W_0\to\dots\to W_t$ is a reduced
computation of $\bf Pr$ with base $Q^1PP^{-1}(Q^1)^{-1}$ or $(Q^2)^{-1}P^{-1}PQ^2$
and $W_0\equiv q^ip^iu(p^i)^{-1}(q^i)^{-1}$ ($i=1,2$) or\\ $W_0\equiv(q^i)^{-1}(p^i)^{-1}vp^{i}q^i$ ($i=1,2$)  for some words $u,v$, then $|W_j|_a\ge |W_0|_a$ for every $j=0,\dots,t$.
\end{lemma}

\proof The statement follows from the projection argument \ref{proj}.
\endproof

\begin{rk} \label{right} Also we will use the right analog ${\bf Pr}^*$ of $\bf Pr$:  now $p^1$ should move right, meet $q^2$ locking
$PQ^2$-sector and turning into $p^2$, and move back towards
$Q^1$.
\end{rk}

\begin{rk} \label{ps} Assume that a standard base has two (or more) subwords of the form $Q^1PQ^2$, for example
$Q^1PQ^2P'Q^3$. Then one can define  {\it parallel} or {\it sequential} \label{composm} composition of two primitive S-machines.

For the
parallel composition, the same rule changes both subwords with bases $Q^1PQ^2$ and $Q^2P'Q^3$. One assumes that the tape alphabet of the sector $Q^2P'$ (of $P'Q^3$) is a copy of the tape
alphabet of the sector $Q^1P$ (of $PQ^2$, resp.) and the rules
of $\bf Pr$ change simultaneously the subwords with bases
$Q^1PQ^2$ and $Q^2P'Q^3$
 (e.g. simultaneously
moves left both $p^1$ and $(p')^1$, and so on). Every rule (e.g $\zeta^1(a)$) is applicable to a word iff it is
applicable to both these subwords.

In case of sequential work, we have a primitive
S-machine working with one of these two bases, say with $Q^1PQ^2$, while the sector $P'Q^3$ is locked
with the state letter $(p')^1$.
The second primitive S-machine can compute with base $Q^2P'Q^3$ when the sector $PQ^2$ is locked (and so the first S-machine stays idle) with state  $p^2$. For this goal, one needs a connecting
rule $\zeta^{21}$; it locks the sector $PQ^2$ and changes the
state $P$-letters for new ones to switch off the first primitive S-machine and to switch on the second one.

It is clear that in the same way one can define a more complex
compositions $\bf P$ of primitive S-machines, with several stages
of parallel and sequential work. We will consider compositions $\bf P$, where every sector can be changed at one stage only.
\end{rk}

\begin{lemma}\label{Hprim} Let us have a composition $\bf P$ of primitive S-machines with parallel or/and sequential work, and ${\cal C}:
W_0\to\dots\to W_t$ be a reduced computation of $\bf P$ with standard
base. Then

(a) $|W_j|_a\le\max (|W_0|_a,|W_t|_a)$ for every configuration $W_j$ of $\cal C$; moreover, $|W_0|_a\le \dots\le |W_t|_a$ if every $P$-letter
neighbors some $Q$-letter in the word $W_0$;

(b) $ t\le ||W_0||+||W_t||-4$, moreover, $t\le 2||W_t||-4$ if every $P$-letter neighbors some $Q$-letter in the word $W_0$.
\end{lemma}

\proof (a) Let $W_r$ be a shortest word of the computation.
Then either $|W_r|_a=|W_{r+1}|_a=\dots= |W_t|_a$, or
$|W_r|_a=|W_{r+1}|_a=\dots= |W_s|<|W_{s+1}|_a$ for
some $s$. It follows that the number of sectors increasing
their lengths by two at the transition $W_s\to W_{s+1}$ is greater than the number of the sectors decreasing the lengths by
$2$. Now it follows from Lemma \ref{prim} (1) that the same
primitive S-machine will continue increasing the lengths of
the whole configurations, i.e., $|W_{s+1}|_a<|W_{s+2}|_a<\dots $.
So for every $j\ge r$, we have $|W_j|_a\le |W_t|_a$.
Similarly, we have $|W_r|_a\le |W_0|_a$ for $j\le r$.
Under the additional assumption about $P$-letters, $W_0$ is the shortest configuration by the projection argument.

(b) If the rules of $\bf P$ do not change the lengths of configurations, then every control letter runs back and forward
only one time, and the inequality follows. (One takes into account that the base has length at least $3$.) If $||W_r||<||W_{r+1}||$
for some $r$, then every next transition keeps increasing the
length by Lemma \ref{prim} (1), and so the inequality holds as well.
\endproof

\end{itemize}

\subsection{S-machine with historical sectors} \label{hs}

To control the space of computations, we endow a given S-machine
with historical sectors. Let us assume that an $S$-machine ${\bf M}_1$ satisfies the
conditions of Lemma \ref{simp} and has hardware $(Q,Y)$, where
$Q=\sqcup_{i=0}^m Q_i$, and the set  of rules $\Theta$. The new
S-machine ${\bf M}_2$ has hardware
$$Q_{0,r}\sqcup Q_{1,\ell}\sqcup Q_{1,r}\sqcup Q_{2,\ell}\sqcup Q_{2,r}\sqcup\dots \sqcup Q_{m,\ell},\;\; Y_h=  Y_1
\sqcup X_1\sqcup Y_2\sqcup\dots \sqcup X_{m-1}\sqcup Y_m,$$
where $Q_{i,\ell}$ and $Q_{i,r}$ are (left and right) copies of $Q_i$,  $X_i$ is a disjoint union of two copies of $\Theta^+$, namely $X_{i,\ell}$ and $X_{i,r}$. (There is neither $Q_{0,\ell}$, nor $Q_{m,r}$, nor $X_0$, nor $X_m$.)
The positive rules of ${\bf M}_2$ are in one-to-one correspondence with the positive rules of ${\bf M}_1$.

If $\theta=[U_0\to V_0,...,U_m\to V_m]$ is a positive rule of ${\bf M}_1$ with parts $U_i\to V_i$ of the form
$v_iq_iu_{i+1}\to v_{i}'q_{i}'u_{i+1}'$, then the corresponding two parts of the rule $\theta_h$ are
$U_{i,\ell}\to V_{i,\ell}$ and $U_{i,r}\to V_{i,r}$, with
$$U_{i,\ell}\equiv v_iq_{i,\ell}a_{\theta,i},  V_{i,\ell}\equiv v_i' q'_{i,\ell}\;\; and \;\;
U_{i,r}\equiv q_{i,r}u_{i+1},  V_{i,r}\equiv b_{\theta,i} q'_{i,r}u_{i+1}',$$ where
$a_{\theta,i}$ (resp., $b_{\theta,i}$) is the copy of $\theta$ in $X_{i,\ell}$ (in $X_{i,r}$).
We also claim that a sector $Q_{i,r}Q_{i+1,\ell}$ is locked by $\theta_h$ if and only
if the sector $Q_{i}Q_{i+1}$ is locked by $\theta$ ($i=1,\dots,m-1$).

\medskip

{\it Comment.} Every computation of the S-machine ${\bf M}_2$ with history $H$ coincides with the computation of ${\bf M}_1$ if
one observes it only in {\it working} sectors $Q_{i,r}Q_{i+1,l}$. In the standard base, the \label{wsect} working
sectors of ${\bf M}_2$ alternate with \label{histsec} {\it historical} sectors $Q_{i,\ell}Q_{i,r}$. Every positive rule
$\theta_h$ multiplies the content of the historical sector $Q_{i,\ell}Q_{i,r}$ by the corresponding letter $b_{\theta, i}$ from the right
and by letter $a_{\theta,i}^{-1}$ from the left.

\medskip

\begin{rk} \label{as} The state letters of
the S-machine ${\bf M}_1$ split when passing to ${\bf M}_2$.
There is a rule $\theta_h$ corresponding to the start (to the accept) rule $\theta$ of ${\bf M}_1$. By definition, {\it the set of letters $Y_h(\theta_h)$
has no letters from the right alphabets $X_{i,r}$ (from the left alphabets
$X_{\ell,i}$) if $\theta$ is the start (resp., the stop) rule of ${\bf M}_1$}.

However we do not define
input/stop configurations of ${\bf M}_2$ since the historical sectors are never locked. By definition, every  $Q_{i-1,r}Q_{i,l}$ is the working sector of ${\bf M}_2$. The input sector
of ${\bf M}_2$ is the working sector corresponding to the input sector
of ${\bf M}_1$.
\end{rk}

\begin{rk} It follows from the definition of ${\bf M}_2$ that only
Properties (1) and (2) of Lemma \ref{simp} hold for ${\bf M}_2$, but not Property (3).
\end{rk}

The sectors of the form $Q_{i,\ell}Q_{i,\ell}^{-1}$ and $Q_{i,r}^{-1}Q_{i,r}$ (in a non-standard base)
are also called {\it historical}. Historical sectors help to give a linear estimate of the
space of every computation $W_0\to\dots\to W_t$ in terms of $||W_0||$ and $||W_t||$.

\begin{lemma}\label{w} Let $W_0\to\dots\to W_t$ be a reduced computation of ${\bf M}_2$ with base $Q_{i,\ell}Q_{i,r}$ and history $H$. Assume that the $a$-letters of $W_0$ belong to one
of the alphabets $X_{i,\ell}$, $X_{i,r}$.
Then $||H||\le |W_t|_a$ and $|W_0|_a\le |W_t|_a$.
\end{lemma}
\proof Let $W_0 \equiv qv_0q'$ and assume that $v_0$ has no letters from $X_{i,r}$. Then $W_t\equiv q''v_tq'''$ with $v_t=uv_0u'$, where $u$ is a copy of $H^{-1}$ in the alphabet $X_{i,\ell}$ and $u'$ is a copy of $H$
in $X_{i,r}$. So no letter of $u'$ is cancelled in the
product $uv_0u'$, Therefore $|W_t|_a\ge ||u'||=||H||$ and
$|W_t|_a\ge |W_0|_a$
\endproof

\begin{lemma} \label{three} If the base of an admissible word of the S-machine ${\bf M}_2$ has length at least 3,
then it contains a historical sector.
\end{lemma}

\proof The base contains a subword of the form $Q'Q''Q'''$ with three letters from $Q^{\pm 1}$. It follows from the definition of admissible word that either $Q'Q''$-sector or $Q''Q'''$-sector is historical
since every non-historical sector of the S-machine ${\bf M}_2$ has to have neighbor  historical sectors.
\endproof

\begin{lemma} \label{wi} Let a reduced computation $W_0\to W_1\to\dots \to W_t$ of the S-machine ${\bf M}_2$ have $2$-letter base
and the history of the form  $H \equiv  H_1H_2^k H_3$ ($k\ge 0$). Then for each tape word $w_i$
between two state letters of $W_i$ ($i=0,1,\dots,t$), we have inequality \\$||w_i||\le ||w_0|+ ||w_t||+2h_1+ 3h_2+2h_3$, where $h_j=||H_j||$ ($j=1,2,3$).
\end{lemma}

\proof By Lemma \ref{simp} (2) and the definition of ${\bf M}_2$, we have $|\;||w_i||-||w_{i-1}||\;|\le 2$
for every $i=1,\dots, t$. Therefore for  $i\le h_1$, we have $||w_i||\le||w_0||+2h_1$.
Similarly, $||w_i||\le||w_t||+2h_3$ for $i\ge t-h_3$. It remains to assume that $h_1<i< t-h_3$.

Denote  by $W_j$ the words $w_i$ with $i=h_1+j h_2$, $j=0,1,\dots,k$.
If $W_1=W_0\cdot  H_2 = uW_0v$ for some words $u$ and $v$ depending on $H_2$ and on the sector, then $W_2= uW_1v= u^2 W_0v^2$ in free group, since the histories of the computations $W_0\to\dots\to W_1$ and
$W_1\to\dots\to W_2$ are both equal to $H_2$. Hence $W_j=u^jW_0v^j$, where both $u$ and $v$ have length
at most $h_2$ by Lemma \ref{simp} (2) and the definition of ${\bf M}_2$.

By Lemma 8.1 from \cite{OS04},
the length of an arbitrary word $W_j$ is not greater
than $||u||+||v||+||W_0||+||W_k||$ provided $0\le j\le k$.

If $|i -jh_2|\le h_2/2 $ for some $j$, then $|\;||w_i||-||W_j||\;|\le h_2$, and therefore
$||w_i||\le ||u||+||v||+||W_0||+||W_k||+h_2$. Since $||W_0||\le||w_0||+2h_1$ and
$||W_k||\le ||w_t||+2h_3$, we obtain $$||w_i||\le ||u||+||v||+||w_0||+||w_t||+2h_1+2h_3+h_2\le ||w_0|+ ||w_t||+2h_1+2h_3+3h_2$$ for every $i$, as required.
\endproof

\begin{lemma} \label{9} For any reduced computation $W_0\to\dots\to W_t$ of ${\bf M}_2$
with base of length at least $3$, we have $|W_i|_a\le 9( |W_0|_a+|W_t|_a) $ ($0\le i\le t$).
\end{lemma}

\proof Let $Q_{i_1}^{\pm 1} \dots Q_{i_m}^{\pm 1} $ be the base of the computation. There are
computations with the same history $H$ and bases $Q_{i_1}^{\pm 1} \dots Q_{i_{m_1}}^{\pm 1} $,
$Q_{i_{m_1}}^{\pm 1} \dots Q_{i_{m_2}}^{\pm 1} $,..., $Q_{i_{m_{s-1}}}^{\pm 1} \dots Q_{i_{m_s}}^{\pm 1},$
where each base has length $3$ or $4$. Hence it suffices to prove the lemma for
any computation with base of the form $Q'Q''Q'''$ or $Q'Q''Q'''Q^{iv}$. By Lemma \ref{three}, every such computation contains a historical sector, say $Q''Q'''$.
Consider two cases.

{\bf 1.} The historical sector has the form $Q_{i,\ell}Q_{i,r}$. By Lemma \ref{gen1}, we have $||H||\le \frac 12 (|W_0|_a+|W_t|_a)$. It follows from Lemma \ref{simp} (2) that $|\;|W_{i+1}|_a-|W_i|_a\;|\le 6$ for every
neighbor admissible words. Therefore $$|W_i|_a\le \max(|W_0|_a, |W_t|_a)+3||H|| \le $$ $$\max(|W_0|_a, |W_t|_a)+
\frac 32 (|W_0|_a +|W_t|_a)\le \frac 52 (|W_0|_a +|W_t|_a).$$

{\bf 2.} The historical sector has form $Q_{i,\ell}Q_{i,\ell}^{-1}$ or $Q_{i,r}^{-1}Q_{i,r}$.
Then one can apply Lemma \ref{gen2} to the sector $Q''Q'''$ and obtain the factorization $H\equiv H_1H_2^kH_3,$ with $k\ge0$, $h_2\le \min(||u||, ||u'||),$ $h_1\le ||u||/2,$
and $h_3\le ||u'||/2,$ where $u$ and $u'$ are the $a$-words of $W_0$ and $W_t$, respectively in
the historical sector and $h_j=||H_j||$ for $j=1,2,3$. Since every $W_i$ has at most three sectors,
applying Lemma \ref{wi} to each of them, we obtain:
$$|W_i|_a\le |W_0|_a+|W_t|_a+3(2h_1+ 3h_2+2h_3) \le $$ $$|W_0|_a+ |W_t|_a +3|W_0|_a+ 9\min(|W_0|_a, |W_t|_a)+3|W_t|\le 9(|W_0|_a +|W_t|_a).$$
\endproof

\subsection{Division S-machine}\label{mcd}

Here we start with an S-machine ${\bf D}_1$. This S-machine has two input
sectors with words $a^k$ and $b^{\ell}$ and checks whether
$2k$ divides $\ell$ or not.

The standard base of ${\bf D}_1$ is $S(1)S(2)T(1)T(2)$. The first input
sector $S(1)S(2)$ has one-letter alphabet $\{a\}$, the second input
sector $T(1)T(2)$ has alphabet $\{b\}$. Also we have one-letter
alphabet $\{a'\}$ for the sector $S(2)T(1)$.
We omit some parts of the rules in the  list below if these parts do not change
configurations (e.g., $s\to s$ for $s\in S(1)$ is a part of $\tau_1$).

\begin{itemize}
\item $\tau_1$:  $[s_1\to a^{-1}s_1a'], \; [t_1\to t_1b^{-1}]$,\; $s_1\in S(2), t_1\in T(1)$

{\it Comment.} The state letter $s_1$ moves left changing letter
$a$ by  its copy $a'$, while $t_1$ erases one letter $b$.

\item $\tau_{12}$: $[s\tool s, s_1\to s_2]$\;\;  $s\in S(1), s_2\in S(2)$

{\it Comment.} The rule $\tau_{12}$ locks the sector $S(1)S(2)$
and replaces the state letter $s_1$ by $s_2$.

\item $\tau_2$: $[s_2\to as_2(a')^{-1}]\;\;[t_1\to t_1b^{-1}]$

{\it Comment.} $s_2$ moves right toward $T(1)$ replacing $a'$ by $a$, while $t_1$ erases the letter $b$.

\item $\tau_{21}$: $[s_2\tool s_1]$

{\it Comment.} This rule locks the sector $S(2)T(1)$ and replaces
$s_2$ by $s_1$.

\item $\tau_{3}$: $[s_1\tool t_1], \;
[t_1\tool t_2]$

{\it Comment.} The state $t_2$ can appear if both sectors $S(2)T(1)$ and $T(1)T(2)$ are empty.

\end{itemize}

We call a transition $W\to W'$ given by the rule $\tau_{1}^{\pm 1}$ or $\tau_2^{\pm 1}$ {\it wrong} if it increases the lengths of both sectors $S(1)S(2)$ and $S(2)T(1)$.

\begin{lemma} \label{wrong} Let $W_0\to\dots\to W_t$ be a reduced computation of
the S-machine ${\bf D}_1$ with standard base and the first transition $W_0\to W_1$ is wrong. Then all subsequent transitions are wrong too.
\end{lemma}

\proof Let the first rule $\theta$ be $\tau_1^{\pm 1}$. Then the first transition
(restricted to the sectors $S(1)S(2)$ and $S(2)T(1)$) has the form $sus_1vt_1\to sua^{\mp 1}s_1(a')^{\pm 1}vt_1$, where the words $u,v, ua^{\mp 1}, (a')^{\pm 1}v$ are reduced.
It follows that the only possible rule for the next transition is
$\theta$ again, $W_2$ contains reduced form $sua^{\mp 2}s_1(a')^{\pm 2}vt$, and so on. The case $\theta = \tau_2^{\pm 1}$ is similar.
\endproof

\begin{lemma} \label{div}(1) Suppose we have a reduced computation $W_0\to\dots\to W_r$ of ${\bf D}_1$, where $W_0\equiv sa^ks_1t_1b^{\ell}t'$, where $t'\in T(2)$, and
$W_r\equiv \dots t_it'$ ($i=1,2$, i.e. the sector $T(1)T(2)$ is empty). Then the exponent $\ell$
is divisible by $2k$.

(2) Conversely, if $2k$ divides $\ell$, then there is a computation
$sa^ks_1t_1b^{\ell}t'\to\dots\to sa^kt_2t'$
of length $|\ell| +|\ell/k|+1$ for $k\ne 0$ and of length $1$ for $k=0$.

\end{lemma}

\proof (1) If $k=0$, then any transition given by $\tau_1^{\pm 1}$ or $\tau_2^{\pm 1}$ would be wrong, and by Lemma \ref{wrong}. one can never
obtain $W_r$, a contradiction. So we have no such transition in the computation. But other rules do not change the exponent $\ell$. Hence $\ell=0$. Thus, we assume further that $k\ne 0$.

By Lemma \ref{wrong}, there are no wrong transitions in
the computation. Therefore if the first rule is $\theta=\tau_1^{\pm 1}$,
it has to move $s_1$ left. Moreover, we have $\theta^k$ as the prefix
of the history, and $W_k\equiv ss_1(a')^kt_1b^{\ell-k}t'$.

The next transition is not wrong, and the only possible next rule
is $\tau_{12}$. Now $s_2$ has to move right, we have $|k|$ such transitions, and obtain $W_{2k+1}\equiv sa^ks_2t_1b^{\ell-2k}t'$. Since the next rule is not wrong,
it has to be $\tau_{21}$. If $\ell-2k\ne 0$, then
the $T(1)T(2)$ sector is not locked, and the rule $\tau_3$
does not apply. Thus the next transition $W_{2k+2}\to W_{2k+3}$
is given by $\theta$ again, and one should repeat the cycle
obtaining $W_{4k+3}\equiv sa^ks_2t_1b^{\ell-4k}t'$, and so on; the rule $\tau_3$ will never apply if $\ell-(2k)m\ne 0$
for every $m\ge 0$.

There is another possibility for the first rule: $\theta= \tau_{21}^{-1}$. Since the second transition cannot be wrong by
Lemma \ref{wrong}, it is given by $\tau_2^{\mp 1}$. Then we
will obtain cycles as above, but having reverse direction. The
rule $\theta_3$ will never apply if $\ell+2k\ne 0$, $\ell+4k\ne 0$,... Thus, Statement (1) is proved by contradiction.

(2) If $\ell=2km$ or $\ell= -2km$ for some $m\ge 0$, the required
computations can be immediately constructed according to the samples from part (1) of the proof.

\endproof

Now we want to modify the S-machine ${\bf D}_1$ as follows. To define the S-machine ${\bf D}_2$ we add one more part $T(3)$ to the standard base.
The sector $T(2)T(3)$ serves to count the number of cycles
of the S-machine ${\bf D}_1$. So it is empty for the start configuration
of ${\bf D}_2$, and the rule $\tau_{21}$ of ${\bf D}_1$ extended to the rule
of ${\bf D}_2$ has one more part: $[t'\to t'c]$, where $t'\in T(2)$.
Clearly, Properties (1) and (2) of Lemma \ref{div} hold for ${\bf D}_2$ as well. Moreover, repeating the proof  of Property (2), we see
that one obtains $c^m$ in sector $T(2)T(3)$ with $m=\frac{\ell}{2k}$ when the sector $T(1)T(2)$ becomes empty.

The further modification is needed since we should check the
divisibility by $(2k)^3$, which, in turn is necessitated by Lemma \ref{S}.  The S-machine ${\bf D}_3$ has the same standard base as ${\bf D}_2$ but
it checks divisibility by $2k$ three times, so its rules are subdivided
in three parts.

The rules of the first part are exactly the rules of the S-machine ${\bf D}_2$. The rule $\tau_3$ of ${\bf D}_2$ serves as a connecting rule between
the rules of the first part and the rules of the second part.
The difference between these two parts is that the sectors
$T(1)T(2)$ and $T(2)T(3)$ interchange their roles: a state letter
from $T(2)$ erases letters in the sector $T(2)T(3)$ when the analogs of $\tau_1$ and $\tau_2$ work, and a state letter from $T(1)$ add one letter
to the sector $T(1)T(2)$ when the analog of $\tau_{21}$ is applied.
(We do not introduce notation for all state letters and all rules
since we do not need them.) The rules of the third part are absolutely similar to the rules of ${\bf D}_1$, i.e.
the sector
$T(2)T(3)$
is locked.

\begin{rk} \label{D3} Thus, one can repeat the argument from the proof of Lemma \ref{div} (1,2) three times
to conclude that starting with an input configuration, the S-machine
${\bf D}_3$ can empty all the sectors, except for sector $S(1)S(2)$ if
and only if the exponent $\ell$ is divisible by $(2k)^3$.

\end{rk}

Finally, we add a rule $\tau$ erasing $a$-letters  of the
sector $S(1)S(2)$ with non-trivial part $[s\to sa^{-1}]$ for $s\in S(1)$ locking other sectors, and if all the sectors become empty, one more
rule $\tau_0$ (the stop rule) locks all the sectors and changes all the state letters for the letters of the stop configuration. Let us denote the obtained
S-machine by ${\bf D}_4$.

\begin{lemma} \label{div3}(1) Suppose we have a reduced accepting computation $W_0\to\dots\to W_r$ of ${\bf D}_4$, where $W_0\equiv sa^ks_1t_1b^{\ell}t't''$ is an input configuration.
Then the exponent $\ell$
is divisible by $(2k)^3$.

(2) Conversely, if $(2k)^3$ divides $\ell$, then there is an accepting computation starting with
$W_0\equiv sa^ks_1t_1b^{\ell}t't''$
of length $\Theta(|\ell| +|k|)$.

\end{lemma}
\proof Properties (1), (2)  follows from similar properties of ${\bf D}_3$ mentioned in Remark \ref{D3}.

\endproof

The next modification is obtained by adding historical sectors
to the standard base of ${\bf D}_4$. The approach is similar to
that described in Subsection \ref{hs}.
The standard base of ${\bf D}_5$ is $$S(1)_rS(2)_{\ell}S(2)_rT(1)_{\ell}T(1)_rT(2)_{\ell}T(2)_rT(3)_{\ell},$$
and the rules of ${\bf D}_5$ are the extensions of the rules of ${\bf D}_4$ to
historical sectors as this was defined in Subsection \ref{hs}.
The following lemma is an analog of Lemma \ref{9}.

\begin{lemma} \label{9D} For any reduced computation $W_0\to\dots\to W_r$ of the S-machine ${\bf D}_5$
with base of length at least $3$, we have $|W_i|_a\le 9( |W_0|_a+|W_r|_a) $ ($0\le i\le r$).
\end{lemma}

\subsection{Control state letters}\label{csl}

The work of the main S-machine will be checked by control state
letters running back and forward along the sectors from time to time. The control letters behave as $p$-letters of the primitive
S-machines $\bf Pr$ or ${\bf Pr}^*$.

Suppose $M$ is an $S$-machine with a standard base $Q_0Q_1\dots Q_s$. We denote by $M_c$ the S-machine
with standard base $$P_0 \sqcup Q_0\sqcup R_0\sqcup P_1\sqcup Q_1\sqcup R_1\sqcup\dots \sqcup P_s \sqcup Q_s\sqcup R_s. $$ For every rule $\theta$ of $M$,
its $i$-th part $[v_{i}q_iu_{i+1}\to v_{i}'q_i'u_{i+1}']$ is replaced in $M_c$ with three parts
\begin{equation}\label{3p}
[v_ip^i\tool v'_ip^i], [q_i\tool q'_i], [r^iu_{i+1}\to r^i u_{i+1}']
\end{equation}
 ($i=0,\dots, s$, $p^i\in P_i, r^i\in R_i$). Here we should use one more $\tool$
if there is $\tool$ in the definition of the  component of $M$.

{\it Comment.} Thus, the sectors $P_iQ_i$ and $Q_iR_i$ are always locked, and three state letters $p^i, q^i, r^i$ work together in $M_c$ as the single $q^i$ in $M$. Of course, such a modification is useless for solo work of $M$. But it will be helpful when one constructs
a composition of $M_c$ with other S-machines, because the {\it control
letters} from the parts $P_i$ and $R_i$ will work when $M_c$ stands idle.

\section{The main S-machine.}

\subsection{Definitions of machines ${\bf M}_3 - {\bf M}_6$} \label{M6}

We use the S-machine ${\bf M}_2$ from Section \ref{hs} and auxiliary S-machines to compose
the main machine needed for this paper.

At first we add control state letters to ${\bf M}_2$
and obtain  S-machine ${\bf M}_3$ as it was defined in Section \ref{csl}. Let the standard
base of ${\bf M}_2$ be $Q_0Q_1\dots Q_s$, where sectors $Q_0Q_1$, $Q_2Q_3$,...,$Q_{s-1}Q_s$
are working sectors, $Q_{s-1}Q_s$ is the input sector, and $Q_1Q_2$, $Q_3Q_4,\dots$ are historical sectors. Then the standard base of ${\bf M}_3$ is $$P_0Q_0R_0P_1Q_1R_1\dots P_sQ_sR_s,$$ where
$P_i$ (resp., $R_i$) contains control $p$-letters ($r$-letters), $i=0,\dots, s$.

Since the rules of ${\bf M}_3$ treat every syllable $P_iQ_iR_i$ as
a single base letter, the working and the historical sectors for ${\bf M}_3$ are of the form $R_{i-1}P_i$. In particular,
every historical sector
has the form $R_{i-1}P_i$ with even $i$.

The rules of the next S-machine ${\bf M}_4$ will be partitioned in subsets corresponding to
ten {\it steps} with auxiliary rules $\theta(ij)$ \label{conrule} {\em connecting} $i$-th and $j$-th steps.
The state letters are also disjoint for different steps. Therefore we need ${\cal Q}_0$,
which is the disjoint union of ten subsets, ${\cal P}_0$, which is the
disjoint union of ten subsets, and so on. Thus, the rules of different steps of a computation on ${\bf M}_4$ must be separated by
connecting rules.

We want to combine the S-machines ${\bf M}_3$, the machine $({\bf D}_5)_c$
(i.e. the S-machine  $({\bf D}_5)$ from subsecion \ref{mcd}
endowed with control state letters),
and compositions of primitive S-machines
introduced in Subsection \ref{pm}. We interbreed the
input sector of ${\bf M}_3$ and the first sectors of
$({\bf D}_5)_c$. Namely, the state letters from $Q_{s-1}$ and from  $S(1)_r$
will be included in the part ${\cal Q}_{s-1}$ of the new S-machine ${\bf M}_4$, $Q_s$ and $S(2)_{\ell}$ will be included in ${\cal Q}_{s}$. The reader will see below that at some steps of computations, the
part of base ${\cal Q}_{s-1}{\cal R}_{s-1}{\cal P}_s {\cal Q}_s{\cal R}_s$ works as ${\bf M}_3$ while at other steps it works as ${\bf D}_5$.

The new S-machine ${\bf M}_4$ repeats
the computation of ${\bf M}_3$ many times and ${\bf D}_5$ bounds the number of such cycles. The standard
base of ${\bf M}_4$ is $${\cal P}_0 {\cal Q}_0{\cal R}_0{\cal P}_1 {\cal Q}_1{\cal R}_1\dots {\cal P}_{s-1} {\cal Q}_{s-1}{\cal R}_{s-1}{\cal P}_s {\cal Q}_s{\cal R}_s\times\;\;{\cal P}{\cal S}(2)_r {\cal R} \times $$ $${\cal P}^{1,\ell}{\cal T}(1)_{\ell}{\cal R}^{1,\ell}{\cal P}^{1,r}{\cal T}(1)_r{\cal R}^{1,r}{\cal P}^{2,\ell}{\cal  T}(2)_{\ell}{\cal R}^{2,\ell}{\cal P}^{2,r}{\cal  T}(2)_r{\cal R}^{2,r}{\cal P}^{3,\ell}{\cal T}(3)_{\ell}{\cal R}^{3,\ell}.
$$
(Starting with ${\cal P}_{s-1}$, this base looks looks like
the base of ${\bf D}_5$ equipped with control $\cal P$- and $\cal R$-parts.)
The historical
sectors of the form ${\cal  R}_{i-1}{\cal P}_i$ with even $i$ are called {\it big
historical} sectors while ${\cal R}^{1,\ell}{\cal P}^{1,r}$ and ${\cal R}^{2,\ell}{\cal P}^{2,r}$ are {\it small
historical} sectors. The sector ${\cal R}_s{\cal P}$ is also small historical one. It corresponds to the sector $S(2)_{\ell}S(2)_r$ of ${\bf D}_5$.

The rules of ${\bf M}_4$ will be partitioned in subsets $\Theta_{i^-}$ and $\Theta_i$ ($i=1,\dots,5$) corresponding to
ten \label{step} {\it Steps}.
We will not  list all state letters here since it would be complicated and not too helpful. It suffices to define
the work of ${\bf M}_4$ at different steps as a composition
of the S-machines defined in Section \ref{am}.

\begin{figure}
\begin{center}
\includegraphics[width=1.0\textwidth]{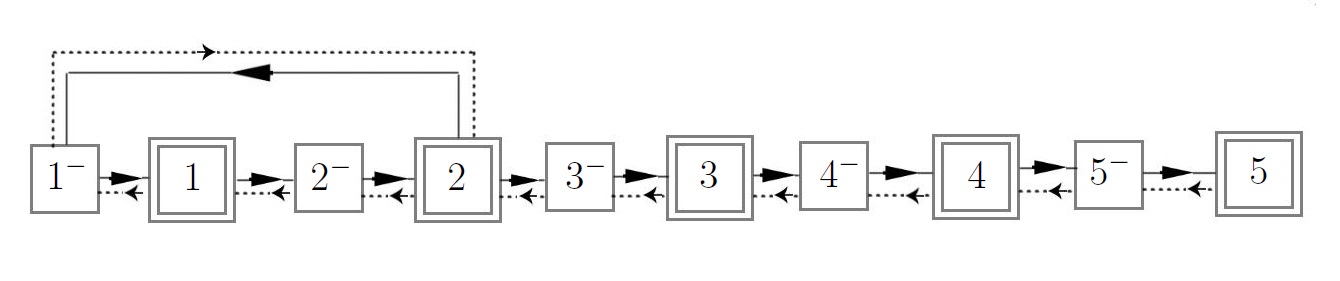}
\end{center}
\caption{Graph of steps of S-machine ${\bf M}_4$.\label{Pic4}}
\end{figure}

The Steps $1^-,2^-,\dots, 5^-$ are \label{cstep} {\it control} steps, where the copies of primitive S-machines work. For example,
we want to put Step $2^-$ between Step 1 and Step 2 (see fig. \ref{Pic4}). So we define the composition ${\bf P}_{12^-}$
of primitive S-machines working after the connecting rule $\theta(12^-)$ and the composition ${\bf P}_{22^-}$ of primitive S-machines working after $\theta(2^-2)^{-1}$, provided the
inverse of the canonical computation of ${\bf P}_{22^-}$ should
follow right after the canonical computation of ${\bf P}_{12^-}$.

Thus, to define the control S-machine of Step $2^-$ below, one should
define the order of the work of primitive components for
${\bf P}_{12^-}$ and ${\bf P}_{22^-}$  and choose either $p$-letters or $r$-letters to be control letters for these primitive
components.

\begin{rk} The control steps are used for double purpose.
If the base of a computation is standard, then the history of a control step restores all the configurations (Lemma \ref{121} (1)).
If the base is not reduced, then the control steps and
the right order of the work of their primitive components redound
to a linear bound of the space of the computation
in terms of the lengths of the first and the last words (Lemma \ref{nonst}).
\end{rk}

By default, every {\it connecting} rule $\theta(ij)$ locks a sector if this sector
is locked by all rules from $\Theta_i$ or if it is locked by all rules from $\Theta_j$.
It also changes all state letters used at Step $i$ by there copies from the disjoint set of state letters used at Step $j$.

{\bf Step $1^-$.} This is a control step between Steps $2$ and $1$. So  we define below the canonical work of the S-machines
${\bf P}_{2^-1}$ and ${\bf P}_{1^-1}$, keeping in mind that the last rule (as $\zeta^{21}$ in Remark \ref{ps}) of the canonical computation of ${\bf P}_{2^-1}$ switches
on the inverse computation for the canonical one of ${\bf P}_{1^-1}$. This conjunction of the S-machines ${\bf P}_{2^-1}$ and ${\bf P}_{1^-1}$ is the S-machine of Step $1^-$ denoted by ${\bf P}_{1^-}$.

Let us define ${\bf P}_{2^-1}$.
At first, we have the parallel work  of primitive S-machines in all
big historical sectors, and control $r$-letters
run forward and back according the rules from Subsection \ref{pm}. The tape
alphabet for every such primitive S-machine is the left alphabet ${\cal X}_{i,\ell}$ of the big historical sector.

The next primitive S-machine (see Remark \ref{ps} for the definition of composition) starts working similarly in the input
subsector of ${\bf M}_3$ after the above mentioned primitive S-machines
stop working.

Then we have parallel work of primitive S-machines
in the small historical sectors. Again, the running control letters are $r$-letters.

Finally, the primitive S-machine is switched on that checks the input word (of ${\bf D}_5$) between
${\cal R}^{1,r}$ and ${\cal P}^{2,\ell}$
(with running state letters from  ${\cal R}_s$).

The running control letters of ${\bf P}_{11^-}$ are $r$-letters again (not $p$-letters), and this S-machine is a copy of ${\bf P}_{21^-}$ with another set of state letters.

The transition rule $\theta(1^-1)$ changes all state letters of Step $1^-$ by their
copies in Step 1, which contain the letters of the start vector $\vec s_1$ of ${\bf M}_1$. It locks all sectors except for historical sectors, the sector ${\cal R}_{s-1}{\cal P}_s$ (we can call it the input since it comes from the input sectors of ${\bf M}_3$) and the sector
${\cal R}^{1,r}{\cal P}^{2,\ell}$.

The $\theta(1^-1)$-admissible words may involve the copies $p^{1,*}_i$ and $r^{1,*}_i$ of the letter
$p^1$ of a primitive S-machine $\bf Pr$, but no copies of $p^2$. They may contains letters from alphabets $ X_{i,\ell}$ but not from $X_{i,r}$.
\medskip

{\bf Step 1.} The rules $\theta({\bf M}_4)$ from $\Theta^+_1$ restricted to the base
${\cal P}_0{\cal Q}_0{\cal R}_0\dots{ \cal P}_s {\cal Q}_s {\cal R}_s $ are just the (copies of the) positive rules $\theta$ of ${\bf M}_3$. They do not change other sectors and lock the sector ${\cal R}^{2,r}{\cal P}^{3,\ell}$.

{\it Comment.}  At Step 1, ${\bf M}_4$ works as the S-machine ${\bf M}_3$.

The connecting rule $\theta(12^-)$ changes the state letters by their copies in disjoint alphabet, in particular, the letters from the accept vector $\vec s_0$  of ${\bf M}_1$ are replaced by their copies. The $\theta(12^-)$-admissible words have no letters from `left'
alphabets $X_{i,\ell}$.

Besides, the rule $\theta(12^-)$ `removes' one letter in  the sector
${\cal R}^{1,r}{\cal P}^{2,\ell}$: $\;\; [bp^1_{2,\ell}\to (p^1_{2,\ell})']$.

\medskip

{\bf Step $2^-$.}
This Step is similar to Step $1^-$, the difference in the definition of the S-machine ${\bf P}_{2^-}$ working at Step $2^-$ is that
alphabets $X_{i,\ell}$
should be replaced by alphabets $X_{i,r}$ ($i=1,\dots, s$)
and the control letters are $p$-letters (not $r$-letters).

{\it Comment.} The copies of primitive S-machine   check several sectors again.

The transition rule $\theta(2^-2)$ replaces all state letters of Step $2^-$ with their
copies in Step $2$, and the letters of  the accept vector $\vec s_1$  of ${\bf M}_1$ are among them. It locks non-historical sectors
except for the sector ${\cal R}^{1,r}{\cal P}^{2,\ell}$.

\medskip

{\bf Step 2.} The positive rules from $\Theta_2$ are just copies of the negative rules from $\Theta_1$.

{\it Comment.} ${\bf M}_4$ works as at Step 1, but reverses the computation procedure.

The connecting rule $\theta(21^{-})$ completes the cycle.

The connecting rule $\theta(23^-)$ makes possible final Steps $3^--6$.

\medskip

{\bf Step $3^-$.} As at Step $1^-$, the S-machine ${\bf P}_{3^-}$
is the conjunction of two S-machines: ${\bf P}_{23^-}$  and
${\bf P}_{33^-}$, where the first one is just a copy of ${\bf P}_{21^-}$ with different set of state letters.

 For the S-machine ${\bf P}_{33^-}$, the running state letters are
 $r$-letters too. Its canonical work is as follows. At first, the primitive S-machines
 simultaneously check the small historical sectors. Then the next primitive S-machine checks
 the sector ${\cal R}^{1,r}{\cal P}^{2,\ell}$, then the input sector ${\cal R}_{s-1}{\cal P}_s$ is checked,
 and finally the big historical sectors are simultaneously checked.

The connecting rule $\theta(3^-3)$ cannot be applied to an admissible word having $a$-letters from right alphabets of  historical sectors.

\medskip

{\bf Step 3.} The rules from $\Theta_3$ extend the rules of S-machine ${\bf D}_5$ as follows.  The rules
of ${\bf D}_5$ on the  configurations with base
\begin{equation}\label{D5}
S(1)_rS(2)_{\ell}S(2)_rT(1)_{\ell}T(1)_rT(2)_{\ell}T(2)_rT(3)_{\ell}
\end{equation}
are now the rules of Step 3 on the base
$${\cal P}_{s-1}{\cal Q}_{s-1}{\cal R}_{s-1}{\cal P}_s {\cal Q}_s{\cal R}_s {\cal P}{\cal S}(2)_r {\cal R} \times $$ $${\cal P}^{1,\ell}{\cal T}(1)_{\ell}{\cal R}^{1,\ell}{\cal P}^{1,r}{\cal T}(1)_r{\cal R}^{1,r}{\cal P}^{2,\ell}{\cal  T}(2)_{\ell}{\cal R}^{2,\ell}{\cal P}^{2,r}{\cal  T}(2)_r{\cal R}^{2,r}{\cal P}^{3,\ell}{\cal T}(3)_{\ell}{\cal R}^{3,\ell}$$
with control letters, according to the definition given in Subsection \ref{csl}, but the control letters do not work at this
step, and so ${\cal P}^{1,r}{\cal T}(1)_r{\cal R}^{1,r}$ in this
base behaves as $T(1)_r$ in (\ref{D5}), and so on.

 The rules of $\Theta_3$  do not change big historical sectors and lock non-historical sectors of ${\bf M}_3$,
except for the input sector ${\cal R}_{s-1}{\cal P}_s$.

{\it Comment.} After standard work with consecutive Steps $1^--2$
and control Step $3^-$, the (copy of) S-machine ${\bf D}_5$ checks if the length of $a$-word in the sector ${\cal R}^{1,\ell}{\cal P}^{2,r}$
divisible by the eight cubes of the $a$-length of the input sector.

The rule $\theta(34^-)$ locks all sectors except for historical ones. It cannot be applied to a word having $a$-letters from left alphabets of small historical sectors.

\medskip

{\bf Step $4^-$.} The tape alphabets of the control S-machine ${\bf P}_{4^-}$ are
right alphabets ${\cal X}_{i,r}$ for  small historical sectors
and left alphabets for big historical sectors. All working
sectors are locked.

The first half of ${\bf P}_{4^-}$ is the S-machine ${\bf P}_{34^-}$. Its
running control letters are $p$-letters, and the
canonical work
starts with  parallel work  of primitive S-machines in all
small historical sectors
followed by the parallel work of primitive S-machines
in the big historical sectors.

The second S-machine ${\bf P}_{44^-}$
starts with  parallel work  of primitive S-machines in all
small historical sectors with control $r$-letters followed by parallel work  of primitive S-machines in all
big historical sectors with control $p$-letters.

\medskip

{\bf Step 4.} The rules of $\Theta_4$ simultaneously erase the letters from
small historical sectors. The corresponding parts of the positive rules are
$r^1_{j,\ell}x\to r^1_{j,\ell}$ ($j=1,2,3$) for every positive letter $x$ from the right alphabet of a small historical sector.

The connecting rule $\theta(45^-)$ locks all sectors except
for big historical sectors.

{\bf Step $5^-$.} We define ${\bf P}_{5^-}= {\bf P}_{45^-}$. The control $p$-letters simultaneously
check the big historical sectors, while all other sectors are locked.

The connecting rule $\theta(5^-5)$ locks all sectors except
for big historical sectors.

{\bf Step 5.}
The rules of $\Theta_5$ simultaneously erase the letters from
big historical sectors. The corresponding parts of the positive rules are
$xp_j^1\to p_j^1$ ($j=1,\dots,s$) for every positive letter $x$ from the left alphabet of a big historical sector.

The accept command $\theta_0$ from $\Theta_5$ can be applied when all the sectors are empty. So it locks
all the sectors, changes the state letters and terminates the work of ${\bf M}_4$.

So ${\bf M}_4$ has a unique accept configuration.

\begin{lemma} \label{lo} Let a computation $\cal C$ of an S-machine ${\bf P}_{i^-}$ ($i=1,\dots, 5$) start with a connecting rule $\theta$ and end with a connecting rule $\theta'\ne\theta^{-1}$. Then for every sector of the standard base,
there is a rule in the history of $\cal C$ locking this sector.
\end{lemma}
\proof We consider ${\bf P}_{1^-}$ only.
Since the computation starts with $(21^-)$ and ends with $(1^-1)$,
all the primitive S-machines listed in the definition of the S-machine ${\bf P}_{1^-}$
have to start and finish their work.
So every sector of the standard base is locked either by  $\theta(21^-)$ or by $\theta(1^-1)$, or by a rule of a primitive S-machine of the form $\zeta^{12}$ - see Subsection \ref{pm} and Lemma \ref{prim} there, because every sector unlocked by
these connecting rules is checked by one of the primitive S-machines
forming the S-machine ${\bf P}_{1^-}$.
\endproof

\begin{rk}
Every cycle of the Steps $1^{-1}, 1, 2^{-}, 2$ just changes
the length of the sector ${\cal R}^{1,r}{\cal P}^{2,\ell}$ by $1$.
(See the definition of the connecting rule $\theta(12^-)$.) If
this length $\ell$ becomes divisible by $8k^3$, where $k$ is the length of the input sector of ${\bf M}_3$, the copy of the S-machine
$({\bf D}_5)_c$ can accept at Step 3, and one can obtain the stop
configuration of ${\bf M}_4$ after Steps $4^-, 4, 5^-, 5$. Hence
the shortest accepting computation has at most $\Theta(k^3)=\Theta(g(n))$ cycles of the Steps $1^{-1}, 1, 2^{-}, 2$ if the length of the sector ${\cal R}^{1,r}{\cal P}^{2,\ell}$
is $\Theta(n)$.

This is an informal answer to the question why the division
S-machine is needed. Indeed, if an auxiliary S-machine just checks if
$\ell$ is equal to $k^3$, then the number of cycles as above could
be $\Theta(n)$, which would lead to at least cubic Dehn function.

Another question: Since we want to repeat the cycle of the Steps $1^{-1}, 1, 2^{-}, 2\;\;\Theta(f(n)^3)$ times,
why does ${\bf M}_4$ recognize the values $f(n)$ instead
of $g(n)=f(n)^3$ ? - Because the Turing machine has to take time at least $\Theta(g(n))$ to recognize $g(n)$. By Lemma \ref{S},
the S-machine ${\bf M}_1$ should work as long as $\Theta(g^3(n))$
or longer for the same goal. If $g(n)$ is `almost' linear
function, then this approach makes the time function of ${\bf M}_1$ almost cubic, and the growth of the Dehn function we are going to construct,
becomes almost biquadratic.
\end{rk}

Let $B$ be the standard base of ${\bf M}_4$ and $B'$ be its disjoint copy. By ${\bf M}_5$ we denote
the S-machine with standard base $B(B')^{-1}$ and rules $\theta({\bf M}_5)=[\theta, \theta]$,
where $\theta\in \Theta$ and $\Theta$ is the set of rules of ${\bf M}_4$. So the rules of $\Theta({\bf M}_5)$ are the same
for ${\bf M}_4$-part of ${\bf M}_5$ and for the mirror copy of ${\bf M}_4$. Therefore we will denote $\Theta({\bf M}_5)$ by
$\Theta$ as well, assuming that the sector between $B$ and $(B')^{-1}$ is locked by any rule
from $\Theta$.

\medskip

Finally, the main S-machine ${\bf M}={\bf M}_6$ is a cyclic analog of ${\bf M}_5$. We add one more base letter $\{t\}$. So the standard base of ${\bf M}_6$ it $\{t\}B(B')^{-1}\{t\}$, where the part $\{t\}$ has only one letter $t$
and the first part $\{t\}$ is identified with the last part. (For example, $\{t\}B(B')^{-1}\{t\}B(B')^{-1}$ can be a base of an admissible word for ${\bf M}_6$. Furthemore,  sectors involving $t^{\pm 1}$ are
 locked by every rule from $\Theta$.
The stop word $W_0(\bf M)$ is defined accordingly: every letter in the standard base $B(B')^{-1}$  of ${\bf M}_5$ is replaced by the corresponding letters from the stop word of ${\bf M}_5$.

The `mirror' symmetry of the base of $\bf M$ will be used in Lemma \ref{001}. For a different purpose mirror symmetry
of Turing machines was used in the papers of W.W.Boon and P.S. Novikov (see \cite{R}).

\subsection{Standard computations of $\bf M$}\label{SC}

The history $H$ of a reduced computation of $\bf M$ can be factorized so that every
factor corresponds to one of the Steps $1^- - 5$. If, for example, $H=H'H''H'''$,
where $H'$ is a product of rules from Step 2, $H''$ has only one rule
$\theta(21^-)$ and $H'''$ is a product of  rules from Step $1^-$, then we say
that the \label{steph} {\it step history} of the computation (or its type) is $(2)(21^-)(1^-)$ or just
$(2)(1^-)$ since the only rule connecting the computations of Steps 2 and $1^-$ is
$\theta(21^{-})$ and for the most asymptotic estimates of the length of steps
(e.g.,  $||H''||$) or of the lengths of their admissible words, it does not matter
to which of the two possible steps the connecting rule is attributed.

There are no computations of some types, say $(1)(3)$, as it immediately follows from
the definition of connecting rules (and from fig. \ref{Pic4}). In this subsection, we eliminate some
other subwords in step histories.

\begin{lemma} \label{212} (1) There are no reduced computations of $\bf M$ with step histories
$(1^-1)(1)(1^-1)$, $(21^-)(2)(21^-)$ and $(23^-)(2)(23^-)$ (with $(12^-)(1)(12^-)$ and $(2^-2)(2)(2^-2)$)) if the base
of the computation contains at least one big historical sector $\dots {\cal P}_*$ (resp., big historical sector ${\cal R}_*\dots$ ).

(2) There are no reduces computations of $\bf M$ with step histories
$(3^-3)(3)(3^-3)$ (with $(34^-)(3)(34^-)$) if the base of the computation contains a small historical sectors $\dots{\cal P}^*$
(resp., a small historical sectors ${\cal R}^*\dots $).
\end{lemma}
\proof
(1) We consider only the type $(1^-1)(1)(1^-1)$ since other variants are similar.
Indeed, if the history $H$ of the part (1) is non-empty,
then a state $p$-letter inserts a copy of $H^{-1}$ in historical letters of the alphabet $X_{i,r}$ (see Subsection \ref{hs}). Recall that
the words with non-empty subwords over $X_{i,r}$ are not $\theta_{1^- 1}$-admissible, but we should
have  $\theta_{1^- 1}$-admissible words both in the beginning and at the end of the subcomputation
with history $H$, a contradiction.

(2) The proof is similar.

\endproof

\begin{lemma} \label{121} (1) There are no reduced computations of $\bf M$ with step histories of the form $(i^-j)(i^-)(i^-j)$ for $i=1,\dots,5$ if the base is standard.

(2) Let the base of a computation  ${\cal C}: W_0\to\dots\to W_t$ be standard and ${\cal C}$ has one of the step histories $(21^-)(1^-)(1^-1)$,
 $(12^-)(2^-)(2^-2)$, $(23^-)(3^-)(3^-3)$,  $(34^-)(4^-)(4^-4)$,  $(45^-)(5^-)(5^-5)$.
 Then all admissible words of $\cal C$ are uniquely defined by the history
 $H$ of $\cal C$,
 $|W_0|_a=|W_1|_a=\dots=|W_t|_a$, and  $||H||\le 4||W_0||$.

\end{lemma}

\proof (1) Consider only the step history $(21^-)(1^-)(21^-)$ and the work of the primitive S-machines switched on by the rule $(21^-)$. If one of the rules $\zeta^1(a)^{\pm 1}$ of this primitive S-machine (see Subsection \ref{pm})
does not move a control state letter right/left, but instead just insert $a'^{\pm 1}$ from the right and
$a^{\mp 1}$ from the left, then the rule $\zeta^{12}$ is not applicable since the
sector is not locked. So the next rule should be $\zeta^1(b)^{\pm1}$ which is not
the inverse one for $\zeta^1(a)^{\pm 1}$. Hence the control state letter has to insert letters from both sides
without cancellations, and neither $\zeta^{12}$ nor $(\zeta^{12})^{-1}$ can be ever applied, a contradiction.

Therefore the primitive S-machine must work canonically,
as it was described in Subsection \ref{pm} (also see Lemma
\ref{prim} (2)). Hence the history of its work uniquely restores the
words in the sectors controlled by this primitive S-machine.

Hence the first  primitive S-machine has to complete
its canonical work and switch on the next primitive S-machine
according to the definition of ${\bf P}_{1^-}$ for Step $(1^-)$, and so on.
Thus, one never obtains a
		$\theta(21^-)^{-1}$-admissible word, a contradiction.

(2) By Lemma \ref{prim} (3), the histories
of the primitive S-machines subsequently restore the tape words in all unlocked sectors.
Lemma \ref{Hprim} applied to $\cal C$ and to the inverse computation, implies equalities $|W_0|_a=|W_1|_a=\dots=|W_t|_a$, and  gives $||H||\le 4||W_0||$ if one takes into account that both the S-machine ${\bf P}_{21^-}$ and ${\bf P}_{21^-}$  control each of the sectors.
\endproof
\begin{lemma} \label{545} There are no reduced computations with standard base and step histories $(5^-)(4)(5^-)$.
\end{lemma}

\proof Assume such a computation exists. Note that the small historical sectors are empty at the both
transitions $(45^-)$. However every rule of Step 4 multiplies
the word in the small historical sector by a letter $x$, and one
obtains no cancellations of these letters since the computation is reduced and different rules multiply by different letters. Hence
the part $(4)$ is empty, and the lemma is proved by contradiction.
\endproof
	
Below we need a rougher subdivision of the history of a reduced computation with {\it standard} base. We say that the
Steps $1-4^-$ are {\it fundamental} steps and the Steps $4,5^-,5$ are {\it erasing} steps. So the \label{blockh} {\it block history} of every computation of $\bf M$ is a subword of $(F)(E)(F)(E)\dots$, where
$(F)$ (where $(E)$) are maximal parts of the history with fundamental (resp., erasing) steps only. The separating rules
for neighbor blocks are $\theta(4^-4)^{\pm 1}$.

\begin{lemma} \label{474} If the block history of a computation is
$(E)$, then its step history is a subword of the word $(4)(5^-)(5)(5^-)(4)$.
\end{lemma}
\proof
Proving by contradiction and taking
into account Lemma \ref{121} (1), we should get a subword
$(5^-)(4)(5^-)$
in the step history,
contrary to Lemma \ref{545}.
\endproof

\begin{lemma}\label{E} Let $W_0\to\dots\to W_t$ be a computation
with block history $(E)$. Then

(1) $|W_j|_a\le \max(|W_0|_a, |W_t|_a)$
for $j=0,1,\dots, t$;

(2) $t\le 10(||W_0||+ ||W_t||)$.
\end{lemma}

\proof
(1) If the history has only one Step 4 or 5, then Statement (1) follows from Lemma \ref{gen} (c). For single Step $5^-$ it follows
from Lemma \ref{Hprim} (a).

If there is Step 5 in the computation, then by Lemma \ref{474}  we have only one maximal subcomputation $W_k\to\dots\to
W_{\ell}$ of Step 5. Here $|W_k|_a\le |W_0|_a$ since $W_k$ has no
non-empty sectors except for big historical sectors, which are
unchanged at Steps $4$, while Steps $5^-$ cannot
decrease the sum of length of these sectors by projection argument (see Remark \ref{proj}). Hence
it suffices to prove Statement (1) for subcomputations
with step histories $H'$ and $H''$, where $H'(5)H''$ is the step history.
Therefore we may
prove Statement (1) under assumption that there are no Steps 5
in the step history.

For the step history $(4)(5^-)$, the $a$-length of the configuration separating two steps cannot be longer than the
final configuration by projection argument, which reduce the proof to one-step histories.

(2) By Lemma \ref{474}, we have at most 5 steps, and by Property (1) it suffices to prove (2) for one-step histories but with coefficient $2$. Indeed such estimates for the lengths of histories are obtained for Steps $4, 5^-, 5$ in Lemmas \ref{gen} (b) and \ref{Hprim} (b).
\endproof

\begin{lemma} \label{histF} Let $W_0\to\dots\to W_t$ be a computation
with block history $(F)$.

(1) Then the step history of this computation
is a subword of the word ${\cal W}(m)\equiv (4^-)(3)(3^-)\big( (1^-)(1)(2^-)(2)\big)^m(3^-)(3)(4^-)$ for some non-zero integer $m$. (Here we define $\big( (1^-)(1)(2^-)(2)\big)^{-1}=(2)(2^-)(1)(1^-)$.)

(2) If the step history is equal to ${\cal W}(m)$ for some
$m\ne 0$, then  the exponent $m$ is divisible by $(2k)^3$, where $k$ is the $a$-length
of the input sector ${\cal R}_{s-1}{\cal P}_s$ after the application
of the connecting rule $\theta(3^-3)$.

(3) If the step history is equal to $\big( (1^-)(1)(2^-)(2)\big)^m(3^-)(3)(4^-)$ and the history starts with a connecting rule, then  the exponent $m$ is
congruent to $l$ modulo $(2k)^3$,
 where $k$ is as in Statement (2)
 and $b^{l}$ is the tape word of the sector ${\cal R}^{1,r}{\cal P}^{2,\ell}$ in the beginning of the computation.
\end{lemma}

\proof (1) Since the block history is $(F)$, Lemmas \ref{121} (1) and \ref{212}   forbid `reverse
moves' $\theta^{\pm 1}H\theta^{\mp 1}$ in the history, where
$H$ is a one-step history and $\theta$ is a connecting rule,
Statement (1) follows from the definition of
connecting rules between steps (see figure \ref{Pic4}).

(2) Let us restrict the subcomputation with step history
$(3^-3)(3)(34^-)$ to the input sector ${\cal R}_{s-1}{\cal P}_s$
and the sector ${\cal R}^{1,r}{\cal P}^{2,\ell}$. Then we can apply
Lemma \ref{div} to conclude that the exponent $l$ of the tape word $b^{l}$ in the sector ${\cal R}^{1,r}{\cal P}^{2,\ell}$ at the beginning of this computation is divisible by $(2k)^3$.
Similarly, for $ b^{l'}$ at the end of the computation
with step history $(34^-)(3)(3^-3)$, we obtain that $l'$ is  divisible by $(2k')^3$, where $k'$ is the $a$-length of the sector ${\cal R}_{s-1}{\cal P}_s$.

Note that a computation with step history $(1^-1)(1)(2^-)(2)(21^-)$
(or \\$(1^-1)(1)(2^-)(2)(23^-)$)
does not change the $a$-length of the sector ${\cal R}_{s-1}{\cal P}_s$
since it is preserved by the S-machine ${\bf P}_{2^-}$ by Lemma \ref{prim} and the history
of Step 2 is inverse (of the copy) of the history of Step 1 here.
The same is true for the subcomputations with step history
$(21^-)(1^-)(1^-1)$ and $(23^-)(3^-)(3^-3)$ by Lemma \ref{prim}.
So $k'=k$ and $l'-l$ is divisible by $(2k)^3$ by Lemma \ref{div}.

On the other hand, every rule $\theta(12^-)$ multiplies
the tape word in the sector ${\cal R}^{1,r}{\cal P}^{2,\ell}$. by $b^{-1}$.
Therefore $\ell'-\ell=m$, whence $m$ is divisible by $(2k)^3$, as required.

(3) The proof is similar to the proof of Statement (2).
\endproof

\subsection {Computations with faulty bases}\label{fau}

\begin{lemma} \label{123} If the step history of a reduced computation $\cal C$ of $\bf M$ is $(21^-)(1^-)(1^-1)$ or $(12^-)(2^-)(2^-2)$,
or $(23^-)(3^-)(3^-3)$,
or $(34^-)(4^-)(4^-4)$, or $(45^-)(5^-)(5^-5)$,
then

(1) the base  of $\cal C$ is a reduced word;

(2) the first letter of the base, the length of the base and the history $H$ of $\cal C$ completely determine the computation $\cal C$.
\end{lemma}
\proof (1) This follows from Lemmas \ref{lo} and \ref{qqiv}.

(2) Consider for example the step history $(21^-)(1^-)(1^-1)$.
By Property (1) and the definition of admissible word, the base
of $\cal C$ is determined by its length and the first letter. Since
every sector unlocked by $\theta(21^-)$ has to be checked by a primitive S-machine, the copy of the content of this sector is contained in the history of the
computation as a product of the letters of the form $\zeta^1(a)^{\pm 1}$ defined for the particular primitive S-machine.
\endproof

\begin{lemma} \label{WV} Suppose that a reduced computation $W_0\to\dots\to W_t$ of Step 1 or 2 (of Step 3) starts  with a connecting rule. Assume that the length of its base $B$ is bounded from above by a constant  $N_0$, and $B$ has
a big historical sector (a small historical sector, resp.) of the form $\cal RP$ (with indices). There is a constant $c=c(N_0)$ such that
$|W_0|_a\le c|W_t|_a$.
\end{lemma}

\proof
 Let $V_0\to\dots\to V_t$
be the restriction of the computation  to
the sector with base ${\cal R}{\cal P}$. By Lemma \ref{w}, we have
$t\le |V_t|_a$ and $|V_0|_a\le |V_t|_a$.

It follows from lemma \ref{simp} (2) that

$$|W_0|_a
\le |W_t|_a+ 2N_0t\le |W_t|_a+ 2N_0|V_t|_a\le (2N_0+1)|W_t|_a.$$

It suffices to choose $c=2N_0+1$.
\endproof

\begin{df}\label{faul}
We call a base of $\bf M$ {\it faulty} if

\begin{itemize}

\item{it starts and ends with the same base letter,}

\item{it has no proper subwords with this property, and}

\item{it is a not a reduced word.}

\end{itemize}

\end{df}

Note that if there is a computation $\cal C$
with a faulty base $U_1\dots U_i\dots U_s$
(where $U_1=U_s$), then one can replace every
admissible word of this computation by the cyclic shift of it with faulty base $U_i\dots U_sU_1\dots U_{i-1}U_i$ and obtain the {\it cyclic
shift} $\cal C'$ of $\cal C$.

The main lemma of this subsection is

\begin{lemma} \label{nonst} There is a constant $C=C(\bf M)$, such that for every reduced computation ${\cal C}:\; W_0\to\dots\to W_t$
of $\bf M$ with a faulty base and every $j=0,1,\dots,t$, we have $|W_j|_a\le C\max(|W_0|_a, |W_t|_a)$.
\end{lemma}

\proof {\bf 0.} If the faulty base is not of the form $(pp^{-1}p)^{\pm 1}$ for some control state letter $p$ (or $r$)),
then the words $W_0,\dots,W_t$ have to contain non-control state letters. Hence we can replace
the computation $\cal C$ by a cyclic shift of it
and suppose further that the first and the last
letters of $W_0,\dots,W_t$ are not control letters.

{\bf 1.} By Lemma \ref{123}, the step history (and the inverse step history) has no subwords $(21^-)(1^-)(1^-1)$, $(12^-)(2^-)(2^-2)$,
$(23^-)(3^-)(3^-3)$, $(34^-)(4^-)(4^-4)$ and $(45^-)(5^-)(5^-5)$.

\medskip

{\bf 2.} Assume that the  history has only one Step $(2^-)$ and
the base is $(pp^{-1}p)^{-1}$ for some control state letter $p$ (or $r$)).
Obviously, the running letter $p$ cannot lock any sector and
so every rule has the same type (either $\zeta^1(*)$ or $\zeta^2(*)$ as in Subsection \ref{prim}).

If $|W_{j+1}|_a>|W_j|_a$ for some $j$, then one of two sectors
(e.g., the first one) increases its length by $2$, while another
one does not decrease the length under the transition $W_{j}\to W_{j+1}$. But then the first sector of $W_{j+2}$
has to increase again by $2$ since the computation is reduced.
It follows that $|W_{j+1}|\le |W_{j+2}|\le\dots\le |W_t|$.
Hence the length of every admissible word does not exceed
$\max(|W_0|_a,|W_t|_a)$.

If the faulty base has no subwords of the form $(pp^{-1}p)^{\pm 1}$,
the same inequality for a computation with history of Step $2^-$ follows from Lemmas \ref{prim} and \ref{ewe}. Thus, one may
assume further that the step history is not $(2^-)$. Similarly,
it is neither $(1^-)$ nor $(3^-)$, nor $(4^-)$, nor $(5^-)$.

\medskip

{\bf 3.} Assume there is $(1)$ in the step history.
Then by item 1, the set of steps is either (a) $\{(1)\}$ or (b) $\{(1), (1^-)\}$, or (c) $\{(1), (2^-)\}$, or (d) we have the subword $(1^-)(1)(2^-)$ in the step history (or in the inverse step history).

{\bf 3a.} In this case, the required inequality follows from
Lemma \ref{9}. (Recall that by definition the connecting rule $\theta(12^-)$ inserts/deletes one letter,
but this small change of length is not an obstacle.)

{\bf 3c.} Assume that the step history is $(2^-)(1)$ (or $(1)(2^-)$) and $W_0\to\dots\to W_j$ is a maximal subcomputation
with step history $(2^-)$. Then $|W_0|_a\ge |W_j|_a$ by the projection
argument. Therefore it suffices to prove the statement of the lemma
for a subcomputations with step histories $(2^-)$ and $(1)$.
For the case $(2^-)$, we refer to item 2. The case $(1)$ is considered in item 3a. Thus, we
assume further that the step history has length at least 3.

Suppose the step history has a subword $(2^-)(1)(2^-)$. Then the
base (or the inverse base) has no big historical sectors ${\cal R}\dots $ (with indices)  by Lemma \ref {212}.

So the only possible bases for
 big historical sectors have form ${\cal P}^{-1}{\cal P}$.
 The state control letters of such sector should start running at Step $(2^-)$ after the application of the connecting rule $\theta(12^-)$
but they cannot ever meet state letter from ${\cal R}$ and
will run forever, and $|W_t|_a\ge |W_j|_a$ by Lemma \ref{ewe},
if $W_j$ is obtained at the application of $\theta(12^-)$.
Hence one can cut the computation at $W_j$ reducing the problem
to items 2 and 3a.
If there are no big historical sectors, then  there are no working sectors of
the S-machine ${\bf M}_3$.
The other sectors of $\bf M$ (which
could come from the base of ${\bf D}_5$) do not work at Steps $(1)$ and
cannot decrease their length starting from a connecting rule at Steps $(2^-)$ by Lemma \ref{prim}. This makes the statement of the lemma obvious for them.

Hence one may assume that there are no subwords $(2^-)(1)(2^-)$, and so the step history is  $(1)(2^-)(1)$.

Suppose the base has a big historical sector ${\cal R}{\cal P}$ (with indices).
Then for the maximal subcomputation $W_r\to\dots\to W_t$ of Step 1, we have $|W_r|_a\le c|W_t|_a$ by Lemma \ref{WV} because the length of a faulty base is bounded by a constant
$N_0$ depending on the S-machine $\bf M$ only.
Hence one can reduce our task to the subcomputations with the
step histories $(1)$ and $(2^-)$. (Of course, the constant in the
desired inequality changes when we pass to step histories involving
more types of steps.) Hence we assume further that the base
has no big historical sectors ${\cal R}{\cal P}$.

 Also there are no big historical sectors ${\cal P}^{-1} {\cal P}$
 because state control letters of such sector should start running after the application of the connecting rule $\theta(12^-)$
but they cannot ever meet state letter from ${\cal R}$ and
will run forever by Lemma \ref{ewe}; the last Step 1 will not be reached.

So all big historical sectors (if any) are of the form ${\cal R}^{-1}{\cal R}$. Recall that the alphabets for
the $\theta(12^-)^{\pm 1}$-admissible words are $X_{i,r}$,
and so the word in this alphabet will be conjugated
at Step 1 by the letters from the disjoint alphabets $X_{i,\ell}$
in sectors with bases ${\cal R}_{i-1}^{-1}{\cal R}_{i-1}$.
Hence after application of $\theta(12^-)^{-1}$, each rule of Step 1
will increase the length of such sector by $2$.
By Lemmas \ref{three} and \ref{simp}, we have
$|W_r|_a\le\dots\le|W_t|_a$ if the last Step $1$ starts with $W_r$. This reduces the problem to the subcomputations with the
step histories $(1)$ and $(2^-)$ again.

If there are no big historical sectors, then we have no working sectors
except for the sectors of the S-machine ${\bf D}_5$ because one may assume that the left-most sector of the standard base of ${\bf M}_2$ is always locked
and because the base is faulty. The other sectors (which
could come from the base of ${\bf D}_5$) do not work at Steps $(1)$ and
cannot decrease their length towards $W_0$ or $W_t$ at Steps $(2^-)$ by Lemma \ref{prim}. This makes the statement of the lemma obvious for them and completes case 3c.

{\bf 3b.} This case is similar to 3c  up to exchange of $\cal R$ with $\cal P$ and $X_{i,\ell}$ with $X_{i,r}$.

{\bf 3d.} Assume that the step history has a subword $(1^-)(1)(2^-)$. Then there are no big historical sectors ${\cal P}^{-1}{\cal P}$ and ${\cal R}{\cal R}^{-1}$ since
the conjugation in free group given by Step $1$ cannot transform a non-trivial
word in the alphabet $X_{i,\ell}$ in a word in the alphabet $X_{i,r}$. So all big historical sectors have base of the form
${\cal R}{\cal P}$ (or inverse one).

Consider the word $W_j$ obtained after the application of the last
connecting rule $\theta(12^-)^{\pm 1}$.
Only big historical sectors of $W_j$ are not locked (except for the
sectors of the S-machine ${\bf D}_5$, which are not touched by Step 1).
If the next step is $1$, then its rules cannot make the historical sector shorter by Lemma \ref{w}. If the next step is $2^-$, then no sector  becomes shorter
by Lemma \ref{prim}. The same is true for Step $1^-$ if it follows Step 1.
Repeating this procedure, we have $|W_j|_a\le |W_t|_a$. Therefore
one can reduce Case 3d to previous cases by subdivision of
the computation along the transitions between Steps $1$ and $2^-$. (The occurrences
of $\theta(12^-)$ and $\theta(12^-)^{-1}$ in the history of the computation have to alternate in
Case 3d inserting/deleting the same letter, and so this does not affect the
desired inequality.)

Thus, we may assume further that there are no Steps $1$ in the computation.

\medskip

{\bf 4.} We may also assume that there are no Steps $2$ in the computation. The proof copies the proof at item 3 with subcases
(a) $\{(2)\}$ (b) $\{(2), (2^-)\}$, (c) $\{(2), (3^-)\}$ and (d), where one considers the subword $(2^-)(2)(3^-)$ and the word $W_j$
provided by the last connecting rule $\theta(2^-2)^{\pm 1}$.

\medskip

{\bf 5.} If there is Step 5 in the step history, then we cannot have
steps except for $5$ and $5^-$ by item 1.
It follows from the definition of Step $5^-$ and
Lemma \ref{prim} (4) that the subword $(5)(5^-)(5)$ in the step
history is not possible if there is a big historical sector with base of type
${\cal R}{\cal P}$. By Lemma \ref{ewe}, historical sectors
with base ${\cal P}^{-1}{\cal P}$ are not possible too.
However the historical sectors with base ${\cal R}{\cal R}^{-1}$ do not change  words at Step 5. It
follows that the word $W_j$ obtained after the transition from $(5^-)$ to $(5)$ is not longer than $W_t$ (see similar argument
at item 3d). So one reduces the task to shorter step histories.

Thus, one may assume that there are no subwords $(5)(5^-)(5)$
in the step history.  Therefore assuming that there is $(5)$ in the step history, we should consider only the history $(5^-)(5)(5^-)$
or its subwords. Again the rules $\theta(5^-5)$ defines a word $W_j$, which is not longer than $W_t$ by Lemma \ref{prim}. This reduce the task to step histories $(5^-)$ and $(5)$. For $(5^-)$, the problem is solved in item 2, and it is solved by Lemma \ref{gen} (c) for $(5)$.

Thus one may assume from now that there are no Steps (5) in the step history.

\medskip

{\bf 6.} If there is $(5^-)$ in the step history, then by items
1 and 5, there are no other steps except for $(4)$ and $(4^-)$.
However the transition $(45^-)$  provides us with the shortest words
in the computations since neither computations of Step $5^-$
nor those from $(4)$ can make big historical sectors shorter by
Lemmas \ref{prim} and \ref{ewe}. So cutting the computation
along such transitions, we can decrease the number of steps.
Since a single Step $5^-$ is eliminated in item 2, we may further
assume that there are no Steps $(5^-)$ in the computation.

\medskip

{\bf 7.} If there  is $(4)$ in the step history, then by items
1 and 6, there are no other steps except for  $(4^-)$. For such histories we will repeat some arguments from item 3c using now
small histories instead of big ones.

The case with a single Step $(4^-)$ is eliminated
by item 2. The brief history $(4)$ is also eliminated by Lemma \ref{gen} (c) (for small historical
sectors $\cal RP$) and Lemma \ref{gen3} (for small
historical sectors ${\cal RR}^{-1}$).
Assume that the step history is $(4^-)(4)$  and $W_0\to\dots\to W_j$ is a maximal subcomputation
with step history $(4^-)$. Then $|W_0|_a\ge |W_j|_a$ by the projection
argument. Therefore the problem is reduced to single step histories Thus, we assume further that the step history has length at least 3.

Assume that there is a subword $(4)(4^-)(4)$ in the step history and
there is the letter ${\cal R}_s$ in the base.
Then it cannot follow
by ${\cal R}_s^{-1}$, because the letter from  ${\cal R}_s$ must start running right after the connecting rule
$\theta(4^-4)^{-1}$ and it cannot ever reach the part ${\cal P}$,
so the next connecting
rule $\theta(4^-4)$ cannot appear, a contradiction. Hence there is
a sector ${\cal R}_s{\cal P}$. It follows that there is ${\cal R}^{1,\ell}$ in the base since the sectors between ${\cal P}$ and
${\cal R}^{1,\ell}$  are locked by $\theta(4^-4)$. Then we obtain ${\cal P}^{1,r}$ and all other parts of small historical sectors.

Also there is the part
${\cal P}_{s-1}$ in the base, becase all the sectors between
${\cal P}_{s-1}$ and  ${\cal R}$ are locked by $\theta(4^-4)$.
Hence right after the sector ${\cal R}_s{\cal P}$
is checked by the first primitive S-machine of ${\bf P}_{44^-}$,
a letter from ${\cal P}_{s-1}$ starts running checking the
big historical sectors, and again, there should be the part ${\cal R}_{s-2}$ in the base, since otherwise the next occurrence of $\theta(4^-4)$
does not happen.

This implies it turn, that there is ${\cal P}_{s-2}$
in the base, and so on, that is we have all the sectors of the standard
base of $\bf M$ and the base of our computation has no cancellation,
contrary to the definition of faulty base.

Hence there are no parts ${\cal R}_s$ in the base,
provided $(4)(4^-)(4)$ is a subword of the step history.  Similar
argument shows that there are no ${\cal R}^{1,\ell}$ and ${\cal R}^{2,\ell}$. (For example, if there is ${\cal R}^{1,\ell}$,
then there is a rule locking the sector ${\cal R}^{1,\ell}{\cal P}^{1,r}$ by the definition of primitive S-machines; this rule has to lock sector ${\cal R}_s{\cal P}$ too, and so the part ${\cal R}_s$
 occurs in the base too.) Hence nothing changes at Steps 4, since
 only $\cal R$-letters can erase the small historical sectors.
The transition $(44^-)$  provides us with the shortest words
in the computations since neither computations of Step $4^-$
nor those from $(4)$ can make small historical sectors shorter by
Lemmas \ref{prim} and \ref{ewe}. So cutting the computation
along such transitions, we can decrease the number of steps.
Therefore we may further
assume that there are no subwords $(4)(4^-)(4)$ in the step
history.

Thus, the step history is $(4^-)(4)(4^-)$, and as above it can be
subdivided in one-step histories. Therefore we may assume from now
that there are no Steps 4 in the step history.

\medskip

{\bf 8.} Now we assume that there is Step $4^-$ in the step history.

Suppose the step history has a subword $(4^-)(3)(4^-)$. Then the
base has no small historical sectors ${\cal R}^{*}\dots$  by Lemma \ref {212} (2).
So the only possible bases for
 small historical sectors are $({\cal P}^{i,r})^{-1}{\cal P}^{i,r}$ or ${\cal P}^{-1}{\cal P}$.
The state control letters of such sectors should start running after the application of the connecting rule $\theta(34^-)$
but they cannot ever meet state $\cal R$-letters
and will run forever, whence $|W_t|_a\ge |W_j|_a$ by Lemma \ref{ewe},
if $W_j$ is obtained at the application of $\theta(34^-)$. So the
whole step history is $(4^-)(3)(4^-)$.
Hence one can cut the computation at $W_j$ reducing the problem
to item 2 and Lemma \ref{9D}.
If there are no small historical sectors, then  there are no sectors of
the S-machine ${\bf D}_5$ changed at Step 3.
The other sectors (which
could come from the base of ${\bf M}_2$) do not work at Steps $3$ and
cannot decrease their length when starting from a connecting rule at Steps $4^-$ by Lemma \ref{prim}. As usual, this allows induct on the
number of steps.

Hence one may assume that there are no subwords $(4^-)(3)(4^-)$. The subwords $(3^-)(3)(3^-)$ are eliminated in the same way.

Assume that the step history has a subword $(3^-)(3)(4^-)$. Then there are no small historical sectors with non-reduced bases $UU^{-1}$
since
the conjugation in free group given by Step $3$ cannot transform a non-trivial
word in the alphabet $X_{i,\ell}$ in a word in the alphabet $X_{i,r}$. So all small historical sectors have base of the form
${\cal R}{\cal P}$
(or inverse ones).

Consider the word $W_j$ obtained after the application of the last
connecting rule $\theta(34^-)^{\pm 1}$.
Only small historical sectors of $W_j$ are not locked (except for the
sectors of S-machine ${\bf M}_3$, which are not touched by Step 3).
If the next step is $(3)$, then its rules cannot make the historical sector shorter by Lemma \ref{w}. If the next Step $(4^-)$, then no sector  becomes shorter
by Lemma \ref{prim}. The same is true for Step $3^-$ if it follows Step 3.
Repeating this procedure, we have $|W_j|_a\le |W_t|_a$. Therefore
one can reduce the problem to shorter step histories.

Assume that the step history is $(3)(4^-)(3)$. This resembles item 3c, but below we consider small historical sectors instead of the big ones.

 Suppose the base has a small historical sector ${\cal R}^{i,\ell}{\cal P}^{i,r}$ or ${\cal R}_s\cal P$. Let $W_r\to\dots\to W_t$
be the maximal subcomputation with step history $(3)$.
Then we obtain inequality $|W_r|_a\le c|W_t|_a$ by Lemma \ref{WV}.
Hence one can reduce our task to  the subcomputations with the
step histories $(3)$ and $(4^-)$. Therefore we assume further that the base
has no small historical sectors of the form ${\cal R}{\cal P}$.

 Also there are no small historical sectors $({\cal P}^{i,r})^{-1}{\cal P}^{i,r}$ (or ${\cal P}^{-1}\cal P$),
 because state control letters of such sector should start running after the application of the connecting rule $\theta(34^-)$
but they cannot ever meet state letter from ${\cal R}$ and
will run forever by Lemma \ref{ewe}; the last Step 3 will not be reached.

So all small historical sectors are of the form
${\cal R}^{i,\ell}({\cal R}^{i,\ell})^{-1}$ (or ${\cal R}_s {\cal R}_s^{-1}$)
But then the word from such a sector in `right' alphabet will be conjugated at Step 3 by the letters from a left alphabet.
Hence after application of $\theta(34)^{-1}$, each rule of Step 3
will increase the length of such sector by $2$.
By Lemmas \ref{three} and \ref{simp}, we have
$|W_r|_a\le\dots\le|W_t|_a$ if the last Step $3$ starts with $W_r$. This reduces the problem to the subcomputations with the
step histories $(3)$ and $(4^-)$ again.

If there are no such sectors, then there are no sectors of
the S-machine ${\bf D}_5$ changed at Step 3  since the base is faulty. The other sectors (which
could come from the base of ${\bf M}_2$) do not work at Steps $3$ and
cannot decrease the length towards $W_0$ or $W_t$ at Steps $3$.
This makes the statement of the lemma obvious.

It remains to consider the brief history $(3)(4^-)$ and apply the projection
argument to the part $(4^-)$. (Compare with case $(2^-)(1)$ in item 3c.)

Hence one may assume from now that there are no Steps $4^-$ in the computation.

\medskip

{\bf 9.} If there is Step 3 in the step history, then there are
no steps except for $3$ and $(3^-)$. As in item 8, one may assume
that the length of step history is at least 3. Then the subwords
$(3^-)(3)(3^-)$ can be eliminated by the same argument we used
in item 8 to eliminate subwords $(4^-)(3)(4^-)$. It remains to consider
computations with step history $(3)(3^-)(3)$. Again, one refer to item 8
since one can eliminate such history in the way the histories
$(3)(4^-)(3)$ were eliminated in item 8.

The lemma is proved.
\endproof

\subsection{Space and time of $\bf M$-computations with standard base}

\begin{lemma}\label{cycle}

Let ${\cal C}: W_0\to\dots\to W_t$ be a computation with standard base and step history $(21^-)(1^-)(1)(2^-)(2)(21^-)$. Then the
configuration $W_t$ is a copy of $W_0$ except for the sector
${\cal R}^{1,r}{\cal P}^{2,l}$ and the mirror copy of it, whose
lengths in $W_0$ and $W_t$ differ by one.

If $H(1^-)$ and $H(2^-)$ are the histories of the subcomputations
${\cal C}(1^-)$ and ${\cal C}(2 ^-)$ of $\cal C$ with step histories $(21^-)(1^-)(1^-1)$ and $(12^-)(2^-)(21^-)$,
respectively, then \\$|||H(1^-)||-||H(2^-)|||= 2$.
\end{lemma}

\proof The  subcomputations ${\cal C}(1^-)$ and ${\cal C}(2 ^-)$  do not change the $a$-words in the historical sectors by Lemma \ref{Hprim}(a), and so the histories
 of Steps 1 and 2 are inverse copies of each other. Taking into account
 that the transition $(12^-)$ changes the length of the sector
 ${\cal R}^{1,r}{\cal P}^{2,l}$ (and the mirror it) by 1, we
 obtain the first statement of the lemma.

 After the sector ${\cal R}^{1,r}{\cal P}^{2,l}$ changes length
 by one, the primitive S-machine checking this sector changes the computation time by $2$, as it follows from Lemma \ref{prim} (3).
 This proves the second statement.
 \endproof

 Recall that the blocks $(E)$ and $(F)$ of a history were defined before Lemma \ref{474}.

\begin{lemma}\label{standard}
(1) If a configuration $W_0$ is accepted and $\theta$-admissible for a rule $\theta$ from block $(E)$, then
there is a reduced
accepting computation $W_0\to\dots\to W_t$ with block
history $(E)$ and $t\le 3||W_0||$;

(2) There is a constant $c_1$ depending on $\bf M$ only, such that for any
 computation  ${\cal C}:\;W_0\to\dots\to W_t$  of $\bf M$, which is the beginning of a reduced accepting computation with block history $(E)$ or $(F)(E)$, we have
$||W_j||\le c_1(||W_0||+s)$ for every $j=0,\dots,t$, where $s$ is the length of the step history.

 \end{lemma}

 \proof
 (1) Note that all the sectors of $W_t$ are empty. It follows that all the words in big historical sectors of $W_0$ are
 copies of the same word since every rule of $\bf M$ multiplies
 the tape words of these sectors by the copies of the same letter or does not change all these words. Similarly, the tape words in the small historical sectors
 (and in their mirror sectors) are copies of each other.

 If the small historical sectors of $W_0$ are non-empty,
 then $\theta$ is a rule of Step 4.  So there is an
 accepting computation erasing all the letters of this sector
 (and of its mirror copy). The length of the next control Step $5^-$ will be at most $2||W_0||$ by the definition of the S-machine
 ${\bf P}_{5^-}$ and Lemma \ref{prim} (3). Then the rules of Step $5$
 can erase all tape letters in the big historical sectors (and their mirror copies).
 This gives the total upper estimate $t\le 2||W_0||+||W_0||$, as required. If $\theta$ is a rule of Step $5^-$ or $5$, then
 the estimate is even better.

 (2) If the computation is accepting and has type $(E)$, then the step history is a suffix of the word $(4)(5^-)(5)$ by Lemma \ref{474}. The rules of Step $5^-$ cannot increase the lengths of configurations by Lemma \ref{prim}.
 Clearly, the rules of Steps $4$ and $5$ cannot insert letters too.
 Hence $||W_j||\le |W_0||$, and so it suffices to prove the
 same inequality under the assumption that $j\le r$, where
 $W_0\to\dots\to W_r$ is the maximal subcomputation with block history  $(F)$.

 By Lemma \ref{histF}, the step history of this subcomputation is a suffix of the word  ${\cal W}(m)\equiv (4^-)(3)(3^-)\big( (1^-)(1)(2^-)(2)\big)^m(3^-)(3)(4^-)$ for some non-zero integer $m$.

 At first we consider the subcomputations ${\cal C}_i$ with step histories \\$\big((21^-)(1^-)(1)(2^-)(2)(21^-)\big)^{\pm 1}$. By Lemma \ref{cycle}, we conclude
 that in the beginning and at the end of  ${\cal C}_i$,  the difference
 of lengths of configurations is equal to $\pm 2$. The number of
 such subcomputations ${\cal C}_i$ does not exceed $s/2$.

 The number of one-step subcomputations, which are not subcomputations of any ${\cal C}_i$, is at most $7$. The transitions of some of them ($1^-, 2^-, 3^-, 4^-$) do not increase the lengths of configurations by Lemma \ref{Hprim} (a). All transitions of each other step ($1,2,3$) can increase
 the $a$-length but Lemma \ref{WV} bounds  possible enlargement from above. This proves Statement (2).\endproof

 \begin{lemma}\label{boundH} Let a history $H$ of a computation ${\cal C}:\;W_0\to\dots\to W_t$ with standard base  have type $(F)$ and end with a connecting rule. Suppose that $\cal C$ is a beginning of a reduced accepting computation and there are at most $10$ steps in $\cal C$. Then for  a constant $c_2=c_2(\bf M)$, we have $||H||\le c_2||W_0||$.

 \end{lemma}

 \proof Recall that the length of the history is linearly bounded for each of the steps in terms of the lengths of their configurations
(Lemma \ref{gen}  for Steps 4 and 5, Lemma \ref{w} for Steps 1, 2 and 3, Lemmas \ref{121}(2) and \ref{prim} (1) for control steps.) The lengths
of these configurations are linearly bounded in Lemma \ref{standard} (2). Taking into account that the number of steps is at most $10$
we come to the desired inequality.
\endproof

 \begin{lemma}\label{stand}

 There is a constant $c_3$ such that for any accepted configuration $W_0$, which is $\theta(2,1^-)^{\pm 1}$-admissible, there exists an accepting computation with block history $(F)(E)$ of length
at most $c_3(k^3+1)(||W_0||+k^3)$, where $k$ is the $a$-length of the input sector ${\cal R}_s{\cal P}$ of the word $W_0$. The number
of steps in this computation is less than $32k^3+4.$
\end{lemma}

\proof
By Lemma \ref{standard} (1), an accepted word $\theta$-admissible
 for a rule $\theta$ from block $(E)$ can be accepted by a computation
 having only one block . Hence it suffices to consider computations with block histories
 $(F)(E)$.

By Lemma \ref{histF}, the step history of block $(F)$ is a  word $\big( (1^-)(1)(2^-)(2)\big)^m(3^-)(3)(4^-)$ for some integer $m$ since $W_0$ is $\theta(21^-)^{\pm 1}$ admissible.

Assume that $|m|>(2k)^3$. Then the rule $\theta(12^-)$ occurs at least $(2k)^3+1$ times in the history. The rule $\theta(12^-)$ permanently changes the length of the sector  ${\cal R}^{1,r}{\cal P}^{2,\ell}$ by one multiplying it by the same letter.
On the other hand, by Lemma \ref{cycle}, the subcomputations
with step histories $\big((21^-)(1^-)(1)(2^-)(2)(21^-)\big)^{\pm 1}$ do  not change the $a$-word of length $k$ in the input sector ${\cal R}_s\cal P$. Therefore there is a transition
(not the last one), where the $a$-length of the sector ${\cal R}^{1,r}{\cal P}^{2,\ell}$
is divisible by $(2k)^3$, and so the transition to Step 3 was
possible earlier by Lemma \ref{div3} (2). Hence it suffices
to prove the lemma under the assumption  $|m|\le (2k)^3$,
and so the number of steps does not exceed $4|m|+3\le 32k^3+3$
by Lemma \ref{474}.

The length of the history is linearly bounded for each of the steps in terms of the lengths of the configurations
(Lemma \ref{gen}  for Steps 4 and 5, Lemma \ref{w} for Steps 1, 2 and 3, Lemmas \ref{121}(2) and \ref{prim} (1) for control steps.) The lengths
of these configurations are linearly bounded in Lemma \ref{standard} (2). Taking into account that the number of steps is linearly bounded in terms
of $m$ and $|m|\le (2k)^3$, we come to the desired inequality.\endproof

\medskip

We will consider a suitable function function $f(n)$ and the functions $g(n)$ and $F(n)$ from Definition \ref{sui} under the assumption that $s=2$ in that definition. (The inequality $s\ge 3$
will appear in the last Subsection \ref{sc}.) So the recognizing
Turing S-machine ${\bf M}_0$ is taken from that definition and ${\bf M}_1$ is given by Lemma \ref{S}.

Without loss of generality, one may assume that the values of
$f(n)$ are greater than some constant, for instance, for every
integer $n\ge 1$, we have

\begin{equation}\label{5}
f(n)\ge 5.
\end{equation}

Besides, it is convenient to enlarge the domain of those functions assuming
that they are defined for every real $x\ge 0$. One may assume that  $f(x)$ is still positive for $x>0$ and non-decreasing.

\begin{lemma}\label{d-x} For every $x>0$ and $d\in [0,x)$, we have
$F(x)-F(x-d)\ge \frac dx F(x)=dxg(x)$.
\end{lemma}
\proof Note that $F(x)/x^2\ge F(x-d)/(x-d)^2$ since the function $g(x)$
is non-decreasing. Hence
$$F(x)-F(x-d)\ge F(x)\big(1-\frac{(x-d)^2}{x^2}\big)\ge F(x)\frac{2dx-d^2}{x^2}\ge F(x)\frac{dx}{x^2}=\frac dx F(x).$$ \endproof

\begin{lemma} \label{k3} Let
${\cal C}: \; W_0\to\dots\to W_t$ be a
 computation  of $\bf M$ with block history $(F)$ and the step history
 of this computation contains a subword $\big((21^-)(1^-)(1)(2^-)(2)(21^-)\big)^{\pm 1}$. Let
 $k$ be the $a$-length of the input sector ${\cal R}_s{\cal P}$ of a $\theta(2^-1)^{\pm 1}$-admissible configuration $W_{j}$
 in $\cal C$. Then (a) $k=f(n)$ for some $n\ge 1$, (b) $k=O(f(||W_0||)$, and (c) $k^3=O(||W_0||)$.
\end{lemma}

\proof Claim (c) follows from (b) since $f(n)^3=O(n)$.

To prove (b) statement, it suffices to prove that $k=O(f(r) )$, where $r$ is the
length of a big historical sector of $W_0$. We have $r'\le r$,
where $r'$ is the
length of a big historical sector of $W_j$ by Lemma \ref{histF};
indeed, the computation $W_j\to\dots\to W_0$ cannot decrease
the length of it by Lemma \ref{Hprim} (a) applied to the control
steps and Lemma \ref{w} applied to Steps 1 and 2. (Step 3 does
not change this sector.) For the same reason $r''\le r'$, and $k= k''$
where $r''$ is the
length of a big historical sector of the first configuration $D_0$
of a subcomputation $\cal D$ of ${\cal C}^{\pm 1}$ with the step history $(1^-1)(1)(12^-)$ and $k''$ is the $a$-length of the input sector ${\cal R}_s{\cal P}$ of $D_0$. Hence it suffices to prove that $k''=O(f(r''))$.

The subcomputation of Step $1$ (restricted to the base of ${\bf M}_3$) is actually the
computation of ${\bf M}_2$. If we ignore the historical sectors,
we have the accepting computation of ${\bf M}_1$ with input sector
$a^{k''}$. By the definition of ${\bf M}_0$ and ${\bf M}_1$,
$k''=f(n)$ for some $n\ge 1$, the accepting computation of ${\bf M}_0$
has length $\Theta(n^{1/3})$, and so the number of rules at
Step 1 is $\Theta(n)$ by Lemma \ref{S}. Hence $r''\ge\Theta(n)$, because   $D_0$ contains the history of Step $1$ written in big historical sectors. Hence
$k''=f(n) = O(f(r''))$. Since $k=k''$, Property (a) is obtained
as well.
\endproof

\begin{lemma} \label{sta} Let $W_0$ be an accepted word and ${\cal C}: \; W_0\to\dots\to W_t$ be a
 reduced  computation  of $\bf M$ with block history $(E)$, where the computation either length-non-increasing or length-non-decreasing,  or $(F)(E)$, where the block $(E)$ is a length-non-increasing subcomputation, or $(F)$.
 Then

(a) $||W_j||\le c_4\max(||W_0||, ||W_t||)$ ($j=0,\dots,t$), where $c_4=c_4({\bf M})$, or

(b) there are accepting computations for $W_0$ and $W_t$ with block histories $(E)$ or $(F)(E)$ and histories
$H_0$ and $H_t$ such that $||H_0||+||H_t||<t/100$,
the history $H$ of $\cal C$ has a factorization $H=H(1)H(2)H(3)$,
where $||H(1)||, ||H(2)||< t/100$, $H(2)$ is of type $\big((21^-)(1^-)(1)(2^-)(2)(21^-)\big)^{m}$ with $|m|\ge c_3$
and the lengths of the subhistories of $H(2)$ of type $\big((21^-)(1^-)(1)(2^-)(2)(21^-)\big)^{\pm 1}$ are less than $t/10$.
\end{lemma}

 \proof If ${\cal C}$ is of type $(E)$, then Property (a)
 follows since ${\cal C}$
 either non-increases or non-decreases the lengths of configurations.

 Let ${\cal C}$ have type $(F)(E)$. Again we obtain Property (a)
 by Lemma \ref{standard} (2) if the number of steps in block $(F)$ is less than $10$.

 Otherwise, by Lemmas \ref{474} and \ref{histF},  block $(F)$ contains
 a subword \\$\big((21^-)(1^-)(1)(2^-)(2)(21^-)\big)^{\pm 1}$ in the step history, and by Lemma \ref{cycle}, the
 corresponding computations ${\cal C}_i$ just multiply the words in the sector
 ${\cal R}^{1,r}{\cal P}^{2,\ell}$  (and in its mirror copy) by a letter
 $b$ or by $b^{-1}$ depending on the sign of the exponent, while the length
 $k$ of the word in the sector ${\cal R}_s\cal P$ is not changed.

 The history corresponding to the block $(F)$  is $H_1H_2H_3$, where the length of the step history at most $4$ for $H_1$ and $H_3$,
 $H_2$ has the form $\big((21^-)(1^-)(1)(2^-)(2)(21^-)\big)^m$, and $H_2$
 is the history of the subcomputation  $W_{s_0}\to \dots \to W_{s_1}\to\dots\to  W_{s_m} $ starting with $\theta(21^-)^{\pm 1}$, and every $W_{s_i}$ is a $\theta(21^-)^{\pm 1}$-admissible word.
 If $w=\max(||W_{s_0}||, ||W_{s_m}||),$ then by Lemma \ref{cycle},
 \begin{equation}\label{sc}
 \max(\sum_{i=0}^{m-1} ||W_{s_i}||, \sum_{i=1}^{m} ||W_{s_i}||)\ge w+|w-2|+\dots +|w-2|m-1||\ge |m|w/2.
 \end{equation}

 Hence the sum of the lengths of the subcomputations ${\cal C}_0,\dots, {\cal C}_{m-1}$ is at least $|m|w$ by Lemma \ref{prim} for the
 Steps $1^-$ and $2^-$ of these subcomputations. Therefore $||H_2||> |m|w$.

 Now if $\max(||W_0||, ||W_t||)=w'\ge w$, then the length of every $W_j$
 with $j\le s_0$ or $j\ge s_m$ is bounded by $c_1(w'+10)$ by Lemma
 \ref{standard} (2), and the same estimate works if $W_j$ is a configuration of some subcomputation ${\cal C}_i$. Hence we have
 the inequality of item (a) of the lemma.

 So we assume now that $w'<w$.
Let $k$ be the $a$-length of the input sector ${\cal R}_s{\cal P}$ of the word $W_{s_0}$.
 If $m\ge c_3^2k^3$, then $t\ge ||H_2||> c_3^2k^3w$. We also have that each of the subcomputations ${\cal C}_i$ with step
 histories $\big((21^-)(1^-)(1)(2^-)(2)(21^-)\big)^{\pm 1}$ have length
 less than $t/10$. This follows from the property, that the difference of the lengths of ${\cal C}_i$ and ${\cal C}_{i+1}$ is at most $4$ (the computational time for control step can change by 2 by Lemma \ref{cycle},
 but the number of such subcomputations is at least $c_3^2>1000$).

Let $W_{s_0}\equiv V_0\to \dots\to V_d$ be a shortest accepting computation
for $W_{s_0}$ with a history $H_0$. To estimate $d$ from above,
we may assume by Lemma \ref{standard} (1) that
its block history is either $(E)$ or $(F)(E)$. The step history of block $(E)$ has length at most $3$ by Lemma \ref{474}.

If the number of steps in $H_0$ is at most $10$, then $||H_0||=O(||W_{0}||)$
by Lemma \ref{boundH}. Otherwise by Lemma \ref{histF}, $H_0=H'H''$, where $H'$ has step
history \\ $\big((21^-)(1^-)(1)(2^-)(2)(21^-)\big)^{\pm 1}$ and $H'H''$ starts
with a $\theta(21^-)^{\pm 1}$-admissible configuration $W_{s_0}$. By Lemma \ref{stand},
\begin{equation}\label{Ow}
||H'H''||\le c_3(k^3+1)(||W_{s_0}||+k^3)\le
c_3(k^3+1)(w+k^3)=c_3(k^3+1)O(w)
\end{equation}
since $k^3=O(w)$ by Lemma \ref{k3}.

Thus by (\ref{Ow}), $||H'H''||\le c_3^2k^3w/400\le ||H_2||/400\le (s-r)/400$ if $c_3$ is large enough.
To estimate the length of the shortest accepting computation
for $W_{0}$, it remains to estimate $s_0=||H_1||$, but this value does not exceed $c_2w\le c_3^2k^3w/800\le t/400$ by Lemma \ref{boundH}.
Therefore $||H_0||\le t/400+t/400=t/200\le ||H_2||/200$. Similarly, we
obtain $||H_t||\le t/200\le ||H_2||/200$.
Now define $H(1)=H_1, H(2)=H_2$ and $H(3)=H_3\bar H$, where $\bar H$ is of type $(E)$ and so has at most $3$ steps, and $H(3)$ has less than $10$ steps. To obtain Property (2) of the lemma if $m\ge c_3^2k^3$, it remains to estimate $H(3)$. Indeed, by Lemmas \ref{standard} (2) and \ref{boundH}, we have
$||H(3)||= ||H_3||+ ||\bar H||\le c_2w +3w \le t/100$.

Assume now that $|m|< c_3^2k^3$. Then the number of steps $s$ in the
subcomputation ${\cal C}$ does not exceed $10c_3^2k^3+20$.
Therefore for every configuration $W_j$ of ${\cal C}$, we obtain
from Lemmas \ref{standard} (2) and \ref{k3} that $||W_j||\le c_1(w'+ 10c_3^2k^3+20)=c_1c_3^2 O(w')\le c_4w',$ if $c_4$ is big enough, and we have Property (a).

The same argument works if the block history is just $(F)$.

\endproof

\begin{lemma}\label{EFE}
Let $W_0$ be an accepted word and ${\cal C}: \; W_0\to\dots\to W_t$ be a
 reduced  computation of $\bf M$ with
 block history of the form $(E)(F)...(F)(E)$, where  the first (the last) block $(E)$ is a length-non-decreasing (resp., length-non-increasing) subcomputation.
 Then
there are accepting computations for $W_0$ and $W_t$ with block histories $(E)$ or $(F)(E)$ and histories
$H_0$ and $H_t$ such that $||H_0||+||H_t||< t/100$.
\end{lemma}

\proof (1) Since the word $W_0$ is $\theta$-admissible for a rule $\theta$
of block $E$, we have $||H_0||\le 3||W_0||$ by Lemma \ref{standard} (1),
and $||W_0||\le ||W_s||$, where ${\cal C'}: W_s\to\dots\to W_r$ is the subcomputation corresponding to the first occurrence of $(F)$ in the
block history of $\cal C$, because the rules of the first block $(E)$ does not decrease the lengths.

Let $k$ and $m$ be the parameters of ${\cal C'}$ defined as in the proof
of Lemma \ref{sta}. Note that $m\ne 0$ by Lemmas \ref{histF} (1) and \ref{121} (1) since $\cal C'$ starts with $\theta(4^-4)^{-1}$.
Due to the maximality of
${\cal C}'$, one can apply Lemmas \ref{histF} (2) and \ref{k3} (a) to ${\cal C}'$
and obtain
$|m|\ge (2k)^3\ge 1000$ since the length
$k$ of any word accepted by ${\bf M}_1$ and ${\bf M}_2$ is at least $5$ by (\ref{5}). Then, as in Lemma
\ref{sta}, we obtain that the length of $\cal C'$ is at
least $2||W_s||+2|||W_s||-2|+\dots +2|||W_s||-999|$, which is at least $1000||W_s||$ and at least $1000||W_0||>300||H_0||$.

Similarly we have that the length of $\cal C''$ is greater than
$300||H_t||$, where the computation $\cal C''$ corresponds to the
last occurrence of $(F)$ in the block history of $\cal C$. It
follows that $||H_0||+||H_t||<t/100$.
\endproof

\begin{lemma}\label{form}
Let $W_0$ be an accepted word,  ${\cal C}: \; W_0\to\dots\to W_t$ be a
 reduced  computation of $\bf M$ and $H_0$, $H_t$ be the histories of the shortest computations accepting $W_0$ and $W_t$, respectively.
 Then

 (1) $\cal C$ is a product of
 at most three subcomputations ${\cal C}_1: \; W_0\to\dots\to W_{n_1}$,
 ${\cal C}_2: \; W_{n_1}\to\dots\to W_{n_1+n_2}$ and ${\cal C}_3: \; W_{n_1+n_2}\to\dots\to W_{n_1+n_2+n_3}$ ($n_1+n_2+n_3=t$), where $\max(||W_{n_1}||,||W_{n_1+n_2}||) \le \max(||W_0||, ||W_t||)$ and for every ${\cal C}_i$ ($i=1,2,3$) either

(a) $||W_j||\le c_4\max(||W_0||, ||W_t||)$, for every configuration $W_j$ of ${\cal C}_i$, where $c_4=c_4(\bf M)$ or

(b) there are accepting computations for the first and the last configuration of ${\cal C}_i$ with block histories $(E)$ or $((F)(E)$ and histories $H'_i$ and $H''_i$ such that $||H'_i||+||H''_i||<n_i/100$.

(2) The sum of lengths of all maximal subcomputations of ${\cal C}$ with block history $(E)$ does not exceed $3(||W||_0+||W_t||)+t/100$.
\end{lemma}

 \proof (1) If the block history of $\cal C$ is $(F)$, then the entire computation $\cal C$ satisfies either (a) or (b) by Lemma \ref{sta}. If the block history of $\cal C$ is $(E)(F)$, we
 consider the subcomputation ${\cal C}': W_0\to\dots\to W_s$ corresponding to the first block $(E)$ of the step history. If
 $W_{n_1}$ has minimal length in ${\cal C}'$, then the subcomputation
 ${\cal C''}: W_{n_1}\to\dots\to W_s$ is length-non-decreasing; this follows
 from Lemma \ref{gen} for Steps $4$ and $5$ and from Lemma \ref{Hprim}
 for Step $4^-$. Similarly, the subcomputation $ W_0\to\dots\to W_{n_1}$ is length-non-increasing. Since $||W_{n_1}||\le ||W_0||$ we can again apply Lemma \ref{sta} to
 the subcomputations $W_0\to\dots\to W_{n_1}$ and $W_{n_1}\to\dots\to W_t$.

 Therefore we assume that there are at least three blocks.

 Consider the subcomputation ${\cal C'}: W_r\to\dots\to W_s$ corresponding to the first block $(E)$ occured in the step history.
 It has a maximal length-non-increasing part $W_r\to\dots\to W_{n_1}$ as in the previous paragraph.

 Observe that $||W_r||\le ||W_0||$. Indeed, only historical sectors
 can be unlocked in $W_r$, but neither control S-machine can increase
 the lengths of these sectors in the computation $W_r\to\dots\to W_0$
 with block history $(F)$ by Lemma \ref{prim}, nor the computations of Steps 1,2,3 can do this by Lemma \ref{w}.

 Since $||W_{n_1}||\le ||W_r||\le ||W_0||$, we can apply Lemma \ref{sta}
 to the subcomputation ${\cal C}_1:\; W_0\to\dots\to W_{n_1}$ and obtain one of the
 properties (a) or (b) for it.
 Similarly, we consider the last block $(E)$ in the block history of $\cal C$ and define the subcomputation ${\cal C}_3$ starting
 with $W_{n_1+n_2}$ whose length does not exceed the length of $W_t$,
 and so either (a) or (b) holds for ${\cal C}_3$. We have $n_2=0$
 if the block history is $(F)(E)(F)$.

 If there are at least two blocks $(E)$, then the middle computation ${\cal C}_2$ satisfies the assumptions of Lemma \ref{EFE}, and
 so Property (b) holds for it.

 (2) The subdivision of each subcomputation  corresponding to
a $(E)$ according the sample of part (1) gives the required estimate. Namely, if
a subcomputation $\cal D$ is a product ${\cal D}_1 {\cal D}_2{\cal D}_3$
with block history  $(E)(F)(E)$ satisfying Lemma \ref{EFE}, then we obtain that
the length of ${\cal D}_1$ plus the length of ${\cal D}_3$ is less than
$0.01$ of the length of ${\cal D}_2$.  If $\cal D$ has brief history $(E)$ or $(F)(E)$,
or $(E)(F)$, the we refer to Lemma \ref{sta}.
\endproof

\endproof

\begin{lemma}\label{gtime} For every accepted word $W_0$ of length at most $n$ there is an accepting computation of length $O(nf(n)^3)$ with number of steps $O(f(n)^3)$. The generalized time function $T'(n)$ of $\bf M$ is equivalent to $\Theta(nf(n)^3)$.
\end{lemma}

\proof By Lemma \ref{standard} (1),  given an accepted word $W_0$
of length $n$, there is a shortest accepting computation
$W_0\to \dots\to W_t$ with block history either $(E)$ or $(F)(E)$ . We denote by $H$ its history. The step history of block $(E)$ has length at most $3||W_0||$ by Lemma \ref{474} and contains at most three steps.

If the number of steps in $H$ is at most $10$, then $||H||=O(||W_0||)=O(n)$
by Lemma \ref{boundH}. Otherwise by Lemma \ref{histF}, $H=H_1H_2H_3$, where $H_1$ has less than $10$ steps and $H_2$ has step
history $\big((1^-)(1)(2^-)(2))^{\pm 1}$ and $H_2H_3$ starts
with a $\theta(21^-)^{\pm 1}$-admissible configuration $W_j$. By Lemma \ref{stand}, $||H_2H_3||\le c_3(k^3+1)(||W_j||+k^3)$,
where $k$ is the length of the sector ${\cal R}_s\cal P$ in $W_j$,
and the number of steps in $H_2H_3$ is $O(k^3)$ (and so the number of steps in the entire $H$ is $O(k^3)$).  Here $||W_j||\le c_1(|W_0|_a+10)=O(n)$ by Lemma
\ref{standard} (2).

Since $k=O(f(||W_j||)=O(n)$ and $k^3=O(||W_j||)$ by Lemma \ref{k3},
we have $O(k^3)=O(f(n)^3)$ for the number of steps and $$||H_2H_3||\le c_3(k^3+1)(||W_j||+k^3)= c_3(f(n)^3)O(n)=O(n)g(n))$$
by the definition of the functions $f(n)$ and $g(n)$.

The length of each of each one-step subhistory of $H_1$ is bounded by $4c_1(|W_0|_a+10)$ (use Lemma \ref{Hprim} (b) for
Steps $1^-,2^-,3^-,4^-,5^-$, Lemma \ref{gen1} (b) for Steps
1, 2 and 3, and Lemma \ref{gen} for Steps 4 and 5).
Hence the length of the whole history $H$ is also $O(n)g(n)$,
as required.

To bound $T'(n)$ from below, we will  construct a series of
accepted words $V(n)$ of length $\Theta(n)$. The base of every $V(n)$ is standard,
and  $V(n)$ is $\theta(2^-1)$-acceptable. The input sector ${\cal R}_s\cal P$ contains
$a^k$, where $k=f(n)>0$, the word in the sector ${\cal R}^{2,\ell}{\cal P}^{3,\ell}$ is
$b^{l}$, where $l=\Theta(n)> 8k^3$ and $l$ congruent to $4k^3$
modulo $8k^3$. (There is such $l$ since $k^3=f(n)^3=O(n)$.)
Each of the big  historical sectors of $V(n)$ contains
the history of an accepting computation for ${\bf M}_2$, written in the  alphabets $X_{i,\ell}$. The length of this history is $O(n)$ by the definition of the suitable function $f(n)$, the definition of the machines ${\bf M}_0-{\bf M}_2$ and Lemma \ref{S}). Each of the small historical  sectors contains the history of the computation of ${\bf D}_5$ (also in left alphabets) that checks that
$l-4k^3$ is divisible by $8k^3$. Since $l=\Theta(n)$, this history has length $O(n)$ by Lemma \ref{div3}.
Thus, we have $||V_n||=\Theta(n)$.

Every word $V(n)$ is accepted. Indeed, the rules of Step $1^-$ can
check all the sectors since the base is standard. Then the rules of
Step 1 can accept $f(n)$, $\theta(12^-)$ replaces $b^{l}$ with
$b^{l-1}$. The rules of Step $2^-$ check the sectors again, the history of Step 2 copies the inverse history of Step 1, it restores
the alphabets $X_{i,\ell}$ in big historical sectors. Then we repeat the cycle  decreasing the exponent at $b$ by one again.
After $4k^3$ such cycles we obtain $b^{l-4k^3}$ in the sector ${\cal R}^{2,\ell}{\cal P}^{3,\ell}$, where $l-4k^3$ is divisible by $8k^3$, and therefore after Step $3^-$,  Step 3 can complete its  work by Lemma \ref{div3} (b). It remains to erase all tape letters using the rules of block $(E)$ and stop computing after the rule $\theta_0$ is applied.

Now let us estimate from below the length of arbitrary (reduced) computation $V(n)\equiv W_0\to\dots\to W_t$. By Lemma \ref{histF}, we have the block ${F}$
of the form \\$\big((1^-)(1)(2^-)(2)\big)^m(3^-)(3)(4^-)$ in the history, where $m\equiv 4k^3$ ($mod \;8k^3$) by Lemma \ref{histF} (2). Hence $|m|\ge 4k^3=4f(n)^3$.

The history of every subcomputation with step history $\big((1^-)(1)(2^-)(2)\big)^{\pm 1}$ has length at least $\Theta(n)$ for the following reason. Every configuration of it has
a word in the sector ${\cal R}^{2,\ell}{\cal P}^{3,\ell}$ of
length $\Theta(n)$ since this length belongs to the segment
$[l -4k^3, l]$. So by Lemma
\ref{prim}, one needs $\Theta(n)$ rules to check
this sector at the control steps $1^-$ and $2^-$.

Therefore the length of the computation $W_0\to\dots\to W_t$ is at least $4k^3\Theta(n)=\Theta(nf(n)^3)$, as desired.
Since we obtain the required lower bound for every $n$ and $||V_n||=\Theta(n)$, the lemma is proved.
\endproof

\begin{rk} \label{ar} A subcomputation with step history
$\big((21^-)(1^-)(1)(2^-)(2)(21^-)\big)^{\pm 1}$ does not change
the length of the sector ${\cal R}^{2,\ell}{\cal P}^{3,\ell}$
by Lemma \ref{Hprim} applied to steps $(1^-)$ and $(2^-)$. Hence
we have the same property for computations $\big((1^-)(1)(2^-)(2)\big)^m$ starting and ending with connecting rules.
Thus, above we obtained $\Theta(nf(n)^3)$ configurations of
length at least $\Theta(n)$ for any computation accepting the word $V(n)$.
\end{rk}

We call a base $B$ of a reduced computation (and the computation itself) \label{revolv} {\it revolving} if
 $B\equiv x v x$ for some letter $x$ and a word $v$, and $B$  has no proper
 subword of this form.

If $v\equiv v_1zv_2$ for some letter $z$, then the word
$zv_2xv_1z$ is also revolving. One can cyclically permute the
sectors of revolving computation with base $xvx$ and obtain
a uniquely defined computation with the base $zv_2xv_1z$,
which is called a cyclic permutation of the original
computation. The history and lengths of configurations do not change
when one cyclically permutes a computation.

\begin{lemma}\label{narrow} There is a constant $c_4=c_4(\bf M)$ such that following holds. For any
computation ${\cal C}: W_0\to\dots\to W_t$ of $\bf M$ with a revolving base  $xvx$  either

(1) we have inequality $||W_j||\le c_4\max(||W_0||, ||W_t||)$, for every word $W_j$ of ${\cal C}$, where $c_4=c_4(\bf M)$, or

(2) we have the following  properties:

(a) the word $xv$ is a cyclic permutation of the standard base $B=B(\bf M)$ or of $B^{-1}$ and

(b) the corresponding cyclic permutations $W'_0$ and $W'_t$ of the words $W_0$ and $W_t$ are  accepted words, and

(c) the step history of $\cal C$ (or of the inverse computation) contains  subwords $(21^-)(1^-)(1^-1)$ and
 $(12^-)(2^-)(2^-2)$,
 and

(d) $\cal C$ is a product of
 at most three subcomputations ${\cal C}_1: \; W'_0\to\dots\to W'_{n_1}$,
 ${\cal C}_2: \; W'_{n_1}\to\dots\to W'_{n_1+n_2}$ and ${\cal C}_3: \; W'_{n_1+n_2}\to\dots\to W'_{n_1+n_2+n_3}$ ($n_1+n_2+n_3=t$), where $\max(||W'_{n_1}||,||W'_{n_1+n_2}||)\le \max(||W_0||,||W_t||)$ and for each ${\cal C}_i$, either

 (d1) $||W'_j||\le c_4\max(||W'_0||, ||W'_t||)$, for every configuration $W_j$ of ${\cal C}_i$ or

(d2) there are accepting computations for the first and the last configuration of ${\cal C}_i$ with histories $H'_i$ and $H''_i$ such that $||H'_i||+||H''_i||<n_i$ and the corresponding block histories are either $(E)$ or $(F)(E)$.
\end{lemma}

\proof If the computation is faulty, then Property (1) is
given by Lemma \ref{nonst} since $c_4>C$. If it is non-faulty, then we have
all sectors of the base in the same order as in the standard base (or its inverse), and we obtain Property (2a). Therefore we may assume now that the base $xv$ is standard and Property (1) does not hold.

If the block history of $\cal C$ is $(E)$, we obtain a contradiction
with Lemma \ref{E} since $c_4>1$.

If the computation has only one step of type $(F)$,
then Property (1) follows from Lemmas \ref{Hprim}, \ref{9} and \ref{9D},
a contradiction again. So there is a connecting rule $\theta$ from block $(F)$ in the history.

Assume there is a block $(F)$ in the block history of $\cal C$,
and this block has at least $8$ steps. Then by Lemma \ref{histF},
the step history of $\cal C$ has a subword $\big((21^-)(1^-)(1)(2^-)(2)(21^-)\big)^{\pm 1}$, and
Property $(c)$ follows.
Moreover the words at the big (small) history sectors are copies of the same word since the subcomputations of Step $(1^-)$ (or $2^-$) have simultaneously controlled these sectors. Therefore after a number of such cycles one can obtain the length of the sector ${\cal R}^{1,r}{\cal P}^{2,\ell}$ divisible by $8k^3$ (where $k$ is the length of the sector ${\cal R}_s\cal P$), which by Lemma \ref{histF} (2), makes possible to accept
after the Steps $3^-, 3, 4^-, 4, 5$. So
one obtains Properties (a), (b) and (c). Then Property (d) follows from
Lemma \ref{form}.

If there are no such blocks $(F)$, then there are no subwords
$(E)(F)(E)$ in the block history by Lemmas \ref{histF} and \ref{121}.
Hence the block history is $(F)(E)(F)$ or a subword of this word.
Let configuration $W'_r$ and $W'_s$ subdivide $\cal C$ in single block
computations. Then $||W'_r||<c_3||W'_0||$, because there are at most $7$ steps
in the subcomputation $W'_r\to\dots\to W'_0$, and each step transition
from $W'_j$ towards $W'_0$ can multiply the length by at most $c$.
(See Lemma \ref{Hprim} for control steps and Lemma \ref{WV} for
Steps $1,2$ and $3$.) Analogously, we have $||W'_s||<c_3||W'_0||$.
Since for every step the lengths of all
configurations are linearly bounded in terms the first and the last
configurations (see Lemmas \ref{9}, \ref{9D} \ref{Hprim} (a))
we have $||W'_j||\le c_3\max(||W'_0||,||W'_t||)$ if $j\le r$ or $j\ge s$.
So to obtain Property (1) (and a contradiction), it suffices to linearly bound the configurations in the subcomputation $W'_r\to\dots\to W'_s$
in terms of $\max(||W'_r||,||W'_s||$. This is done in Lemma \ref{E} (1).
Thus, the proof is complete.
\endproof

\subsection{Two more properties of standard computations}\label{long}

Here we prove two lemmas needed for the estimates in Subsection \ref{ub}. The first one says (due to Lemma \ref{121} (2)) that if a standard computation $\cal C$
is very long in comparison with the lengths of the first and
the last configuration, then it can be completely restored
if one knows the history of $\cal C$, and the same is true for
the long subcomputations of $\cal C$. This makes the auxiliary
parameter $\sigma_{\lambda}(\Delta)$ useful for some estimates
of areas of diagrams $\Delta$. The second lemma is helpful
for the proof of Lemma \ref{led} in Section \ref{ub}.

\begin{lemma} \label{B} Let ${\cal C}: \;W_0\to\dots\to W_t$ be a reduced computation with standard base,
where $t\ge  c_5 \max(||W_0||, ||W_t||)$ for sufficiently large constant $c_5=c_5(\bf M)$. Suppose the word $W_0$ is accepted. Then any subcomputation ${\cal D}: W_r\to\dots\to W_s$ of $\cal C$ (or the inverse for $\cal D$) of  length at least $0.4t$ contains one of the words
$(21^-)(1^-)(1^-1)$, $(12^-)(2^-)(2^-2)$, $(23^-)(3^-)(3^-3)$, $(34^-)(4^-)(4^-4)$ in the step history.
\end{lemma}

\proof If the block history of $\cal C$ is $(F)$,  we refer
to Lemma \ref{sta} as follows.

Assume that Property (a) of that lemma holds. Then every step
of the computation $\cal D$ has length at most $4c_4\max(||W_0||, ||W_t||)$
by Lemma \ref{w} for Steps 1,2 and 3 and by Lemma \ref{Hprim} (b) for other steps.
Hence the number of steps in $\cal D$ has to be at least $10$ since the length
of its history is at least $0.4 c_5 \max(||W_0||, ||W_t||)$ and $c_5$ can be chosen large
enough. It follows from Lemma \ref{histF} (1) that the step history of $\cal D$
contains subwords $(21^-)(1^-)(1^-1)$ and $(12^-)(2^-)(2^-2)$, as required.

If Property (a) of  Lemma \ref{sta} fails, then by
Property (b), we have a subcomputation of length $>0.98t$
with step history $\big((1^-)(1)(2^-)(2)\big)^{\pm m}$,
where every cycle with block history $\big((1^-)(1)(2^-)(2)\big)^{\pm 1}$ has length $<t/10$. Then the subcomputation $\cal D$ of length
$\ge 0.4t$ has to contain such a cycle, and so the step history of $\cal D$ contains $(12^-)(2^-)(2^-2)$, as required. Thus, we may assume
that the block history of  $\cal C$ is not $(F)$.

If the block history of $\cal D$ contains a subword $(E)(F)(E)$,
then the statement follows from Lemma \ref{histF} (where Lemma \ref{121} eliminates the case $m=0$).
So the block history of $\cal D$ is a subword of $(F)(E)(F)$.
By Lemma \ref{form} (2), the length of the $(E)$-subcomputation of $\cal D$ is less that $t/100+3(||W_0||+ ||W_t||)\le t/50$. So one of the
$(F)$-subcomputations of $\cal D$ has length $ >(0.4-0.02)t/2=0.19t$.

{\bf Case 1.} Assume that there is a block $(E)$ in the block history of $\cal D$, and without loss of generality, we may assume that the
computation ${\cal D}$ has a subcomputation
${\cal D}': W_r\to\dots\to W_j$ of type $(F)$ with $j-r>0.19t$
and the subcomputation of type $(E)$ occurs after ${\cal D}'$ in $\cal D$. Proving by contradiction, we conclude that the step history of
${\cal D}'$ is $(4^-)$ since a longer step history
would provide us with the subword $(34^-)(4^-)(4^-4)$ in the step history of $\cal D$.

Suppose the subcomputation $W_0\to\dots\to W_r$ also has a block $(E)$. Then $\cal C$ has a subcomputation with block history $(E)(F)(E)$. Let the subcomputation $\bar{\cal C}$ correspond to
the middle block $(F)$. Then the first and the last configurations
of $\bar{\cal C}$ are admissible for some rules of type $(E)$.
Therefore by Lemma \ref{standard} (1), one can construct an
auxiliary computation $\tilde{\cal C}={\cal C}'\bar{\cal C}{\cal C}''$, where
the first factor (the third one) starts (resp. ends) with Step (5).
Then by Lemmas \ref{histF} and \ref{k3} for $\tilde{\cal C}$, the subcomputation $\bar{\cal C}$ of $\cal C$ has a subcomputation of type $(F)$
with step history $\big( (1^-)(1)(2^-)(2)\big)^{m_1}(3^-)(3)(4^-)$.  $|m_1|\ge 8k^3\ge 1000$, and the
subcomputation corresponding to the last Step $(4^-)$ has length $>0.19t$. It follows that at least $1000$ control  steps of the form $(1^-)$ or $(2^-)$ should have length $>0.19t$ since their control S-machines have to
check the big and small historical sectors too (and the length of
the historical sectors are unchanged by the rules of  $\big( (1^-)(1)(2^-)(2)\big)^{m_1}(3^-)(3)(4^-)$). We obtain a contradiction
since $1000\times 0.19t>t$.

Thus, the computation ${\cal E}:\;W_0\to\dots\to W_j$ is of type $(F)$.
Hence $||W_j||\le ||W_0||$, because only historical sectors
 can be unlocked in $W_j$, but neither control S-machine can decrease
 the lengths of these sectors in the computation $W_j\to\dots\to W_0$
 with block history $(F)$ by Lemma \ref{prim}, nor the computations of Steps 1,2,3 can do this by Lemma \ref{w}.

 If the step history of $\cal E $ ends with $(34^-)(4^-)(4^-4)$,
 then $||W_0||\ge ||W_j||\ge 0.19t /4>0.04t$ by Lemma \ref{Hprim},
 which contradicts to the assumption of the lemma. Hence the step
 history of ${\cal E}$ is $(4^-)$, and so $0.19t \le 4||W_0||$ by Lemma
\ref{Hprim} (b), a contradiction again.

{\bf Case 2.} The block history of $\cal D$ is $(F)$. Since the block
history of $\cal C$ is not $(F)$ but a subword of $(F)(E)\dots $, we
conclude without loss of generality, that $\cal C$ begins with a
maximal subcomputation ${\cal E}: W_0\to\dots\to W_u$ of type $(F)$, where $r\le u<t$. Then as in Case 1, we have $||W_u||\le ||W_0||$.

Now consider the options (a) and (b) provided by Lemma \ref{sta}
for $\cal E$.
The option (b) is eliminated exactly as in the second paragraph
of the proof of the current lemma, where $t$ can be replaced by $0.4t$
since $u\ge r-s\ge 0.4t$.  Hence we have by (a) that every
configuration $W_j$ of $\cal E$ satisfies the inequality $||W_j||\le c_4||W_0||$
since $||W_u||\le ||W_0||$.
Then the length of every single step of $\cal E$
cannot exceed $4c_4||W_0||$ (see Lemma \ref{Hprim} (b) for control
steps and \ref{gen1} for Steps 1,2 and 3). Here we have $4c_4<c_5/10$ by the choice of $c_5$.

Since $s-r \ge 0.4t\ge 0.4 c_5 ||W_0||$, we see that the length
of the step history of $\cal D$ is at least $4$. It follows from
Lemma \ref{histF} that the step history of $\cal D$ contains one of
the words mentioned in the formulation of Lemma \ref{B}.
\endproof

\begin{lemma} \label{pol} Let a reduced computation $W_0\to\dots\to W_t$ start with an accepted word $W_0$, have  standard base, and have step history of
length 1. Assume that for some index $j$, we have $|W_j|_a>3|W_0|_a$.
 Then there is a sector $QQ'$ such that
a state letter  from $Q$ or from $Q'$ inserts a letter increasing
the length of this sector after any transition of the subcomputation
$W_j\to\dots\to W_t$.
\end{lemma}

\proof
Let  the step history be $(1)$. Note that all big historical sectors of any configuration $W_i$
have the same content (up to taking a copy) since the word $W_0$ is accepted. Assume that no rule of the subcomputation ${\cal D}:\; W_0\to\dots\to W_j$ increases the length of big historical sectors. Then by Lemma \ref{gen1}
(b) the length of the history of $\cal D$ does not exceed $h$, where
$h$ is the $a$-length of such sectors.

Every rule of the subcomputation $\cal D$ can change the length of any working sector at most by $1$. (See Lemma \ref{simp} (3)). Hence if its
length in $W_0$ is $\ell$, its length in $W_j$ is at most $\ell+h$. It follows that $|W_{j}|_a\le 3 |W_0|_a$, because
the working sectors of ${\bf M}_2$ and its historical sectors alternate in the standard base. This contradicts to the assumption of the lemma.

Thus, there is a rule in the history of $\cal D$ increasing the
length of a big historical sector $QQ'$. It has to insert a letter from
$X_{i,\ell}$ from the left and a letter from $X_{i,r}$ from the right.
Since the obtained word is not a word over one of these alphabets,
Step 1 is not over, and the next rule has to increase the length
of the sector again in the same manner since the computation is reduced.
This procedure will repeat until one gets $W_t$. This proves the statement. The same proof works for Steps 2 and 3. (In the later case,
one will consider small historical sectors.)

It follows from the definition of Step 4 (of Step 5) that every
rule either increase or decrease the length of small (resp., of big)
history sectors. If any rule increases it, then all the next rules will
increase the lengths of these sectors too. Hence the argument of the previous paragraph works for Steps 4 and 5 as well.

For the control Steps $1^--5^-$, the statement of the lemma follows from Lemma \ref{prim} (1):
if we have a transition of a primitive S-machine, where the control state letter increases
the length of a sector, then it will keep increasing it in any reduced computation.

\endproof

\section{Groups and diagrams}\label{gd}

\subsection{The groups}\label{MG}

Every S-machine can be simulated by finitely presented group (see \cite{SBR} and
also \cite{OS06}, \cite{OS01}). Here we apply such a construction to the S-machine $\bf M.$
To simplify formulas, it is convenient to change the notation. From now on we shall denote by $N$ the length of
the standard base of $\bf M$.

Thus the set of state letters is $Q=\sqcup_{i=0}^{N}Q_i$ (where $Q_{N}=Q_0=\{t\}$) $Y= \sqcup_{i=1}^{N} Y_i,$ and $\Theta$ is the set of rules of the S-machine $\bf M.$

The finite set of generators of the \label{groupM} group $M$ consists of \label{qletter}{\em $q$-letters} corresponding to the states
$Q$, \label{aletter}{\em $a$-letters} corresponding to the tape letters from $Y,$ and
\label{thetal}{\em $\theta$-letters} corresponding to the rules from the positive part $\Theta^+$ of $\Theta$ (the same letter as for the S-machine).

The \label{relations} relations of the group $M$ correspond to the rules of the S-machine $\bf M$.
Recall that the cyclic S-machine $\bf M$ satisfies Property
(1) of Lemma \ref{simp}, and so every rule $\theta\in \Theta^+$ of it has the form $\theta: [U_0\to V_0,\dots U_{N}\to V_{N}]$, where $U_0\equiv U_N$ and $V_0\equiv V_N$.  For every
such rule $\theta$, we introduce the following relation
of the group $M$.
\begin{equation}\label{rel1}
U_i\theta_{i+1}=\theta_i V_i,\,\,\,\, \qquad \theta_j a=a\theta_j, \,\,\,\, i,j=0,...,N
\end{equation}
for all $a\in Y_j(\theta)$. (Here $\theta_{N}\equiv\theta_0.$)
The first type of relations will be
called \label{thetaqr} $(\theta,q)$-{\em relations}, the second type - \label{thetaar}
$(\theta,a)$-{\em relations}.

Finally, the required \label{groupG} group $G$ is given by the generators and
relations of the group $M$ and by one more additional
relation, namely the \label{hubr} {\it hub}-relation
\begin{equation}\label{rel3}
(W_M)^L=1,
\end{equation}
where
$W_M$
is  the accept word (of length $N$) of the S-machine $\bf M$ and
the exponent $L$ is a large enough integer. (It depends on $\bf M$
and will be made more precise later.) The corresponding cells in van Kampen diagrams
looks like hubs in the net of $q$-bands (see
pictures in \cite{R}, \cite{SBR}, \cite{O97}).

\subsection{Van Kampen diagrams}\label{md}

Recall that a van Kampen \label{diagram} {\it diagram} $\Delta $ over a presentation
$P=\langle A\; | \; \mathcal R\rangle$ (or just over the group $P$)
is a finite oriented connected and simply--connected planar 2--complex endowed with a
\label{Lab} labeling function $\Lab : E(\Delta )\to A^{\pm 1}$, where $E(\Delta
) $ denotes the set of oriented edges of $\Delta $, such that $\Lab
(e^{-1})\equiv \Lab (e)^{-1}$. Given a \label{cell} cell (that is a 2-cell) $\Pi $ of $\Delta $,
we denote by $\partial \Pi$ the boundary of $\Pi $; similarly, \label{partial}
$\partial \Delta $ denotes the boundary of $\Delta $. The labels of
$\partial \Pi $ and $\partial \Delta $ are defined up to cyclic
permutations. An additional requirement is that the label of any
cell $\Pi $ of $\Delta $ is equal to (a cyclic permutation of) a
word $R^{\pm 1}$, where $R\in \mathcal R$. The label and the \label{clength} combinatorial length $||\bf p||$ of
a path $\bf p$ are defined as for Cayley graphs.

The van Kampen Lemma states that a word $W$ over the alphabet $A^{\pm 1}$
represents the identity in the group $P$ if and only
if there exists a diagram $\Delta
$ over $P$ such that
$\Lab (\partial \Delta )\equiv W,$ in particular, the combinatorial perimeter $||\partial\Delta||$ of $\Delta$ equals $||W||.$
(\cite{LS}, Ch. 5, Theorem 1.1; our formulation is closer to
Lemma 11.1 of \cite{book}). The word $W$ representing $1$ in $P$ is freely equal
to a product of conjugates to the words from $R^{\pm 1}$. The minimal number
of factors in such products is called the \label{areaw} {\em area} of the word $W.$ The \label{aread}{\it area}
of a diagram $\Delta$ is the number of cells in it.
By van Kampen Lemma, $\area(W)$ is equal
to the area of a diagram having
the smallest number of cells among all diagrams with  boundary label $\Lab (\partial \Delta )\equiv W.$

We will study diagrams over the groups $M$ and $G$. The edges labeled by state
letters ( = $q$-{\it letters}) will be called \label{qedge} $q$-{\it edges}, the edges labeled by tape
letters (= $a$-{\it letters}) will be called \label{aedge} $a$-{\it edges}, and the edges labeled by
$\theta$-letters are \label{thedge} $\theta$-{\it edges}.

We denote by $|\bf p|_a$ (by $|\bf p|_{\theta}$, by
$|\bf p|_q$)
the \label{alength} $a$-{\it length} (resp., the \label{thlength} $\theta$-{\it length}, the \label{qlength} $q$-length) of a path/word $\bf p,$ i.e., the number of
$a$-edges/letters (the number of $\theta$-edges/letters, the number of $q$-edges/letters) in $\bf p.$

 The cells corresponding
to relation (\ref{rel3})  are called \label{hubs} {\it hubs}, the cells corresponding
to $(\theta,q)$-relations are called \label{tq} $(\theta,q)$-{\it cells},
and the cells are called \label{ta} $(\theta,a)$-{\it cells} if they correspond to $(\theta,a)$-relations.

A van Kampen diagram is \label{reducedd}{\em reduced}, if
it does not contain two cells (= closed $2$-cells) that have a
common edge $e$ such that the boundary labels of these two cells are
equal if one reads them starting with $e$
(if such pairs of cells exist, they can be removed to obtain a  diagram of smaller area and with the same boundary label).
To study (van Kampen) diagrams
over the group $G$ we shall use their simpler subdiagrams such as bands and trapezia, as in \cite{O97},  \cite{SBR}, \cite{BORS}, etc.
 Here we repeat one more necessary definition.

\label{band} \begin{df}Let $\cal Z$ be a subset of the set of letters in the set of generators of the group $M$. A
$\cal Z$-band $\bb$ is a sequence of cells $\pi_1,...,\pi_n$ in a reduced \vk
diagram $\Delta$ such that

\begin{itemize}
\item Every two consecutive cells $\pi_i$ and $\pi_{i+1}$ in this
sequence have a common boundary edge ${\bf e}_i$ labeled by a letter from ${\cal Z}^{\pm 1}$.
\item Each cell $\pi_i$, $i=1,...,n$ has exactly two $\cal Z$-edges in the boundary $\partial \pi_i$,
${\bf e}_{i-1}^{-1}$ and ${\bf e}_i$ (i.e. edges labeled by a letter from ${\cal Z}^{\pm 1}$) with the requirement that either
both $\Lab(e_{i-1})$ and $\Lab(e_i)$ are positive letters or both
are negative ones.

\item If $n=0$, then $\bb$ is just a $\cal Z$-edge.
\end{itemize}
\end{df}

The counter-clockwise boundary of the subdiagram formed by the
cells $\pi_1,...,\pi_n$ of $\bb$ has the factorization ${\bf e}\iv {\bf q}_1{\bf f} {\bf q}_2\iv$
where ${\bf e}={\bf e}_0$ is a $\cal Z$-edge of $\pi_1$ and ${\bf f}={\bf e}_n$ is an $\cal Z$-edge of
$\pi_n$. We call ${\bf q}_1$ the \label{bottomb}{\em bottom} of $\bb$ and ${\bf q}_2$ the
\label{topb}{\em top} of $\bb$, denoted \label{bott} $\bott(\bb)$ and \label{topp} $\topp(\bb)$.
Top/bottom paths and their inverses are also called the \label{sideb}{\em
sides} of the band. The $\cal Z$-edges ${\bf e}$
and ${\bf f}$ are called the \label{seedgesb}{\em start} and {\em end} edges of the
band. If $n\ge 1$ but ${\bf e}={\bf f},$ then the $\cal Z$-band is called a \label{annulus} $\cal Z$-{\it annulus}.

We will consider \label{qband} $q$-{\it bands}, where $\cal Z$ is one of the sets $Q_i$ of state letters
for the S-machine $\bf M$, \label{thband}
$\theta$-{\it bands} for every $\theta\in\Theta$, and \label{aband} $a$-{\it bands}, where
${\cal Z}=\{a\}\subseteq Y$.
 The convention is that $a$-bands do not
contain $(\theta,q)$-cells, and so they consist of $(\theta,a)$-cells  only.

\begin{rk} \label{tb} To construct the top (or bottom) path of a band $\cal B$, at the beginning
one can just form a product ${\bf x}_1\dots {\bf x}_n$ of the top paths ${\bf x}_i$-s of the cells $\pi_1,\dots,\pi_n$ (where each $\pi_i$ is a $\cal Z$-bands of length $1$).
No  $\theta$-letter is being canceled in the word
$W\equiv\Lab({\bf x}_1)\dots\Lab({\bf x}_n)$ if $\cal B$ is  a $q$- or $a$-band since
otherwise two neighbor cells of the band would make the diagram non-reduced. For similar reason, there are no cancellations of
$q$-letters in
$W$ if $\cal B$ is  a $\theta$-band

If $\cal B$ is a $\theta$-band then a few cancellations of $a$-letters (but not $q$-letters) are possible in $W.$ (This can happen if one of $\pi_i, \pi_{i+1}$
is a $(\theta,q)$-cell and another one is a $(\theta,a)$-cell.) We will always assume
that the top/bottom label of a $\theta$-band is a reduced form of the word $W$.
This property  is easy to achieve: by folding edges
with the same labels having the same initial vertex, one can make
the boundary label of a subdiagram in a \vk diagram reduced (e.g., see \cite{book} or
\cite{SBR}).
\end{rk}

If the path $({\bf e}\iv {\bf q}_1{\bf f})^{\pm 1} $ or the path $({\bf f} {\bf q}_2\iv {\bf e}\iv)^{\pm 1}$
 is the subpath of the boundary path of $\Delta$ then the band is called
 a \label{rimb}{\it rim} band of $\Delta.$
We shall call a $\cal Z$-band \label{maxb}{\em maximal} if it is not contained in
any other $\cal Z$-band.
Counting the number of maximal $\cal Z$-bands
 in a diagram we will not distinguish the bands with boundaries
 ${\bf e}\iv {\bf q}_1{\bf f} {\bf q}_2\iv$ and ${\bf f} {\bf q}_2\iv {\bf e}\iv {\bf q}_1,$ and
 so every $\cal Z$-edge belongs to a unique maximal $\cal Z$-band.

We say that a ${\cal Z}_1$-band and a ${\cal Z}_2$-band \label{cross}{\em cross} if
they have a common cell and ${\cal Z}_1\cap {\cal Z}_2=\emptyset.$

Sometimes we specify the types of bands as follows.
A $q$-band corresponding to one
of the letters $t$ of the base is called a \label{tband} $t$-{\it band}.

The papers \cite{O97}, \cite{BORS}, \cite{OS04} contain the proof of the
following lemma in a more general setting. (In contrast to Lemmas 6.1 \cite{O97} and
 3.11 \cite{OS04}, we have no $x$-cells here.)

\begin{lemma}\label{NoAnnul}
A reduced van Kampen diagram $\Delta$ over $M$ has no
$q$-annuli, no $\theta$-annuli, and no  $a$-annuli.
Every $\theta$-band of $\Delta$ shares at most one cell with any
$q$-band and with any $a$-band.
\end{lemma} $\Box$

If $W\equiv x_1...x_n$ is a word in an alphabet $X$, $X'$ is another
alphabet, and $\phi\colon X\to X'\cup\{1\}$ (where $1$ is the empty
word) is a map, then $\phi(W)\equiv\phi(x_1)...\phi(x_n)$ is called the
\label{projectw}{\em projection} of $W$ onto $X'$. We shall consider the
projections of words in the generators of $M$ onto
$\Theta$ (all $\theta$-letters map to the
corresponding element of $\Theta$,
all other letters map to $1$), and the projection onto the
alphabet $\{Q_0\sqcup \dots \sqcup Q_{N-1}\}$ (every
$q$-letter maps to the corresponding $Q_i$, all other
letters map to $1$).

\begin{df}\label{dfsides}
{\rm  The projection of the label
of a side of a $q$-band onto the alphabet $\Theta$ is
called the \label{historyb}{\em history} of the band. The step history of this projection
is the \label{stephb}{\it step history} of the $q$-band. The projection of the label
of a side of a $\theta$-band onto the alphabet $\{Q_0,...,Q_{N-1}\}$
is called the \label{baseb} {\em base} of the band, i.e., the base of a $\theta$-band
is equal to the base of the label of its top or bottom.}
\end{df}
As for words, we may  use representatives of
$Q_j$-s in base words.

\begin{df}\label{dftrap}
{\rm Let $\Delta$ be a reduced  diagram over $M$,
which has  boundary path of the form ${\bf p}_1\iv {\bf q}_1{\bf p}_2{\bf q}_2\iv,$ where
${\bf p}_1$ and ${\bf p}_2$ are sides of $q$-bands, and
${\bf q}_1$, ${\bf q}_2$ are maximal parts of the sides of
$\theta$-bands such that $\Lab({\bf q}_1)$, $\Lab({\bf q}_2)$ start and end
with $q$-letters.

\begin{figure}
\begin{center}
\includegraphics[width=1.0\textwidth]{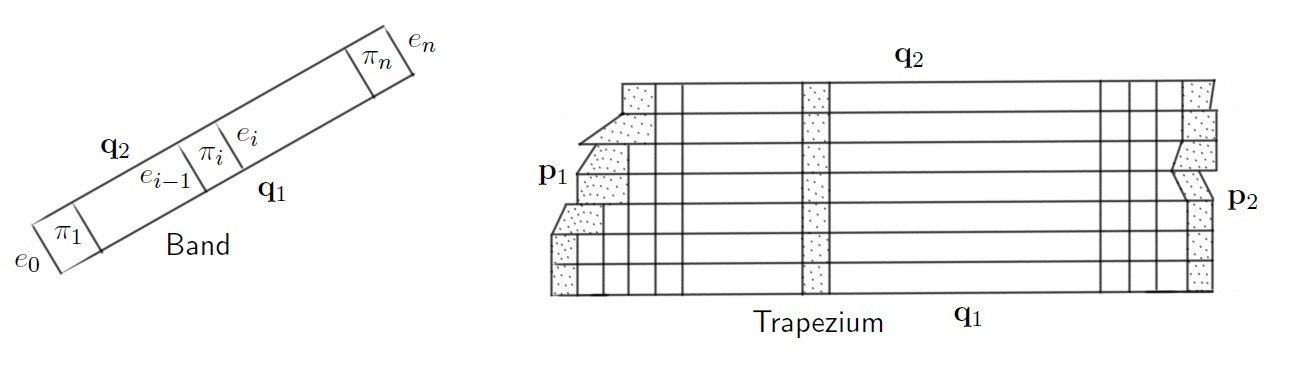}
\end{center}
\caption{Band and Trapezium}\label{bt}
\end{figure}

Then $\Delta$ is called a \label{trapez}{\em trapezium}. The path ${\bf q}_1$ is
called the \label{bottomt}{\em bottom}, the path ${\bf q}_2$ is called the \label{topt}{\em top} of
the trapezium, the paths ${\bf p}_1$ and ${\bf p}_2$ are called the \label{lrsidest}{\em left
and right sides} of the trapezium. The history (step history) of the $q$-band
whose side is ${\bf p}_2$ is called the \label{historyt}{\em history} (resp., \label{stepht} step history) of the trapezium;
the length of the history is called the \label{heightt}{\em height}  of the
trapezium. The base of $\Lab ({\bf q}_1)$ is called the \label{baset}{\em base} of the
trapezium.}
\end{df}

\begin{rk} Notice that the top (bottom) side of a
$\theta$-band $\ttt$ does not necessarily coincide with the top
(bottom) side ${\bf q}_2$ (side ${\bf q}_1$) of the corresponding trapezium of height $1$, and ${\bf q}_2$
(${\bf q}_1$) is
obtained from $\topp(\ttt)$ (resp. $\bott(\ttt)$) by trimming the
first and the last $a$-edges
if these paths start and/or end with $a$-edges.
We shall denote the
\label{trim}{\it trimmed} top and bottom sides of $\ttt$ by \label{ttopp} $\ttopp(\ttt)$ and
\label{tbott} $\tbott(\ttt)$. By definition, for arbitrary $\theta$-band $\cal T,$ $\ttopp(\cal T)$
is obtained by such a trimming only if $\cal T$ starts and/or ends with a
$(\theta,q)$-cell; otherwise $\ttopp(\cal T)=\topp(\cal T).$
The definition of $\tbott(\cal T)$ is similar.
\end{rk}

By Lemma \ref{NoAnnul}, any trapezium $\Delta$ of height $h\ge 1$
can be decomposed into $\theta$-bands $\ttt_1,...,\ttt_h$ connecting
the left and the right sides of the trapezium. The word written on
the trimmed top side of one of the bands $\ttt_i$ is the same as the
word written on the trimmed bottom side of $\ttt_{i+1}$,
$i=1,...,h$.
Moreover, the following lemma claims that every trapezium
simulates the work of $\bf M$. It summarizes the assertions  of Lemmas
6.1, 6.3,
6.9, and 6.16 from \cite{OS04}. For the formulation (1) below, it is important
that $\bf M$  is an $S$-machine. The analog of
this statement is false for Turing machines. (See \cite{OS01} for a discussion.)

\begin{lemma}\label{simul} (1) Let $\Delta$ be a trapezium
with history $\theta_1\dots\theta_d$ ($d\ge 1$).
Assume that $\Delta$
has consecutive maximal $\theta$-bands  ${\cal T}_1,\dots
{\cal T}_d$, and the words
$U_j$
and $V_j$
are the  trimmed bottom and the
trimmed top labels of ${\cal T}_j,$ ($j=1,\dots,d$).
Then the history of $\Delta$ is a reduced word, $U_j$, $V_j$ are admissible
words for $M,$ and
$$V_1\equiv U_1\cdot \theta_1, U_2\equiv V_1,\dots, U_d \equiv V_{d-1}, V_d\equiv U_d\cdot \theta_d.$$

(2) For every reduced computation $U\to\dots\to U\cdot H \equiv V$ of $\bf M$
with $||H||\ge 1$
there exists
a trapezium $\Delta$ with bottom label $U$, top label $V$, and with history $H.$
\end{lemma}

Using Lemma \ref{simul}, one can immediately derive properties of trapezia from the properties of computations obtained earlier.

If $H'\equiv \theta_i\dots\theta_j$ is a subword of the history $\theta_1\dots\theta_d$
from Lemma \ref{simul} (1), then the bands ${\cal T}_i,\dots, {\cal T}_j$ form a subtrapezium
$\Delta'$ of the trapezium $\Delta.$ This subtrapezium is uniquely defined by the
subword $H'$ (more precisely, by the occurrence of $H'$ in the word $\theta_1\dots\theta_d$), and $\Delta'$ is called the \label{H'partt} $H'$-{\it part} of $\Delta.$
\medskip

  We say that a trapezium $\Delta$ is \label{sttrap} {\it standard} if the base of $\Delta$
  is the standard base $\bf B$ of $\bf M$ or ${\bf B}^{-1}$, and the
  step history of $\Delta$ (or the inverse one) contains one of the words
  $(21^-)(1^-)(1^-1)$, $(12^-)(2^-)(2^-2)$, $(23^-)(3^-)(3^-3)$, $(34^-)(4^-)(4^-4)$, $(45^-)(5^-)(5^-5)$.

   \begin{rk}\label{restore} By Lemmas \ref{simul} and \ref{123} (2), given the history $H$, one can reconstruct the entire standard trapezium $\Delta$.
\end{rk}

\begin{df} \label{bigt} We say that a trapezium $\Gamma$ is $\it big$ if

(1) the base of $\Delta$ or the inverse word has the form $xvx$, where $xv$ a cyclic shift of the $L$-s power of the standard base;

(2) the diagram $\Gamma$ contains a standard trapezium .

\end{df}

\begin{lemma}\label{or}

Let $\Delta$ be a trapezium  whose base is
$xvx$, where $x$ occurs in $v$ exactly $L-1$ times and other letters occur $<L$ times each (where $L$ is as in (\ref{rel3})). Then either $\Delta$ is big or the length of a side of every $\theta$-band of $\Delta$ does not exceed $c_5(||W||+||W'||)$, where $W, W'$ are the labels of its top and bottom, respectively.
\end{lemma}

\proof The diagram $\Delta$ is covered by $L$  subtrapezia $\Gamma_i$ ($i=1,\dots,L$) with bases $xv_ix$.

Assume that the the step history of $\Delta$ (or inverse step history) contains one of the subwords $(21^-)(1^-)(1^-1)$, $(12^-)(2^-)(2^-2)$, $(23^-)(3^-)(3^-3)$, $(34^-)(4^-)(4^-4)$. Then
by Lemma \ref{121} (2) (and \ref{simul}), the base of $\Delta$ has the form $(xu)^Lx$,
where $xu$ is a cyclic shift of the standard base (or the inverse one) and the diagrams $\Gamma_i$-s ($i=1,\dots,L$) are just copies of each other.
Since $\Delta$ contains a standard subtrapezia, it is big.

Now, under the assumption that the step history has no subwords mentioned in the previous paragraph, it suffices to bound the
the length of a side of every $\theta$-band of arbitrary $\Gamma_i$ by $ \le c_4(||V|_a+||V'||)$, where $V$ and $V'$ are the labels of the top and the bottom of $\Gamma_i$.

Assume that the word $xv_ix$ has a proper subword $yuy$, where
$u$ has no letters $y$, and any other letter occurs in $u$ at most once.
Then the word $yuy$ is faulty since $v_i$ has no letters $x$. By Lemma
\ref{nonst}, we have $|U_j|_a\le C\max(|U_0|_a,|U_t|_a)$ for every
configuration $U_j$ of the computation given by Lemma \ref{simul} restricted to the base $yuy$. Since $c_4>C$, it suffices to obtain
the desired estimate for the computation whose base is obtained
by deleting the subword $yu$ from $xv_ix$. Hence inducting on the
length of the base of $\Gamma_i$, one may assume that it has no proper
subwords $yuy$, and so the base of $\Gamma_i$ is revolving.
Now the required upper estimate for $\Gamma_i$ follows from Lemma \ref{narrow} (see (1) and (2c) there).
\endproof

\label{papam} \subsection{Parameters} \label{param}

  The following constants will be used for the proofs in this paper.

$$\lambda^{-1},\;  N <<c_0<<\dots<< c_5<<L_0<<L<<K<<J<<\delta^{-1}<< c_6<<c_7<<$$
\begin{equation}\label{const}
N_1<<N_2<<N_3<<N_4.
  \end{equation}

  For each of the inequalities of this paper, one can find the highest
  constant (with respect to the order $<<$) involved in the inequality
  and see that for fixed lower constants, the inequality is correct as soon
  as the value of the highest one is sufficiently large. This principle
  makes the system of all inequalities used in this paper consistent.

\section{Diagrams without hubs}

\subsection{A modified length function} \label{lf}
\label{length}

Let us modify the length function on the words and paths. The
standard length of a word (a path) will be called the {\em
combinatorial length} of it. From now on we use the word {\em `length'} for
the modified length.
We set the length of every
$q$-letter equal 1, and the length of every $a$-letter equal to a small
enough number $\delta$ so that

\begin{equation}
J\delta<1 \label{param}.
\end{equation}

We also set to 1 the length of every word of length $\le 2$ which
contains exactly one $\theta$-letter and no $q$-letters (such
words are called $(\theta,a)$-{\em syllables}). The length of a
decomposition of an arbitrary word in a product of letters and
$(\theta,a)$-syllables is the sum of the lengths of the factors.
\label{modiflf} {\em The length $|w|$ of a word} $w$ is the smallest length of such
decompositions. {\em The length} $|\bf p|$ {\em of a path} in a diagram is the
length of its label. The {\em perimeter} $|\partial\Delta|$ of a
van Kampen diagram is similarly defined by cyclic decompositions
of the boundary $\partial\Delta$.

The next statement follows from the above definitions and
from the property of $(\theta,q)$-relations and their cyclic shifts: The
subword between two $q$-letters (between $\theta$-letters) in arbitrary $(\theta,q)$-relation is a syllable (has at most one $q$-letter and at most two $a$-letters).

\begin{lemma} \label{ochev}
Let ${\bf s}$ be a path in a diagram $\Delta$ having $c$ $\theta$-edges
and $d$ $a$-edges. Then

(a) $ {\bf |s|}\ge \max(c, c+(d-c)\delta)$;

(b) ${\bf |s|}=c$ if ${\bf s}$ is a top or a bottom of a $q$-band.

(c) For any product ${\bf s=s_1s_2}$ of two paths in a  diagram, we
have

\begin{equation}
{\bf |s_1|+|s_2|\ge |s|\ge |s_1|+|s_2|}-\delta. \label{delta}
\end{equation}

(d) Let $\cal T$ be a $\theta$-band with base
of length $l_b$.  Let $l_a$ be the number of $a$-edges in the top
path ${\bf topp}(\cal T)$. Then the length of $\cal T$ (i.e., the number of cells in $\cal T$) is between
$l_a-l_b$ and $l_a+3l_b$.

\end{lemma} $\Box$

\begin{lemma} \label{rim} Let $\Delta$ be a van Kampen diagram
whose rim $\theta$-band $\ttt$ has base with at most $K$ letters. Denote by
$\Delta'$ the subdiagram $\Delta\backslash \ttt$. Then
$|\partial\Delta|-|\partial\Delta'|>1$.
\end{lemma}

\proof Let ${\bf s}$ be the top side of $\ttt$ and
${\bf s}\subset\partial\Delta$. Note that
the difference between the number of $a$-edges in the bottom
${\bf s'}$ of
$\ttt$ and the number of $a$-edges in $s$ cannot be greater than
$2K$, because every $(\theta,q)$-relator has at most two $a$-letters. Hence ${\bf |s'|-|s|}\le 2K\delta$. However, $\Delta'$ is obtained
by cutting off $\ttt$ along ${\bf s'}$, and its boundary contains two
$\theta$-edges fewer than $\Delta$. Hence we have ${\bf |s_0|-|s'_0|}\ge
2-2\delta$ for the complements ${\bf s}_0$ and ${\bf s'}_0$
of $\bf s$ and ${\bf s}'$, respectively, in the boundaries $\partial\Delta$
and $\partial\Delta'$. Finally,
$$|\partial\Delta|-|\partial\Delta'|\ge
2-2\delta-2K\delta-4\delta>1 $$ by (\ref{param}) and
(\ref{delta}).
\endproof

We call a base word $w$ \label{tightb} {\it tight} if

(1) for some letter $x$ the word $w$  has the form $uxvx$,
where the letter $x$ does not occur in $u$ and $x$ occurs in $v$ exactly $L-1$ times,

(2) every proper prefix $w'$  of $w$
does not
satisfy Property (1).

\begin{lemma}\label{width}
If a base $w$ of a $\theta$-band has no tight prefixes, then
$||w||\le K_0$, where $K_0=2LN$.
\end{lemma}
\proof The hub base includes
every base letter $L$ times. Hence every word in this group alphabet of length $\ge K_0+1$ includes one of the letters $L+1$ times. \endproof

From now on we shall fix a constant $K$ such that
\begin{equation}\label{kk}K>2K_0=4LN.\end{equation}

\bigskip

\begin{df} \label{com} {\rm We say that a reduced diagram $\Gamma$ is a {\em comb} if it has a
maximal $q$-band $\cal Q$ (\label{handle} the {\em handle} of the comb), such that

\begin{enumerate}
\item[$(C_1)$]${\bf bott}(\cal Q)$ is a part of $\partial \Gamma$, and every maximal
$\theta$-band of $\Gamma$ ends at a cell in $\cal Q$.
\end{enumerate}

If in addition the following properties hold:
\begin{enumerate}
\item[$(C_2)$] one of the maximal $\theta$-bands $\cal T$ in
$\Gamma$ has a tight base (if one reads the base towards the handle) and

\item[$(C_3)$] other maximal
$\theta$-bands in $\Gamma$ have tight  bases or bases
without tight prefixes
\end{enumerate}
then the comb is called \label{tightc} {\em tight}.

The
number of cells in the handle $\cal Q$ is the \label{heightt} {\em height} of the
comb, and the
maximal length of the bases of the $\theta$-bands of a
comb is called the \label{basicw} {\em basic width} of the comb. }
\end{df}

\begin{figure}
\begin{center}
\includegraphics[width=0.4\textwidth]{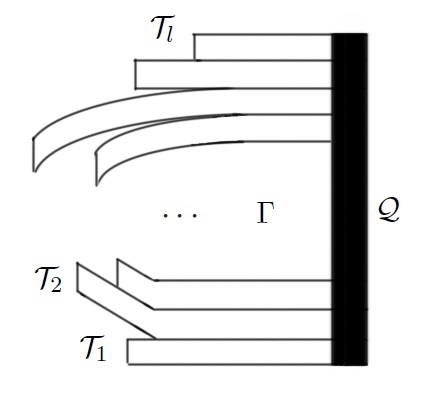}
\end{center}
\caption{Comb.} \label{Pic6}
\end{figure}

Notice that every trapezium is a comb.

\begin{lemma} \label{comb} (\cite{OS06}, Lemma 4.10) Let $l$ and $b$ be the height and
the basic width of a comb $\Gamma$ and let ${\cal T}_1,\dots {\cal T}_l$ be
consecutive $\theta$-bands of $\Gamma$ (as in Figure \ref{Pic6}). We
can assume that ${\bf bot}({\cal T}_1)$ and ${\bf top}({\cal T}_l)$ are contained in
$\partial\Gamma$. Denote by $\alpha=|\partial\Gamma|_a$ the number
of $a$-edges in the boundary of $\Gamma$, and by $\alpha_1$ the
number of $a$-edges on ${\bf bot}({\cal T}_1)$. Then $\alpha + 2lb\ge
2\alpha_1$, and the area of $\Gamma$ does not exceed $c_0bl^2+2\alpha
l$ for some constant $c_0=c_0(\bf M)$.
\end{lemma}

We say that a subdiagram $\Gamma$ of a diagram $\Delta$ is a
\label{subc} {\it subcomb} of $\Delta$ if $\Gamma$ is a comb, the handle of $\Gamma$ divides $\Delta$ in two parts, and $\Gamma$ is one of these parts.

\begin{lemma} \label{est'takaya}
Let $\Delta$ be a reduced diagram
over $G$
with non-zero area, where
every rim $\theta$-band has base
of length at least $K$.
Assume that

(1) $\Delta$ is a diagram over the group $M$ or

(2) $\Delta$ has a subcomb
   of basic width at least $K_0$.

   Then there exists a maximal $q$-band $\cal Q$ dividing
$\Delta$ in two parts, where one of the parts is a tight subcomb with handle $\cal Q$.
\end{lemma}

\proof Let $\ttt_0$ be a rim band of $\Delta$ (fig.\ref{Pic4}). Its base $w$ is of
length at least $K$, and therefore $w$ has disjoint prefix and
suffix of length $K_0$ since $K>2K_0$ by (\ref{kk}). The prefix of
this base word must have its own tight subprefix $w_1$, by
Lemma \ref{width} and the definition of tight words.
A $q$-edge of $\ttt_0$ corresponding to the last $q$-letter of
$w_1$ is the start edge of a maximal $q$-band $\qq'$ which bounds a
subdiagram $\Gamma'$ containing a band $\ttt$ (a subband of
$\ttt_0$) satisfying Property ($C_2$). It is useful to note that a
minimal suffix $w_2$ of $w$, such that $w_2\iv$ is tight, allows us
to construct another band $\qq''$ and a subdiagram $\Gamma''$ which
satisfies ($C_2$) and has no cells in common with $\Gamma'$.

\begin{figure}
\begin{center}
\includegraphics[width=0.7\textwidth]{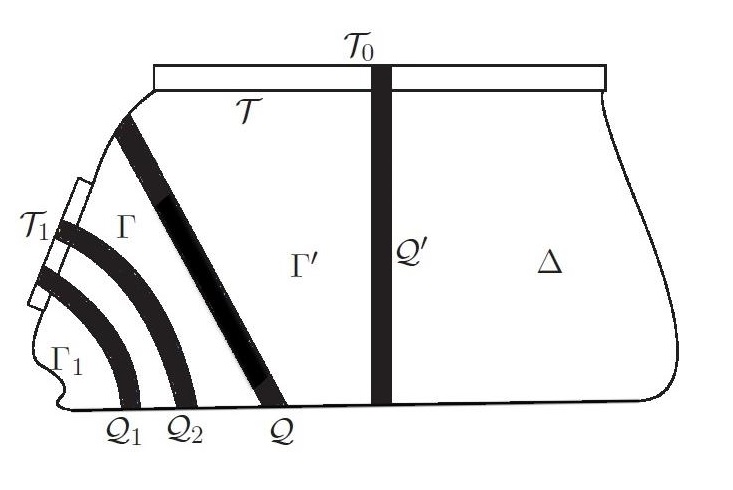}
\end{center}
\caption{Lemma \ref{est'takaya}.} \label{Pic7}
\end{figure}

Thus, there are $\qq$ and $\Gamma$ satisfying ($C_2$). Let us choose
such a pair with minimal $\area(\Gamma)$. Assume that there is a
$\theta$-band in $\Gamma$ which does not cross $\qq$. Then there
must exist a rim band $\ttt_1$ which does not cross $\qq$ in
$\Gamma$. Hence one can apply the construction from the previous
paragraph to $\ttt_1$ and construct two bands $\qq_1$ and $\qq_2$
and two disjoint subdiagrams $\Gamma_1$ and $\Gamma_2$ satisfying
the requirement ($C_2$) for $\Gamma$. Since $\Gamma_1$ and
$\Gamma_2$ are disjoint, one of them, say $\Gamma_1$, is inside
$\Gamma$. But the area of $\Gamma_1$ is smaller than the area of
$\Gamma$, and we come to a contradiction. Hence $\Gamma$ is a comb
and condition ($C_1$) is satisfied.

Assume that the base of a maximal $\theta$-band $\ttt$ of $\Gamma$
has a tight proper prefix (we may
assume that $\ttt$ terminates on $\qq$), and again one obtain a
$q$-band $\qq'$ in $\Gamma$, which provides us with a smaller
subdiagram $\Gamma'$ of $\Delta$, satisfying ($C_2$), a
contradiction. Hence $\Gamma$ satisfies Property ($C_3$) as well.

(2) The proof is shorter since a comb is given in the very beginning. \endproof

\subsection{Mixture on the boundaries of diagrams}\label{mix}

We will need a parameter of diagrams introduced in \cite{O12}. It was called mixture.

Let $O$ be a circle with two-colored
finite set of points (or vertices) on it,
more precisely, let any vertex of this finite set be either black or white. We call $O$ a \label{neckl}{\it necklace} with black and white \label{bead}{\it beads} on it. We want to introduce the {\em mixture} of this finite set of beads.

Assume that there are $n$ white beads and $n'$ black ones on $O$. We define sets \label{Pk}
${\bf P}_j$ of
ordered pairs of distinct white beads as follows. A pair $(o_1,o_2)$
($o_1\ne o_2$) belongs to the set ${\bf P}_j$ if the simple arc of $O$
drawn from $o_1$ to $o_2$ in clockwise direction has at least $j$
black beads. We denote by \label{muKO}$\mu_J(O)$ the sum $\sum_{j=1}^J \#
 {\bf P}_j$ (the \label{Kmix}$J$-{\it mixture} on $O$). Below similar sets for
another necklace $O'$ are denoted by ${\bf P'}_J$. In this subsection, $J\ge
1$, but later on it will be a fixed large enough number $J$ from the list (\ref{const}).

\begin{lemma}\label{mixture} (\cite{O12}, Lemma 6.1) (a) $\mu_J(O)\le J(n^2-n)$.

 (b) Suppose a
necklace $O'$ is obtained from $O$ after removal of a   white bead
$v$. Then $ \# {\bf P'}_j \le \#{\bf
P}_j$ for every $j$, and $\mu_J(O)-Jn<\mu_J(O')\le \mu_J(O).$

(c) Suppose a necklace $O'$ is obtained from $O$ after removal of a
black bead $v$. Then  $\# {\bf P'}_j \le \# {\bf P}_j$ for
every $j,$ and $\mu_J(O')\le \mu_J(O).$

(d) Assume that there are three black beads $v_1, v_2, v_3$ of a necklace
$O,$ such that the clockwise arc $v_1 - v_3$ contains $v_2$ and has
at most $J$ black beads (excluding $v_1$ and $v_3$), and the arcs
$v_1-v_2$ and $v_2-v_3$ have $m_1$ and $m_2$ white beads,
respectively. If $O'$ is obtained from $O$ by removal of $v_2$, then
$\mu_J(O')\le\mu_J(O)-m_1m_2.$
\end{lemma} $\Box$

For any diagram $\Delta$ over $G$, we introduce the following invariant
$\mu(\Delta)=\mu_{J} (\partial\Delta)$
depending on the boundary only.
 To define it, we
consider the boundary $\partial(\Delta),$ as a \label{mu1neckl} {\it necklace}, i.e.,
we consider a circle $O$ with $||\partial\Delta||$ edges labeled as the
boundary path of $\Delta.$ By definition, the
white beads are the mid-points  of the $\theta$-edges of
$O$ and black beads are the mid-points of the $q$-edges
$O$.
Then, by definition, the \label{muK1mix} ${\it mixture}$
on $\partial\Delta$ is $\mu(\Delta)=\mu_J(O).$

\subsection{Quadratic upper bound for quasi-areas of  diagrams over $M$.}\label{qub}

The Dehn function of the group $M$ is greater that the required function $F(n)=n^2f(n)^3$. For example, it is cubic if $f(n)=const.$ However we
are going to find the Dehn function of $G$, and
first we want to bound the areas of the words vanishing in $M$
with respect to the presentation of $G$. For this goal we artificially introduce the concept of $G$-{\it area}. The $G$-area of a big trapezia can be much less that the real area of it in $M$. This concept will be justified at the end of this paper, where some big trapezia are replaced by diagrams with hubs, but having lesser areas.

\begin{df}\label{abt} The $G$-area $\area_G(\Gamma)$ of a big trapezium $\Gamma$ is, by definition, the minimum of the half of its area (i.e., the number of cells) and the product $$c_5h(||{\bf top}(\Gamma)||+||{\bf bot}(\Gamma)||),$$ where $h$ is the height of $\Gamma$ and $c_5$ is the constant from (\ref{const}).

To define the $G$-area of a diagram $\Delta$ over $M$, we consider a family $\bf S$ of big subtrapezia (i.e. subdiagrams, which are trapezia) and single
cells of $\Delta$ such that every cell of $\Delta$ belongs to
a member $\Sigma$ of this family, and if a cell $\Pi$ belongs to different $\Sigma_1$ and $\Sigma_2$ from $\Sigma$, then both $\Sigma_1$ and $\Sigma_2$ are big subtrapezia of $\Delta$ with bases $xv_1x$, $xv_2x$, and $\Pi$ is a $(\theta,x)$-cell.
(In the later case, the intersection $\Sigma_1\cap\Sigma_2$
must be an $x$-band.)
There is
such a family `covering' $\Delta$, e.g. just the family of all cells of $\Delta$.

The $G$-area of $\bf S$ is the sum of $G$-areas of all big
trapezia from $\bf S$ plus the number of single cells from $\bf S$
(i.e. the $G$-area of a cell $\Pi$ is $area(\Pi)=1$). Finally,
the \label{areaGd} $G$-{\it area} $\area_G(\Delta)$ is the minimum of thea $G$-areas
of all `coverings' $\bf S$ as above.

\end{df}

It follows from the definition that $\area_G(\Delta)\le \area(\Delta)$ since the $G$-area of a big trapezium does not
exceed a half of its area.

\begin{lemma} \label{GA} Let $\Delta$ be a reduced diagram, and every cell $\pi$ of $\Delta$ belongs in one of subdiagrams $\Delta_1,\dots,\Delta_m$, where any intersection $\Delta_i\cap\Delta_j$ either has no cells or it is a $q$-band,
Then $\area_G(\Delta)\le \sum_{i=1}^m \area_G(\Delta_i)$.
\end{lemma}

\proof Consider the families ${\bf S}_1,\dots, {\bf S}_m$ given
by the definition of $G$-areas for the diagrams $\Delta_1,\dots,\Delta_m$. Then the set ${\bf S}={\bf S}_1\cup\dots\cup {\bf S}_m$
'covers' the entire $\Delta$ according to the above definition. This implies the required inequality for $G$-areas.
\endproof

We will show that for some constant $N_2$ and $N_1$ the $G$-area of any reduced diagram
$\Delta$ over $M$ with perimeter $n$ does not exceed $N_2n^2+N_1\mu(\Delta)$. (Using the quadratic upper bound for
$\mu(\Delta)$ from Lemma \ref{mixture} (a), one deduces that the $G$-area is bounded by $N'n^2$ for some constant $N'$.) Roughly speaking, we are doing the
following. We use induction on the perimeter of the diagram. First
we remove rim $\theta$-bands (those with one side and both ends on the boundary of
the diagram) with short bases. This operation decreases the
perimeter and preserves the sign of $N_2n^2+N_1\mu(\Delta)-\area_G(\Delta)$, so we can assume that the diagram
does not have such bands. Then we use Lemma \ref{est'takaya} and
find a tight comb inside the diagram with a handle $\cal C$. We also
find a long enough $q$-band $\cal C'$ that is close to $\cal C$. We use
a surgery which amounts to removing a part of the diagram between
$\cal C'$ and $\cal C$ and then gluing the two remaining parts of
$\Delta$ together. The main difficulty is to show that, as a result
of this surgery, the perimeter decreases and the measure and the
mixture change in such a way that the expression $N_2n^2+N_1\mu(\Delta)-\area_G(\Delta)$ does not change its sign. In the
proof, we need to consider several cases depending on the shape of
the subdiagram between $\cal C'$ and $\cal C$. Note that neither
$N_2n^2$ nor $N_1\mu(\Delta)$ nor $\area_G(\Delta)$ alone behave in
the appropriate way as a result of the surgery, but the expression
$N_2n^2+N_1\mu(\Delta)-\area_G(\Delta)$ behaves as needed.

\begin{rk} We introduced the surgery and used induction mentioned
above in \cite{OS06} (Lemma 6.2) to obtain a worse upper bound $n^2\log n$ for the area. But there were neither mixture, nor $G$ (just $M$), nor $G$-area in \cite{OS06}, and a different definition
for length $|*|$ was used there. Besides, we will use an auxiliary function $\Phi(x)$ in the proof to be able to repeat in part our argument
later, for diagrams over $G$. So we shall prove Lemma \ref{main} anew to obtain the better estimate.
\end{rk}

So, $N_1$ and $N_2$ are big enough constants from the list (\ref{const}). Here ``big enough'' means that
they satisfy the inequalities used in the proof of Lemma \ref{main}
(such that as (\ref{param3}),(\ref{param4}), (\ref{param5}), (\ref{param6}),
(\ref{param7}), (\ref{param8}), (\ref{param9}), (\ref{param10})).
Each of them has the form $N_i>* $ ($i=1,2$), where the right-hand side $*$ does not
depend on $N_i$ (but depends on the constants introduced earlier). Since the number of inequalities is finite,  the right choice of $N_1, N_2$ is possible.

\medskip

Let \label{Phif} $\Phi(x)$ be an arbitrary function defined for real $x\ge 0$ such
that

$\Phi(x) = x^2\phi(x)$ for a non-decreasing function $\phi(x)>0$ with
$\phi(1)\ge 1$ and

\begin{equation}\label{xy}
 \Phi(x)-\Phi(x-y)\ge xy\phi(x)\;\; for\;\; 0\le y \le x.
\end{equation}

\begin{rk} For this section, it suffices to take quadratic $\Phi(x)$ and $\phi(x)=const$, but to estimate the $G$-area of
diagrams with hubs, we will take the functions $\Phi(x)=F(x)$ and $\phi(x)=g(x)$, satisfying inequality (\ref{xy}) by Lemma \ref{d-x}.
\end{rk}

We are going to prove that the $G$-area of a reduced diagram $\Delta$ over $M$ does not exceed $N_2\Phi(n) +N_1\phi(n)\mu(\Delta)$, where
$n=|\partial\Delta|$.
Arguing by contradiction in the remaining part of this section, we consider a {\bf counter-example}
$\Delta$ with minimal perimeter $n$. Of course, its $G$-area is
positive, and, by Lemma \ref{NoAnnul}, we have at least two
$\theta$-edges on the boundary $\partial\Delta$, and so $n\ge 2$.

If $\Gamma$ is a comb with handle $\cal C$ and $\cal B$ is another
maximal $q$-band in $\Gamma$, then $\cal B$ cuts up $\Gamma$ in two parts,
where the part that does not contain $\cal C$ is a comb $\Gamma'$ with
handle $\cal B$. It follows from the definition of comb, that every
maximal $\theta$-band of $\Gamma$ crossing $\cal B$ connects $\cal B$
with $\cal C$. If  $\cal B$ and $\cal C$ can be connected by a $\theta$-band
containing no $(\theta,q)$-cells, then $\Gamma'$ is called the \label{ders} {\it derivative
subcomb} of $\Gamma$. Note that no maximal $\theta$-band of $\Gamma$ can
cross the handles of two derivative subcombs.

\begin{lemma} \label{notwo} (1) The diagram $\Delta$ has no two disjoint subcombs $\Gamma_1$ and $\Gamma_2$ of basic widths at most $K$ with handles ${\cal B}_1$ and ${\cal B}_2$ such that some ends of these handles are
connected by a subpath $\bf x$ of the boundary path of $\Delta$, where $\bf x$ has at most $N$ $q$-edges.

(2) The boundary of every subcomb $\Gamma$ with basic width $s\le K$ has $2s$ $q$-edges.
\end{lemma}

\proof We will prove the statements (1) and (2) using simultaneous induction on $A=\area(\Gamma_1)+\area(\Gamma_2)$
(resp., on $A=\area(\Gamma)$). Arguing by contradiction, we consider a counter-example with minimal $A$.

\begin{figure}
  \centering
  \includegraphics[width=1.0\textwidth]{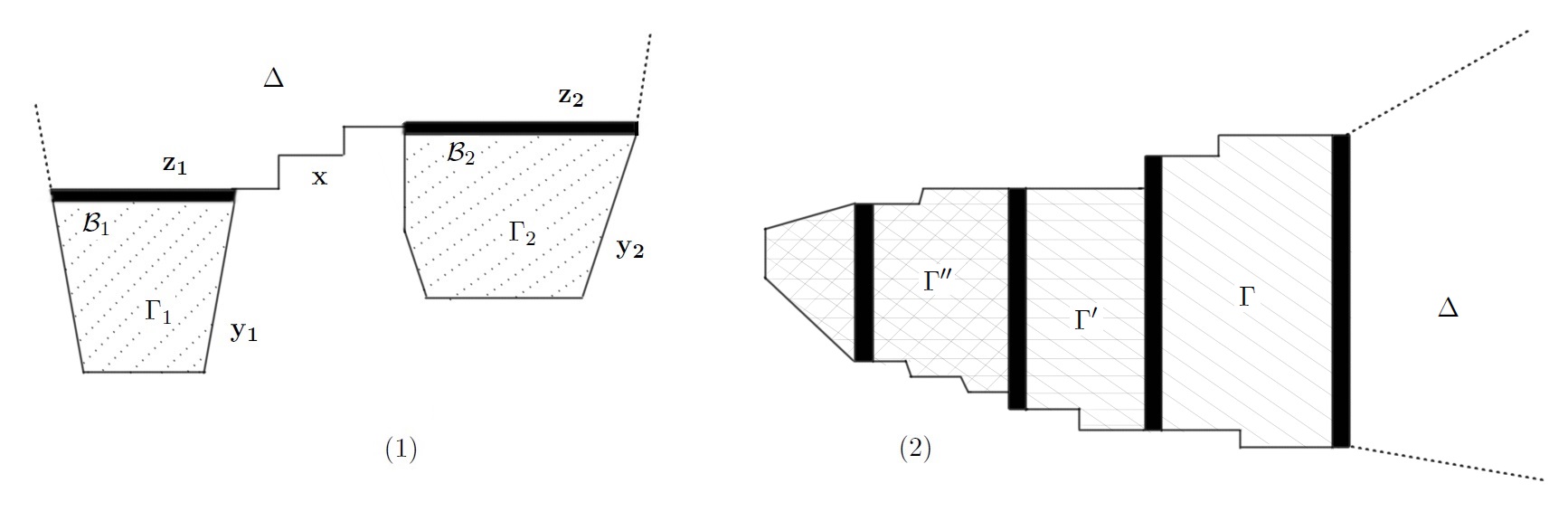}
  \caption{Lemma \ref{notwo}}\label{Pic8}
\end{figure}

(1) Since the area of $\Gamma_i$
($i=1,2$) is less than $A$, we may use Statement (2), and so we have at most $2K$ $q$ edges in $\partial\Gamma_i$.

Let $h_1$ and $h_2$ be the lengths of the handles ${\cal B}_1$ and ${\cal B}_2$ of $\Gamma_1$ and $\Gamma_2$, resp.
Without loss of generality, we assume that $h_1\le h_2$.
Denote by ${\bf y}_i{\bf z}_i$ the boundaries of $\Gamma_i$ ($i=1,2$), where ${\bf y}_i$ is the part of $\partial\Delta$
and ${\bf z}_i$ is the side of the handle of $\Gamma_i$ (so ${\bf y}_1{\bf x}{\bf y}_2$ is the part of the boundary path of $\Delta$, see fig. \ref{Pic8} (1)). Then each of the $\theta$-edges $\bf e$
of ${\bf y}_1$ is separated in $\partial\Delta$ from every $\theta$-edge $\bf f$ of ${\bf y}_2$ by
less than $4K+N < J$ $q$-edges.  Hence every such pair $({\bf e,f})$ (or the pair of white beads on these edges) makes a contribution to $\mu(\Delta)$.

Let $\Delta'$ be the diagram obtained by deleting the subdiagram
$\Gamma_1$ from $\Delta$. When passing from $\partial\Delta$ to $\partial\Delta'$,
  one replaces the $\theta$-edges  from ${\bf y}_1$ by the $\theta$-edge
  of ${\bf z}_1$ belonging to the same maximal  $\theta$-band. The same is true for white beads.

  But each of the $h_1h_2$ pairs in the corresponding set $P'$ of white beads is
  separated in $\partial\Delta'$ by less number of black
  beads than the pair defined by $\Delta$.
   Indeed, since the handle of $\Gamma_1$ is removed when one replaces $\partial\Delta$ by $\partial\Delta'$, two black bead at the ends of this handle are removed, and  therefore
  \begin{equation}\label{mudd}
 \mu (\Delta)-\mu (\Delta')\ge h_1h_2
 \end{equation}
  by Lemma \ref{mixture} (d).

  Let $\alpha$ be the number of $a$-edges in $\partial\Gamma_1$.
  It follows from Lemma \ref{comb} that the area, and so the $G$-area of $\Gamma_1$, does not exceed $
C_1(h_1)^2+2\alpha h_1$, where $C_1=c_0K$.

\begin{rk} The constants $C_1, C_2, C_{12}, C_3$ are not included in the list
(\ref{const}) since their values chosen here make sense only in the present subsection.
\end{rk}

  Since the boundary of $\Delta'$ has at
least two $q$-edges fewer than $\Delta$ and $|{\bf z}_1|=h_1\le |{\bf y}_1|$, we have $|\partial
\Delta'|\le |\partial\Delta|-2$. Moreover, we have from Lemma
\ref{ochev} (a) and Lemma \ref{NoAnnul} that
\begin{equation}\label{deriv}
|\partial \Delta|-|\partial\Delta'|\ge\gamma =\max(2,
\delta(\alpha-2h_1)),
\end{equation}
because the top/the bottom of ${\cal B}_1$ has at most $h_1$
$a$-edges.

This inequality, inequality (\ref{mudd}), and the inductive
assumption related to  $\Delta'$, imply that the $G$-area of
$\Delta'$ is not greater than
$$N_2\Phi(n-\gamma)+ N_1\phi(n)\mu(\Delta)-N_1\phi(n)h_1h_2.$$
Adding the $G$-area of $\Gamma_1$ and using inequality (\ref{xy}), we  see that by Lemma \ref{GA}, the $G$-area of $\Delta$ does
not exceed $$N_2\Phi(n)-N_2\gamma n +N_1\phi(n)\mu(\Delta)
-N_1\phi(n)h_1h_2+C_1h_1^2+2\alpha h_1.$$
Since $h_1\le h_2$ and $\phi(n)\ge 1$, this will contradict the choice of the
counter-example $\Delta$ when we prove that
\begin{equation}\label{nado1}
- N_2\gamma n -N_1h_1^2+C_1h_1^2+2\alpha h_1<0.
\end{equation}

If $\alpha\le 4h_1$, then inequality (\ref{nado1}) follows from
the inequalities $\gamma\ge 2$ and
\begin{equation}\label{param3i}
N_1\ge C_1+ 8.
\end{equation}

 Assume that $\alpha> 4h_1$. Then by (\ref{deriv}), we have
$\gamma\ge \frac12 \delta\alpha$ and $N_2\gamma n >2\alpha h_1$
since $n\ge 2h_1$ by Lemma \ref{NoAnnul}, and

\begin{equation}\label{param4}
N_2> 2\delta^{-1}.
\end{equation}
Since $N_1h_1^2>C_1h_1^2$ by (\ref{param3i}), the
inequality (\ref{nado1}) follows.

(2) If there are at least two derivative subcombs of $\Gamma$, then one can
find two of them satisfying the assumptions of Statement (1) (moreover, with $|x|_q=0$), and $\area(\Gamma_1)+\area(\Gamma_2)<\area(\Gamma)=A$, a contradiction. Therefore
there is a most one derivative subcomb $\Gamma'$ in $\Gamma$ (fig \ref{Pic8} (2)). In turn, $\Gamma'$
has at most one derivative subcomb $\Gamma''$, and so one. It follows that
there are no maximal $q$-bands in $\Gamma$ except for the handles of
$\Gamma', \Gamma'',\dots $. Since the basic width of $\Gamma$ is $s$, we have
$s$ maximal $q$-bands in $\Gamma$, and the lemma is proved.
\endproof

\begin{lemma} \label{twocombs} There are no pair of subcombs $\Gamma$ and $\Gamma'$
in $\Delta$ with handles $\cal X$ and $\cal X'$ of length $\ell$ and $\ell'$ such that $\Gamma'$
is a subcomb of $\Gamma$, the basic width of $\Gamma$ does not exceed $K_0$ and
$\ell'\le \ell/2$.
\end{lemma}

\begin{figure}
\begin{center}
\includegraphics[width=0.9\textwidth]{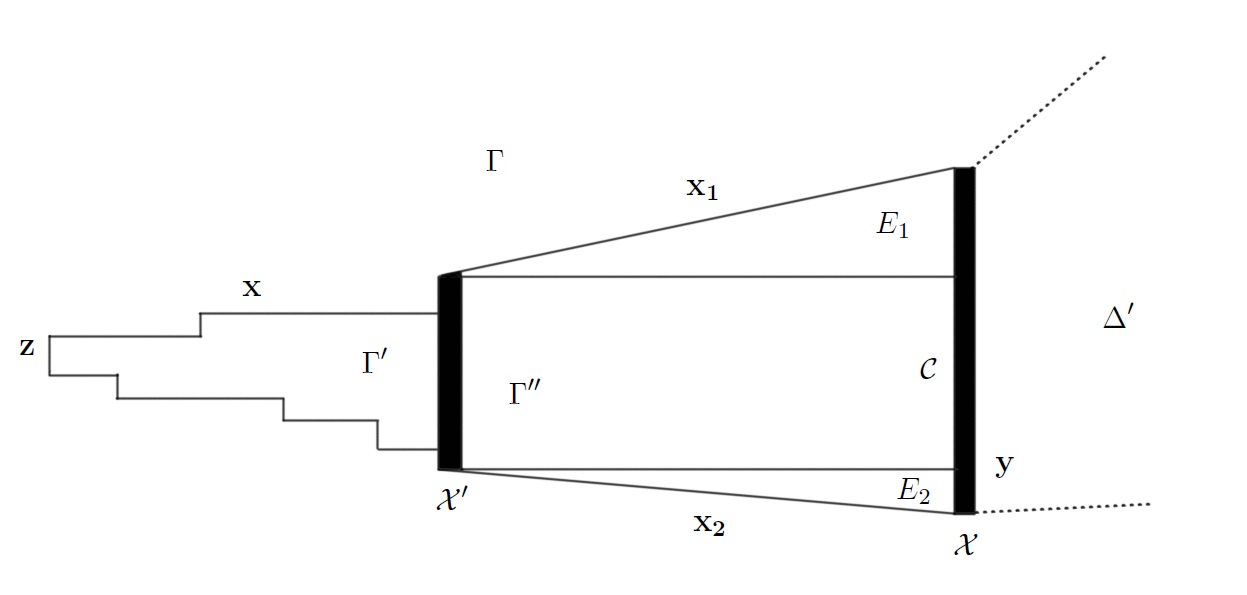}
\end{center}
\caption{Lemma \ref{twocombs}}\label{Pic9}
\end{figure}

\proof
Proving by contradiction, one can choose $\Gamma'$ so that $\ell'$ is minimal for
all subcombs in $\Gamma$ and so $\Gamma'$  has no
proper subcombs, i.e. its basic width is $1$ (fig. \ref{Pic9}).
It follows from Lemma \ref{comb}
that for $\alpha =|\Gamma'|_a$, we have
\begin{equation}\label{G'}
\area_G(\Gamma')\le\area(\Gamma')\le c_0(l')^2+2\alpha l'.
\end{equation}

Let $\Delta'$ be the diagram obtained after removing the subdiagram
$\Gamma'$ from $\Delta$. The following inequality is the analog of (\ref{deriv})
(where $h_1$ is replaced by $\ell'$)
\begin{equation}\label{eq679}
|\partial \Delta|-|\partial\Delta'|\ge\gamma =\max(2,
\delta(\alpha-2l')).
\end{equation}

The $q$-band $\cal X$ contains a subband $\cal C$ of length
$l'$. Moreover one can choose $\cal C$ so that all maximal $\theta$-bands of
$\Gamma$  crossing the handle $\cal X'$ of $\Gamma'$, start from $\cal C$.
These $\theta$-bands form a comb $\Gamma''$ contained in $\Gamma$, and in turn,
$\Gamma''$ contains $\Gamma'$. The two parts of the
complement ${\cal X}\backslash{\cal C}$ are the handles of two subcombs $E_1$ and $E_2$
formed by maximal $\theta$-bands of $\Gamma$, which do not cross $\cal X'$.
Let the length of these two handles be $\ell_1$ and $\ell_2$, respectively, and so we have $\ell_1+\ell_2=l-l'>l'$. ($E_1$ or $E_2$ can be empty;
then $\ell_1$ or $\ell_2$ equals $0$.)

It will be convenient to assume that $\Gamma$ is drawn from the left of the vertical handle $\cal X$. Denote by ${\bf yz}$ the boundary path
of $\Gamma$, where ${\bf y}$ is the right side of the band $\cal X$. Thus, there are $l_1$ (resp., $l_2$)
$\theta$-edges on the common subpath ${\bf x}_1$ (subpath ${\bf x}_2$) of ${\bf z}$ and $\partial E_1$ (and $\partial E_2$).

By Lemma \ref{notwo} (2),
the path ${\bf z}$ contains at most $2K_0$ $q$-edges, because the basic width of $\Gamma$ is at most $K_0$.

Consider the factorization ${\bf z=x}_2{\bf xx}_1$, where ${\bf x}$ is a subpath
of $\partial \Gamma'$.
It follows that between every white bead on ${\bf x}_1$ (i.e. the middle point of the $\theta$-edges on ${\bf x}_1$) and a white bead on $\bf x$ we
have at most $2K_0$ black beads (i.e. the middle points of the
$q$-edges of the path $\bf x$). Since $J$ is greater than $2K_0$,
every pair of white beads, where one bead belongs to $\bf x$
and another one belongs to ${\bf x}_1$ (or, similarly,  to ${\bf x}_2$) contributes $1$ to $\mu(\Delta)$. Let $P$ denote the set
of such pairs. By the definition of $E_1$
and $E_2$, we have $\# P =l'(\ell_1+\ell_2)=l'(l-l')>(l')^2$.

  When passing from $\partial\Delta$ to $\partial\Delta'$,
  one replaces the left-most $\theta$-edges of every maximal $\theta$-band from $\Gamma'$ with the right-most $\theta$-edges
  lying on the right side of ${\cal X}'$. The same is true for white beads.
  But each of the $l'(l-l')$ pairs in the corresponding set $P'$  of white beads is
  separated in $\partial\Delta'$ by less number of black
  beads since the $q$-band $\cal X'$ is removed. Therefore every
  pair from $P'$ gives less by $1$ contribution to the mixture, as it follows from the definition of mixture. Hence $\mu (\Delta)-\mu (\Delta')\ge l'(l-l')\ge(l')^2$.
This inequality, inequality (\ref{eq679}), and the inductive
assumption related to $\Delta'$, imply that the $G$-area of
$\Delta'$ is not greater than
$$N_2\Phi(n-\gamma)+ N_1\phi(n)\mu(\Delta)-N_1\phi(n)(\ell')^2.$$
Adding the $G$-area of $\Gamma'$ (\ref{G'}) and applying inequality (\ref{xy}), we see that the $G$-area of $\Delta$ does
not exceed $$N_2\Phi(n) +N_1\phi(n)\mu(\Delta)- N_2\gamma n
-N_1\phi(n)(l')^2+c_0(l')^2+2\alpha l'.$$
This will contradict the choice of the
counter-example $\Delta$ when we prove that
\begin{equation}\label{nado}
- N_2\gamma n -N_1(l')^2+c_0(l')^2+2\alpha l'<0,
\end{equation}
 because $\phi(n)\ge 1$. Consider two cases.

(a) Let $\alpha\le 4l'$. Then inequality (\ref{nado}) follows from
the inequalities $\gamma\ge 2$ and
\begin{equation}\label{param3}
 N_1\ge c_0+ 8.
\end{equation}

(b) Assume that $\alpha> 4l'$. Then by (\ref{eq679}) we have
$\gamma\ge \frac12 \delta\alpha$ and $N_2\gamma n >2\alpha l'$
since $n\ge 2l\ge 4l'$ by Lemma \ref{NoAnnul}, and

\begin{equation}\label{param4}
N_2> \delta^{-1}.
\end{equation}
Since $N_1(l')^2>c_0(l')^2$ by (\ref{param3}), the
inequality (\ref{nado}) follows.

Thus, the lemma is proved by contradiction.
\endproof

\begin{lemma} \label{nori} $\Delta$ has no rim $\theta$-band
whose base has  $s\le K$ letters.
\end{lemma}

\begin{figure}
\begin{center}
\includegraphics[width=0.5\textwidth]{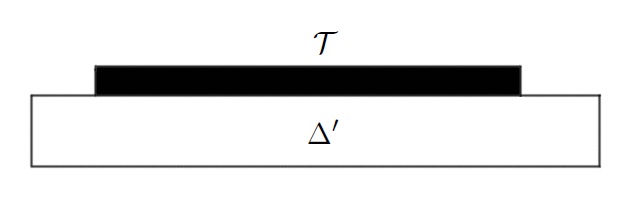}
\end{center}
\caption{Rim $\theta$-band}\label{Pic10}
\end{figure}

\proof Assume by contradiction that such a rim $\theta$-band $\cal T$ exists,
and ${\bf top}(\cal T)$ belongs in
$\partial(\Delta)$ (fig.\ref{Pic10}). When deleting $\cal T$, we obtain, by Lemma
\ref{rim}, a diagram $\Delta'$ with $ |\partial\Delta'|\le
n-1$. Since ${\bf top}(\cal T)$ lies on $\partial\Delta$, we have
from the definition of the length , that the
number of $a$-edges in ${\bf top}(\cal T)$ is less than
$\delta^{-1}(n-s)$. By Lemma \ref{band}, the length of $\cal T$ is at
most $3s+\delta^{-1}(n-s)< \delta^{-1}n$. Thus, by applying the inductive
hypothesis to $\Delta'$, we have that $G$-area of $\Delta$ is not
greater than $N_2\Phi(n-1) +N_1\phi(n)\mu(\Delta)
+\delta^{-1}n$ because $\mu(\Delta')\le\mu(\Delta)$ by Lemma
\ref{mixture} (b). But the first term of this sum does not exceed $N_2\Phi(n)-N_2n$ by (\ref{xy}), and so the entire sum is bounded by
$N_2\Phi(n) +N_1\phi(n)\mu(\Delta)$
 provided

\begin{equation}\label{param5}
N_2\ge\delta^{-1}.
\end{equation}

This contradicts to the choice of $\Delta$, and the lemma is proved.
\endproof

\begin{lemma}\label{main} The $G$-area of a reduced diagram $\Delta$ over $M$ does not exceed  $N_2\Phi(n) +N_1\phi(n)\mu(\Delta)$, where
$n=|\partial\Delta|$.
\end{lemma}

\proof We continue studying the hypothetical minimal counter-example $\Delta$. By Lemma \ref{nori}, now we can apply Lemma \ref{est'takaya} (1). By that
lemma, there exists a tight subcomb $\Gamma\subset\Delta$. Let $\cal T$ be a
$\theta$-band of $\Gamma$ with a tight base.

\begin{figure}
\begin{center}
\includegraphics[width=1.0\textwidth]{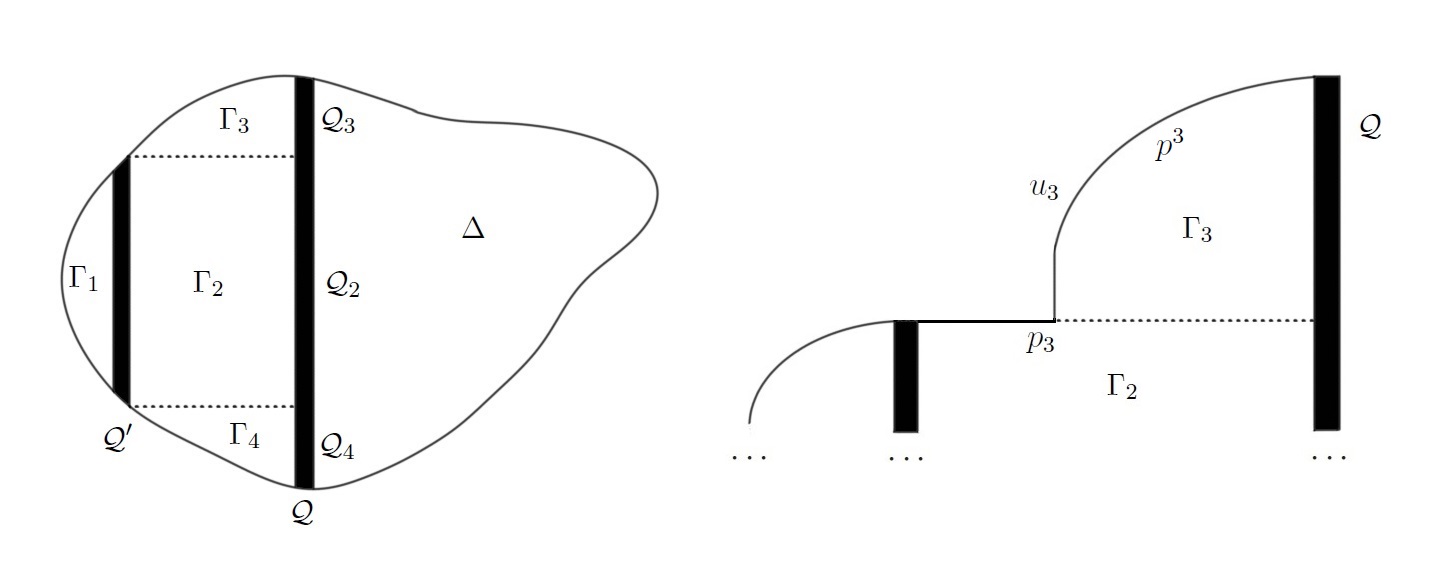}
\end{center}
\caption{Lemma \ref{main}.} \label{Pic18}
\end{figure}

The basic width of $\Gamma$ is less than $K_0$ by Lemma \ref{width}.
Since the base of $\Gamma$ is tight, it is equal to $uxvx$ for
some $x$, where the last occurrence of $x$ corresponds to the handle $\cal Q$ of $\Gamma$,
 the word $u$ does not contain $x$, and $v$ has exactly $L-1$  occurrences of $x$.
Let $\cal Q'$ be the maximal
$x$-band of $\Gamma$ crossing $\cal T$ at the cell corresponding to
the first occurrence of $x$ in $uxvx$ (fig. \ref{Pic18} (a)).

We consider the smallest subdiagram $\Gamma'$ of $\Delta$ containing
all the $\theta$-bands of $\Gamma$ crossing the $x$-band $\cal Q'$. It
is a comb with handle ${\cal Q}_2\subset{\cal Q}$. The comb $\Gamma'$ is
covered by a trapezium $\Gamma_2$ placed between $\cal Q'$ and $\cal Q$,
and the comb $\Gamma_1$ with handle $\cal Q'$. The band $\cal Q'$
belongs to both $\Gamma_1$ and $\Gamma_2$. The remaining part of
$\Gamma$ is a disjoint union of two combs $\Gamma_3$ and $\Gamma_4$
whose handles ${\cal Q}_3$ and ${\cal Q}_4$ contain the cells of $\cal Q$ that do
not belong to the trapezium $\Gamma_2$. The handle of $\Gamma$ is
the composition of handles ${\cal Q}_3$, ${\cal Q}_2$, ${\cal Q}_4$ of $\Gamma_3$,
$\Gamma'$ and $\Gamma_4$ in that order.

Let the lengths of ${\cal Q}_3$ and ${\cal Q}_4$ be $l_3$ and $l_4$,
respectively. Let $l'$ be the length of the handle of $\Gamma'$.
Then by Lemma \ref{twocombs}, we have
\begin{equation}\label{ll'}
l'>l/2 \;\;\; and \;\;\;l=l'+l_3+l_4.
\end{equation}

For $i\in \{3,4\}$ and $\alpha_i=|\partial\Gamma_i|_a$, Lemma
\ref{comb} gives inequalities
\begin{equation}\label{gamma34}
A_i\le C_1l_i^2+2\alpha_i l_i,
\end{equation}
where $A_i$ is the $G$-area of $\Gamma_i$. (We take into account that $G$-area cannot exceed area.)

Let
${\bf p}_3, {\bf p}_4$ be the top and the bottom of the trapezium $\Gamma_2$.
Here ${\bf p}_3^{-1}$ (resp. ${\bf p}_4^{-1}$) shares some initial edges with $\partial
\Gamma_3$ (with $\partial \Gamma_4$), the rest of these paths belong
to the boundary of $\Delta$. We denote by $d_3$ the number of
$a$-edges of ${\bf p}_3$ and by $d'_3$ the number of the $a$ edges of ${\bf p}_3$ which do
not belong to $\Gamma_3$. Similarly, we introduce $d_4$ and $d'_4$.

Let $A_2$ be the $G$-area of $\Gamma_2$. Then by Lemma \ref{or} and the definition of the $G$-area for big trapezia,
we have

\begin{equation}\label{gamma2}
A_2\le C_2l'(d_3+d_4+1)
\end{equation}
for some constant $C_2<\delta^{-1}$, because the basic width of $\Gamma_2$ is less than $K$.

Now we observe that the handle ${\cal Q}_2$ of $\Gamma'$ is a copy of $\cal Q'$
because both maximal $q$-bands of the trapezium $\Gamma_2$
correspond to the same basic letter $x$.

This makes the following
surgery possible. The diagram $\Delta$ is covered by two
subdiagrams: $\Gamma$ and another subdiagram $\Delta_1$, having only
the band $\cal Q$ in common. We construct a new auxiliary diagram by
attaching $\Gamma_1$ to $\Delta_1$
	with identification of the of the band
$\cal Q'$ of $\Gamma_1$ and the band ${\cal Q}_2$.
We denote the constructed diagram by $\Delta_0$.

$\Delta_0$ is a reduced diagram because every pair of
its cells having a common edge, has a copy either in $\Gamma_1$ or
in $\Delta_1$. Now we need the auxiliary

\begin{lemma} \label{A00} The $G$-area $A_0$ of $\Delta_0$ is at least the
sum of the $G$-areas of $\Gamma_1$ and $\Delta_1$ minus $l'$.
\end{lemma}

\proof Consider a minimal `covering' $\bf S$ of $\Delta_0$ from the definition of $G$-area, and assume that there is a big trapezium $E\in \bf S$, such that neither $\Gamma_1$ nor $\Delta_1$ contains it. Then $E$ has a base $ywy$, where
$(yw)^{\pm 1}$ is a cyclic permutation of the $L$-th power of the standard base,
and the first $y$-band of $E$ is in $\Gamma_1$, but it is not a subband of $\cal Q'$.

Since the history $H$ of the big trapezium $E$ is a subhistory
of the history of $\Gamma_2$,
we conclude that $\Gamma_2$ is a big trapezium itself, and therefore $(xv)^{\pm 1}$ is an $L$-th power of a  cyclic shift of the standard base (or of the inverse of it). Since the first $y$ occurs in uxvx
before the first $x$ it follows that we have the $(L+1)-th$
occurrence of $y$ before the last occurrence of $x$ in the word $uxvx$. But this contradicts to the definition of tight comb $\Gamma$.

Hence every big trapezium from $\bf S$ entirely belongs
either to $\Gamma_1$ or to $\Delta_1$. Therefore one can obtain
'coverings' $\bf S'$ and $\bf S''$ of these two diagrams
if (1) every $\Sigma$ from $\bf S$ is assigned either to $\bf S'$ or to $\bf S''$ and then (2) one add at most $l'$ single cells since the common band ${\cal Q}'$ in $\Delta_0$ should
be covered twice in disjoint diagrams $\Gamma_1$ and $\Delta_1$.
These construction complete the proof of the lemma.
\endproof

By Lemma \ref{GA}, the $G$-area of $\Delta$ does not exceed the sum
of $G$-areas of the five subdiagrams
$\Gamma_1$, $\Gamma_2$, $\Gamma_3$, $\Gamma_4$ and $\Delta_1$. But
the direct estimate of each of these values is not efficient.
Therefore we will use Lemma \ref{A00} to bound the $G$-area of the
auxiliary diagram $\Delta_0$ built of two pieces $\Gamma_1$ and $\Delta_1$.

It follows from our constructions and Lemmas \ref{GA}, \ref{A00},  that
\begin{equation}\label{AGD}
 \area_G(\Delta)\le A_2+A_3+A_4+A_0+l'.
 \end{equation}

Now we continue proving Lemma \ref{main}.

\medskip

Let ${\bf p}^3$ be the segment of the boundary $\partial\Gamma_3$ that
joins $\cal Q$ and $\Gamma_2$ along the boundary of $\Delta$ (fig. \ref{Pic18} (b)). It
follows from the definition of $d_3$, $d'_3$, $l_3$ and $\alpha_3$,
that the number of $a$-edges lying on ${\bf p}^3$ is at least $\alpha_3
-(d_3-d'_3)-l_3$.

Let ${\bf u}_3$ be the part of $\partial\Delta$ that contains ${\bf p}^3$ and
connects $\cal Q$ with $\cal Q'$. It has $l_3$ $\theta$-edges. Hence we
have, by Lemma \ref{ochev}, that
at least
$$|{\bf u}_3|\ge \max(l_3, l_3+\delta(|p^3|_a-l_3))\ge
\max(l_3, l_3+\delta(\alpha_3 -(d_3-d'_3)-2l_3)).$$ Since ${\bf u}_3$ includes a subpath of length $d'_3$ having no $\theta$-edges, we also have by Lemma \ref{ochev} (c) that $|{\bf u}_3|\ge
l_3+\delta(d'_3-1)$.

One can similarly define ${\bf p}^4$ and ${\bf u}_4$ for $\Gamma_4$. When
passing from $\partial\Delta$ to $\partial\Delta_0$ we replace the
end edges of $\cal Q'$, ${\bf u}_3$ and ${\bf u}_4$ by two subpaths of
$\partial\cal Q$ having lengths $l_3$ and $l_4$. Let
$n_0=|\partial\Delta_0|$. Then it follows from the previous paragraph
that

\begin{equation}\label{nbezn0}
n-n_0\ge2+\delta(\max(0,d'_3-1, \alpha_3-(d_3-d'_3)-2l_3)+\max(0,
d'_4-1,\alpha_4-(d_4-d'_4)-2l_4)).
\end{equation}

In particular, $n_0\le n-2$. By the inductive hypothesis,
\begin{equation}\label{A0}
A_0\le N_2\Phi(n_0) +N_1\phi(n_0)\mu(\Delta_0).
\end{equation}

We note that the mixture $\mu(\Delta_0)$ of $\Delta_0$ is not
greater than $\mu(\Delta)- l'(l-l')$ . Indeed, by Lemma \ref{twocombs} (2), one can use the same trick as in Lemma \ref{twocombs} as follows.
For every pair of white beads, where
one bead corresponds to a $\theta$-band of $\Gamma_2$ and
another one to a $\theta$-band of $\Gamma_3$ or $\Gamma_4$,
the contribution of this pair to $\mu(\Delta_0)$ is less than the contribution to $\Delta$. It remains to count
the number of such pairs: $l'(l_3+l_4)=l(l-l')$.

Therefore, by inequalities (\ref{A0}) and (\ref{xy}), the $G$-area of $\Delta$ is not
greater than
\begin{equation}\label{ADelta}
N_2\Phi(n) +N_1\phi(n)\mu(\Delta)- N_2 n(n-n_0) -N_1
\phi(n)l'(l-l')+A_2+A_3+A_4+l'.
\end{equation}

In view of inequalities (\ref{gamma2}), (\ref{gamma34}) and (\ref{AGD}), to
obtain the desired contradiction, it suffices to prove that
\begin{equation}\label{tsel}
N_2 n(n-n_0)+ N_1 l'(l-l')\ge  C_{12}l'(d_3+d_4+1)+C_{12}(l_3^2+l_4^2)+2\alpha_3 l_3 + 2\alpha_4l_4+l',
\end{equation}
where $C_{12}=\max(C_1,C_2)$.

Since $l=l'+l_3+l_4$, it suffices to prove that
\begin{equation}\label{tsel'}
N_2n(n-n_0)+N_1l'(l-l')\ge  C_3l'(d_3+d_4+1)+C_3(l_3^2+l_4^2)+2\alpha_3 l_3 + 2\alpha_4l_4,
\end{equation}
where $C_3=C_{12}+1$.

Note that we can assume that
\begin{equation}\label{c3}
C_3>>1.
\end{equation}

First we can choose $N_1$ big enough so that
$N_1 l'(l-l')/3\ge C_3(l_3+l_4)^2\ge C_3(l_3^2+l_4^2)$.
Indeed, by (\ref{ll'}), we obtain
$\frac{N_1}{3} l'(l-l')\ge \frac{N_1}{3}
(l_3+l_4)(l_3+l_4)$, so it is enough to assume that
\begin{equation}
\label{param6} N_1>3C_3.
\end{equation}

We also have that
\begin{equation}\label{a4}
\frac{N_2}2 n(n-n_0)\ge C_3l',
\end{equation}
because $n-n_0\ge 2$, $n\ge 2l'$ and $N_2\ge C_3$ by (\ref{param6}).

It remains to prove that \begin{equation}\label{t1} \frac{N_2}2
n(n-n_0)+\frac{2N_1}{3}l'(l-l')>
C_3l'(d_3+d_4)+2\alpha_3l_3+2\alpha_4l_4.
\end{equation}

We assume without loss of generality that $\alpha_3\ge\alpha_4$, and
consider two cases.

\medskip

(a) Suppose $\alpha_3\le 2C_3(l-l')$.

Since $d_i\le \alpha_i+d'_i$ for $i=3,4$,  by inequality
(\ref{nbezn0}), we have $$d_3+d_4\le
\alpha_3+\alpha_4+d_3'+d_4'<4C_3(l-l')+\delta^{-1}(n-n_0)+2-2\delta\iv<
4C_3(l-l')+\delta\iv(n-n_0).
$$

Therefore

\begin{equation}\label{a1}
\frac{N_1}{3} l'(l-l')+\frac{N_2}2 n(n-n_0)\ge
4C_3^2l'(l-l')+C_3\delta\iv(n-n_0)l'>
  C_3l'(d_3+d_4)
\end{equation}
since we can assume that

\begin{equation}\label{param7}
N_1> 12C_3^2,\qquad N_2/2> C_3\delta^{-1}.
\end{equation}

We also have  by (\ref{ll'}):
\begin{equation}\label{a3}
\frac{N_1}{3} l'(l-l')\ge \frac{N_1}{3}(l_3+l_4)(l_3+l_4)\ge
\frac{N_1}{3}\frac{\alpha_3+\alpha_4}{4C_3}(l_3+l_4)>
2\alpha_3l_3+2\alpha_4l_4
\end{equation}
since we can assume that
\begin{equation}\label{param8}
N_1>24C_3.
\end{equation}

 The sum of inequalities (\ref{a1}) and (\ref{a3}) gives us the
desired inequality (\ref{t1}).

\medskip

(b) Assume now that $\alpha_3>2C_3(l-l')$. Then, applying Lemma
\ref{comb} to the comb $\Gamma_3$, we obtain
\begin{equation}\label{dd'}
d_3-d'_3<\frac12\alpha_3+K_0l_3\le\frac56\alpha_3
\end{equation}
since $l_3\le
l-l'<\frac{\alpha_3}{2C_3}$ and
\begin{equation}\label{param9} C_3>3K_0.\end{equation} We also have
$d_4-d'_4<\frac12\alpha_4+K_0l_4\le\frac56 \alpha_3$. These two
inequalities and inequality (\ref{nbezn0}) lead to
\begin{equation}\label{d34}
d_3+d_4\le \frac53\alpha_3+\delta^{-1}(n-n_0).
\end{equation}

It follows from (\ref{dd'}) that
$$\alpha_3-(d_3-d'_3)-2l_3\ge \frac16 \alpha_3 -
\frac{2}{2C_3}\alpha_3\ge\frac17\alpha_3,$$ since $l_3\le
l-l'<\frac{\alpha_3}{2C_3}$ and $C_3>42$ by (\ref{c3}). Therefore,
by (\ref{nbezn0}),
\begin{equation}\label{raznitsa}
n-n_0\ge \frac17 \delta\alpha_3.
\end{equation}
Thus, by (\ref{d34}),

\begin{equation}
\label{d10} d_3+d_4<13 \delta\iv(n-n_0).
\end{equation}

Since $2l'<n$ and $n-n_0\ge 2$, inequality (\ref{d10}) implies

\begin{equation}\label{b1}
\frac{N_2}3 n(n-n_0)> C_3l'(d_3+d_4),
\end{equation}
because we can assume that \begin{equation}\label{param10}
N_2>>C_3\delta\iv
\end{equation}
($N_2>21C_3\delta\iv$ is enough).

Inequalities (\ref{raznitsa}), (\ref{param10}),
$\alpha_3\ge\alpha_4$, and $4(l_3+l_4)\le n$ give us
\begin{equation}\label{b2}
\frac{N_2}6 n(n-n_0)\ge \frac72 C_3 \delta^{-1}(n-n_0)n\ge
2\alpha_3(l_3+l_4)\ge 2\alpha_3l_3+2\alpha_4 l_4.
\end{equation}

The inequality (\ref{t1}) follows now from inequalities (\ref{b1}),
and (\ref{b2}).
\endproof

\section{Minimal diagrams over $G$}\label{midi}

Given a reduced diagram $\Delta$ over the group $G,$ one can construct a planar graph whose vertices
are the hubs of this diagram plus one improper vertex outside $\Delta,$ and the
edges are maximal $t$-bands of $\Delta.$

Let us consider two hubs $\Pi_1$ and $\Pi_2$ in a reduced diagram,
connected by two neighbor
$t$-bands
${\cal C}_i$ and  ${\cal C}_{i+1}$, where
there are no other hubs between these $t$-bands.
These bands, together with parts of
\begin{figure}
\begin{center}
\includegraphics[width=0.8\textwidth]{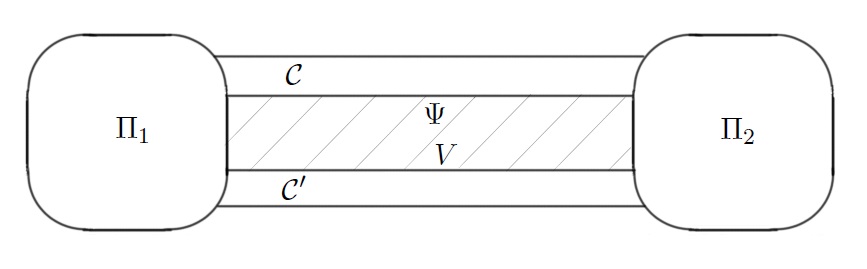}
\end{center}
\caption{Cancellation of two hubs}\label{2h}
\end{figure}
$\partial\Pi_1$ and $\partial\Pi_2,$ bound either a subdiagram
having no cells,
or  a trapezium $\Psi$ of height $\ge 1$ (fig. \ref{2h}). The former case is impossible since in this case
the hubs have a common $t$-edge and, the diagram is not reduced
since all cyclic permutations of the word $(W_M)^L$ starting with $t$ are equal.
We want to show that the latter case is not possible either if the diagram $\Delta$ is
chosen with minimal number of hubs among the diagrams with the same boundary label.

Indeed,  both the top and the bottom labels of $\Psi$ are equal to the word $(W_Mt)^{\pm 1},$ and removing ${\cal C}_{i+1}$ from $\Psi$,
we get a subdiagram $\Psi'$ with top/bottom label $(W_M)^{\pm 1}$ and the same label $V$  of its sides. It follows that $W_M$ and $V$ commute in
the group $M$. Hence the word $U\equiv (W_M)^{1-L}V(W_M)^{L-1}V^{-1}$ is equal to $1$ in $M$. But $U$
is the boundary label of a subdiagram $\Gamma$ containing $\Psi'$
and both $\Pi_1$ and $\Pi_2$. Hence one can replace $\Gamma$ with
a diagram over $M$, decreasing the nubmer of hubs in $\Delta$,
a contradiction.

If $W$ is an $\bf M$-accepted word, then the word $(W)^L$ is equal to $1$ in $G$. To see this, one can glue up $L$
copies $\Delta_1,\dots, \Delta_L$ ot the trapezia corresponding to the accepting computation of $W$, identifying
the right side of each $\Delta_i$ and the left side of $\Delta_{i+1}$ (indices are taken modulo $L$). The obtained
annulus has inner boundary labeled by the hub word $(W_M)^L$,
and so the hole can be glued up by a hub cell.

As in \cite{SBR} and \cite{O97}, we will increase the set of relations of $G$ by adding
the (infinite) set of {\it disk relators} $(W)^L$ for
every accepted word $W$. So we will consider diagrams
with \label{diskr} {\it disks}, where every disk cell (or just {\it disk}) is labeled by
such a word  $(W)^{\pm L}$. (In particular, the hub is a disk.)

Again, if two disks are connected by two $t$-bands and there are no other disks between these $t$-bands, then one can reduce the number of disks in the diagram. For this aid,
it suffices to replace the disks with hubs and the cells
corresponding to the defining relations of $M$, and apply
the trick exploited above.

\begin{df}\label{minimald}
We will call a reduced diagram $\Delta$ {\it minimal} if

(1) the number of disks is minimal for the diagrams with the same boundary label and

(2) $\Delta$ has minimal number of $(\theta,t)$-cells among
the diagrams with the same boundary label and with minimal number of  disks.

Clearly, a subdiagram of a minimal diagram is minimal itself.

\end{df}

Thus, no two disks of a minimal diagram are connected by  two $t$-bands,
such that the subdiagram bounded by them contains no other disks.
This property makes the disk graph of a reduced diagram
hyperbolic, in a sense, if the degree $L$ of every proper vertex (=disk) is high ($L>> 1$).
Below we give a more precise formulation (proved for diagrams with such a disk graph, in particular,
in \cite{SBR}, Lemma 11.4 and in  \cite{O97}, Lemma 3.2).

\begin{figure}
\begin{center}
\includegraphics[width=0.7\textwidth]{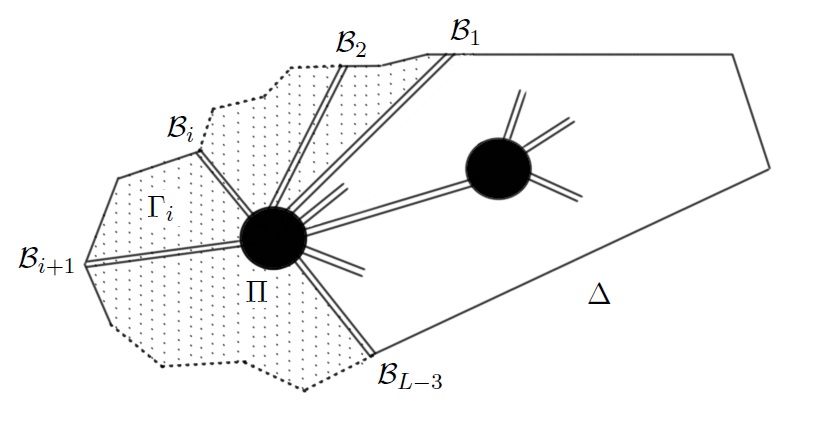}
\end{center}
\caption{Lemma \ref{extdisc}}\label{extd}
\end{figure}

\begin{lemma} \label{extdisc} If a a minimal diagram  contains at least one disk,
then there is a disk $\Pi$ in $\Delta$ such that $L-3$ consecutive maximal $t$-bands ${\cal B}_1,\dots
{\cal B}_{L-3} $ start on $\partial\Pi$ , end on the boundary $\partial\Delta$, and for any $i\in [1,L-4]$,
there are no disks in the subdiagram $\Gamma_i$ bounded by ${\cal B}_i$, ${\cal B}_{i+1},$ $\partial\Pi,$ and $\partial\Delta$ (fig. \ref{extd}).
\end{lemma} $\Box$

A maximal $q$-band starting on a disk of  a diagram is called a \label{spoke}
{\it spoke}.

Lemma \ref{extdisc} implies by induction on the number of hubs:

\begin{lemma} \label{mnogospits} (see \cite{O12}, Lemma 5.19) If a reduced diagram $\Delta$ has $m\ge 1$ hubs, then the number of spokes of $\Delta$ ending on the boundary $\partial\Delta$, and therefore
the number of $q$-edges in the boundary path of $\Delta$, is greater than  $mLN/2$ .
\end{lemma} $\Box$

Recall the following transformation for diagrams
with disks, exploited earlier in \cite{O97}, \cite{SBR}. Assume there is a disk $\Pi$ and a $\theta$-band $\cal T$ subsequently crossing some spokes ${\cal B}_1,\dots, {\cal B}_k$ which start (say, counter-clockwise\cal) from $\Pi$.
Assume that $k\ge 2$ and there are no other cells between $\Pi$
and the bottom of $\cal T$, and so there is  a subdiagram $\Gamma$
formed by $\Pi$ and $\cal T$.

\begin{figure}
\begin{center}
\includegraphics[width=0.8\textwidth]{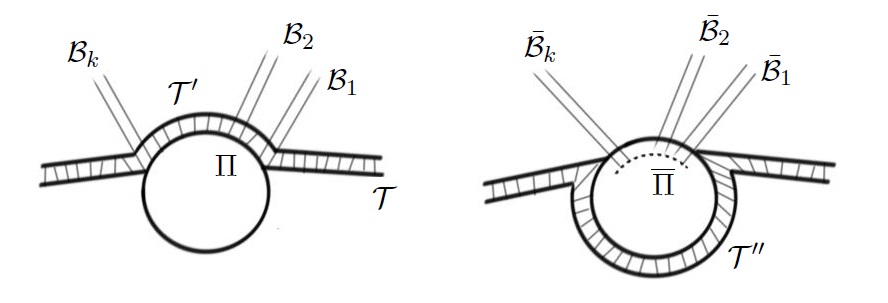}
\end{center}
\caption{Transposition of a $\theta$-band and a disk}\label{tran}
\end{figure}

We describe the \label{transpo} {\it transposition} (band moving construction in the terms of \cite{SBR}) of the disk and the band as follows.
We have a word $V\equiv tW_1t\dots tW_{k-1}t$
written on the top of the subband ${\cal T'}$ of
$\cal T$, that starts on ${\cal B}_1$ and ends on ${\cal B}_k$. The bottom ${\bf q}_2$ of $\cal T'$ is the subpath
of the boundary path ${\bf q}_2{\bf q}_3$ of $\Pi$ (fig. \ref{tran}).

Note that $W_1\equiv W_2\equiv \dots\equiv W_{k-1}$ and $tW_1$ is an accepted word by Lemma \ref{simul}. Therefore one can construct a new disk $\bar{\Pi}$
with boundary label $(tW_1)^L$ and boundary ${\bf s}_1{\bf s}_2$,
where $\Lab({\bf s}_1)\equiv V$. Also one construct
an auxiliary band $\cal T''$ with top label $W_1^{-1}t^{-1}\dots t^{-1}W_1^{-1}$, where the number
of occurrences of $t^{-1}$ is $L-k$, and attach it
to ${\bf s}^{-1}_2$, which has the same label. Finally
we replace the subband $\cal T'$ by ${\cal T''}$
(and make cancellations in the new $\theta$-band $\bar{\cal T}$ if any appear). The new diagram $\bar\Gamma$
formed by $\bar{\Pi}$ and $\bar{\cal T}$ has the same
boundary label as $\Gamma$.

\begin{rk}\label{less} After the transposition, the first
$(\theta,t)$-cells of $t$-spokes ${\cal B}_1,\dots, {\cal B}_k$ are removed and the total number of common $(\theta,t)$-cells of the new spokes $\bar{\cal B}_1,\dots, \bar{\cal B}_k$ of $\bar\Pi$ and $\bar{\cal T}$ is less
than the number of common $(\theta,t)$-cells of ${\cal B}_1,\dots, {\cal B}_k$ and $\cal T$ at least by $k$.
In particular, if $k>L-k$, then
the number of $(\theta,t)$-cells in $\bar\Gamma$ is less than the number of $(\theta,t)$-cells in $\Gamma$.
This observatiom implies
\end{rk}

\begin{lemma} \label{withd} Let $\Delta$ be a minimal diagram.

(1) Assume that a $\theta$-band ${\cal T}_0$ crosses $k$
$t$-spokes ${\cal B}_1,\dots, {\cal B}_k$ starting on a disk $\Pi$, and there are no
disks in the subdiagram $\Delta_0$, bounded by these spokes, by ${\cal T}_0$ and by $\Pi$. Then $k\le L/2$.

(2) Assume that there are two disjoint $\theta$-bands $\cal T$ and ${\cal S}$ whose bottoms are parts of the boundary of a disk $\Pi$ and these bands correspond to the same rule $\theta$ (if their histories are read towards the disk).
Suppose $\cal T$ crosses $k\ge 2$  $t$-spokes starting
on $\partial\Pi$ and ${\cal S}$ crosses $\ell\ge 2$ $t$-spokes
starting on $\partial\Pi$. Then $k+\ell\le L/2$.

(3) $\Delta$ contains no $\theta$-annuli.

(4) A $\theta$-band
cannot cross a maximal $q$-band (in particular, a spoke) twice.
\end{lemma}
\proof
(1) Since every cell, except for disks, belongs to a
maximal $\theta$-band, it follows from Lemma \ref{NoAnnul} that there is a $\theta$-band $\cal T$
such that $\cal T$ crosses all ${\cal B}_1,\dots, {\cal B}_k$ and $\Delta_0$ has no cells between $\cal T$ and $\Pi$. If $k>L/2$, then by Remark \ref{less}, the transposition of $\Pi$
and $\cal T$ would decrease the number of $(\theta, t)$-cells in $\Delta$, a contradiction, since $\Delta$ is a minimal diagram.

(2) As above, let us transpose $\cal T$ and $\Pi$. This operation
removes $k$ $(\theta,t)$-cells but add $L-k$ new $(\theta,t)$-cells in $\bar{\cal T}$.
However $\ell$ $(\theta,t)$-cells of $\cal S$ and $\ell$ $(\theta,t)$-cells of $\bar{\cal T}$
will form mirror pairs, and so after cancellations one will have at most $L-k-2\ell$ new
$(\theta,t)$-cells. This number is less than $k$ if $k+\ell>L/2$ contrary to the minimality
of the original diagram. Therefore $k+\ell\le L/2$.

(3) Proving by contradiction, consider the
subdiagram $\Delta'$ bounded by a $\theta$-annulus.
It has to contain disks by Lemma \ref{NoAnnul}. Hence
it must contain spokes ${\cal B}_1,\dots, {\cal B}_{L-3}$ introduced in Lemma \ref{extdisc}. But this
contradits to item (1) of the lemma since $L-3>L/2$.

(4) The argument of item (3) works if there is
a subdiagram $\Delta'$ of $\Delta$ bounded by an $q$-band and a $\theta$-band.
\endproof

The transposition transformation will be used for extrusion of disks from quasi-trapezia. The definition
of a \label{quasitr} {\it quasi-trapezium} sounds as the definition of trapezium, but quasi-trapezia may contain disks. (So a quasi-trapezium without disks is a trapezium.)

\begin{lemma} \label{qt} Let a minimal diagram $\Gamma$ be a quasi-trapezium with standard factorization of the boundary as ${\bf p_1^{-1}q_1p_2q_2^{-1} }$. Then there is a
diagram $\Gamma'$ such that

(1) the boundary of $\Gamma'$ is ${\bf (p'_1)^{-1}q'_1p'_2(q'_2)^{-1} }$, where $\Lab({\bf p'}_j)\equiv\Lab({\bf p}_j)$ and $\Lab({\bf q'}_j)\equiv\Lab({\bf q}_j)$ for $j=1,2$;

(2) the numbers of hubs and
$(\theta, q)$-cells in $\Gamma'$ are the same as in $\Gamma$;

(3) the vertices $({\bf p'}_1)_-$ and $({\bf p'}_2)_-$ (the vertices $({\bf p}_1')_+$ and $({\bf p'}_2)_+$) are connected by
a simple path ${\bf s}_1$ (by ${\bf s}_2$, resp.)
such that we have three subdiagrams $\Gamma_1,\Gamma_2,\Gamma_3$ of $\Gamma'$,
where $\Gamma_2$ is a trapezium
with standard factorization of the boundary ${\bf (p'_1)^{-1}s_1p'_2s_2^{-1} }$ and all cells
of the subdiagrams $\Gamma_1$ and $\Gamma_3$
with boundaries $\bf q'_1s_1^{-1}$ and $\bf s_2 (q'_1)^{-1}$ are disks;

(4) All maximal $\theta$-bands of $\Gamma$ and all
maximal $\theta$-bands of $\Gamma_2$ have the same number
ot $(\theta,t)$-cells (equal for $\Gamma$ and $\Gamma_2$) .
\end{lemma}

\proof Every maximal $\theta$-band of $\Gamma$ must
connect an edge of $\bf p_1$ with an edge of ${\bf p}_2$; this follows from Lemma \ref{withd} (3). Hence
we can enumerate these bands from bottom to top:
${\cal T}_1,\dots,{\cal T}_h$, where $h=|{\bf p}_1|=|{\bf p}_2|$.

If $\Gamma$ has a disk, then by Lemma
\ref{extdisc}, there is a disk $\Pi$ such that at
least $L-3$  $t$-spokes of it end on ${\bf q}_1$ and ${\bf q}_2$, and there are no disks between the spokes
ending on ${\bf q}_1$ (and on ${\bf q}_2$). By Lemma \ref{withd} (2), at least $L-3 -L/2\ge 2$ of these spokes must end on  ${\bf q}_1$ (resp., on ${\bf q}_2$).

If $\Pi$ lies
between ${\cal T}_j$ and ${\cal T}_{j+1}$,
then the number of its $t$-spokes crossing
${\cal T}_j$ (crossing ${\cal T}_{j+1}$) is at least $2$.
So one can make a transposition of $\Pi$ with each of these two $\theta$-bands. So we can move the disk toward $\bf q_1$ (or toward $\bf q_2$) until the disk
is extruded from the quasi-trapezium. (We use that
if $k$ $t$-spokes ${\cal B}_1,\dots,{\cal B}_k$ of $\Pi$ end on $\bf q_1$, then after transposition, we again have $k$
$t$-spokes $\bar{\cal B}_1,\dots,\bar{\cal B}_k$ of $\bar\Pi$ ending
on $\bf q_1$. - See the notation of Remark \ref{less}.)

No pair ${\cal T}_j$ and ${\cal T}_{j+1}$ corresponds to two mutual inverse letters of the history. This is clear if there are no discs between these $\theta$-bands. If there is a disk, then this is
impossible too by Lemma \ref{qt} (2) since one could choose a disk
$\Pi$ as in the previous paragraph. So the projection of the
label of ${\bf p}_1$  on the history is reduced.

Let us choose $i$ such that the number $m$ of $(\theta,t)$-cells in ${\cal T}_i$ is minimal.
It follows that $\Gamma$ has at least $hm$ $(\theta,t)$-cells.

If the disk $\Pi$ lies above ${\cal T}_i$,
we will move it upwards using transpositions. So after
a number of transpositions all such (modified) disks
will be placed above the $\theta$-band number $h$
and form the subdiagram $\Gamma_1$. Similarly
we can form $\Gamma_3$ moving other disks downwards.

In the obtaining diagram
$\Gamma_2$ lying between $\Gamma_1$ and $\Gamma_3$,
every $\theta$-band is reduced by the definition of transpositions. The neighbor maximal $\theta$-band
of $\Gamma_2$ cannot correspond to mutual inverse
letters of the history since the labels of ${\bf p}_1$ and
${\bf p}'_1$ are equal.
It follows that the diagram $\Gamma_2$ (without disks)
is a reduced diagram, and so it is a trapezium of height
$h$.

The $\theta$-band ${\cal T}_i$ did not participate
in the  transpositions. Therefore it is a maximal $\theta$-band of $\Gamma_2$. Hence the trapezium $\Gamma_2$ contains exactly $mh$ $(\theta,t)$-cells, which does not
exceed the number of $(\theta,t)$-cells in $\Gamma$. In fact these two numbers are equal since $\Gamma$ is
a minimal diagram. So every maximal $\theta$-band
of $\Gamma$ and every maximal $\theta$-band of $\Gamma_2$ has $m$ $(\theta,t)$-cells.

\endproof

We say that a history $H$ is \label{standh} {\it standard} if there is a standard trapezium with history $H$.

\begin{df}\label{shf}
Suppose we have a disk $\Pi$ with boundary label $(tW)^L$ and $\cal B$ be a $t$-spoke starting on $\Pi$.
Suppose there is a subband $\cal C$ of $\cal B$, which also starts on $\Pi$ and has a standard history $H$, for which the word $tW$ is $H$-admissible.
Then we call the $t$-band $\cal C$ a {\it shaft}.

For a constant $\lambda\in [0;1/2)$ we also define
a stronger concept of $\lambda$-shaft at $\Pi$ as follows.
A shaft $\cal C$ with history $H$ is a $\lambda$-shaft    if
for every factorization of  the history $H\equiv H_1H_2H_3$, where $||H_1||+||H_3||<\lambda ||H||$, the middle part $H_2$ is still a standard history.
(So a shaft is a $0$-shaft).
\end{df}

\begin{lemma} \label{str} Let $\Pi$ be a disk in a minimal diagram
$\Delta$ and $\cal C$ be a $\lambda$-shaft at
$\Pi$ with history $H$. Then $\cal C$ has no factorizations  ${\cal C}=
{\cal C}_1{\cal C}_2{\cal C}_3$ such that

(a) the sum of lengths of ${\cal C}_1$ and ${\cal C}_3$ do not exceed $\lambda ||H||$ and

(b) $\Delta$ has a quasi-trapezium $\Gamma$
such that top (or bottom) label of $\Gamma$ has $L+1$ occurrences of $t$ and ${\cal C}_2$ starts on  the bottom and ends on the top of $\Gamma$.
\end{lemma}

\proof Proving by contradiction, we first replace $\Gamma$ by a trapezium $\Gamma'$ according to Lemma \ref{qt}. The transpositions used for this goal do not affect neither
$\Pi$ nor $\cal C$ since $\cal C$ crosses all the
maximal $\theta$-bands of $\Gamma$. Also one can replace $\Gamma'$ by a trapezium with shorter base and so we assume that the base of it starts and ends with letter $t$.

\begin{figure}
\begin{center}
\includegraphics[width=0.9\textwidth]{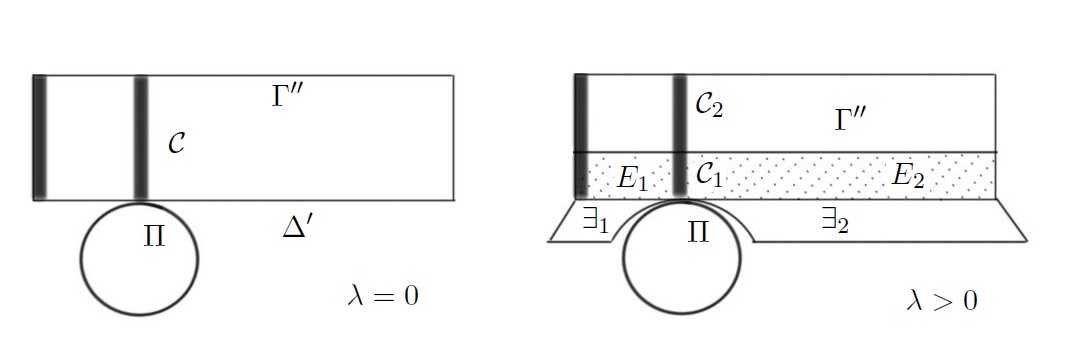}
\end{center}
\caption{Lemma \ref{str}.}\label{shft}
\end{figure}

For the beginning,  we assume that $\cal C$ is a shaft
(i.e. $\lambda=0$). Then it follows from the definition
of shaft and Remark \ref{restore} that ${\bf bot} (\Gamma')$ is labeled by $(tW)^Lt$, where $ (tW)^L$ is the boundary label of $\Pi$.  One can remove the first or the last maximal
$t$-band from $\Gamma'$ and obtain a subtrapezium
$\Gamma''$ whose boundary label coincides with
the label of $\partial\Pi$ (up to cyclic permutation),
and $\partial\Gamma''$ shares a $t$-edge with $\partial\Pi$ (fig.\ref{shft} with $\lambda=0$). It follows that
the subdiagram $\Delta'=\Pi\cup\Gamma''$ has boundary label freely equal to $\Lab({\bf top}(\Gamma''))$.
However $\Lab({\bf top}(\Gamma'')\equiv (tW')^L$,
where $tW'=(tW)\cdot H$ by Lemma \ref{simul}, and so
there is a disk $\Pi'$ with boundary label $(tW')^L$.
Therefore the subdiagram $\Delta'$ can be replaced
by a single disk. So we decrease the number of $(\theta,t)$-cells contrary to the minimality
of $\Delta$.

Now we consider the general case, where ${\cal C}=
{\cal C}_1{\cal C}_2{\cal C}_3$. As above, we  replace
$\Gamma$ by a trapezium $\Gamma'$ and obtain a trapezium $\Gamma''$
after removing of one $t$-band in $\Gamma'$.
To obtain a contradiction, it suffices to consider
the diagram $\Delta'=\Pi\cup {\cal C}_1{\cal C}_2\cup\Gamma''$ (forgetting of the complement of $\Delta'$ in $\Delta$) and find another diagram $\Delta''$ with one disk and fewer $(\theta,t)$-cells such that
$\Lab(\partial\Delta'')=\Lab(\partial\Delta')$ in the free group.

Since both histories $H$ and $H_2$ (and so $H_1H_2$)
are standard, one can enlarge $\Gamma''$ and construct
a trapezium $\Gamma'''$ with history $H_1H_2$.
(The added parts $E_1$ and $E_2$  are dashed in figure \ref{shft} with $\lambda >0$).
Note that we add $<\lambda ||H||L$ new $(\theta,t)$-cells since every maximal $\theta$-band of $\Gamma'''$
has $L$ such cells. As in case $\lambda=0$, this trapezium $\Gamma'''$ and the disk $\Pi$ can be replaced
by one disk $\Pi'$. However to obtain the boundary label equal to $\Lab(\partial\Delta')$, we should
attach the mirror copies $\exists_1$ and $\exists_2$ of $E_1$ and $E_2$ to $\Pi'$.
The obtained diagram $\Delta''$ has at most
$\lambda||H_1||L$ $(\theta,t)$-cells, while $\Delta'$
has at least $||H_2||L\ge (1-\lambda)||H||$ $(\theta,t)$-cells. Since $\lambda<1-\lambda$,
we have the desired contradiction.
\endproof

Lemma \ref{str} will be used to obtain a linear bound, in terms
of perimeter $|\partial\Delta|$, for the sum of lengths $\sigma=\sigma_{\lambda}(\Delta)$ of all
$\lambda$-shafts in a minimal diagram $\Delta$, which make possible
to exploit $\sigma_{\lambda}$ as an inductive parameter along with $|\partial\Delta|$. One more tool needed to linearly bound $\sigma_{\lambda}$,
is a combinatorial proposition of two finite systems
of disjoint segments on Euclidean plane proved in the next section.

\section{Designs in topological disk}\label{des}

\begin{figure}
\begin{center}
\includegraphics[width=0.45\textwidth]{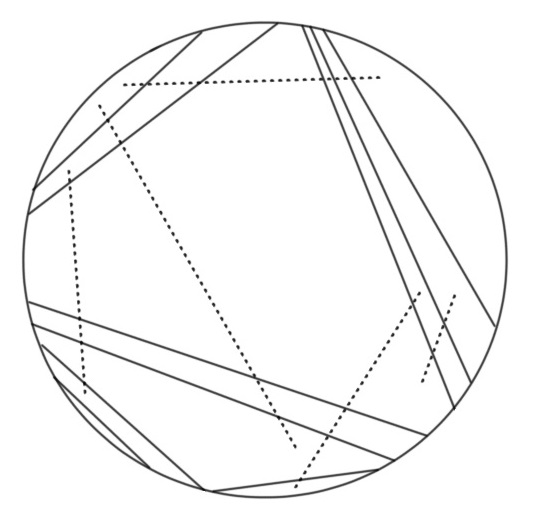}
\end{center}
\caption{Design}\label{desig}
\end{figure}

By $\cal D$, we denote Euclidean closed disk of radius $1$.
Let $\bf T$ be a finite set of disjoint \label{chord} chords (plain lines in fig. \ref{desig}) and $\bf Q$ a finite set of disjoint simple curves in $\cal D$ (dotted lines in fig. \ref{desig}). One may think of a curve as a non-oriented broken line, i.e. it is built from finitely many finite segments. To distinguish
the elements from $\bf T$ and $\bf Q$, we will say that the elements of $\bf Q$ are \label{arc} {\it arcs}.

We shall assume that the arcs belong to the open disk $D^o$, an arc may cross a  chord transversally  at most once, and the intersection point cannot coincide with one of the two ends of an arc.

Under these assumptions, we shall say that the pair
$(\bf T, Q)$ is a \label{desn} {\it design} .

By definition, the \label{lena} {\it length} $|C|$ of an arc $C$
is the number of the chords crossing $C$. The term
{\it subarc} will be used in natural way; obviously
one has $|D|\le |C|$ if $D$ is a subarc of an arc $C$.

We say that an arc $C_1$ is \label{paralla} {\it parallel} to an arc $C_2$
and write $C_1\parallel C_2$ if every chord (from $\bf T$) crossing $C_1$ also crosses $C_2$. So the relation
$\parallel$ is transitive. (The arc of length $2$ is parallel to
the arc of length $5$ in fig. \ref{desig}.)

\begin{df} \label{propP} Given  $\lambda\in (0;1)$ and an integer
$n\ge 1$, Property $P(\lambda,n)$ of a design
says that for any $n$ different arcs $C_1,\dots, C_n$,
there exist no subarcs $D_1,\dots, D_n$, respectively,
such that $|D_i|>(1-\lambda)|C_i|$ for every $i=1,\dots, n$ and $D_1\parallel D_2\parallel\dots\parallel D_n$.
\end{df}

By definition, the length $\ell(\bf Q)$ of the set
of arcs $\bf Q$ is defined by the equality
\begin{equation}
\ell({\bf Q})= \sum_{C\in \bf Q} |C|.
\end{equation}

The number of chords will be denoted by $\#\bf T$.
The goal of this subsection is to prove the following

\begin{theorem} \label{design} There is a constant $c=c(\lambda,n)$
such that for any design $(\bf T,Q)$ with Property
$P(\lambda,n)$, we have

\begin{equation}\label{QT}
\ell({\bf Q})\le c(\#\bf T).
\end{equation}
\end{theorem}

To prove Theorem \ref{design}, we may assume that $\bf Q$ has no arcs of length $0$ and that every chord is crossed by an arc. Also we may assume that
$\#\bf T>1$ since otherwise all the arcs are parallel,
and Property $P(\lambda, n)$ implies that the number
of arcs is at most $n-1$; therefore one can take $c=n-1$.

Every chord $T$ divides the disk $\cal D$ in two half-disks. If one of these half-disks contains no
other chords, we call the chord $T$ \label{peric} {\it peripheral}
and denote the \label{perif} {\it peripheral} half-disk (without chord) by $O_T$.

An arc $D$ is called an \label{extena} {\it extension} of an arc $C$
if $C$ is a subarc of $D$. (An extension need not be an element of $\bf Q$.) We will consider only extensions
of $C\in \bf Q$ such that replacing $C$ by $D$ we
again obtain a design $(\bf T, Q')$ (but  Property
$P(\lambda, n)$ can be violated for the new design).

An arc $C$ of a design is called {\it maximal} if there exists no extension $D$ of $C$ with $|D|>|C|$.

\begin{lemma} \label{extarc} Let $(\bf T,Q)$ be a design with $\# \bf T\ge 1$. Then
all the arcs $C_1, C_2,\dots$ from $\bf Q$ have maximal extensions $D_1, D_2,\dots $ forming a set of arcs $\bf Q'$ such that the design $(\bf T, Q')$ has the following property: for every arc $D_i$, its ends
belong to two different peripheral half-disks.
\end{lemma}

\proof Since no arc can cross a chord twice and the set of chords is finite, there is a system of maximal arcs $D_1, D_2,\dots $ such that every $D_i$ is an extension of $C_i$. It suffices to prove
that the ends of every $D_i$ belong to peripheral half-disks.

Assume that we have an end $o$ of an arc $D=D_i$
and $o$  belongs to no peripheral half-disk. Let us choose the direction for $D$ toward $o$, and assume that $T$ is the last chord crossed by $D$. Let $H$
be the half-disk defined by $T$, where the point $o$ belongs to.
The half-disk $H$ is not peripheral, and so it contains
a chord $T'\ne T$. None of such $T'$ could be crossed by $D$ because otherwise $D$
had to cross the chord $T$ at least twice. We may assume that $T'$ is the closest to $o$ in the sense that one can connect $o$ and $T'$ by a path ${\bf p}$
inside $H$, ${\bf p}$ consequently intersect $\ell$ arcs $D_{i_1},\dots, D_{i_{\ell}}$ from the set $\{D_1, D_2,\dots\}$, and the number $\ell$ is minimal.

\begin{figure}
\begin{center}
\includegraphics[width=0.9\textwidth]{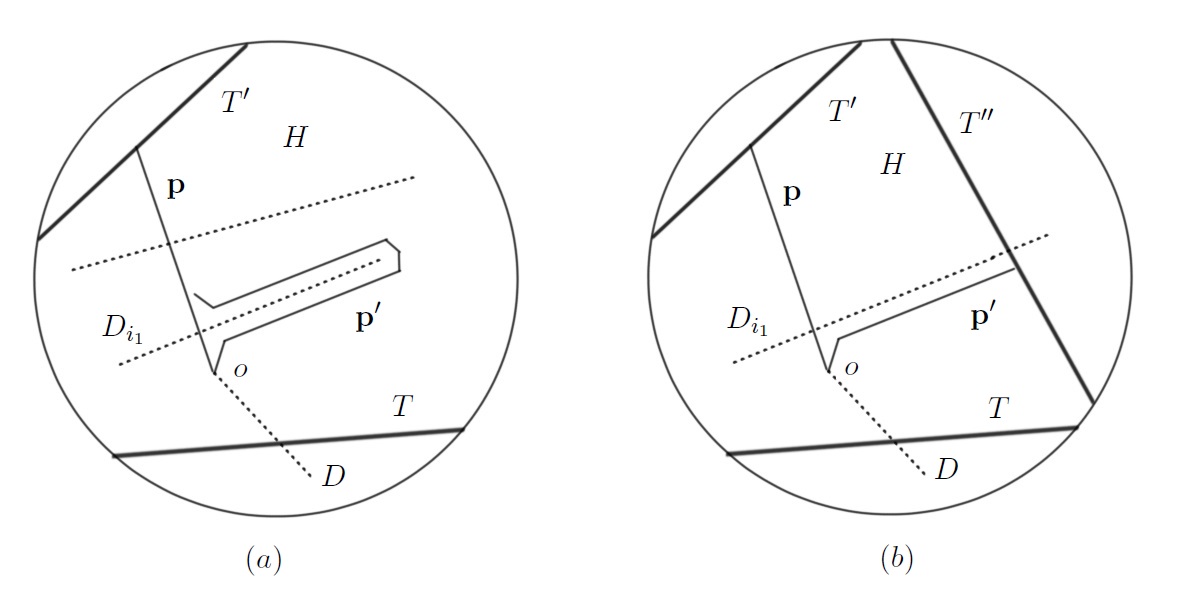}
\end{center}
\caption{Lemma \ref{extarc}}\label{Pic13}
\end{figure}

If $\ell=0$, then using $\bf p$ one could extend $D$
so that the extension crosses $T'$. If $\ell>0$ we
come to a contradiction too. Indeed, let us consider a small neighborhood $U$ of $D_{i_1}$, which contains
neither points of other arcs, nor the boundary points
of $\cal D$. If one can bypath  $D_{i_1}$ in $U$ moving around it, thereby replacing $\bf p$ by a path ${\bf p}'$
having $\ell-1$ intersections with arcs (see fig. \ref{Pic13} (a)), we get a contradiction. Otherwise going around $D_{i_1}$ in $U$ clockwise or counter-clockwise one will cross
an interval of some chord $T''\ne T$ from $H$ (fig \ref{Pic13} (b)). Then one can use $U$ to connect $o$
with  the chord $T''$ and an extention of $D$ crosses $T''$, a contradiction again.

\endproof

To continue the proof of Theorem \ref{design}, we
modify the number $\#\bf T$, taking every chord $T$
with a weight $\nu(T)$.  By definition,

$\nu(T)=1$, if $T$ crosses exactly one arc from $\bf Q$,

$\nu(T)=2$, if $T$ crosses exactly two arcs,

...

$\nu(T)=2n-2$, if $T$ crosses exactly $2n-2$ arcs,

$\nu(T)=2n-1$, if $T$ crosses {\it at least} $2n-1$ arcs.

Clearly, we have

$$\#{\bf T} \le \nu({\bf T})\; \stackrel{def}{=} \;  \sum_{T\in {\bf T}}\nu(T)\le(2n-1)\#\bf T. $$

Therefore, instead of (\ref{QT}), it  suffices to prove the following inequality

\begin{equation}\label{nu}
\ell({\bf Q})\le  d\;\nu (\bf T)
\end{equation}
for some $d=d(\lambda,n)>0$. We will prove (\ref{nu}) for any $d\ge \frac{1}{\lambda}$ by induction on the number of arcs in $\bf Q$. If there is only one arc, then it is nothing to prove since $\ell({\bf Q})\le
\#\bf T $ in this case. So we will assume that there
are at least two arcs.

Assume that there is an arc $C$, which, after one choose a direction, can be factorized as $C_1C_2C_3$,
where

(a) $|C_1|+|C_3|\ge\lambda |C|$ and

(b) every chord crossing $C_1$ or $C_3$ has weight at most $2n-2$.

Let us remove $C$ from $\bf Q$. We obtain a new
design $(\bf T',Q')$. Here $\bf T'$ has the same chords as $\bf T$, but their weights change. Obviously, Property $P(\lambda,n)$ holds for the
design $(\bf T', Q')$. Hence by inductive hypothesis,
the inequalty $\ell({\bf Q'})\le d\;\nu(\bf T')$ is true.

Also we have $\ell(\bf Q)=\ell({\bf Q'})+|C|$ and
$$\nu({\bf T})\ge \nu({\bf T'})+|C_1|+|C_3|\ge \nu({\bf T'})+\lambda|C|$$
since all chords crossing $C_1$ and $C_3$ loss their weight by $1$. It follows that
$$\ell({\bf Q})=\ell({\bf Q'})+|C|\le d\;\nu({\bf T'})+|C|
\le d(\nu({\bf T'})+\lambda|C|)\le d\;\nu(\bf T),$$
as desired, since $d\ge\lambda^{-1}$.

It remains to obtain a contradiction assuming that no arc $C\in \bf Q$ has a factorization with Properties (a) and (b). In other words, every arc $C\in \bf Q$ has a subarc $D$ maximal with respect to  the following properties:

(A) $|D|>(1-\lambda)|C|$ and

(B) the first and the last chords crossing $D$ have weight $2n-1$.

We denote by $(\bf\bar T,\bar Q)$ the design obtained
after the transition $C\to D$ for every arc $C$. Here we assume that $\bf\bar T$ contains the chords from $\bf T$, which cross some arcs from $\bf\bar Q$.
Observe that all chords of weight $2n-1$ from $\bf T$ have the same weight in $\bf\bar T$, as it follows
from the definition of $\bf\bar Q$. (We do not claim
Property $P(\lambda,n)$ for $(\bf\bar T,\bar Q)$.)

Let $\bar T$ be a peripheral chord from $\bf\bar T$.
It is crossed by (at least) $2n-1$ arcs $D_1,...,D_{2n-1}$ since $\bar T$ is the first/last
chord crossing the arcs.

By Lemma \ref{extarc}, there are maximal extensions
$\tilde D_1,\dots, \tilde D_{2n-1}$ of $D_1,...,D_{2n-1}$, respectively. Moreover, such
extensions can be constructed for {\it every} peripheral
chord, and the two ends of every extension must belong
to different peripheral half-disks of the design $(\bf\bar T,\bar Q)$.

Suppose one can choose $n$ extension, say $\tilde D_1,\dots, \tilde D_n$ starting in $O_{\bar T}$ and
ending in the same half-disk $O_{T'}$ ($T'\in\bf\bar T$).
Then every chord of $\bf\bar T$ crossing $\tilde D_i$
has to cross $\tilde D_j$ for $1\le i,j\le n$. The same is true for the chords of $\bf T$ disposed between $\bar T$ and $T'$. Since
$\tilde D_i$ starts with $D_i$ and $\tilde D_j$ starts
with $D_j$, the inequality $|D_i|\le |D_j|$ implies that
every chord of $\bf T$ crossing $D_i$ has to cross $D_j$ too. Therefore assuming that $|D_1|\le\dots\le |D_n|$, we have $D_1\parallel\dots\parallel D_n$.
However this violates Property $P(\lambda,n)$ for
the arcs $C_1,\dots C_n$ since $|D_i|> (1-\lambda)|C_i|$
for every $i=1,\dots,n$.

Thus, there are no $n$ arcs among $\tilde D_1,\dots,\tilde D_{2n-1}$ connecting the half-disk
$O_{\bar T}$ with the same peripheral half-disk. We see
that for every peripheral half-disk $O_{\bar T}$ one
can find three arcs, say $\tilde D_i, \tilde D_j,\tilde D_k$ connecting $O_{\bar T}$ with three {\it different}
peripheral half-disks.

Now let us choose one vertex in every peripheral half-disk (e.g, on the boundary of the disk $\cal D$)
and connect it with three other vertices  using the triples of arcs obtained in the previous paragraph.
We obtain an outerplanar graph with at least four vertices, where every vertex
has degree at least $3$. However there exist no such
graphs (\cite{Har}, Corollary 11.9). The obtained contradiction completes the proof of Theorem \ref{design}.  $\Box$

\begin{rk} One may allow the ends of arcs to belong
to the boundary of $\cal D$, and then the same inequalty (\ref{QT}) holds since one can cut off the
ends of every arc $C$ preserving the length of $C$.
\end{rk}

\medskip

Let us have a parameter $\lambda\in [0,1/2)$. For every
$t$-spoke $\cal B$ of a minimal diagram $\Delta$, we choose the
$\lambda$-shaft of maximal length in it (if a $\lambda$-shaft exists). If $\cal B$
connects two disks $\Pi_1$ and $\Pi_2$, then there can
be two maximal $\lambda$-shafts: at $\Pi_1$ and  at $\Pi_2$.
We denote by \label{slam} $\sigma_{\lambda}(\Delta)$ the sum of lengths of all $\lambda$-shafts in this family.

\begin{lemma} \label{clam}There is a constant $c=c(\lambda)$ such that
$\sigma_{\lambda}(\Delta)\le c |\partial\Delta|$ for
every minimal diagram $\Delta$ over the group $G$.
\end{lemma}

\proof Let us associate the following design with $\Delta$. We say that the middle lines of the maximal $\theta$-bands (they cross $\theta$-edges of the bands
in the middle points) are the chords and the middle
lines of the maximal $\lambda$-shafts are the arcs.
Here we use two disjoint middle lines if two maximal
$\lambda$-shafts share a $(\theta,t)$-cell. By Lemma \ref{withd}
(3), (4), we obtain a design, indeed.

Observe that the length $|C|$ of an arc is the
number of cells in the $\lambda$-shaft and $\#{\bf T}\le |\partial\Delta|/2$ since every maximal $\theta$-band
has two $\theta$-edges on $\partial\Delta$.

Thus, by Theorem \ref{design}, it suffices to show
that the constructed design satisfies  the condition
$P(\lambda,n)$, where $n$ does not depend on $\Delta$.

Let $n=2L+1$. If Property $P(\lambda,n)$ is violated, then we have $n$ maximal $\lambda$-shafts
${\cal C}_1,\dots, {\cal C}_n$ and a subband $\cal D$ of
${\cal C}_1$, such that $|{\cal D}|>(1-\lambda)|{\cal C}_1|$, and every maximal $\theta$-band crossing $\cal D$ must cross each of
${\cal C}_2,\dots,{\cal C}_n$. (Here $|\cal B|$ is the length of a $t$-band $\cal B$.) It follows that  each of these $\theta$-bands crosses at least $L+1$
maximal $t$-bands. (See Lemma \ref{withd} (3,4). Here we take
into account that the same $t$-spoke can generate two
arcs in the design.) Hence using the $\lambda$-shaft
${\cal C}_1$ one can construct a quasi-trapezium of height
$|\cal D|$, which contradicts  the statement of Lemma
\ref{str}.
\endproof

\section{Upper bound for $G$-areas of diagrams over the group $G$.}\label{ub}

By definition, the \label{areaGdi} $G$-{\it area of a disk} $\Pi$ is just the minimum of areas of diagrams over the presentation (\ref{rel1},\ref{rel3}) of $G$ having the same
boundary label as $\Pi$.

\begin{lemma} \label{disk} There is a constant $c_6$ such that the $G$-area of any disk does not exceed
$c_6F(|\partial\Pi|)$.
\end{lemma}

\proof The disk $\Pi$ can be built of a hub and $L$ standard accepting trapezia over $M$.
By Lemma \ref{gtime}, and the definition of the functions $f(n), g(n)$, there are such trapezia of  height $O(||\partial\Pi||/L)g(||\partial\Pi||/L)=O(|\partial\Pi|)g(|\partial\Pi|)$.
The step history of these trapezia has length $O(f(|\partial\Pi|)^3)=O(|\partial\Pi|)$ by Lemma \ref{gtime}.
Therefore the length of every $\theta$-band in
it is $O(|\partial\Pi|)$ by Lemmas \ref{stand}
and \ref{ochev} (a, d). The statement of the lemma follows.
\endproof

By definition, the \label{areaGmd} $G$-area of a minimal diagram $\Delta$ over $G$ is the sum of $G$-areas of its disks plus the $G$-area of the complement. For the complement, as in Subsection \ref{qub}, we consider a family $\bf S$ of big subtrapezia and single
cells of $\Delta$ such that every cell of $\Delta$ belongs to
a member $\Sigma$ of this family, and if a cell $\Pi$ belongs to different $\Sigma_1$ and $\Sigma_2$ from $\Sigma$, then both $\Sigma_1$ and $\Sigma_2$ are big subtrapezia of $\Delta$ with bases $xv_1x$, $xv_2x$, and $\Pi$ is an $(\theta,x)$-cell.)
Hence the statement of Lemma \ref{GA} holds for minimal diagrams
over $G$ as well.

\medskip

We want to prove that for big enough constants $N_3$
and $N_4$, $\area(\Delta)\le N_4 F(n+\sigma_{\lambda}(\Delta))
+N_3\mu(\Delta)g(n)$ for every minimal diagram $\Delta$
with perimeter $n$. To prove this property by induction,  we have to consider a large class of diagrams, called weakly minimal.

Let $\cal C$ be a cutting $q$-band of a reduced diagram $\Delta$ with disks, i.e. it starts and ends on $\partial\Delta$ and cut up the diagram. We call $\cal C$ a \label{stb} {\em stem band}, if it either
a rim band of $\Delta$ or both components of $\Delta\backslash \cal C$ contain disks.
The (unique) maximal subdiagram  of $\Delta$, where every cutting $q$-band is stem, is called
{\it the stem} \label{stD}  $\Delta^*$ of $\Delta$. It is obtained by removing all \label{crown} {\it crown} cells from
$\Delta$,  where a cell $\pi$ is called crown, if it belongs to the component $\Gamma$ defined by a cutting $q$-band $\cal B$, where $\Gamma$ contains no disks and $\pi$ is not in $\cal B$. In particular, all the disks and $q$-spokes of $\Delta$ belong in the stem $\Delta^*$.
The stem of a diagram without disks is empty.

\begin{df} \label{wminimald} A reduced diagram $\Delta$ (with disks) is called weakly minimal
if the stem
$\Delta^{*}$ is a minimal diagram.
\end{df}

\begin{lemma} \label{wmin} (a) If $\Delta_1$ is a subdiagram of weakly minimal diagram $\Delta$,
then $\Delta_1$ is weakly minimal and $\Delta_1^*\subset\Delta^*$;

(b) under the same assumption, we have $\sigma_{\lambda}(\Delta_1^*)\le \sigma_{\lambda}(\Delta^*)$;

(c) There is a constant $c=c(\lambda)$ such that
$\sigma_{\lambda}(\Delta^*)\le c |\partial\Delta|$ for
every weakly minimal diagram $\Delta$ over the group $G$;

(d) If a diagram $\Delta$ has a cutting $q$-band $\cal C$ and two components
$\Delta_1$ and $\Delta_2$ of the compliment of $\cal C$ such that $\Delta_1\cup\cal C$
is a reduced diagram without disks and ${\cal C}\cup \Delta_2$ is a weakly minimal
diagram, then $\Delta$ is weakly minimal itself;

(e) a weakly minimal diagram $\Delta$ contains no $\theta$-annuli, and
a $\theta$-band cannot cross a $q$-band of $\Delta$ twice.

\end{lemma}

\proof (a) Every crown cell $\pi$ of $\Delta$ belonging is $\Delta_1$ is crown in $\Delta_1$
since the cutting $q$-band $\cal B$ separating $\pi$ from all the disks of $\Delta$
separates (itself or the subbands of $\cal B$ in the intersection of $\cal B$ and  $\Delta_1$)
$\pi$ from $\Delta_1^*$. Therefore we have $\Delta_1^*\subset\Delta^*$, and so  $\Delta_1^*$
is minimal being a subdiagram of a minimal diagram.

(b) Now it follows from the definition of shaft, that every $\lambda$-shaft of $\Delta_1^*$ is
a $\lambda$-shaft in $\Delta^*$, which implies  inequality $\sigma_{\lambda}(\Delta_1^*)\le \sigma_{\lambda}(\Delta^*)$.

(c) If a cutting  $q$-band $\cal C$ of a reduced diagram $\Delta$ gives a decomposition $\Delta=\Gamma_1\cup{\cal C}\cup \Gamma_2$, where $\Delta_1=\Gamma_1\cup{\cal C}$ has no disks,
then every maximal $\theta$-band starting in the subdiagram $\Delta_1$ with $\cal C$
cannot ends on $\partial\Gamma_1$ by Lemma \ref{NoAnnul}. Hence $|\partial\Delta_2|\le|\partial\Delta|$ by Lemma \ref{ochev}. So removing subdiagrams as
$\Gamma_1$ from $\Delta$, we obtain by induction that $|\partial\Delta^*|\le
|\partial\Delta|$. Now the property (c) follows from Lemma  \ref{clam} applied to the minimal
subdiagram $\Delta^*$.

(d) The diagram $\Delta$ is reduced since both $\Delta_1\cup\cal C$ and $\Delta_2\cup\cal C$
are reduced subdiagrams sharing the cutting band $\cal C$. Since $\Delta_1$ has no disks,
we have $\Delta^*=(\Delta_2\cup\cal C)^*$ by the definition of stem. Therefore
the stem $\Delta^*$ is a minimal diagram and $\Delta$ is weakly minimal.

(e) The statement follows from Lemma \ref{withd} (3, 4) if the bands belong in the
stem $\Delta^*$. By the same reason, a $\theta$-band cannot cross a rim $q$-band
of $\Delta^*$ twice. It remains to assume that the bands belong to the crown of $\Delta$,
and in this case, the statement follows from Lemma \ref{NoAnnul} since the crown
is a union of disjoint reduced subdiagrams over the group $M$.
\endproof

\begin{rk} The statement (d) of Lemma \ref{wmin} fails if one replaces the words
``weakly minimal'' with ``minimal''.
\end{rk}

We  will prove that $\area(\Delta)\le N_4 F(n+\sigma_{\lambda}(\Delta^*))
+N_3\mu(\Delta)g(n)$ for every weakly minimal diagram $\Delta$
with perimeter $n$. For this goal, we will argue by contradiction in this section and study a
weakly minimal {\bf counter-example} $\Delta$ giving opposite inequality
\begin{equation}\label{ce}
\area_G(\Delta)>  N_4 F(n+\sigma_{\lambda}(\Delta^*))
+N_3\mu(\Delta)g(n)
\end{equation}
with minimal possible sum $n+\sigma_{\lambda}(\Delta^*)$.

\begin{lemma} \label{norim} The diagram $\Delta$ has no rim
$\theta$-bands with base of length at most $K$.
\end{lemma}

\proof  The functions $F(x)$ and $g(x)$ satisfy the definition given
for $\Phi(x)$ and $\phi(x)$, and the inequality (\ref{xy}) by Lemma \ref{d-x}. Hence the proof of Lemma
\ref{nori} works for the minimal counter-example over $G$ as follows.
It suffices
to replace $N_2$ and $N_1$ with $N_4$ and $N_3$, resp., replace $n$ with $n+\sigma_{\lambda}(\Delta^*)$, and notice that the subdiagram $(\Delta')^*$ is weakly minimal and $\sigma_{\lambda}((\Delta')^*)\le \sigma_{\lambda}(\Delta^*)$ by Lemma \ref{wmin}
(a,b).
\endproof

By Lemma \ref{main},
$\Delta$ has at least one disk. Applying Lemma \ref{extdisc} to the stem $\Delta^*$, we fix a disk $\Pi$ in $\Delta$ such that $L-3$ consecutive maximal $t$-bands \label{b1bL3}${\cal B}_1,\dots
{\cal B}_{L-3} $ start on $\partial\Delta$ , end on the boundary $\partial\Pi$, and for any $i\in [1,L-4]$,
there are no disks in the subdiagram  bounded by ${\cal B}_i$, ${\cal B}_{i+1},$ $\partial\Pi,$ and $\partial\Delta.$
(See fig. \ref{extd}.)

We denote by
\label{cl.} $\Psi=cl(\Pi,{\cal B}_1,{\cal B}_{L-3})$ the subdiagram without disks bounded by the spokes ${\cal B}_1$, ${\cal B}_{L-3}$
(and including them) and by subpaths of the boundaries of $\Delta$
and $\Pi,$ and call this subdiagram a \label{clove} {\it clove}. Similarly one can defined the cloves \label{cloves}
$\Psi_{ij}=cl(\Pi,{\cal B}_i,{\cal B}_j)$ if $1\le i<j\le L-3$.

\begin{lemma} \label{nocomb} The clove $\Psi=cl(\Pi,{\cal B}_1,{\cal B}_{L-3})$ has no subcombs of basic width at least $K_0$.
\end{lemma}

\proof Proving by contradiction, we may assume that there is a tight
subcomb $\Gamma$ by Lemma \ref{est'takaya} (2). Then contradiction
appears as in  Lemmas \ref{notwo}--\ref{main}, since Lemma \ref{d-x} allows us to define $\Phi(x)=F(x)$:
We may assume that there is a tight
subcomb $\Gamma$ by Lemma \ref{est'takaya} (2). Then contradiction
for the counter-example $\Delta$, which is a weakly minimal diagram over the group $G$, appears exactly as in  the proofs of Lemma \ref{main} with the following modifications of few constants and references.

We should
replace $N_2$ and $N_1$ with $N_4$ and $N_3$, replace $n$ with $n+\sigma_{\lambda}(\Delta^*)$, and notice that the value of $\sigma_{\lambda}$ does nor increase when passing from $\Delta$ to a subdiagram by Lemma \ref{wmin} (b). We should use Lemma \ref{wmin} (e) instead of Lemma
\ref{NoAnnul} used in the proofs of Lemmas \ref{notwo} - \ref{main}. The diagram $\Delta_0$
is weakly minimal because it is constructed from the reduced diagram $\Gamma_1\cup\cal Q$ over $M$ and the weakly minimal diagram $\Delta_1\cup\cal Q$ according to the assumption of Lemma \ref{wmin} (d).
\endproof

The statements of auxiliary Lemmas
\ref{notwo}, \ref{twocombs} and \ref{A0} holds as well
for the minimal counter-example over $G$. Below we use the following
analog of Lemma \ref{notwo}:

\begin{lemma} \label{notw} (1) The counter-example $\Delta$ has no two disjoint subcombs $\Gamma_1$ and $\Gamma_2$ of basic widths at most $K$ with handles ${\cal C}_1$ and ${\cal C}_2$ such that some ends of these handles are
connected by a subpath ${\bf x}$ of the boundary path of $\Delta$ with $|{\bf x}|_q\le N$.

(2) The boundary of every subcomb $\Gamma$ with basic width $s\le K$ has $2s$ $q$-edges.
\end{lemma}
$\Box$

\begin{lemma}\label{psi1}
(1) Every maximal
$\theta$-band of $\Psi$ crosses either ${\cal B}_1$ or ${\cal B}_{L-1}$.

(2)
There exists $r$, $L/2-3\le r \le L/2$,   such that the $\theta$-bands of $\Psi$
crossing \label{Br} ${\cal B}_{L-3}$
do not cross ${\cal B}_r$, and the $\theta$-bands of $\Psi$ crossing ${\cal B}_1$
do not cross ${\cal B}_{r+1}$;
\end{lemma}

\begin{figure}
\begin{center}
\includegraphics[width=1.0\textwidth]{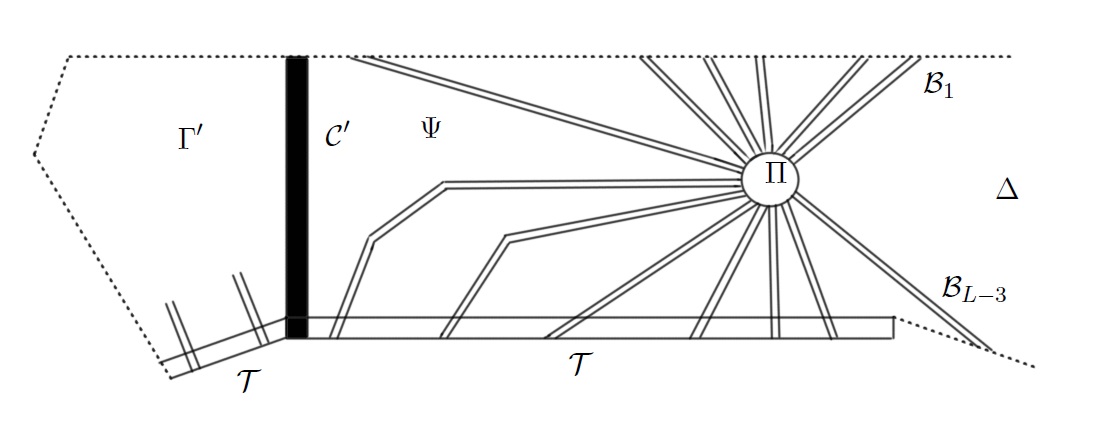}
\end{center}
\caption{Lemma \ref{psi1}}\label{lempsi}
\end{figure}

\proof (1) If the claim were wrong, then one could find a rim $\theta$-band $\cal T$ in $\Psi$, which
crosses neither ${\cal B}_1$ nor ${\cal B}_{L-3}$. By Lemma
\ref{norim}, the basic width of $\cal T$ is greater than $K$.
Since (1) a disk has $LN$ spokes, (2) no $q$-band of $\Psi$ intersects $\cal T$ twice by Lemma
  \ref{NoAnnul} (3), $\cal T$ has at least $K$  $q$-cells, and (4) $K>2K_0+LN$, there exists a maximal $q$-band $\cal C'$
 such that a subdiagram $\Gamma'$ separated from $\Psi$ by $\cal C'$ contains no edges of the spokes
  of $\Pi$ and the part of $\cal T$ belonging to $\Gamma'$ has at least $K_0$  $q$-cells (fig. \ref{lempsi}).

  If $\Gamma'$ is not a comb, and so a maximal $\theta$-band of it does not cross $\cal C',$
  then $\Gamma'$ must contain another rim band $\cal T'$ having at least $K$
  $q$-cells. This makes possible to find a subdiagram $\Gamma''$ of
  $\Gamma'$ such that a part of $\cal T'$ is a rim band of $\Gamma''$ containing at least $K_0$
  $q$-cells, and $\Gamma''$ does not contain $\cal C'$.
  Since $\area(\Gamma')>\area(\Gamma'')>\dots$ , such a procedure must stop. Hence, for some $i$, we
  obtain a subcomb $\Gamma^{(i)}$ of basic width $\ge K_0$,
  contrary to Lemma \ref{nocomb}.

  (2) Assume there is a maximal $\theta$-band $\cal T$ of $\Psi$ crossing
  the spoke ${\cal B}_1$.   Then assume that $\cal T$ is the closest
  to the disk $\Pi$, i.e. the intersection of $\cal T$ and ${\cal B}_1$ is
  the first cell of the spoke ${\cal B}_1$. If ${\cal B}_1,\dots, {\cal B}_r$ are all the spokes crossed by $\cal T$, then $r\le L/2$ by Lemma \ref{withd}.
  Since the band $\cal T$ does not cross the spoke ${\cal B}_{r+1}$, no other
  $\theta$-band of $\Psi$ crossing ${\cal B}_1$ can cross ${\cal B}_{r+1}$
  and no $\theta$-band crossing the spoke ${\cal B}_{L-3}$ can cross ${\cal B}_r$.
  The same argument shows that $r+1\ge L/2 -2$ if there is a $\theta$-band of $\Psi$
  crossing the spoke ${\cal B}_{L-3}$.
  \endproof

  For the clove $\Psi=cl(\pi,{\cal B}_1,{\cal B}_{L-3})$ in  $\Delta$,
we denote by \label{pPsi} ${\bf p}={\bf p}(\Psi)$ the common subpath of $\partial\Psi$ and
$\partial\Delta$ starting with the $t$-edge of ${\cal B}_1$ and ending
with the $t$-edge of ${\cal B}_{L-3}.$ Similarly we define the (outer)
path \label{pij} ${\bf p}_{ij}={\bf p}(\Psi)_{ij}$ for every smaller clove $\Psi_{ij}$.

  \begin{lemma} \label{2K0} Every path ${\bf p}_{i,i+1}$ ($i=1,\dots, L-4$) has less than $3K_0$ $q$-edges.
\end{lemma}

\proof Let a maximal $q$-band $\cal C$ of $\Psi$ start on ${\bf p}_{i,i+1}$ and do not end on $\Pi$. Then it has to end on ${\bf p}_{i,i+1}$ too.

If $\Gamma$ is the subdiagram
without disks separated by $\cal C$, then every maximal $\theta$-band of $\Gamma$ has to
cross the $q$-band $\cal C$ since its extension in $\Psi$ must cross either ${\cal B}_1$
or ${\cal B}_{L-3}$ by Lemma \ref{psi1}. Therefore $\Gamma$ is a comb with handle $\cal C$.

Consider the $q$-bands of this kind defining maximal subcombs $\Gamma_1,\Gamma_2,\dots\Gamma_k$
in $\Psi_{i,i+1}$. The basic width of each of them is less than $K_0$ by Lemma \ref{nocomb}.
Therefore $k\le 1$ since otherwise one can get two subcombs contradicting to Lemma \ref{notw} (1), because there are at most $N+1$ maximal $q$-bands
starting on $\partial\Pi$ in $\Psi_{i,i+1}$.
By Lemma \ref{notw} (2),
such a subcomb has at most $2K_0$ $q$-edges in the boundary. Hence there are at most $2K_0+N<3K_0$ $q$-edges in the path ${\bf p}_{i,i+1}$.
\endproof

\begin{figure}
\begin{center}
\includegraphics[width=0.8\textwidth]{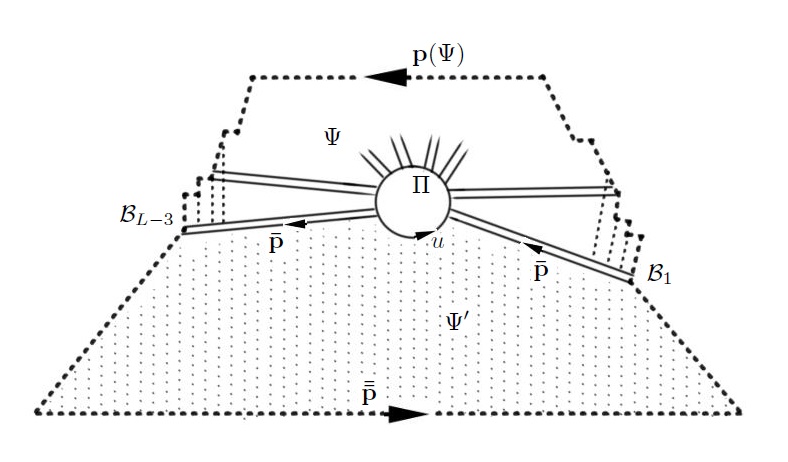}
\end{center}
\caption{Boundaries of $\Psi$ and $\Psi'$}\label{bou}
\end{figure}

We denote
by \label{bdel} $\bar\Delta$ the subdiagram formed by $\Pi$ and $\Psi$, and
denote by \label {barp} $\bf \bar p$ the path $\topp({\cal B}_1) {\bf u}^{-1}
\bott({\cal B})_{L-3}^{-1},$ where $\bf u$ is a subpath of $\partial\Pi,$ such that
$\bf \bar p$ separates $\bar\Delta$ from the remaining subdiagram \label{psi'} $\Psi'$
of $\Delta$ (fig. \ref{bou}).

Similarly we define subdiagrams \label{barij} $\bar\Delta_{ij}$, paths \label{barpij}
${\bf\bar p}_{i,j}={\bf top}({\cal B}_i) {\bf u}_{ij}^{-1}
\bott({\cal B})_{j}^{-1},$ where ${\bf u}_{ij}$ is a subpath of $\partial\Pi$, and
the subdiagrams \label{Psiij'} $\Psi'_{ij}$.

We denote by \label{histL3} $H_1,\dots, H_{L-3}$ the  histories  of the spokes ${\cal B}_1,\dots,
{\cal B}_{L-3}$ (read starting from the disk $\Pi$) and by \label{hidtl} $h_1,\dots,h_{L-3}$ their lengths, i.e. the numbers of $(\theta,t)$-cells. By Lemma \ref{psi1},
these lengths non-increase and then non-decrease as follows:

\begin{equation}\label{downup}
h_1\ge h_2\ge\dots \ge h_r;\;\; h_{r+1}\le\dots \le h_{L-3}\;\;(L/2-3\le r \le L/2),
\end{equation}
and therefore $H_{i+1}$ is a prefix of $H_{i}$ ($H_j$ is a prefix $H_{j+1}$) for $i=1,\dots, r-1$ (resp., for $j=r+1,\dots, L-4$).

Recall that the boundary label of $\partial\Pi$ has the form $W^L$,
i.e., it is the $L$-th power of an accepted word $W$.

\begin{lemma} \label{ppbar} We have the following inequalities
$$ |{\bf \bar p}_{ij}|\le h_i+h_j+ (L-j+i)|W|-1$$
and, if $i\le r$ and $j\ge r+1$, then $$|{\bf p}_{ij}|\ge |{\bf p}_{ij}|_{\theta}+ |{\bf p}_{ij}|_q\ge h_i+h_j+(j-i)N+1.$$
\end{lemma}

\proof The first iequality follows from Lemma \ref{ochev} (b) since the
path ${\bf u}_{ij}$ has $L-j+i-1$ $t$-edges. To prove the second inequality,
we observe that the path $|{\bf p}_{ij}|$ has $(j-i)N+1$ $q$-edges and it has
$h_i+h_j$ $\theta$-edges by Lemma \ref{psi1}.
\endproof

The large constants $L$ and $L_0$ are chosen so that
\begin{equation}\label{LL0}
L_0^3\le L.
\end{equation}

\begin{lemma} \label{muJ} If $j-i>L/2$, then we have $\mu(\Delta)-\mu(\Psi'_{ij})>
-2Jn(h_i+h_j)\ge -2Jn|{\bf p}_{ij}|$.
\end{lemma}
\proof
The number of $q$-edges in the path ${\bf\bar p}_{ij}$ (or in the path ${\bf u}_{ij}$)
does not exceed the similar number for ${\bf p}_{ij}$ provided $j-i\ge L/2$. Therefore any two white beads $o, o'$ of the necklace
on $\partial\Delta$, provided they both  do not belong to ${\bf p}_{ij},$ are separated by at least the same
number of black beads in the necklace for $\Delta$ as in the necklace for $\Psi'_{ij}$  (either the clockwise arc $o-o'$ includes ${\bf p}_{ij}$ or not). So such a pair  contributes
to $\mu(\Delta)$ at least the amount it contributes to $\mu(\Psi'_{ij})$.
Thus, to estimate $\mu(\Delta)- \mu(\Psi'_{ij})$ from below, it suffices to consider
the contribution to $\mu(\Psi')$ for the pairs $o, o', $ where one of the two beads
lies on ${\bf p}_{ij}$. The number of such (unordered) pairs is bounded by $n(h_{i}+h_{j})$
by Lemma \ref{psi1}. Taking into account the definition of $\mu $ of diagrams and inequalities (\ref{downup}), we get the required inequality.
\endproof

\begin{lemma} \label{epsi} If $j-i>L/2$, then  the following inequality holds: $|{\bf p}_{ij}|<(1+\varepsilon)|{\bf \bar p}_{ij}|$, where $\varepsilon=N_4^{-\frac 12}$.
Moreover, we have $|{\bf p}_{ij}|+\sigma_{\lambda}(\bar\Delta_{ij}^*)<(1+\varepsilon)|{\bf \bar p}_{ij}|$.
\end{lemma}

\proof It suffices to prove the second statement. Let $d$ be the difference \\$|{\bf p}_{ij}|+\sigma_{\lambda}(\bar\Delta_{ij}^*)-|{\bf \bar p}_{ij}|$ and
assume by contradiction that $d\ge\varepsilon |{\bf \bar p}_{ij}|$. Then \\ $d\ge
|{\bf p}_{ij}|+\sigma_{\lambda}(\bar\Delta_{ij}^*)-\varepsilon^{-1}d $,
whence
\begin{equation}\label{de}
d\ge (1+\varepsilon^{-1})^{-1}(|{\bf p}_{ij}|+\sigma_{\lambda}(\bar\Delta_{ij}^*))
\ge\frac{\varepsilon}{2} (|{\bf p}_{ij}|+\sigma_{\lambda}(\bar\Delta_{ij}^*))\ge \frac{\varepsilon y}{2},
\end{equation}
where by definition, $y=|{\bf p}_{ij}|+\sigma_{\lambda}(\bar\Delta_{ij}^*)$.

We have
$(|\partial\Delta|+\sigma_{\lambda}(\Delta^*))-(|\partial\Psi'_{ij}|+\sigma_{\lambda}((\Psi'_{ij})^*))\ge d >0$, because $|\partial\Delta|-|\partial\Psi'_{ij}|\ge |{\bf p}_{ij}|-|{\bf \bar p}_{ij}|$ and $\sigma_{\lambda}(\bar\Delta_{ij}^*)+\sigma_{\lambda}((\bar\Psi'_{ij})^*)\le \sigma_{\lambda}(\Delta^*)$. Therefore  for $x=n+\sigma_{\lambda}(\Delta^*)$, we obtain from  the minimality of the counter-example $\Delta$ that

$$\area_G(\Psi'_{ij})\le N_4 F(x-d)
+N_3\mu(\Psi'_{ij})g(n)\le N_4 F(x)- N_4F(x)x^{-1}d +N_3\mu(\Delta)g(n)$$
\begin{equation}\label{DP}
+ 2N_3Jn|{\bf  p}_{ij}|g(n)
\le N_4 F(x)
+N_3\mu(\Delta)g(n)- N_4F(x)x^{-1}d+ 2N_3Jnyg(n)
\end{equation}
by Lemma \ref{muJ}, inequality
$\sigma_{\lambda}((\Psi'_{ij})^*)\le \sigma_{\lambda}(\Delta^*)$,
and Lemma \ref{d-x}.
By Lemma \ref{ppbar}, $|{\bf\bar p}_{ij}|< |{\bf p}_{ij}|+|\partial\Pi|$, and so the perimeter  $|\partial\Psi_{ij}|$ is less
than $2|{\bf p}_{ij}|+|\partial\Pi|$. Since $|\partial\Pi|\le
L|{\bf\bar p}_{ij}|<L(|{\bf p}_{ij}|+\sigma_{\lambda}(\bar\Delta_{ij}^*))$, we have
\begin{equation}\label{Lp2}
|\partial\Psi_{ij}|<(2+L)|{\bf p}_{ij}|+L\sigma_{\lambda}(\bar\Delta_{ij}^*)\le (L+2)y.
\end{equation}

By the inequalities $N_2>N_1$, (\ref{Lp2}), Lemmas \ref{main} and \ref{mixture} (a), the $G$-area of $\Psi_{ij}$ does not exceed
\begin{equation}\label{aGP}
N_2(2+L)^2y^2+N_1\mu(\Psi_{ij})\le N_2(J+1)(2+L)^2y^2.
\end{equation}
By Lemma \ref{disk}, the $G$-area of $\Pi$ does not exceed
$c_6F(|\partial\Pi|)\le c_6 F((L+2)y)$, and by definition of the functions $f$ and $F$, there is
a constant $c_7=c_7(L)$ such that $\area_G(\Pi)\le c_7 F(y)$.

This estimate and (\ref{aGP}) give the inequality $\area_G(\bar\Delta_{ij})\le N_2(J+1)(2+L)^2y^2+c_7 F(y)$, and we
obtain with (\ref{DP}) that
$$\area_G(\Delta)\le  N_4 F(x)
+N_3\mu(\Delta)g(n) $$ $$-N_4F(x)x^{-1}d+ 2N_3Jnyg(n)
+N_2(J+1)(2+L)^2y^2+c_7 F(y).$$

To obtain the desired contradiction, it suffices to show that here, the number $T=N_4F(x)x^{-1}d/3$ is greater than each of the last three summands.
Recall that $F(x)x^{-1}=xg(x)\ge ng(n)$, $d>\varepsilon y/2$ by (\ref{de}), $\varepsilon=N_4^{-1/2}$, and so $T>2N_3Jnyg(n)$ if $N_4$ is large enough in comparison with $N_3$
and other constants chosen earlier. Also we have $T> N_2(J+1)(2+L)^2y^2$, because $x=n+\sigma_{\lambda}(\Delta^*)>
|{\bf p}_{ij}|+\sigma_{\lambda}(\bar\Delta_{ij}^*)=y$, and so
$F(x)x^{-1}d>xg(x)\varepsilon y/2\ge \varepsilon y^2/2$.
Finally, $T>c_7 F(y)$ since $$F(x)x^{-1}d>xg(x)\varepsilon y/2\ge y^2g(y)\varepsilon/2=\varepsilon F(y)/2.$$
\endproof

For every path ${\bf p}_{i,i+1}$ we will fix a shortest path ${\bf q}_{i,i+1}$ homotopic
to ${\bf p}_{i, i+1}$ in the subdiagram $\Psi_{ij}$, such that the first and the last $t$-edges of ${\bf q}_{i,i+1}$ coincide with the first  and the last $t$-edges of ${\bf p}_{i,i+1}$. For $j>i+1$ the path \label{qij} ${\bf q}_{i,j}$
is formed by ${\bf q}_{i,i+1},\dots {\bf q}_{j-1,j} $.
The following lemma is similar to the second part of Lemma
\ref{ppbar}.

\begin{lemma}\label{qq}
If $i\le r$ and $j\ge r+1$, then $$|{\bf q}_{ij}|\ge |{\bf q}_{ij}|_{\theta}+ |{\bf q}_{ij}|_q\ge h_i+h_j+(j-i)N+1.$$
\end{lemma} $\Box$

Let \label{Psi0ij} $\Psi_{ij}^0$ (let $\Psi^0$, $\Delta^0$) be the subdiagram of $\Psi_{ij}$ (of $\Psi$, of $\Delta$) obtained after replacement of the subpath ${\bf p}_{ij}$ (of $\bf p$ ) by ${\bf q}_{ij}$ (by ${\bf q}={\bf q}_{1,L-3}$, resp.) in the boundary.

\begin{lemma} \label{noq} (1) The subdiagram  $\Psi_{ij}^0$ has no maximal $q$-bands except for
the $q$-spokes
starting from $\partial\Pi$. (2) Every $\theta$-band of $\Psi_{ij}^0$ is crossed
by the path ${\bf q}_{ij}$ at most once.
\end{lemma}

\proof (1) Assume there is a $q$-band $\cal Q$ of $\Psi_{ij}^0$ starting and ending on ${\bf q}_{ij}$. Then $j=i+1$ and ${\bf q}_{i,i+1}=\bf uevfw$, where $\cal Q$ starts with the $q$-edge $\bf e$ and ends with the $q$-edge $\bf f$. Let $\cal Q$ have length $\ell$. Then $|v|\ge \ell$
since every maximal $\theta$-band of $\Psi_{i,i+1}^0$ crossing $\cal Q$ has to end on the
subpath  $\bf v$. So one has $|{\bf evf}|=\ell+2$, and replacing the subpath ${\bf evf}$
by a side of $\cal Q$ of length $\ell$ one replaces the path ${\bf q}_{i,i+1}$ with
a shorter homotopic path by Lemma \ref{ochev}. This contradicts to the choice
of ${\bf q}_{i,i+1}$, and so the first statement is proved. The prove of the second
statement is similar.

(2) Assume there is a $\theta$-band $\mathcal T$ of $\Psi_{i, i+1}^0$ starting and ending on ${\bf q}_{i,i+1}$. Then ${\bf q}_{i,i+1}=\bf uevfw$, where $\mathcal T$ starts with the $\theta$-edge $\bf e$ and ends with the $\theta$-edge $\bf f$. Moreover, one can chose $\cal T$ such that
$v$ is a side of this $\theta$--band. By Statement (1) the band $\cal T$ has less than $N$
$(\theta, q)$-cells. Therefore if $v'$ is another side of $\cal T$, we have
$|v'|_Y-|v|_Y \le 2N$. It follows from the definition of length in Subsection \ref{lf}
that $|evf|-|v'|\ge 2 - 2\delta N>1+2\delta$. Therefore, by Lemma \ref{ochev} (c), replacing the subpath $evf$ with $v'$ we decrease the length of  ${\bf q}_{i,i+1}$ at least by $1$,
a contradiction.
\endproof

It follows from Lemma \ref{psi1} that between the spokes ${\cal B}_j$ and ${\cal B}_{j+1}$ ($1\le j\le r-1$), there is a trapezium \label{Gammaj}
$\Gamma_j$ of height $h_{j+1}$ with the side $t$-bands ${\cal B}_{j+1}$ and ${\cal B}'_j$, where ${\cal B}'_j$ is the beginning of length $h_{j+1}$ of the $t$-spoke ${\cal B}_{j}$. Similarly, we have trapezia $\Gamma_j$ for $r+1\le j\le L-4$.
By Lemma \ref{noq} (2),  every trapezium $\Gamma_j$ is contained in both $\Psi_{j,j+1}$ and $\Psi_{j,j+1}^0$.
The bottoms \label{yj} ${\bf y}_j$ of all trapezia $\Gamma_j$ belong to $\partial\Pi$ and have the same label $Wt$. We will use \label{zj} ${\bf z}_j$ for the tops of these trapezia.  Since $\Gamma_j$ and $\Gamma_{j-1}$ ($2\le j\le r-1$) have the same bottom labels and the history $H_j$ is a prefix of $H_{j-1}$,
by Lemma \ref{simul}, $h_j$ different $\theta$-bands of $\Gamma_{j-1}$ form the copy \label{Gammaj'} $\Gamma'_{j}$ of the trapezium $\Gamma_j$ with top and bottom paths ${\bf z}'_j$ and ${\bf y}'_j={\bf y}_{j-1}$.

\begin{figure}
\begin{center}
\includegraphics[width=1.0\textwidth]{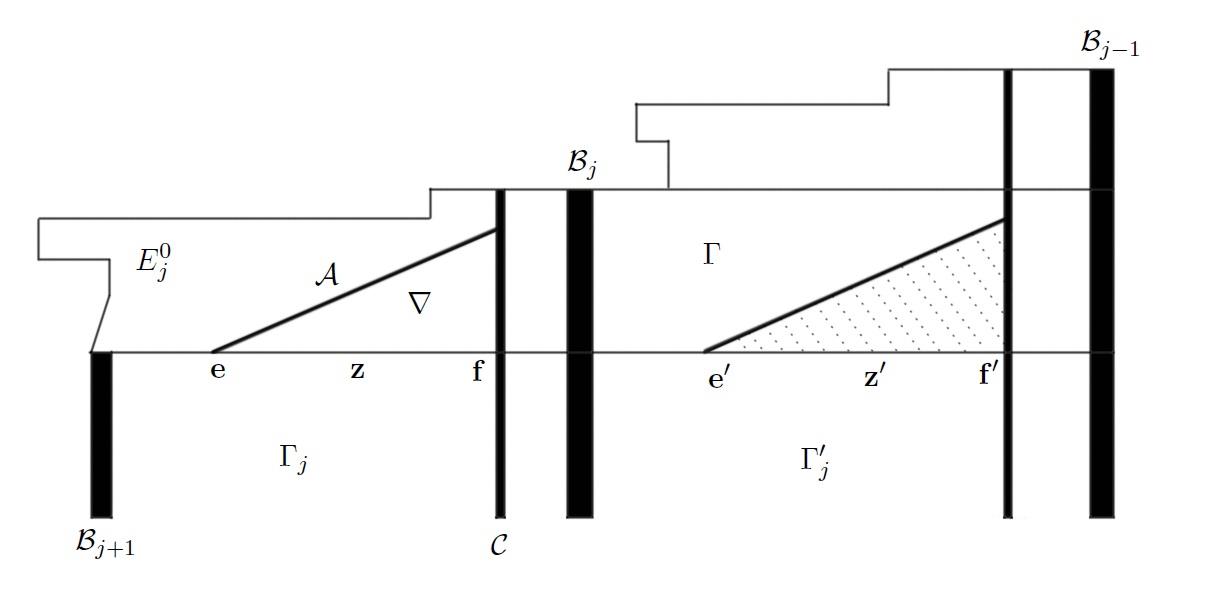}
\end{center}
\caption{Lemma \ref{copy}}\label{Pic15}
\end{figure}

We denote by \label{Ej} $E_j$
(by $E_j^0$ )
the comb formed by the maximal $\theta$-bands of $\Psi_{j,j+1}$
(of $\Psi_{j,j+1}^0$, resp.)
crossing the $t$-spoke ${\cal B}_j$  but not crossing ${\cal B}_{j+1}$ ($1\le j\le r-1$, see fig. \ref{Pic15}). Its handle ${\cal C}_j$ of
height $h_j-h_{j+1}$ is contained in ${\cal B}_j$. The boundary $\partial E_j$ (resp., $\partial E_j^0$)   consists of the
side of this handle, the path ${\bf z}_j$ and the path ${\bf p}_{j,j+1}$
(the path ${\bf q}_{j,j+1}$, respectively).

Assume that a maximal $a$-band $\cal A$ of $E_j^0$ ($2\le j\le r-1$) starts on the path ${\bf z}_j$ and ends
on a side $a$-edge of a maximal $q$-band $\cal C$ of $E_j^0$. Then $\cal A$, a part of $\cal C$ and a part ${\bf z}$ of ${\bf z}_j$
bound a comb $\nabla$.

\begin{lemma} \label{copy} There is a copy of the comb $\nabla$ in the trapezium $\Gamma=\Gamma_{j-1}\backslash\Gamma'_j$.
\end{lemma}

\proof The subpath ${\bf z}$ of ${\bf z}_j$ starts with an $a$-edge $\bf e$ and ends with a $q$-edge $\bf f$.
There is a copy ${\bf z}'$ of $\bf z$ in ${\bf z}'_j$ starting with ${\bf e'}$ and ending with $\bf f'$. Note that the $\theta$-cells $\pi$ and $\pi'$ attached to $\bf f$ and to $\bf f'$ in $\nabla$ and in $\Gamma$ are copies of each other since they correspond to the same letter of the history. Now moving from $\bf f$ to $\bf e$, we see that
the whole maximal $\theta$-band ${\cal T}_1$ of $\nabla$ containing $\pi$ has a copy in $\Gamma$. Similarly
we obtain a copy of the next maximal $\theta$-band ${\cal T}_2$ of $\nabla$, and so on.
\endproof

\begin{lemma} \label{le6} At most $N$ $a$-bands starting on the path ${\bf y}_j$ can end on
the $(\theta,q)$-cells of the same $\theta$-band. This property holds
for the $a$-bands starting on ${\bf z}_j$ too.
\end{lemma}

\proof We will prove the second claim only since the proof of the first
one is similar. Assume that the $a$-bands ${\cal A}_1,\dots, {\cal A}_s$
start from ${\bf z}_j$ and end on some $(\theta,q)$-cells of a $\theta$-band $\cal T$. Let ${\cal T}_0$ be the minimal subband of $\cal T$, where
the $a$-bands ${\cal A}_2,\dots, {\cal A}_{s-1}$ end and ${\bf \bar z}_j$ be the minimal subpath
of ${\bf z}_j$, where they start. Then by Lemma \ref{NoAnnul}, every maximal $q$-band starting
on ${\bf\bar z}_j$ has to cross the band ${\cal T}_0$ and vice versa. Hence the base of ${\cal T}_0$
is a subbase of the standard base (or of its inverse). Since every rule of $\bf M$ can
change at most $N-2$ $a$-letters  in a word with standard base, all $(\theta,q)$-cells of ${\cal T}_0$ have at most $N-2$ $a$-edges,
and the statement of the lemma follows.
\endproof

Without loss of generality, we assume that
\begin{equation}\label{L0}
h=h_{L_0+1}\ge h_{L-L_0-3}.
\end{equation}

\begin{lemma} \label{02} If $h\le L_0^2|W|_a$, then the number of trapezia $\Gamma_j$ with the
properties $|{\bf z}_j|_a\ge |W|_a/c_5N$ for $j\in [L_0+1,r-1]$ or $j\in [r+1, L-L_0-5]$,
is less than $L/5$.
\end{lemma}
\proof
Consider $\Gamma_j$ as in the assumption of the lemma with $j\in [L_0+1,r-1]$.
The subcomb $E_j^0$ has at most $N$ maximal $q$-bands by Lemma \ref{noq}.
So there are at most $N$ maximal $a$-bands starting on ${\bf z}_j$ and ending
on each of the $\theta$-bands of $E_j^0$. Proving by contradiction, we
have at least $L|W|_a/5c_5N$ such $a$-bands for all $j\in S$, where $S$ the set of integers in $[L_0+1,r-1]\cup [r+1, L-L_0-5]$;
denote this set of $a$-bands by ${\bf A}$. But the number of maximal $\theta$-bands in all such
subcombs $E_j^0$ does not exceed $2h$. Therefore at least $L|W|_a/5c_5N-2hN$ bands from
$\bf A$ end on the subpaths ${\bf q}_{j,j+1}$ for $j\in S$. Therefore by Lemmas \ref{qq} and \ref{ochev}, we
have $$|{\bf p}_{L_0+1, L-L_0-5}|\ge |{\bf q}_{L_0+1, L-L_0-5}|\ge h_{L_0+1}+h_{L-L_0-5}
+LN/2+ \delta(L|W|_a/5c_5N-2hN)$$
\begin{equation}\label{pp}
\ge h_{L_0+1}+h_{L-L_0-5}
+LN/2+ \delta L|W|_a/10c_5N
\end{equation}
since  $2hN\le 2L_0^2N|W|_a<L_0^3|W|_a/10c_5N\le L|W|_a/10c_5N$
by the choice of $L_0$ and $L$ (\ref{LL0}).

Also by Lemma \ref{ppbar}, we have
$$|{\bf \bar p}_{L_0+1, L-L_0-5}|\le  h_{L_0+1}+h_{L-L_0-5}
+3L_0N+ 3L_0\delta|W|_a$$
\begin{equation}\label{barbar}
\le h_{L_0+1}+h_{L-L_0-5}
+3L_0N+ \delta L|W|_a/20c_5N,
\end{equation}
because by the choice of $L$, $3L_0< L/20c_5N$. Since $h_{L_0+1}+h_{L-L_0-5}\le 2h\le 2L_0^2|W|_a<L|W|_a$,
$L$ is chosen after $c_5N$, and $\varepsilon =N_4^{-1/2}$ is chosen after $L$, the
inequality $$\frac{|{\bf p}_{L_0+1, L-L_0-5}|}{|{\bf\bar p}_{L_0+1, L-L_0-5}|}\ge 1+\frac{\delta}{20c_5N}> 1+\varepsilon$$
follows from (\ref{pp}, \ref{barbar}), contrary to Lemma \ref{epsi}. The lemma is proved by contradiction.
\endproof

\begin{lemma} \label{hh} If $h\le L_0^2|W_a|$, then the histories $H_1$ and $H_{L-3}$ have different first letters.
\end{lemma}
\proof Let $\cal T$ and $\cal S$ be the maximal $\theta$-bands of $\Psi$
crossing ${\cal B}_1$ and ${\cal B}_{L-3}$, respectively, and the closest to the disk $\Pi$. Let they cross $k$ and $\ell$ spokes
of $\Pi$, respectively. By Lemma \ref{02}, $k+\ell>L-L/5-3L_0>2L/3$,
and also $k,\ell\ge 2$ since $L/2-3\le r \le L/2$. It follows from  Lemma \ref{withd} that the first letters of $H_1$ and
$H_{L-3}$ are different.
\endproof

\begin{lemma}\label{arG} We have $h> L_0^2|W|_a$.
\end{lemma}
\proof If this inequality is wrong, then by Lemma \ref{02}, there are at least $L-L/5-3L_0> 0.7L$ trapezia
$\Gamma_j$ with $|{\bf z}_j|_a<|W|_a/c_5N$, and one can choose two such trapesia
$\Gamma_k$ and $\Gamma_{\ell}$ such that  $k<r$, $\ell\ge r+1$ and $\ell-k>0.6L$.
Since $H_{k+1}$ (resp. $H_{\ell}$) is a prefix of $H_1$ (of $H_{L-3}$), it follows from
Lemma \ref{hh} that the first letters of $H_{k+1}$ and $H_{\ell}$ are different.

Since the bottoms of $\Gamma_k$ and $\Gamma_{\ell}$ (which belong to $\partial\Delta$) have the same label,
one can construct an auxiliary trapezium $E$ identifying the bottom of a copy of $\Gamma_k$
and the bottom of a mirror copy of $\Gamma_{\ell}$. The history of $E$ is $H_{\ell}^{-1}H_{k+1}$, which is a reduced word since the first letters of $H_k$ and $H_{\ell}$ are different, i.e. $E$ is a trapezium indeed
by Lemma \ref{simul}.

The top and the bottom of $E$  have $a$-lengths less than $|W|_a/c_5N$. Without loss of generality, one may assume that
$h_{k+1}\ge h_{\ell}$, and so $h_{k+1}\ge t/2$, where $t$ is the
height of $E$.

Note that the difference of $a$-lengths $|W|_a-|W|_a/c_5N>|W|_a/2$, and so
\begin{equation}\label{2N}
h_{k+1}, h_{\ell}> |W|_a/2N
\end{equation}
 since the difference of $a$-lengths for the top and the bottom of every maximal $\theta$-band of $E$ does not exceed $N$. Therefore $t>|W|_a/N$,
and the computation corresponding $E$ satisfies the assumption of Lemma \ref{B}.

So for every factorization $H'H''H'''$ of the
history of $\Gamma_k$, where $||H'||+||H''||\le \lambda ||H'H''H'''||$, we have $||H''||>0.4t$, since $\lambda< 1/5$.
Therefore by Lemma \ref{B}, the spoke ${\cal B}_{k+1}$
is a $\lambda$-shaft.

Using Lemma \ref{ppbar}, we obtain:
\begin{equation}\label{pb}
|{\bf  p}_{k+1, \ell}|+\sigma_{\lambda}(\bar \Delta_{k+1,\ell}) \ge h_{k+1}+h_{\ell}+0.6 LN+h_{k+1}.
\end{equation}

By inequality (\ref{2N}), we have  $\delta L|W|_a\le 2LN\delta h_{k+1}<h_{k+1}$ by the choice of $\delta$ and by Lemma \ref{ppbar},

\begin{equation}\label{bp}
|{\bf \bar p}_{k+1,\ell}|\le h_{k+1}+h_{\ell}+0.4 LN+0.4L\delta|W|_a\le h_{k+1}+h_{\ell}+h_{k+1}/2+0.4LN
\end{equation}

The right-hand side of the inequality (\ref{pb}) divided  by the right-hand side of (\ref{bp}) is
greater than $1.1$
(because $h_{k+1}\ge h_{\ell}$),
which contradicts Lemma \ref{epsi}. Thus, the lemma is proved.
\endproof

\begin{lemma} \label{hi} We have $h_i> \delta^{-1}$ for every $i=1,\dots, L_0$.
\end{lemma}
\proof By inequalities (\ref{L0}) and (\ref{downup}), we have $h_i\ge h_{L-L_0-3}$.
Proving by contradiction, we obtain $|W|_a<h_i
\le \delta^{-1}$
by Lemma \ref{arG}.
Then $$|{\bf \bar p}_{i, L-L_0-3}|< h_{i}+h_{L-L_0-3}+3L_0(N+\delta^{-1} \delta)\le
h_{i}+h_{L-L_0-3}+4L_0N$$ by Lemma \ref{ppbar}, and
$|{\bf p}_{i, L-L_0-3}|\ge h_{i}+h_{L-L_0-3}+LN/2$. Since $h_{i}+h_{L-L_0-3}\le 2\delta^{-1}$
and $4L_0N< LN/4$, we see that $\frac{|{\bf p}_{i, L-L_0-3}|}{|{\bf \bar p}_{i, L-L_0-3}|}> 1+\delta>1+\varepsilon$ contrary to Lemma \ref{epsi}. The lemma is proved by contradiction.
\endproof.

\begin{lemma} \label{shaft} None of the spokes ${\cal B}_1,...,{\cal B}_{L_0}$ contains
a $\lambda$-shaft at $\Pi$ of  length at least $\delta h$.
\end{lemma}

\proof On the one hand, by Lemmas \ref{ppbar} and \ref{arG},
\begin{equation}\label{pbp}
|{\bf \bar p}_{L_0+1, L-L_0-3}|< h_{L_0+1}+h_{L-L_0-3}+3L_0(N+\delta|W|_a)<h_{L_0+1}+h_{L-L_0-3}+3L_0(N+\delta L_0^{-2}h).
\end{equation}
On the other hand, by Lemma \ref{ppbar},

\begin{equation}\label{bpb}
|{\bf p}_{L_0+1,L-L_0-3}|> h_{L_0+1}+h_{L-L_0-3}+(L-3L_0)N.
\end{equation}

If the statement of the lemma were wrong, then we would have $\sigma_{\lambda}(\bar\Delta)\ge \delta h$, and  inequalities (\ref{pbp}) and (\ref{bpb}) would imply that
$$|{\bf p}_{L_0+1,L-L_0-3}|-|{\bf \bar p}_{L_0+1, L-L_0-3}|+\sigma_{\lambda}(\bar\Delta)\ge (L-6L_0)N - 3L_0^{-1}\delta h +\delta h\ge LN/2 +\delta h/2. $$
The right-hand side of the last inequality divided  by the right-hand side of (\ref{pbp}) is greater than $\varepsilon=N_4^{-\frac12}$,
because $h\ge h_{L_0+1}, h_{L-L_0-3}$,
which would contradict to Lemma \ref{epsi}. Thus, the lemma is proved.
\endproof

\begin{lemma}\label{zgh}
For every $j\in [1,L_0-1]$, we have
$|{\bf z}_{j}|_a> h_{j+1}/c_5$.
\end{lemma}

\proof If $|{\bf z}_{j}|_a\le h_{j+1}/c_5$, then the computation ${\cal C}: \;W_0\to\dots\to W_t$
corresponding to the trapezium $\Gamma_j$ satisfies the assumption of Lemma \ref{B},
since $t=h_{j+1}>c_5 |W_t|_a=c_5|{\bf z}_j|_a$ and by Lemma \ref{arG}, $t= h_{j+1} \ge L_0^2 |W_0|_a \ge c_5|W|_a$ since $L_0>c_5$.
Hence ${\cal B}_{j+1}$ is a $\lambda$-shaft by Lemma \ref{B} since $\lambda<1/2.$
We obtain a contradiction with Lemma \ref{shaft} since $\delta h\le h\le h_{j+1}$.
Thus, the lemma is proved.
\endproof

\begin{lemma} \label{2000N} For every $j\in [1,L_0-1]$, we have $h_{j+1}< (1-\frac{1}{10c_5N})h_j$.
\end{lemma}
\proof By Lemma \ref{zgh}, we have $|{\bf z}_{j}|_a\ge h_{j+1}/c_5$.
Let us assume that $h_{j+1}\ge (1-\frac{1}{10c_5N})h_j$, that is the handle ${\cal C}_j$ of $E_j$ has height at most $h_j/10c_5N$.
By Lemma \ref{le6}, at most $h_{j}/10c_5$ maximal $a$-bands of $E_j$ starting on ${\bf z}_j$ can end on
the $(\theta,q)$-cells of $E_j$. Hence at least $$|{\bf z}_j|_a-h_{j}/10c_5\ge |{\bf z}_j|_a-2h_{j+1}/10c_5\ge 0.8h_{j+1}/c_5> 0.7h_{j}/c_5$$ of them have to end on the path ${\bf p}_{j,j+1}$.

The path ${\bf p}_{j,j+1}$ has at most $\frac{h_j}{10c_5N}\;\;$ $\theta$-edges. Hence by Lemma \ref{ochev},
$$|{\bf p}_{j,j+1}|\ge h_j-h_{j+1}+\delta(0.7h_{j}/c_5-h_j/10c_5N)\ge h_j-h_{j+1}+0.6\delta h_{j}/c_5,$$ and therefore by Lemma \ref{ppbar},
$|{\bf p}_{j, L-L_0-3}|\ge LN/2+h_j + h_{L-L_0-3}+0.6\delta h_{j}/c_5$. On the other hand by Lemma \ref{ppbar}, we have
$$|{\bf\bar p}_{j, L-L_0-3}|\le h_j + h_{L-L_0-3}+ 3NL_0+3L_0\delta|W|_a\le  h_j + h_{L-L_0-3}+ 3NL_0+3L_0^{-1}\delta h_{j+1}$$
by Lemma \ref{arG} and inequality $h\le h_{j+1}$. Hence
$\frac{|{\bf p}_{j, L-L_0-3}|}{|{\bf\bar p}_{j, L-L_0-3}|}\ge (1+\delta/10c_5)$
since $h_{L-L_0-3}\le h_{L_0+1}\le h_{j+1}\le h_j$ and $L_0>>c_5$. We have a contradiction with Lemma \ref{epsi} since
$\delta/10c_5>\varepsilon$.
The lemma is proved by contradiction.
\endproof

The proof of the next lemma is similar.

\begin{lemma} \label{zlh} For every $j\in [2,L_0-1]$,
we have $|{\bf z}_j|_a \le 2Nh_j$,
\end{lemma}
\proof Assume that $|{\bf z}_j|_a \ge 2Nh_j$.
By Lemma \ref{le6}, at most $Nh_j$ maximal $a$-bands of $E_j$ starting on ${\bf z}_j$ can end on
the $(\theta,q)$-cells of $E_j$. Hence at least $|{\bf z}_j|_a-Nh_j\ge Nh_j$ of them has to end on the path ${\bf p}_{j,j+1}$.
The path ${\bf p}_{j,j+1}$ has at most $h_j$ $\theta$-edges. Hence by Lemma \ref{ochev},
$|{\bf p}_{j,j+1}|\ge h_j-h_{j+1}+\delta(Nh_j-h_j)= h_j-h_{j+1}+\delta(N-1)h_j$ and therefore by Lemma \ref{ppbar},
$|{\bf p}_{j, L-L_0-3}|\ge LN/2+h_j + h_{L-L_0-3}+\delta(N-1)h_j$.
On the other hand by Lemmas \ref{ppbar} and \ref{arG}, we have
$$|{\bf\bar p}_{j, L-L_0-3}|\le h_j + h_{L-L_0-3}+ 3NL_0+3L_0\delta|W|_a\le  h_j + h_{L-L_0-3}+ 3NL_0+\frac{3\delta h_j}{L_0},$$ because $h\le h_j$. Since $h_j \ge h\ge  h_{L-L_0-3} $, we have
$\frac{|{\bf p}_{j, L-L_0-3}|}{|{\bf\bar p}_{j, L-L_0-3}|}\ge 1+\varepsilon$,
a contradiction by Lemma \ref{epsi}.
\endproof

\begin{lemma} \label{001} There is no $i\in [2,L_0-3]$ such that the histories
$H_{i-1}=H_iH'=H_{i+1}H''H'=H_{i+2}H'''H''H'$ and the computation $\cal C$ with history $H_i$ corresponding to the trapezium $\Gamma_{i-1}$ satisfy the following condition:

(*) The history $H'''H''H'$ has only one step, and for the subcomputation $\cal D$
with this history,
there is a sectors $Q'Q$ such that
a state letter  from $Q$ or from $Q'$ inserts a letter increasing
the length of this sector after every transition of $\cal D$.

\end{lemma}
\proof

Recall that the standard base of $\bf M$ is built of the standard base $B$ of ${\bf M}_4$
and its inverse copy $(B')^{-1}$ (plus letter $t$).
Due to this mirror symmetry of the standard base, we have mirror symmetry for any accepting computation, in particular, for $\cal C$ and $\cal D$. Therefore proving by contradiction,
we may assume that the $a$-letters are inserted from the left of $Q$.

Let  $\cal Q$ be the maximal $q$-spoke of the subdiagram
$E_i^0\subset \Gamma_i$ corresponding to the base letter $Q$.
If ${\cal Q'}$ is the neighbor from the left $q$-spoke for $\cal Q$
(the spokes are directed from the disk $\Pi$),
then the subpath $\bf x$ of ${\bf z}_i$ between these two $q$-spokes has at least
$h_{i+1}-h_{i+2}=||H'''||$ $a$-letters. Indeed, $\Gamma_i$ contains
a copy $\Gamma'_{i+1}$ of $\Gamma_{i+1}$, the bottom of the trapezium
$\Gamma_i\backslash\Gamma'_{i+1}$ is the copy ${\bf z}'_{i+1}$ of ${\bf z}_{i+1}$
and the top of it iz ${\bf z}_i$,
and so the subcomputation with history $H'''$ has already increased
the length of the $Q'Q$-sector. Thus,  by Lemmas \ref{2000N}, \ref{arG} and the choice of $L_0>100c_5N$, we have
\begin{equation}\label{xa}
|{\bf x}|_a\ge  h_{i+1}-h_{i+2}\ge \frac{1}{10c_5N}h_{i+1}\ge 10L_0 |W|_a.
\end{equation}

Note that an $a$-band $\cal A$ starting on ${\bf x} $ cannot end on a $(\theta,q)$-cell from $\cal Q$.
Indeed, otherwise by Lemma \ref{copy}, there is a copy of this configuration in
the diagram $\Gamma_{i-1}$, i.e. the copy of $\cal A$ ends on the copy of $\cal Q$ contrary
the assumption that the rules of computation with history $H'''H''H'$ do not delete $a$-letters.

Let us consider the comb bounded by $\cal Q$, ${\cal Q}'$, $\bf x$
and the boundary path of $\Delta^0$ (without the cells from ${\cal Q}'$). If the lengths of $\cal Q$ and ${\cal Q}'$ are
$s$ and $s'$, respectively, then there are $|{\bf x}|+s$ maximal
$a$-bands starting on $\bf x$ and $\cal Q$ and ending either on
${\cal Q}'$ or on $\partial\Delta^0$ since the comb has no
maximal $q$-bands by Lemma \ref{noq}. At most $s'<s$ of these
$a$-bands can end on ${\cal Q}'$. Therefore at least $|\bf x|+s-s'$
of them end on the segment of the boundary path of $\Delta^0$ lying between the ends of ${\cal Q}'$ and $\cal Q$.

Since
this segment has $s-s'$ $\theta$-edges, its length is at least $s-s'+\delta|{\bf x}|_a$ by Lemma \ref{ochev}.
This inequality, Lemma \ref{qq} and inequality (\ref{xa})  imply
$$|{\bf p}_{i,L-L_0-3}|\ge |{\bf q}_{i,L-L_0-3}|\ge LN/2+h_{i}+h_{L-L_0-3}+\frac{\delta}{10c_5N}h_{i+1}$$ $$ \ge LN/2+h_{i}+h_{L-L_0-3}+10\delta L_0|W|_a.$$
Therefore
\begin{equation}\label{pgp}
|{\bf p}_{i,L-L_0-3}|-\frac{7\delta}{100c_5N}h_{i+1}> 3L_0N+h_{i}+h_{L-L_0-3}+3\delta L_0|W|_a\ge |{\bf \bar p}_{i,L-L_0-3}|,
\end{equation}
by Lemma \ref{ppbar}, and since $\Delta$ is a minimal counter-example, we obtain by the definition of
$F(x)$, $g(x)$ and inequality (\ref{xy}) that
$$\area_G(\Psi'_{i+1,L-L_0-3})\le N_4F(n+\sigma_{\lambda}(\Delta^*)-\frac{7\delta}{100c_5N}h_{i+1})+N_3g(n)\mu(\Psi'_{i,L-L_0-3})$$
\begin{equation}\label{areai}
\le N_4F(n+\sigma_{\lambda}(\Delta^*))-N_4\frac{7\delta n}{100c_5N}h_{i+1}g(n)+N_3g(n)\mu(\Psi'_{i,L-L_0-3}).
\end{equation}

By Lemma \ref{arG}, $|W|_a \le L_0^{-2}h_i$, and by Lemma \ref{hi}, $h_i>\delta^{-1}>100L_0N$,  whence $|{\bf \bar p}_{i,L-L_0-3}|\le 2h_i+3L_0N+3\delta L_0|W|_a\le(2+0.03+\frac{3\delta}{L_0})h_i \le 2.1h_i$ and by Lemma \ref{epsi},
we have
\begin{equation}\label{22}
|{\bf  p}_{i,L-L_0-3}|\le (1+\varepsilon)|{\bf \bar p}_{i,L-L_0-3}|<2.2h_i.
\end{equation}

By Lemmas \ref{main} and \ref{mixture} (a) and inequalities (\ref{22}) and (\ref{pgp}), the $G$-area of $\Psi_{i,L-L_0-3}$ does not exceed
\begin{equation}\label{areaP}
N_2(2|{\bf p}_{i,L-L_0-3}|)^2+N_1\mu(\Psi_{i, L-L_0-3})\le N_2(4J+4)|{\bf p}_{i,L-L_0-3}|^2\le 5N_2(4J+4)h_i^2.
\end{equation}

By Lemma \ref{disk}, the $G$-area of $\Pi$ is bounded by $c_6F(|\partial\Pi|)$.
The inequalities (\ref{pgp}) and (\ref{22}) imply the inequality
$|\partial\Pi|<L|{\bf\bar p}_{i,L-L_0-3}|<L|{\bf p}_{i,L-L_0-3}|<3Lh_i$.
Therefore one may assume that the constant $c_7$ is chosen so that
\begin{equation}\label{c7}
\area_G(\Pi)< c_6F(|\partial\Pi|)< c_7F(h_i)=c_7h_i^2g(h_i)\le c_7h_i^2g(n).
\end{equation}
(Recall that $h_i\le n/2$ here since $h_i$ is the number $\theta$-bands crossing ${\cal B}_i$;
they start and end on $\partial\Delta$.)
It follows from (\ref{areaP}) and (\ref{c7}) that
\begin{equation}\label{Dij}
\area_G(\bar\Delta_{i,L-L_0-3})\le 5N_2(4J+4)h_i^2+ c_7h_i^2g(n).
\end{equation}

\begin{figure}
\begin{center}
\includegraphics[width=0.8\textwidth]{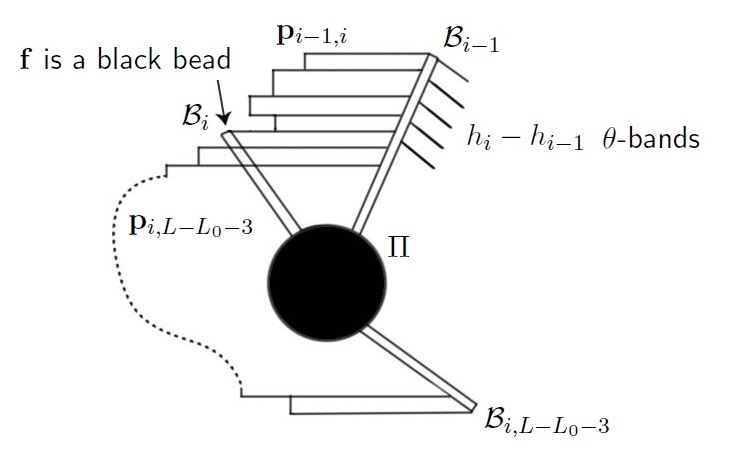}
\end{center}
\caption{$\mu(\Psi'_{i+1,L-L_0-3})-\mu(\Psi'_{i,L-L_0-3})$}\label{Pic17}
\end{figure}

We need an estimate for $\mu(\Psi'_{i+1,L-L_0-3})-\mu(\Psi'_{i,L-L_0-3})$ now. To obtain it,
we observe that by Lemma \ref{psi1}, the common $q$-edge $\bf f$ of the spoke
${\cal B}_{i}$ and $\partial\Delta$
 separates at least
$h_{i-1}-h_{i}=m_1\;\;\;\;\theta$-edges of the path ${\bf p}_{i-1,i}$
and $m_2$ ones lying on ${\bf p}_{i,L-L_0-3}$,
where $m_2=h_{i}+h_{i,L-L_0-3}$ (see fig. \ref{Pic17}).
Since the number of $q$-edges of ${\bf p}$ is less than $3K_0L<J$ by Lemma \ref{2K0},  one
decreases $\mu(\Psi'_{i+1,L-L_0-3})$ at least by $m_1m_2$ when erasing the black
bead on $\bf f$ in the necklace on $\partial\Psi'_{i+1,L-L_0-3}$ by Lemma \ref{mixture} (d,b,c). Hence
$$\mu(\Psi'_{i+1,L-L_0-3}))-\mu(\Psi'_{i,L-L_0-3})\ge m_1m_2$$ $$=(h_{i-1}-h_{i})(h_{i}+h_{L-L_0-3})\ge\frac{1}{c_5N}h_{i-1}(h_{i}+h_{L-L_0-3})$$
by Lemma \ref{2000N}.
This inequality and Lemma \ref{muJ} applied to $\Psi'_{i+1, L-L_0-3}$, imply
$$\mu(\Delta)-\mu(\Psi'_{i,L-L_0-3})\ge -2Jn(h_{i+1}+h_{L-L_0-3})+\frac{1}{10c_5N}h_{i-1}(h_{i}+h_{L-L_0-3}).$$
Note that $(h_{i+1}+h_{L-L_0-3})\le 2h_{i+1}$ by (\ref{downup}) and (\ref{L0}). Hence
\begin{equation}\label{mumu}
N_3\mu(\Delta)-N_3\mu(\Psi'_{i,L-L_0-3})\ge -4N_3Jn h_{i+1}+\frac{N_3}{10c_5N}h_{i-1}(h_{i}+h_{L-L_0-3}).
\end{equation}
It follows from (\ref{Dij}, \ref{mumu}, \ref{areai}) that $$\area_G(\Delta)\le \area_G(\Psi'_{i, L-L_0-3})+\area_G(\bar\Delta_{i, L-L_0-3})\le $$
$$\le N_4F(n+\sigma_{\lambda}(\Delta^*))-N_4\frac{7\delta n}{100c_5N}h_{i+1}g(n)+N_3g(n)\mu(\Psi'_{i,L-L_0-3})+ 5N_2(4J+4)h_i^2+ c_7h_i^2g(n)$$

$$\le N_4F(n+\sigma_{\lambda}(\Delta^*)+N_3\mu(\Delta)g(n)-N_4\frac{7\delta n}{100c_5N}h_{i+1}g(n)-
\frac{N_3}{10c_5N}h_{i-1}(h_{i}+h_{L-L_0-3})g(n)$$ $$+4N_3Jnh_{i+1}g(n)+ 5N_2(4J+4)h_i^2+c_7h_i^2g(n).$$

We come to a contradiction since we obtain inequality $\area_G(\Delta)\le N_4F(n+\sigma_{\lambda}(\Delta^*))+N_3\mu(\Delta)g(n)$, because $N_4\frac{7\delta}{100c_5N}>4N_3J$
and $\frac{N_3}{10c_5N}>5N_2(4J+4)+c_7$.
\endproof

\begin{lemma} \label{led} There exists no counter-example $\Delta$, and therefore $\area_G(\Delta)\le N_4F(n+\sigma_\lambda(\Delta^*))+N_3\mu(\Delta)g(n)$ for any minimal diagram $\Delta$ with $|\partial\Delta|=n$.
\end{lemma}
\proof Recall that
for $j=1,\dots,L_0-1$, we have
$h_{j+1}< (1-\frac{1}{10c_5N})h_j$ by Lemma \ref{2000N},
and by Lemmas \ref{zgh} and \ref{zlh}, we have inequalities
$|{\bf z}_j|_a\ge h_{j+1}/c_5$ and $|{\bf z}_k|_a\le 2Nh_k$ for $2\le k\le L_0-1$. One can choose an integer
$\rho=\rho(\bf M)$ (it depends on the S-machine $\bf M$ only as $c_5$ and $N$) so that $(1-\frac{1}{10c_5N})^{\rho}<\frac{1}{6Nc_5}$, and so $h_{j+1}>6Nc_5h_k$
if $k-j-1\ge \rho $. Hence
$$|{\bf z}_j|_a\ge h_{j+1}/c_5\ge 6Nh_k >3|{\bf z}_k|_a.$$

If $L_0$ is large enough, say $L_0>2000\rho$, one can obtain $1000$ indices
$j_1<j_2<\dots<j_{1000}<L_0$ such that for $i=2,\dots, 1000$, one obtains inequalities

\begin{equation}\label{zh}
|{\bf z}_{j_{i-1}}|_a>3|{\bf z}_{j_i}|_a\;\;and\;\; h_{j_{i-1}}\ge h_{j_{i-1}+1}> 6c_5N h_{j_i}.
\end{equation}

Let ${\cal C}:\;W\equiv W_0\to\dots\to W_t$  be the computation corresponding to the trapezium $\Gamma_{j_2}$. Since it contains the copy $\Gamma_{{j_2}+1}'$ of $\Gamma_{{j_2}+1}$, which in turn
contains a copy of $\Gamma_{{j_2}+2}$ and so on, we have
some configurations $W(k)$ in $\cal C$ ($k=1,\dots, 999$), that are the labels of some ${\bf z}_{i_k}$ and $|W(k+1)|_a>3|W(k)|_a$ for $k=1,\dots,998$.
If for some $k$ we were obtain one-step
subcomputation $W(k)\to\dots\to W(k+4)$, then the statement of Lemma \ref{pol} would give
a subcomputation $W(k+1)\to\dots\to W(k+4)$
contradicting  to the statement of Lemma \ref{001}. Hence no five consecutive words $W(k)$-s are configuration
of a one-step subcomputation, and so the number of steps in $W(1)\to\dots \to W(999)$ in at least $100$.

It follows now from Lemmas \ref{histF} and \ref{474} that the step
history of $\Gamma_{j_2}\backslash \Gamma$, where $\Gamma$ is the copy of $\Gamma_{L_0}$
in $\Gamma_{j_2}$,  has a subword $\big((21^-)(1^-)(1)(2^-)(2)(21^-)\big)^{\pm 1}$. Without loss of generality we assume that the exponent is $+1$. Therefore the
history $H_{j_2+1}$ of $\Gamma_{j_2}$ can be decomposed as $H'H''H'''$, where $H''$ has step history $(12^-)(2^-)(2^-2)$,
$||H'||\ge ||H''||$ by Lemma \ref{cycle} and $||H'||\ge h$ since
$H_{L_0}$ is a prefix of $H'$.

Since $h_{j_1+1} > 2h_{j_2}$ by (\ref{zh}), the history $H_{j_1+1}$
of $\Gamma_{j_1}$ has a prefix $H'H''H^*$, where $||H^*||=||H'||\ge ||H''||$, and so the $t$-spoke ${\cal B}_{j_1+1}$ has a $t$-subband $\cal C$ starting with $\partial\Pi$
and having the history $H'H''H^*$.

For any factorization ${\cal C}={\cal C}_1{\cal C}_2{\cal C}_3$
with $||{\cal C}_1||+||{\cal C}_2||\le||{\cal C}||/3$, the history
of ${\cal C}_2$ contains the subhistory $H''$, since $||H^*||=||H'|||\ge ||H''||$. It follows that ${\cal C}$ is
a $\lambda$-shaft, because $\lambda<1/3$. The shaft has length
at least $||H'||\ge h$ contrary to Lemma \ref{shaft}. We come to the final contradiction in this section.

\endproof

\section{Proof of Theorem \ref{beta}}\label{end}

\subsection{Dehn function of the group $G$}\label{end1}

\begin{lemma} \label{big} For every big trapezia $\Delta$, there is a diagram
$\tilde\Delta$ over $G$ with the same boundary label, such that
the area of $\tilde\Delta$ does not exceed $2\area_G(\Delta)$.
\end{lemma}

\proof Consider the computation ${\cal C}:\; V_0\to\dots\to V_t$ corresponding to $\Delta$. According to Definition \ref{abt},
one may assume that $\area_G(\Delta)=c_5h(||V_0||+||V_t||)$ since
otherwise $\tilde\Delta=\Delta$.

$\Delta$ is covered by $L$ trapezia $\Delta_1,\dots,\Delta_L$
with base $xvx$, where $xv$ (or the inverse word) is a cyclic shift ot the standard base
of $\bf M$. By Remark \ref{restore} all $\Delta_1,\dots,\Delta_L$ are copies of each other.
Let us apply Lemma \ref{narrow} to any of them, say to $\Delta_1$,
whose top and bottom have labels $W_0$ and $W_t$. If we have Property (1) of that lemma, then the area of $\Delta_1$ does not exceed
$c_4h(||W_0||+||W_t||)$ since every maximal $\theta$-band of $\Delta_1$ has at most $c_4(||W_0||+||W_t||)$ cells in this case.
Hence area of $\Delta$ does not exceed $Lc_4h(||W_0||+||W_t||)\le 2c_4h(||V_0||+||V_t||)<c_5h(||V_0||+||V_t||)=\area_G(\Delta)$, i.e.
$\tilde\Delta=\Delta$ in this case too.

Hence one may assume that Property (2) of Lemma \ref{narrow} holds
for $\Delta_1$. By that Lemma, items (b,d), the corresponding cyclic shifts $W'_0$ and $W'_t$
are accepted, and there is a factorization of ${\cal C}={\cal C}_1{\cal C}_2{\cal C}_3$ (we use the same letter for the computations corresponding to $\Delta$, to $\Delta_1$, and for the revolving computation with standard base),
where ${\cal C}_1: \; W'_0\to\dots\to W'_{n_1}$,
 ${\cal C}_2: \; W'_{n_1}\to\dots\to W'_{n_1+n_2}$ and ${\cal C}_3: \; W'_{n_1+n_2}\to\dots\to W'_{n_1+n_2+n_3}$ ($n_1+n_2+n_3=t$), where $\max(||W'_{n_1}||,||W'_{n_1+n_2}||)\le \max(||W_0||,||W_t||)$ and for each ${\cal C}_i$, either

 (d1) $||W'_j||\le c_4\max(||W'_0||, ||W'_t||)$, for every configuration $W_j$ of ${\cal C}_i$ or

(d2) there are accepting computations for the first and the last configuration of ${\cal C}_i$ with histories $H'_i$ and $H''_i$ such that $||H'_i||+||H''_i||<n_i$

So $\Delta$ is built of at most three trapezia, where  $\Delta(i)$ ($i=1,2,3$) corresponds to the computations ${\cal C}_i$.
Since their tops and bottoms have lengths at most $\max(||V||_0, ||V_t||)$, it suffices to estimate the area of $\tilde \Delta(i)$ for $i=1,2,3$.
Again, we have $\tilde \Delta(i)$= $\Delta(i)$ in the case (d1).

Assume that we have Property (d2) for ${\cal C}_i$. Denote by
$U_1$ and $U_2$ the first and the last configurations of ${\cal C}_i$ with standard base. By (d2), $||U_1||, ||U_2||\le max(||W_0||, ||W_t||)$.

By Property (d2), there is an accepting
computations $\cal D$ of length $\le n_i$ for $U_1$, and we may assume that $\cal D$
is the shortest such computation. Then case (b) of Lemma \ref{sta}
gives a contradiction of the form $\ell<\ell/100$ for the length $\ell$ of the computation $\cal D$. Hence we should have case (a),
and so every configuration of $\cal D$ has length at most $c_4||U_1||$.
If $\Gamma$
is the trapezium corresponding  to $\cal D$ with bottom (top) label $U_1$ (resp., $W_M$), then the lengths
of its $\theta$-bands are less than $2c_4||U_1||$ by Lemma \ref{ochev}
and therefore the area of $\Gamma$ is less than $2n_ic_4||U_1||-1$.
Therefore $L$ copies
of $\Gamma$ can be attached to an auxiliary hub so that one gets an auxiliary disk $\Pi_1$ of area $\le 2Ln_ic_4||U_1||\le 3 c_4 n_i ||V(1)||$, where
$V(1)t$ is the label of the bottom of $\Delta(i)$ up to cyclic permutations. Thus, the word $V(1)$
is equal to the boundary label of $\Pi_1$. Similarly, one can construct a disk of area $\le 3c_4 n_i ||V(2)||$ for the top of
$\Delta(2)$.

Denote by $\Delta_-$ the diagram $\Delta$ without maximal rim
$x$-band. So $\Delta_-$ has the boundary ${\bf p}_1{\bf q}_1
{\bf p}_2^{-1}{\bf q}_2^{-1}$, where $\Lab({\bf p}_1)$ and
$\Lab({\bf p}_2)$ are the boundary labels of the disks
$\Pi_1$ and $\Pi_2$ (up to cyclic permutations) and $\Lab(\bf q)_1)\equiv\Lab(\bf q)_2)$ since two $x$-bands
with the same history have the same boundary labels.

If we attach disks $\Pi_1$ and $\Pi_2$ along their boundaries
to the top and the bottom of $\Delta_-$, we obtain a diagram,
whose boundary label is trivial in the free group. Hence there
is a diagram $E$ with two disks whose boundary label is equal
to the boundary label of $\Delta_-$ and the area is less than
$\le 3c_4 n_i (||V(1)||+||V(2)||$. If we attach one $x$-band
of length $n_i$ to $E$, we construct the required diagram
$\tilde\Delta(i)$ of area at most $\le 3c_4 n_i (||V(1)||+||V(2)||)+n_i< c_5n_i (||V(1)||+||V(2)||)$
\endproof

The proof of Lemma \ref{big} shows that the area of a minimal diagram with
some boundary label $V$ can be much greater than the area of $V$, which is  equal to the minimal number of
cells in all diagrams with boundary label $V$ over the presentation (\ref{rel1}-\ref{rel3}).
So to obtain the lower bound for the Dehn function of $G$, we prove in the next
lemma that these two areas `almost equal' for the words having no $\theta$-letters.

\begin{lemma}\label{notheta}
Let $(tW)^L=1$ in $G$, where
the reduced word $W$ has no $\theta$-letters and no letters $t^{\pm 1}$. Then
there exists a reduced diagram $\Delta$ over the presentation (\ref{rel1}-\ref{rel3}) such that it has exactly one hub, has  boundary label $V\equiv (tW)^L$ and
$\area(\Delta)\le 2 \area(V)$.
\end{lemma}

\proof Let $\Delta_0$ be a diagram over the presentation (\ref{rel1}-\ref{rel3}) of $G$ with  boundary label $(tW)^L$, where $\area(\Delta_0)=\area(tW)^L$.  We say that $\Gamma$ is a disk subdiagram of $\Delta_0$
if it has reduced boundary, has exactly one hub and every $\theta$-cell of $\Gamma$ (if any) belongs in a
$\theta$-annulus surrounding this hub. The diagram $\Delta_0$ can be covered by a family of
subdiagrams $\Gamma$, where each $\Gamma$ is either a disk subdiagram or a $\theta$-cell
and different subdiagram of this covering ${\bf S}_0$ have no cells in common. Let $A({\bf S}_0)$ be the sum of the areas
of all disk subdiagram from ${\bf S}_0$ plus doubled number of the single $\theta$-cells from ${\bf S}_0$.
By $A(\Delta_0)$ we denote the minimum of the numbers $A({\bf S}_0)$ over all such coverings ${\bf S}_0$.
Clearly $A(\Delta_0)\le 2\area(\Delta_0)=2\area(V)$, and so it suffices to prove that
a reduced diagram $\Delta$ with boundary label  $(tW)^L$ and minimal possible value of $A(\Delta)$ has exactly one disk, because $\area(\Delta)\le A(\Delta)\le A(\Delta_0)$.

Below we fix the covering $\bf S$ of $\Delta$ such that $A(\Delta)=A(\bf S)$. Note that every $\theta$-annulus
of $\Delta$ surrounds at least one disk by Lemma \ref{NoAnnul} since the diagram $\Delta$ is reduced.

By induction on the number of $\theta$-annuli in a disk subdiagram $\Gamma\in \bf S$, we see
that the boundary label of $\Gamma$ has the form  $(tU)^{L}$. Therefore there is only
one cyclic shift of  the word $(tU)^L$ starting with $tU$.
Note that there are no
two distinct disk subdiagrams
$\Gamma_1$ and $\Gamma_2$ in $\bf S$ whose boundaries share at least two $t$-edges, provided there are no other
disk subdiagrams between $\Gamma_1$ and $\Gamma_2$, because such pair of subdiagrams could be
canceled  out, which would decrease the value of $A(\Delta)$.

Assume that $\bf S$ has at least one single $\theta$-cell. Let $\cal T$ be the maximal $\theta$-band
of $\Delta$ containing this cell. It has to be a $\theta$-annulus, where every cell is
a member of $\bf S$ since $\cal T$ can end neither on $\partial\Delta$ nor on
the boundary of a disk subdiagram from $\bf S$. So one can choose a minimal $\theta$-annulus
$\cal T$ whose cells do not belong to the disk subdiagrams from the family $\bf S$,
and  $\cal T$ surrounds a subdiagram $E$ having no single $\theta$-cells from $\bf S$.

The reduced diagram $E$ must contain disk subdiagrams by Lemma \ref{NoAnnul}. Hence as
in Lemma \ref{extdisc}, we have a disk graph, where there are no two different edges connecting
neighbor disk subdiagrams in $E$ (and crossing the $t$-edges on the boundaries of these
subdiagrams) provided there are no other disk subdiagrams between these two edges of
the disk graph. Hence there is a disk subdiagram $\Gamma$ in $E$ sharing a
boundary subpath $\bf q$ with a side of $\cal T$, where $\Lab({\bf q}) = (tU)^{L-4}t$.
After the transposition of $\cal T$ and $\Gamma$ we can obtain a new disk subdiagram
$\Gamma'$ with $\area(\Gamma')\le\area(\Gamma)+L(u+1)$, where $u$ is the number of cells
between two neighbor $(\theta,t)$-cells in $\cal T$. However, the transposition
removes $(L-4)u+L-3$ cells from $\cal T$ and add at most $4u+3$ new cells.  Since\\
$(L-4)u+L-3 - 4u-3 > \frac{L}{2} (u+1)$, we have a new diagram and new covering $\bf S'$
after the transposition, where $A({\bf S}')<A({\bf S})$, because the single $\theta$-cell
is taken with coefficient $2$ in the above definition of $A(\cdot)$; a contradiction.

Thus, the covering $\bf S$ has no single $(\theta,q)$-cells. Then the standard argument implies that
$\bf S$ has at most one disk subdiagram (see Lemma \ref{mnogospits}). The diagram $\Delta$ cannot be
a diagram over $M$ since all $t$-letters occur in the boundary label $(tV)^L$ with exponent $+1$,
and so the $t$-edges of $\partial\Delta$ cannot be connected in $\Delta$ by a $t$-band. Thus, the number
of hubs of $\Delta$ is 1.
\endproof

\begin{lemma} The Dehn function $d(n)$ of the group $G$ is equivalent to $F(n)$.
\end{lemma}
\proof To obtain the upper bound for $d(n)$ (with respect to
the finite presentation of $G$ given in Section \ref{gd}),
 it suffices, for every word $W$ vanishing in $G$ with $||W||\le n$, to find a diagram over $G$ of area $O(F(n))$ with boundary label $W$. Since $|W|\le ||W||$, van Kampen's lemma and Lemma \ref{arG} provide us with a minimal diagram $\Delta$ such that
 $\area_G(\Delta)\le N_4F(n+\sigma_\lambda(\Delta^*))+N_3\mu(\Delta)g(n)$ for some costants $N_3$ and $N_4$ depending on the presentation of $G$.
 By Lemmas \ref{clam}, \ref{mixture} (a) and the definition of $\mu(\Delta)$, the right-hand side does not exceed
 $N_4F((1+c)n)+N_3J n^2g(n)$. Since $F(O(n))=O(F(n))$ and
 $n^2g(n)=F(n)$, we conclude that $\area_G(\Delta)\le C_0F(n)$
 for some constant $C_0$.

 Recall that in the definition of $G$-area, the subdiagrams, which are big trapezia $\Gamma, \Gamma', \dots,$ can have common cells in their rim $q$-bands only.  By Lemma \ref{big}, any big trapezia $\Gamma$ from this list with top ${\bf p}_1$ and bottom ${\bf p}_2$ can be replaced
 by a diagram $\tilde \Gamma$ with (combinatorial) area at most $2\area_G(\Gamma)$.
 When replacing all big trapezia $\Gamma, \Gamma', \dots$ in this way, we should add $q$-bands for the possible intersection of
 big trapezia, but for every $\Gamma$ of height $h$, we add at most $2h$ new cells. So the area of the modified diagram $E$ is at most
 $3\area_G(\Delta)\le 3C_0F(n)$. Hence the required diagram is found for given word $W$.

 To obtain the lower bound for $d(n)$, we will use the series of
 $\bf M$-accepted words $V(n)$ of (combinatorial) length $\Theta(n)$ constructed in the
 proof of Lemma \ref{gtime}. Since $V(n)^L=1$ in $G$, it will be suffice to bound from below
the areas of the diagrams $\Delta(n)$ given by Lemma \ref{notheta}: $\Delta(n)$ has boundary label
$V(n)^L$, exactly one disk and the area equal to $\area(V(n)^L)$ up to a multiplicative factor
from the segment $[1,2]$.

 A $q$-band starting on the hub $\pi$ of $\Delta(n)$ cannot end on it
 since all occurrences of a particular $q$-letter in the hub relation have the same exponent.
 Hence the spokes of $\pi$ end on $LN$ $q$-edges of $\partial\Delta$.
 Hence $\Delta$ has $L$ trapezia corresponding to an accepting
 reduced computation $\cal C$ for $V(n)$, and it suffices to
 get a lower bound for the area of one trapezium $\Gamma$.

 By Remark \ref{ar}, $\Gamma$ has at least $\Theta(ng(n))$
maximal $\theta$-bands of length at least $\Theta(n)$. Therefore
the area of $\Gamma$ is at least $\Theta(F(n))$. Since $||\partial(\Delta)||=\Theta(n)$, the  Dehn function $d(n)$ is bounded from below by a function equivalent to $F(n)$, as required.
\endproof

\subsection{Supercubic Dehn functions}\label{sc}

Here we show that for the Dehn functions $F(n)$ obtained
earlier, one can construct a finitely presented group with
Dehn function $nF(n)$. This will complete the proof of Theorem
\ref{beta}. For this goal we modify the control S-machines
used in the definition of the S-machine $\bf M$. The unnecessary extra-control will just slow down the work of $\bf M$. The construction
resembles the one from Subsection 4.2 of \cite{O12}.
We will modify only the  S-machine
 $\bf P$ defined in Subsection \ref{pm}.
The copies of the auxiliary primitive S-machine $Z(A)$  will  work
between the applications of the (copies of) the rules of $\bf P$.

\medskip

For every set of letters
$A,$ let $A'$, $A''$ and $A'''$ be  disjoint  copies of $A$, the maps $a\mapsto a'$, $a\mapsto a''$  and $a\mapsto a'''$ identify $A$ with $A'$, $A''$ and $A''',$ resp. Let \label{overZ} $\overleftarrow Z$ be the $S$-machines
with tape alphabet $A'\sqcup A''\sqcup A'''$, state alphabet $\{L\}\sqcup K\sqcup P\sqcup \{R\},$ with
the following positive rules.

$$\chi_1(a)=[L\tool L, k(1)\to k(1), p(1)\to (a'')^{-1}p(1)a''', R\to R], a\in A;$$

$$\chi_2:=[L\tool L, k(1)\tool k(2) ,  p(1)\to p(2), R\to R];$$

$$\chi_3(a)=[L\tool L, k(2)\to k(2), p(2)\to a'' p(2)(a''')^{-1}, R\to R];$$

$$\chi_4=[L\tool L, k(2)\to k(3), p(2)\tool p(3), R\to R].$$

The rules of $\overrightarrow Z$ are similar, but the moving
base letter is $K$, while the sector $PR$ is locked.

To define \label{M3} the composition ${\bf P}\circ\{\overrightarrow Z,\overleftarrow Z \} $   we insert the base of  $\overleftarrow Z$ (and $\overrightarrow Z$) between every two consecutive state letters of $\bf P$.
In this subection, we assume that $\bf P$ has the standard base $Q_0Q_1...Q_N$
(and forget more detailed earlier notation).

For every $i=1,...,N$, we make  copies $Y'_i$, $Y''_i$ and  $Y'''_i$ of the alphabet
$Y_i$ of $\bf P$ ($i=1,...,N$).
Let $\Theta$ be the set of positive commands of $\bf P.$
The set of state letters of ${\bf P}\circ\{\overrightarrow Z,\overleftarrow Z \} $   is $$S_0\sqcup K_1\sqcup P_1\sqcup
S_1\sqcup K_2\sqcup P_2\sqcup...\sqcup P_{N}\cup S_N,$$
where $P_i=\{p^{(i)}, p^{(i,1)}, p^{(i,0)},
p^{(\theta,i)}(1), p^{(\theta,i)}(2),  p^{(\theta,i)}(3)\mid
\theta\in\Theta\}$, $K_i$ is defined similarly for
$i=1,...,N$, $S_i=Q_i\sqcup (Q_i\times \Theta)$.
Thus the state letters $L$ and $R$ of the copies of the S-machines $\overrightarrow Z$ and $\overleftarrow Z$ are identified with the corresponding $S$-letters. We shall call the state letters from $P_i$-s and $K_i$-s the  \label{pletter} $p$-{\it letters},
and the other state letters (i.e. the copies of the state letters of $\bf P$), the {\em basic state letters}.

The set of tape letters of ${\bf P}\circ\{\overrightarrow Z,\overleftarrow Z \} $  is $Y=Y_1\sqcup \dots \sqcup Y_{3N}= Y'_{1}\sqcup Y''_{1}\sqcup Y'''_{1}\sqcup
Y'_{2}\sqcup Y''_{2}\sqcup Y'''_{2}\sqcup...\sqcup
Y'_{N}\sqcup Y''_{N} \sqcup Y'''_{N}.$

Assume $\theta$ is a positive $\bf P$-rule
of the form
$$[s_0u_1\to s_0'u'_1, v_1s_1u_2\to
v_1's_1'u'_2,...,v_{N}s_N\to v_{N}'s_{N}'],$$
where $s_i, s_i'\in
Q_i$, and $v_i$-s, $u_i$-s are words in $Y.$ Then this rule is replaced in ${\bf P}\circ\{\overrightarrow Z,\overleftarrow Z \} $
by positive $$\bar\theta=\left[\begin{array}{l}s^{(\theta,0)}\tool (s')^{(\theta,0)},\; k^{(\theta,1)}(3)u_1
\to k^{(\theta,1)}(1)u'_1,\; v_1p^{(\theta,1)}(3)\tool v'_1p^{(\theta,1)}(1),\\ s^{(\theta,1)}\tool (s')^{(\theta,1)},\;k^{(\theta,2)}(3)u_2\to k^{(\theta,2)}(1)u'_2,\; v_2p^{(\theta,2)}(3)
\tool v'_2p^{(\theta,2)}(1), ...,
\end{array}\right]$$
 with $Y_{3i-2}(\bar\theta)=\emptyset $, $Y_{3i-1}(\bar\theta)=  Y_{i}''(\theta)$ and $Y_{3i}(\bar\theta)=Y_{i}'''(\theta).$

 Now we want to describe the alternating work of the auxiliary S-machines \label{Ztilr}
 $\overleftarrow Z^{(\theta,i)}$ and $\overrightarrow Z^{(\theta,i)}$. Normally
 each of them is switched on exactly once in the frame of the rule $\theta,$
 but the sequence of their turning on depends on $\theta.$
First, we need the following {\em transition}
\label{zmt} rule $\chi_-(\theta).$
This rule adds $\theta$ to all state letters
and turns all $k^{(j)}$ and $p^{(j)}$ into $k^{(\theta,j)}(1)$ and $p^{(\theta,j)}(1)$:

$$[s_i\tool s^{(\theta,i)}, k^{(j)}\to
k^{(\theta,j)}(1), p^{(j)}\tool
p^{(\theta,j)}(1), i=0,...,N, j=1,...,N].$$

Then the S-machines $\overrightarrow Z^{(\theta,1)},\dots,$ $ \overrightarrow Z^{(\theta,N)}$ and $\overleftarrow Z^{(\theta,1)},\dots,$ $ \overleftarrow Z^{(\theta,N)}$ are switched on in a specific
order (defined below) after the rule $\chi_-(\theta)$ is applied.
So the state letters $k^{(\theta,j)}(1)$, $p^{(\theta,j)}(1)$ ($j=1,\dots, N$) successively turn
into $k^{(\theta,j)}(3)$, $p^{(\theta,j)}(3),$ find themselves just after $s_{i-1}$- and before $s_{i}$-letters, respectively, and the rule $\bar\theta$ can be applicable.

After an application of  $\bar\theta,$ the S-machines
$\overrightarrow Z^{(\theta,1)},\dots,$ $ \overrightarrow Z^{(\theta,N)}$ and $\overleftarrow Z^{(\theta,1)},\dots,$ $ \overleftarrow Z^{(\theta,N)}$ are switched on again in the
following order.

Assume that the rule $\theta$
is a rule of a primitive S-machine $\cal P$. (Recall that $\bf P$ is composed
from primitive machines.) The S-machine $\cal P$ can work in several sectors.
(For example, $\theta$ can be the
control rule checking all the big historical sectors simultaneously.)
Let $i$ be the minimal index such that $S_i$ has a control running state
letter of $\cal P$. Then the rule $\bar\theta$
first switches on the S-machines $\overleftarrow Z^{(\theta,j)}$
for $j=i$ and simultaneuosly for all other $j$-s, where $S_j$ is also
has a control running state letter of $\cal P$.
The last rule $\chi_4$ of this S-machine switches on the S-machine
$\overleftarrow Z^{(\theta,i-1)}$ and similar S-machines in similar
sectors (e.g. in all small historical sectors if the sectors $S_{i-2}S_{i-1}$ is a small historical sector). The next S-machine is
$\overleftarrow Z^{(\theta,i-2)}$,
if it did not work ealier,
and so on. Then the S-machines $\overrightarrow Z^{(\theta,i+1)},\dots$ subsequently work, except for the sectors, where the auxiliary $Z$-machines worked earlier.

The same S-machines work after the application of the rule
$\chi_-(\theta)$ but they are switched on in the inverse order.

Finally, the transition
rule \label{zmt} $\chi_+(\theta)$ removes the index $\theta$ from all state letters, and turns all  $p^{(\theta,j)}(3)$ into $p^{(j)}$:

$$[s^{(\theta,i)}\tool s_i, k^{(\theta,j)}(3)\to
k^{(j)},  p^{(\theta,j)}(3)\to
p^{(j)}, i=0,...,N, j=1,...,N].$$

 For every admissible word $W$ of $\bf P$ with standard base, let
$\iota(W)$ be the admissible word of ${\bf P}\circ\{\overrightarrow Z,\overleftarrow Z \} $ obtained by inserting state letters $k^{(i)}$ and $p^{(i)}$ next to the
right of each $s_{i-1}$  and next to the
left of each $s_{i}$, $i\le N$.

Assume that $W\to W\cdot\theta$ is a computation of the S-machine $\bf P$
with standard base and a positive rule $\theta$.
Then, by the definition of ${\bf P}\circ\{\overrightarrow Z,\overleftarrow Z \}$,
we have the canonically defined reduced computation $\dots \to\iota(W)\to \iota(W)\cdot\bar\theta\to\dots$ starting and ending with words whose state
letters have no $\theta$-indices and all other words do have $\theta$-indices.
The computation of ${\bf P}\circ\{\overrightarrow Z,\overleftarrow Z \} $  with
these properties is unique by Lemma \ref{prim} (3) since
the S-machines $\overrightarrow Z$ and $\overleftarrow Z  $ are primitive.
Thus the following claim is true.

\begin{lemma}\label{pi23} (similar to Lemma 4.24 from \cite{O12}). For every computation $W\to W\cdot\theta$ of the S-machine
$\bf P$  with standard base and a positive rule $\theta,$ there is a unique
reduced ${\bf P}\circ\{\overrightarrow Z,\overleftarrow Z \}$ -computation
$\dots \to \iota(W)\to \iota(W)\cdot\bar\theta\to\dots$
starting and ending with words whose state
letters have no $\theta$-indices and all other words have $\theta$-indices. The history
of this computation starts with $\chi_-(\theta)$ and ends with $\chi_+(\theta).$
\end{lemma}

For every admissible word $W$ of ${\bf P}\circ\{\overrightarrow Z,\overleftarrow Z \} $
with
the standard base,
let \label{pi32} $\pi(W)$ be the word obtained by removing
state $k$- and $p$-letters, $\theta$-indices of state letters, and the indices that
distinguish $a$-letters from the left and from the right of $k$- and $p$-letters.
After possible cancellations of $a$-letters, we obtain an admissible
word of $\bf P$. Note that we have

Given a computation $\cal C$ of the  S-machine ${\bf P}\circ\{\overrightarrow Z,\overleftarrow Z \} $ with standard base and
history $H$ involving a $\bar\theta$-rule ,
we define the {\it projection} $\pi({\cal C})$ of it, which is a
computation of $\bf P$. To obtain it, one removes all transitions given by $\chi$-rules and replaces the configurations by their projections.
Note that this operation
makes sense since $\chi$-rules do not change the projection of
the word  onto a word in $Y_i$.
The projection $\pi(H)$ of the history is defined in obvious way:
one forgets the $\chi$-rules and removes bars over $\theta$-rules.

\begin{lemma} \label{red} (also see Lemma 4.28 in \cite{O12}). If $\cal C$ is a reduced computation of ${\bf P}\circ\{\overrightarrow Z,\overleftarrow Z \} $ with standard base, then $\pi(\cal C)$ is a reduced computation of ${\bf P}$.
\end{lemma}
\proof
Assume that we have a subword $\bar\theta H'\bar\theta^{-1}$ in the history $H$ of ${\cal C}$,
where $H'$ has no
rules of the form $\bar\theta$. If $H'$ is non-empty, then one obtains a contradiction
by Lemma \ref{prim} (4) applied to the work of the primitive
S-machines $\overrightarrow Z^{(\theta,i)}$ and $\overleftarrow Z^{(\theta,i)}$ . Hence $\pi(H)$ is a reduced history.
\endproof

Below we change some formulations of Subsection \ref{prim} as applied to ${\bf P}\circ\{\overrightarrow Z,\overleftarrow Z \} $.

\begin{lemma}\label{primm} (duplicate of \ref{prim}). Let ${\cal C}:\; C_0\to\dots\to C_t$ be a reduced computation of the S-machine ${\bf P}\circ\{\overrightarrow Z,\overleftarrow Z \} $ with the standard
base and  with $t\ge 1$. Then the following properties hold.

(1) If $|C_i|_a>|C_{i-1}|_a$ for some $i=1,\dots,t-1$, then $|C_{i}|_a\le |C_{i+1}|_a\le |C_{i+2}|_a,\le\dots$

(2) $|C_i|_a\le\max(|C_0|_a, |C_t|_a)$ for every $i=0,1,\dots t$.

(3) Assume that the words $C_0$ and $C_t$ have $a$-letters only from the
subalphabets $Y''_i$ ($i=1,\dots, N$) and that for the primitive S-machines $\cal P$ forming $\bf P$, all their subwords in $\pi(C_0)$ and $\pi(C_t)$ look like in Lemma \ref{prim} (3),
i.e. as $q^1up^{1}q^2$ and $q^1vp^{2}q^2$ for some words $u,v$ in the notation of Lemma \ref{prim}.
Then  $a$-words in the corresponding sectors of $C_0$ and $C_t$ are equal,    $|C_i|_a=|C_0|_a$ for every $i=0,\dots,t$ and $t=\Theta(s^2)$, where $s=|C_0|_a$.

(4) Assume that the words $C_0$ and $C_t$ have $a$-letters only from the
subalphabets $Y''_i$ ($i=1,\dots, N$) and that for the primitive S-machines $\cal P$ forming $\bf P$, all their subwords in $\pi(C_0)$ and $\pi(C_t)$ look like in Lemma \ref{prim} (4). Then
it is not possible that
the configurations $C_0$ and $C_t$ have the same set of state letters.

(5)
If $C_0$ (or $C_t$) satisfies the assumptions of item (3), then $|C_i|_a\ge |C_0|_a$
(respectively,  $|C_i|_a\ge |C_t|_a$) for every $i=0,\dots,t$.

\end{lemma}

\proof Let us start with Property (1).
If ${\cal C}_i={\cal C}_{i-1}\chi_1^{(\theta,i)}(a)^{\pm1}$ (or ${\cal C}_i={\cal C}_{i-1}\chi_2^{(\theta,i)}(a)^{\pm1}$), then $p(1)^{(\theta,i)}$
inserts letters from both sides and the next rule of the computation must be again $(\chi_1^{(\theta,i)})^{\pm 1}(c)$ for some $c$. It again must increase the length of the configuration by two, and so on.

If ${\cal C}_i={\cal C}_{i-1}\bar\theta$ for some $\bar\theta$-rule,
then the transition $\pi(C_{i-1})\to\pi({\cal C}_i)$ increases the length
by Lemmas \ref{red} and \ref{prim} (1). The work of $\overrightarrow Z^{(\theta,i)}$,$\overleftarrow Z^{(\theta,i)}$ cannot decrease configuration length by Lemma \ref{prim} (5) for these primitive S-machines.
Therefore Statement (1) is true and Statement (2)
is also true since one can choose a shortest $C_j$ and consider the subcomputation $C_j\to\dots\to C_t$ and inverse
subcomputation $C_j\to\dots\to C_0$.

To prove equalities $|C_i|_a=|C_0|_a$ in Statement (3), one just
apply Lemma \ref{prim} (3) first, to $\pi({\cal C})$ and second,
to the maximal subcomputations of $Z$-machines. Since
the length of $\pi({\cal C})$ is $2k+1$ by that Lemma, and if $k>0$
the maximal subcomputations of $Z$-machines work
with configurations of lengths, $k-1,k-2,\dots 1$, we obtain
the later claim of (3) by summation.

Property (4)  follows from Lemma \ref{prim} (4) applied to the
computation $\pi(\cal C)$ and the maximal subcomputations of  $Z$-machines. Property (5) follows from the projection argument as
in Lemma \ref{prim} (5).
\endproof

\begin{lemma}\label{ewee} (duplicate of \ref{ewe}) If $C_0\to\dots\to C_t$ is a reduced
computation of ${\bf P}$ with base $S_{i-1}K_iK_i^{-1}S_{i-1}^{-1}$ or $S_i^{-1}P_i^{-1}P_iS_i$ and $C_0$ has $a$-letters from the alphabet $Y_i''$ only.
Then $|C_i|_a\ge |C_0|_a$ for every $i=0,\dots,t$.
\end{lemma}

\proof The statement follows from the projection argument \ref{proj} as in Lemma \ref{ewe}.
\endproof
Let us call the constructed S-machine ${\bf P}\circ\{\overrightarrow Z,\overleftarrow Z \} $
{\it biprimitive}.

To define the modified S-machine ${\bf M}_4'$ we insert two more base
letters in each pair ${\cal R}_{i-1}$ and ${\cal P}_i$ of the standard base,
i.e. this base has the subwords ${\cal Q}_{i-1}{\cal R}_{i-1}{\cal R}'_{i-1}{\cal P}'_i{\cal P}_i{\cal Q}_i$;  now the first big historical
sector is  ${\cal R}'_0{\cal P}'_1$ (instead of ${\cal R}_0{\cal P}_1$). In the definition of Step $1^-$, we now replace the rules of the primitive control S-machines
with the rules of the corresponding biprimitive machines. By definition, at all other
steps
the control S-machines are
just primitive, i.e. the sectors ${\cal R}_{i-1}{\cal R}'_{i-1}$
and ${\cal P}'_i{\cal P}_i$ are locked, there are no
$\chi$-rules, and two base letters of the locked sectors
work as one letter of a primitive S-machine.

Thus, only at Step $1^-$, the biprimitive S-machine ${\bf P}'$ works.

\begin{lemma}\label{Hprimm}  (duplicate of \ref{Hprim}) Let ${\cal C}:\;
W_0\to\dots\to W_t$ be a reduced computation of ${\bf P}'$ with standard base. Then

(a) $|W_j|_a\le\max (|W_0|_a,|W_t|_a)$ for every configuration of $\cal C$; moreover, $|W_0|_a\le \dots\le |W_t|_a$ if every
control $P$- and $R$-letter
neighbors some $Q$-letter and $P'$-letter ($R'$-letter) neighbors a $P$-letter (resp., an $R$-letter) in the word $W_0$;

(b) we have $t= O(||W_0||^2+||W_t||^2)$,
moreover, $t= O(||W_t||^2)$
if every $p$-letters and control letters has a neighbor in the word $W_0$ as in item (a).
\end{lemma}

\proof (a)
The first property is given by Lemma \ref{primm} (2).
Under the additional assumption for control letters and $p$-letters, $W_0$ is the shortest configuration by the projection argument.

(b) Let us say that a subcomputation  without $\bar\theta$-rules in the history is a $Z$-subcomputation.
If the computation is a $Z$-subcomputation, then the statement
(b) follows from Lemma \ref{Hprim} (b) for the (composition of) primitive S-machines.
Otherwise we chose a $\bar\theta$-transition $W_r\to W_{r+1}$
with minimal $\min(||W_r||, ||W_{r+1}||)$. Without loss of generality, one may assume that this minimum is $||W_r||$. Lemma
\ref{Hprim} (a) for $\pi(\cal C)$ and for $Z$-subcomputations
implies that $||W||_r\le ||W_i||$  for every $i$,
and therefore it suffices to bound the histories of length-non-decreasing
computations ${\cal C}':\;W_r\to\dots\to W_t$ and ${\cal C}'':\;W_r\to\dots\to W_0$.

The length of $\pi({\cal C}')$ is at most $2||W_t||$ by Lemma
\ref{Hprim} (b) for the S-machine $\bf P$. By the same Lemma for maximal $Z$-subcomputations, the length of every maximal
$Z$-subcomputation is at most $2||W_t||$ too.
Since the number of such maximal subcomputations is at most
$2||W_t||$, we have that $t-r=O(||W_t||)^2)$.
Similarly, we obtain $r=O(||W_0||)^2)$, and the
the first estimate of (b) is obtained. For the second estimate,
one can choose $r=0$, since the whole computation $\cal C$ is length-non decreasing by Lemma \ref{primm} (5,1).
\endproof

We see that Lemmas \ref{primm} (3) and \ref{Hprimm} (b) provide
us with quadratic estimate of the computation time for the biprimitive S-machine against the linear time for S-primitive machines.
This have a few consequences mentioned below.

For the same function $f(x)$ recognized by the original
Turing machine and  $F(x)=x^2g(x)$,  we define the function $g(x)$ to be equivalent to $xf(x)^3$ now.

The extra sectors of the biprimitive S-machine do not affect the
work of all other steps except for $1^-$ since they are locked therein.

The formulation of Lemma \ref{121}(1) is unchanged but the
proof is now based on Lemma \ref{primm} (4). For the step history
$(21^-)(1^-)(1^-1)$ in the formulation of Lemma \ref{121}(2),
we have now $||H||=O(||W_0||^2)$, which follows from
Lemma \ref{Hprimm} (b). The estimates of $||H||$ for other step
histories mentioned in Lemma \ref{121}(1) remain unchanged.
Also we add an item to the formulation of Lemma  \ref{121}:

(3) {\it Let the history of a reduced computation ${\cal C}: W_0\to\dots\to W_t$ with standard base
have a subword $\bar\theta_1H\bar\theta_2$, where $\bar\theta_1$ and $\bar\theta_2$ are
$\bar\theta$-rules of the S-machine ${\bf P}'$ and $H$ has only $\chi$-rules of ${\bf P}'$.
Then all configurations of $\cal C$ are uniquely determined by
$H$, $|W_1|_a=\dots=|W_{t-1}|_a$ and $||W_j||=\Theta(||H||)$ for $j\le t$.}

The proof of this statement follows the proof of Lemma \ref{121} (2), but
now one refers to Lemma \ref{pi23} instead of Lemma \ref{prim} (3).

To the assumption of Lemma \ref{123}, we add: {\it "or the history of the computation $\cal C$ has a subword $\theta_1H\theta_2$, where $\theta_1$ and $\theta_2$ are
$\bar\theta$-rules of the S-machine ${\bf P}'$ and $H$ has only $\chi$-rules of ${\bf P}'$"}.
In the proof of Lemma \ref{123} (2), one should use the following
property. The computation of the biprimitive S-machine does not
change the length of configurations by Lemma \ref{Hprimm} (a).
Hence the reduced computation of it is canonical by Lemma \ref{pi23}, and the history
restores the tape words.

Consider the single Step $1^-$ in item 2 of the proof of Lemma
\ref{nonst}. If we have the work of primitive S-machines only,
then the proof is unchanged. If there is a $\bar\theta$-rule of the
biprimitive S-machine in the computation $\cal C$, we condider
the projection $\pi(\cal C)$, where the $a$-lengths of all configurations are at most $C(|W_0|_a+|W_t|_a)$ by Lemma
\ref{nonst}. It remains to consider maximal subcomputation
$W_r\to\dots\to W_s$ of $Z$-machines, where
$|W_r|_a, |W_s|_a\le C(|W_0|_a+|W_t|)_a$. By Lemma \ref{Hprim} (a) we have $|W_i|_a \le \max(|W_r|_a, |W_s|_a)\le C(|W_0|_a+|W_t|_a)$
for $r\le i\le s$, as required.

In item 3 of the proof of Lemma \ref{nonst}, the base letters of
the first big historical sector should be replaced by their dashed duplicates.

The estimate $||H||\le c_2||W_0||$ of Lemma \ref{boundH} changes now by the quadratic
estimate $||H||\le c_2||W_0||^2$ due to the application of Lemma
\ref{Hprimm} (b). Respectively, the upper bound $c_3(k^3+1)(||W_0||+k^3)$ of Lemma \ref{stand} is now replaced
by $c_3(k^3+1)(||W_0||^2+k^3)$.

The formulation of Lemma \ref{sta} does not change since comparing
the lengths of histories, we now increase all of them. In particular, we have now in equation (\ref{sc}) that $2w^2+2(w-1)^2+\dots+2(w-(m-1))^2\ge |m|w^2/2$, we have that
the difference of lengths of subcomputations ${\cal C}_{i-1}$ and
${\cal C}_i$ does not exceed $10w$, we have $||H_0||=O(w^2)$
and obtain $||H'H''||\le c_2(k^3+1)O(w^2)$, which leads to the same estimates for $||H_0||$ and $||H_t||$ since the constant $c_4$
can be chosen large enough.

The change of linear estimates by quadratic ones in the proof
of Lemma \ref{EFE} just sharpens the required inequalities.

Lemma \ref{gtime} claims now that the generalized time function $T'(n)$ of the S[machine $M'$ is equivalent to $n^2f^3(n)$. Thus, we multiply
the generalized time function of $\bf M$ by $n$. This is sufficient for
the upper bound. Indeed the replacements of the form $||W_j||\to ||W_j||^2$ in the proof of the modified Lemma \ref{gtime} can
multiply the generalized time function at most by $n$ since it
is shown there that $||W_j||=O(n)$ for every configuration $W_j$.
The lower bound obtained in the original Lemma \ref{gtime} must also
be multiplied by $n$ now. Indeed, the time of the constructed
subcomputations with step history $\big((21^-)(1^-)(1)(2^-)(2)(21^-)\big)^{\pm 1}$
will be at least $\Theta(n^2)$ now (instead of $\Theta(n)$) since by Lemma \ref{primm} (3), we have such lower bound for the Step $1^-$.

Similar replacement should be made in Remark \ref{ar}: There are $\Theta(n^2f(n)^3)$ configurations of
length at least $\Theta(n)$ for any computation accepting the word $V(n)$.

The statement of Lemma \ref{B} is modified now by adding the words
"{\it or the history of $\cal C$ has a subword $\bar\theta_1H\bar\theta_2$, where $\bar\theta_1$ and $\bar\theta_2$ are
$\bar\theta$-rules of the S-machine ${\bf P}'$ and $H$ has only $\chi$-rules of ${\bf P}'$}".
If the step history of $\cal C$ is $(F)$ and Property (a) of Lemma \ref{sta} holds
we have the same proof as in the original Lemma \ref{B} since the length of a subcomputation of Step $1^-$ is at most $5\max (||W_0||_a, ||W_t||_a)$ provided
the history of $\cal D$ has no subwords $\bar\theta_1H\bar\theta_2$. Indeed, in this case,
the history of any subcomputation of Step $1^-$ is a subword of $H'\bar\theta H''$,
where $H'$ and $H''$ have only $\chi$-rules of ${\bf P}'$, and one can bound
each of  $||H'||$ and $||H''||$ by $2c_4(||W_0||_a, ||W_t||_a)$ applying Lemma
\ref{Hprim} (b).

The formulation of Lemma \ref{pol} remains unchanged for all steps except for
Step $1^-$. If the step is $1^-$, we add the assumption that {\it the history has
no $\bar\theta$-rules of the S-machine} ${\bf P}'$. Since the biprimitive S-machine ${\bf P}'$
works as a primitive one in the later case, the proof does not change.

The set of {\it standard trapezia} is enlarged now, namely, a trapezium with standard
base having a subword $\bar\theta_1H\bar\theta_2$ in the history, where $\bar\theta_1$ and $\bar\theta_2$ are $\bar\theta$-rules of the S-machine ${\bf P}'$ and $H$ has only $\chi$-rules of ${\bf P}'$,
is also standard by definition.

The above modification of Lemma \ref{pol} changes the proof of Lemma \ref{led}
as follows. One obtains a contradiction if
for some $k$, there is one-step
subcomputation $W(k)\to\dots\to W(k+4)$, where the step differs from $1^-$
{\it or the step is $1^-$ and the history of this subcomputation has no $\bar\theta$-rules
of the S-machine ${\bf P}'$}. Hence either every subcomputation $W(k)\to\dots\to W(k+8)$
has at least two steps, and so the computation $W(1)\to\dots \to W(999)$ has at least $100$ steps, which leads to a contradiction as in the original proof of Lemma \ref{led},
or for some $k$, one obtains at least two $\bar\theta$-rules in the history of
$W(k)\to\dots\to W(k+8)$, i.e., the history of this subcomputation has a subword $\bar\theta_1H\bar\theta_2$, where $\bar\theta_1$ and $\bar\theta_2$ are
$\bar\theta$-rules of the S-machine ${\bf P}'$ and $H$ has only $\chi$-rules of ${\bf P}'$.
In the later case one should use the modified Lemma \ref{123} (instead of
the original Lemma \ref{123}) to complete the proof.

\medskip

Lemma \ref{led} proves Theorem \ref{beta} for $s=2$. In the
present subsection, we have modified the main S-machine and the corresponding groups $M$ and $G$ so that the modified Lemma \ref{led} provides us with the statement of Theorem \ref{beta} for $s=3$. One can make further modifications, which similarly slow
down the work of the previously modified S-machines. This will
give the proofs of Theorem \ref{beta} for $s=4, 5,\dots$. However
we can leave the details to the reader taken into account that for $s\ge 4$, Corollary \ref{22m} is obtained in \cite{SBR}. Thus, the proof of Theorem \ref{beta} is complete.

\medskip

{\bf Acknowledgment.}

The author is grateful to Mark Sapir for useful discussions.

\twocolumn

\noindent {\bf Subject index} \label{sind}
\bigskip

\noindent $a$-band \pr{aband}\\
\noindent accepted configuration \pr{acceptedc}\\
\noindent accept configuration \pr{acceptc}\\
\noindent admissible word \pr{admissiblew}\\
\noindent $a$-edge \pr{aedge}\\
\noindent $a$-length \pr{alen}\\
\noindent $a$-letter \pr{aletter}\\
\noindent arc \pr{arc}\\
\noindent application of a rule \pr{application}\\
\noindent area of diagram \pr{aread}\\
\noindent area of word \pr{areaw}\\
\noindent band \pr{band}\\
\noindent base of $\theta$-band \pr{baseb}\\
\noindent base of word \pr{basew}\\
\noindent basic width of comb \pr{basicw}\\
\noindent bead \pr{bead}\\
\noindent big trapezium \pr{bigt}\\
\noindent block history \pr{blockh}\\
\noindent bottom of band \pr{bottomb}\\
\noindent bottom of trapezium \pr{bottomt}\\
\noindent chord \pr{chord}\\
\noindent clove \pr{clove}\\
\noindent comb \pr{com}\\
\noindent combinatorial length $||\cdot||$\pr{lengthc}\\
\noindent composition of S-machines \pr{composm}\\
\noindent computable in time $T(n)$ \pr{compT}\\
\noindent computation \pr{computation}\\
\noindent configuration \pr{config}\\
\noindent connecting rule \pr{conrule}\\
\noindent control state letter \pr{csl}\\
\noindent control step \pr{cstep}\\
\noindent copy of word \pr{copyw}\\
\noindent crown \pr{crown}\\
\noindent Dehn function \pr{Dehnf}\\
\noindent derivative subcomb \pr{ders}\\
\noindent design \pr{desn}\\
\noindent diagram \pr{diagram}\\
\noindent disk \pr{diskr}\\
\noindent equivalence of functions \pr{equivf}\\
\noindent equivalence of machines \pr{equivM}\\
\noindent extension of arc \pr{extena}\\
\noindent faulty base \pr{faul}\\
\noindent functions $\Phi(x)$ and $\phi(x)$ \pr{Phif}\\
\noindent $G$-area of diagram over $G$ \pr{areaGmd}\\
\noindent $G$-area of diagram over $M$ \pr{areaGd}\\
\noindent $G$-area of big trapezium \pr{abt}\\
\noindent $G$-area of disk \pr{areaGdi}\\
\noindent generalized time function \pr{gentf}\\
\noindent groups $M$ and $G$ \pr{MG}\\
\noindent handle of comb \pr{handle}\\
\noindent height of trapezium \pr{heightt}\\
\noindent historical sector \pr{histsec}\\
\noindent history of $q$-band \pr{historyb}\\
\noindent history of computtion \pr{historyc}\\
\noindent history of trapezium \pr{historyt}\\
\noindent $H'$-part of trapezium \pr{H'partt}\\
\noindent hub cell \pr{hubs}\\
\noindent hub relation \pr{hubr}\\
\noindent input configuration \pr{inpt}\\
\noindent $J$-mixture on a necklace \pr{Kmix}\\
\noindent length of arc \pr{lena}\\
\noindent locked sector \pr{locks}\\
\noindent S-machine checking disibility \pr{mcd}\\
\noindent main S-machine \pr{M6}\\
\noindent maximal band \pr{maxb}\\
\noindent minimal diagram \pr{minimald}\\
\noindent mixture on the boundary \pr{muK1mix}\\
\noindent modified length function $|\cdot|$ \pr{modiflf}\\
\noindent necklace \pr{neckl}\\
\noindent parallel arc \pr{paralla}\\
\noindent parallel and sequential work of primitive S-machines \pr{ps}\\
\noindent parameter \pr{param}\\
\noindent part of rule \pr{part}\\
\noindent peripheral chord \pr{peric}\\
\noindent peripheral half-disc \pr{perif}\\
\noindent positive, negative rules \pr{positiver}\\
\noindent primitive S-machine \pr{pm}\\
\noindent projection of word \pr{projectw}\\
\noindent property $P(\lambda,n)$ \pr{propP}\\
\noindent $q$-band \pr{qband}\\
\noindent $q$-edge \pr{qedge}\\
\noindent $q$-letter \pr{qletter}\\
\noindent quasi-trapezium \pr{quasitr}\\
\noindent reduced computation \pr{reducedc}\\
\noindent reduced diagram \pr{reducedd}\\
\noindent revolving base and computation \pr{revolv}\\
\noindent rim band \pr{rimb}\\
\noindent rule \pr{rule}\\
\noindent S-machine \pr{Smachine}\\
\noindent sector \pr{sectorw}\\
\noindent shaft, $\lambda$-shaft \pr{shf}\\
\noindent side of band \pr{sideb}\\
\noindent sides of trapezium \pr{lrsidest}\\
\noindent spoke \pr{spoke}\\
\noindent standard base \pr{standardb}\\
\noindent standard history \pr{standh}\\
\noindent standard trapezium \pr{sttrap}\\
\noindent start/end edges \pr{seedgesb}\\
\noindent state letter \pr{statel}\\
\noindent stem band \pr{stb}\\
\noindent stem $\Delta^*$ of a diagram $\Delta$ \pr{stD}\\
\noindent step history \pr{steph}\\
\noindent step history of band \pr{stephb}\\
\noindent step history of trapezium \pr{stepht}\\
\noindent steps of ${\bf M}_4$, ${\bf M}_5$ and ${\bf M}_6=\bf M$ \pr{step}\\
\noindent subcomb \pr{subc}\\
\noindent suitable function \pr{sui}\\
\noindent tape letter \pr{tapel}\\
\noindent $t$-band \pr{tband}\\
\noindent $\theta$-band \pr{thband}\\
\noindent $(\theta,a)$-cell \pr{ta}\\
\noindent $(\theta,q)$-cell \pr{tq}\\
\noindent $\theta$-edge \pr{thedge}\\
\noindent $\Theta$-equivalence \pr{Tequiv}\\
\noindent $\theta$-letter \pr{thetal}\\
\noindent $(\theta,a)$-relation \pr{thetaar}\\
\noindent $(\theta,q)$-relation \pr{thetaqr}\\
\noindent tight base \pr{tightb}\\
\noindent tight comb \pr{tightc}\\
\noindent time function \pr{timef}\\
\noindent top of band \pr{topb}\\
\noindent top of trapezium \pr{topt}\\
\noindent transposition of disk and $\theta$-band \pr{transpo}\\
\noindent trapezium \pr{dftrap}\\
\noindent weakly minimal diagram \pr{wminimald} \\
\noindent working sector \pr{wsect}\\
\noindent $\cal Z$-annulus \pr{annulus}\\
\noindent ${\cal B}_1,\dots, {\cal B}_{L-3}$ \pr{b1bL3}\\
\noindent ${\cal B}_r$ \pr{Br}\\
\noindent $E_j, E_j^0$ \pr{Ej}\\
\noindent $H_1,\dots, H_{L-3}$ \pr{histL3}\\
\noindent $h_1,\dots, h_{L-3}$ \pr{hidtl}\\
\noindent $\Lab(\cdot)$ \pr{Lab}\\
\noindent ${\bar\bf p}$ \pr{barp}\\
\noindent ${\bf p}_{ij}$ \pr{pij}\\
\noindent ${\bar\bf p}_{ij}$ \pr{barpij}\\
\noindent ${\bf p}={\bf p}(\Psi)$ \pr{pPsi}\\
\noindent ${\bf q}_{ij}$ \pr{qij}\\
\noindent ${\bf y}_j$ \pr{yj}\\
\noindent  ${\bf z}_j$ \pr{zj}\\
\noindent $\Gamma_j$ \pr{Gammaj}\\
\noindent $\Gamma'_j$ \pr{Gammaj'}\\
\noindent $\bar\Delta$ \pr{bdel}\\
\noindent $\bar\Delta_{ij}$ \pr{barij}\\
\noindent $\bar\Delta^0, \Psi^0, \Psi^0_{ij}$ \pr{Psi0ij}\\
\noindent $\Delta^*$ \pr{stD}\\
\noindent $\Psi, \Psi_{ij}$ \pr{cloves}\\
\noindent $\Psi'$ \pr{psi'}\\
\noindent $\Psi'_{ij}$ \pr{Psiij'}\\
\noindent $\sigma_{\lambda}$ \pr{slam}\\
\noindent $\equiv$ \pr{sequiv}\\
\noindent $\sim$ \pr{equivf}\\
\noindent $\tool$ \pr{tool}\\

\addtocontents{toc}{\contentsline {section}{\numberline { }Subject
index \hbox {}}{\pageref{sind}}}

\bigskip

\begin{minipage}[t]{3 in}
\noindent Alexander Yu. Ol'shanskii\\ Department of Mathematics\\
Vanderbilt University \\ alexander.olshanskiy@vanderbilt.edu\\

\end{minipage}
\begin{minipage}[t]{3 in}
\noindent
and Department of
Higher Algebra,\\ MEHMAT,
 Moscow State University\\

\end{minipage}

\end{document}